\documentclass[a4paper,12pt]{article}
\usepackage{amsfonts}
\usepackage{amssymb}
\usepackage{amsmath}
\usepackage{geophysics}
\usepackage{latexsym}
\usepackage{graphicx}
\usepackage{caption}

\newdimen \dummy
\dummy=\oddsidemargin  \marginparwidth=.675\dummy
\newcounter{cor}
\newcounter{lem}

\newtheorem{theorem}{Theorem}

\newtheorem{corollary}[cor]{Corollary}

\newtheorem{lemma}[lem]{Lemma}

\input{tcilatex}
\begin{document}

\title{Identification and estimation of Structural VARMA models using higher
order dynamics\thanks{%
Thanks to J.C. Escanciano, F.J. Hidalgo, I.N. Lobato, G. Sucarrat and
seminar participants at Nuffield College, LSE, QMUL, BI Norwegian BS and at
Workshop on Time Series Econometrics 2018, Waseda International Symposium
2018, CFE-CM Statistics Conference 2019, EC$^2$ Conference on Identification
in Macroeconomics 2019 for helpful discussions and comments. Financial
support from the Ministerio Econom\'{\i}a y Competitividad (Spain) grant
ECO2017-86009-P is gratefully acknowledged. } \\
}
\author{Carlos Velasco \\
Universidad Carlos III de Madrid\\
carlos.velasco@uc3m.es \  \  \  \  \ }
\date{August 28, 2020}
\maketitle

\begin{abstract}
We use information from higher order moments to achieve identification of
non-Gaussian structural vector autoregressive moving average (SVARMA)
models, possibly non-fundamental or non-causal, through a frequency
domain criterion based on a new representation of the higher order
spectral density arrays of vector linear processes. This allows to identify
the location of the roots of the determinantal lag matrix polynomials based
on higher order cumulants dynamics and to identify the rotation of the model
errors leading to the structural shocks up to sign and permutation. We
describe sufficient conditions for global and local parameter identification
that rely on simple rank assumptions on the linear dynamics and on finite
order serial and component independence conditions for the structural
innovations. We generalize previous univariate analysis to develop
asymptotically normal and efficient estimates exploiting second and non-Gaussian higher
order dynamics given a particular structural shocks ordering without
assumptions on causality or invertibility. Bootstrap approximations to
finite sample distributions and the properties of numerical methods are
explored with real and simulated data.\bigskip

\bigskip


\end{abstract}
\bigskip

\section{1. Introduction}

There is an increasing literature on the application of Structural VARMA
(SVARMA) models for the analysis of economic data which tries to solve the
identification problem of these models by incorporating information from the
distribution of non-Gaussian structural shocks. This information recovered
from the data can substitute, at least in part, the restrictions provided by
economic theory on the impulse response functions (IRF) of endogenous
variables to given shocks. Furthermore, there are many examples where it is
not possible to discard non-fundamental solutions to a wide class of dynamic
macroeconomic models with features affecting the flow of information used by agents to make
decisions, see e.g. the surveys in Alessi, Barigozzi and Capasso (2011) and
Gouri\'{e}roux, Monfort and Renne (2019). In parallel, the identification of
noncausal structural VAR models have been also investigated, see e.g. Lanne
and Saikkonen (2013).

The identification analysis of SVARMA models has to account for both static
(Structural) and dynamic (VARMA) aspects. The dynamic problem
is related to the location of the determinantal roots of the VAR and VMA
polynomials that lead to causal/noncausal or invertible/noninvertible
solutions, respectively, guaranteeing that model errors are unpredictable,
not just a serially uncorrelated white noise sequence. The static
identification relates to the choice of the particular rotation of the
reduced form errors that delivers the true vector of structural shocks with
proper economic interpretation. Therefore, these shocks must additionally
satisfy some mutual independence condition strengthening the uncorrelation
achieved by any square root transformation of the covariance matrix of the
reduced form errors.

Under Gaussianity, SVARMA identification is not possible in absence of
further restrictions provided by economic theory because uncorrelation is
equivalent to independence and therefore all the infinite sequences obtained
by different versions of the lag polynomials obtained by flipping roots
and/or rotating the different shocks through Blaschke orthogonal matrices
(see e.g. Lippi and Reichlin, 1994) would be admissible. However, under
non-Gaussianity and the independence component assumption (ICA) of the
structural shocks it is known that static rotations can be identified up to
permutation and sign, i.e. up to labeling of the shocks obtained, if at most
one of the innovations components is Gaussian (see Comon (1994), Hyv\"{a}%
rinen, Zhang, Shimizu and Hoyer (2010)), while a condition on higher order
cumulants and moments of serially
independent errors guarantees the dynamic identification (see Chan and Ho (2004), Chan, Ho and Tong (2006), Gouri%
\'{e}roux et al. (2019)). However, these results do not lead to specific
methods for designing parameter estimates and inference rules that rely on
the identifying assumptions and are easy to interpret and check for
particular models. Instead, most of available methods are based instead on
moment estimates for which local rank conditions are assumed after a basic
order condition is guaranteed or on (Pseudo) ML procedures which have to be
further justified. This is precisely the aim of this paper, to provide neat
inference methods exploiting efficiently global identification conditions
based on a minimal finite number of moments of the marginal and joint
distributions of the sequence of structural non-Gaussian errors.

Typically, non-Gaussianity is exploited for identification of dynamic models
through conditions on higher order cumulants (Gouri\'{e}roux et al., 2019;
Lanne and Luoto, 2019) or spectral densities (Lii and Rosenblatt, 1982;
Kumon, 1992), but it has also been imposed through particular probability
distribution assumptions on the shocks (Lanne and L\"{u}tkepohl, 2010) or
with conditional (Normandin and Phaneuf, 2004) and unconditional
heteroskedasticity conditions (Rigobon, 2003; Lanne and L\"{u}tkepohl,
2008), possibly with Markov switching dynamics (Lanne, L\"{u}tkepohl and
Maciejowska, 2010; L\"{u}tkepohl and Net\^{s}unajev, 2017). Then, estimation
is performed using ML or approximate versions of it (e.g. Lii and
Rosenblatt, 1992, 1996, for ARMA; Gouri\'{e}roux, Monfort and Renne, 2017,
and Lanne, Meitz and Saikkonen, 2017, for SVAR; Gouri\'{e}roux et al., 2019,
for SVARMA models) or non-Gaussian criteria like LAD or ranks (Breidt, Davis
and Trindade, 2001, and Andrews, Davis and Breidt, 2007, in the univariate
case). Methods based on higher order moments have been also developed, first
for the univariate case in the frequency and time domains (Lii and
Rosenblatt, 1982, Gospodinov and Ng, 2015, respectively). For multivariate
models, Gouri\'{e}roux et al. (2019) proposed a semiparametric 2-step
method, where first the VAR parameters are estimated using a 2SLS approach
under causality, and then the VMA parameters are estimated using moment
conditions on linear combinations of the residuals or using a PMLE
approximation to some prespecified non-Gaussian distribution. The
restrictions on moments of order 2, 3 and 4 are derived from ICA to improve
efficiency and raise the chances that they provide sufficient information to
guarantee usual local identification rank conditions. A similar approach is
pursued in Lanne and Luoto (2019) to achieve local identification of a SVAR
model by imposing a certain set of co-kurtosis conditions.

In this paper we study the problem of SVARMA identification extending the
frequency domain approach of Velasco and Lobato (2018), henceforth VL, to
the multivariate and structural case. VL showed that identification of a
possible noncausal or noninvertible ARMA model can be achieved by checking
that higher order spectral densities are sensitive to the location of the
roots of the lags polynomials unlike the usual second order spectral
density, i.e. they can achieve phase identification as noted in Lii and Rosenblatt
(1992). They also investigated model estimation using a minimum distance
criterion between the higher order periodograms and the parametric
specification of the corresponding higher order spectral densities of the
ARMA model that accounts efficiently for all moment conditions of a given
order at all lags. This approach provides a comprehensive method for dealing
with the problem of the location of roots of the lag polynomials and the
characterization of the non-Gaussian information through higher order
spectral densities to develop robust and efficient estimates, see also
Lobato and Velasco (2018).

To extend these ideas to the SVARMA setting, we first develop a new
representation of higher order spectral densities arrays for linear vector
processes and show that our identifying frequency domain criterion can
indeed discriminate processes that have observationally equivalent linear
covariance dynamics, but whose different IRFs are reflected on their higher
order dynamics. We are able to reproduce the previous dynamic and static
identification results found in the literature for non-Gaussian vector
models assuming only ICA and serial independence up to a given order (third
and/or fourth) and providing some extensions when some non-zero (i.e.
non-Gaussian) cumulant condition is violated 
or when no version of ICA holds but we impose a rank condition on the
innovations third order cumulant array. These results rely on a simple
non-singularity condition on the transfer function of the VARMA system so
that our criterion can evaluate all versions of the model up to a Blaschke
factor and the value of higher order cumulants of structural errors. 

This nonparametric global identification provides a constructive method for
designing minimum distance parameter estimates in the frequency domain which
can exploit efficiently all information contained in the dynamics of moments
of order 2, 3 and 4 without distributional assumptions or factorizations of
the matrix lag polynomials to deal with the simultaneous presence of roots
inside and outside the unit circle. Despite SVARMA identification up to a
signed permutation is enough for IRF and decomposition variance analyses
(given that a particular labeling can be attached to each shock), to obtain
standard asymptotic results for our parameter estimates we fix a unique
identified version of the model using a particular ordering and sign
structure on the innovations. These restrictions could be replaced by
alternative statistical conditions or economic information, which then would
become overidentification restrictions that could be tested in our
framework. We also develop bootstrap approximations for the asymptotic
distribution of parameter estimates and for the computation of efficient
estimates exploiting all moment conditions available. The finite sample
properties of a numerical algorithm to implement these identification and
inference methods are explored with real and simulated data.

The rest of the paper is organized as follows. Section~2 sets the
identification problem and introduces the main concepts and tools. Section~3
provides the basic identification results. Section~4 deals with parameter
identification and Section~5 with parameter minimum distance estimation.
Section~6 analyzes GMM efficient estimates exploiting information from
moments of several orders and bootstrap approximations to the distribution
of estimates. Section~7 presents the numerical methods and the simulation
experiment. Section~8 reanalyses Blanchard and Quah (1989) identification of
a bivariate system for US GNP growth and unemployment. A series of
appendices include additional discussion of concepts used in the paper,
together with proofs and auxiliary results.\bigskip

\section{2. Identification problem and assumptions}

We consider the SVARMA$\left( p,q\right) $ system 
\begin{equation*}
\Phi \left( L\right) Y_{t}=\mu +\Theta \left( L\right) \mathbf{\varepsilon }%
_{t},
\end{equation*}%
where the $d$-vector $\mathbf{\varepsilon }_{t}$ behaves as an independent
identically distributed (\textit{iid}) sequence up to a finite number $k$ of
moments, $k\geq 3,$ with zero mean and covariance matrix $\mathbf{I}_{d},$
the $d$-dimensional identity matrix, but with components not necessarily
mutually independent. 
The vector $\mu $ is an unknown level parameter and the lag polynomials with
matrix coefficients $\Phi \left( L\right) =\mathbf{I}_{p}+\Phi _{1}L+\cdots
+\Phi _{p}L^{p}$ and $\Theta \left( L\right) =\Theta _{0}+\Theta
_{1}L+\cdots +\Theta _{q}L^{q}$, $\Theta _{0}$ nonsingular, satisfy det$%
\left( \Phi \left( z\right) \right) $det$\left( \Theta \left( z\right)
\right) \neq 0$ for $\left \vert z\right \vert =1$. These conditions
guarantee the existence of a stationary solution for $Y_{t}.$ Note that we
allow the roots of the determinants of $\Theta \left( z\right) $ or $\Phi
\left( z\right) $ to be inside or outside the unit circle so that the
expansions of $\Psi \left( z\right) :=\Phi ^{-1}\left( z\right) \Theta
\left( z\right) $ and $\Psi ^{-1}\left( z\right) $ could include powers of $%
z $ and $z^{-1}$ simultaneously, accounting for noncausal or noninvertible
systems.

To investigate the identification problems on the location of the roots of
the matrix polynomials $\Phi \left( z\right) \ $and$\  \Theta \left( z\right) 
$ and on the components of $\mathbf{\varepsilon }_{t}$ determined by $\Theta
_{0}$, we use the device of Blaschke matrices (BM) that are generalized
orthogonal matrices. Following Lippi and Reichlin (1994), and denoting by $%
\ast $ simultaneous transposition and complex conjugation, a $d\times d$
matrix $A\left( z\right) $ is a BM if \newline

\textbf{1}. $A\left( z\right) $ has no poles of modulus smaller or equal to
unity and

\textbf{2}. $A\left( z\right) ^{-1}=A^{\ast }\left( z^{-1}\right) ,$ i.e. $%
A\left( z\right) A^{\ast }\left( z^{-1}\right) =\mathbf{I}_{d}.$ \newline

Further, for any BM, there exists an integer $r$ and complex numbers $a_{j},$
$j=1,\ldots ,r,$ $\left \vert a_{j}\right \vert <1,$ such that%
\begin{equation}
A\left( z\right) =K_{0}R\left( a_{1},z\right) K_{1}R\left( a_{2},z\right)
K_{2}\cdots K_{r-1}R\left( a_{r},z\right) K_{r},  \label{BM}
\end{equation}%
where $K_{j}$ are orthogonal matrices\footnote{%
Notice that in Theorem~1 of Lippi and Reichlin (1995) it is fixed that $%
K_{0}=\mathbf{I}_{d},$ but in general $K_{0}$ needs to be different from
identity to complete their proof as can be seen for the BM $A\left( z\right)
=diag(1,g_{a}\left( z\right) )$, }, $K_{j}K_{j}^{\prime }=\mathbf{I}_{d},$ and 
\begin{equation*}
R\left( a,z\right) =\left( 
\begin{array}{cc}
g_{a}\left( z\right) & 0 \\ 
0 & \mathbf{I}_{d-1}%
\end{array}%
\right) ,\  \  \ g_{a}\left( z\right) =\frac{z-a}{1-a^{\ast }z},
\end{equation*}%
see also the discussion in Hannan (1970, pp.~65-67).

For any BM $A\left( z\right) ,$ we can write\  \ 
\begin{equation}
Y_{t}=\Psi \left( L\right) A\left( L\right) \mathbf{u}_{t},\  \  \  \ 
\label{u}
\end{equation}%
where $\mathbf{u}_{t}=A\left( L\right) ^{-1}\mathbf{\varepsilon }_{t}\ $is a
serially uncorrelated all-pass process though not independent, and, because
its spectral density matrix is constant, $f_{\mathbf{u}}\left( \lambda
\right) =\left( 2\pi \right) ^{-1}A^{-1}\left( e^{-i\lambda }\right)
A^{-1\ast }\left( e^{i\lambda }\right) =\left( 2\pi \right) ^{-1}\mathbf{I}%
_{d},$ we conclude that the spectral density implied by the representation ($%
\ref{u}$) for any BM $A\left( z\right) $ and any $\mathbf{u}_{t}$ is always
the same,%
\begin{equation*}
f\left( \lambda \right) =\Psi \left( e^{-i\lambda }\right) A\left(
e^{-i\lambda }\right) f_{\mathbf{u}}\left( \lambda \right) A^{\ast }\left(
e^{i\lambda }\right) \Psi ^{\ast }\left( e^{i\lambda }\right) =\frac{1}{2\pi 
}\Psi \left( e^{-i\lambda }\right) \Psi ^{\ast }\left( e^{i\lambda }\right) .
\end{equation*}%
The same conclusion arises if $A$ is the inverse of a BM, and in particular
when in representation $\left( \ref{BM}\right) $ it holds that $%
1/a_{j}^{\ast }$ equals an actual root of det$\left( \Psi \left( z\right)
\right) ,$ irrespectively of being inside or outside the complex unit
circle, in a process of flipping the roots of det$\left( \Psi \left(
z\right) \right) $. 
These facts imply at once that using only second order information we can
not identify the location of these roots with respect to the unit circle,
and, even with knowledge of $p$ and $q,$ there are infinite VARMA
representations with the same second order properties but different IRF $%
\Psi \left( L\right) A\left( L\right) $ and error sequence $A(L)^{-1}\mathbf{%
\varepsilon }_{t}$. These alternative IRFs and errors are associated to
invertible and noninvertible representations when $q>0$ and to causal and
noncausal representations when $q=0$ as in this case $\left( \Psi \left(
L\right) A\left( L\right) \right) ^{-1}=A^{-1}\left( L\right) \Theta
_{0}^{-1}\Phi \left( L\right) $ and the (inverse of the) roots of $%
A^{-1}\left( L\right) $ can match those of $\Phi \left( L\right) $ in our
generalized setup. Traditional estimation methods based on Gaussian PML,
like Whittle approximation, only consider causal and invertible
representations, but still have to deal with the static problem that arises
for $A\left( L\right) $ constant.

The static identification problem refers to the well known lack of
identification of standard Structural VAR(MA) models with respect to
orthogonal rotations $\mathbf{u }_{t}=K \mathbf{\varepsilon }_{t}$ of the
structural errors in absence of further identifying assumptions on the IRF
provided by economic theory and/or further model structure (see e.g.
Rubio-Ramirez, Waggoner and Zha, 2010, for equality restrictions, and
Granziera, Moon and Schorfheide, 2018, for sign restrictions). However, it
is possible to consider this static problem within the same framework by
allowing BM which are constant and equal to an orthogonal matrix. Then, when
dynamics are known or sufficient conditions for their identification are
imposed (and e.g. causality and invertibility), we can identify \textit{%
statistically} the structural shocks in SVARMA models by higher order moment
conditions implied by ICA under non-Gaussianity without further restrictions.

To consider all these situations when trying to identify a SVARMA model we
extend the concept of Blaschke Matrix (BM) to any matrix $A\left( z\right) $
that satisfies the orthogonality condition 2. and\newline

\textbf{1}$^{\ast }$. $A\left( z\right) $ has no poles of modulus equal to
unity, but could have some with modulus larger or smaller than unity. 
\newline

Then, in the representation $\left( \ref{BM}\right) $ for a BM $A\left(
z\right) $ we allow for $\left \vert a_{j}\right \vert >1$ as well as $%
\left
\vert a_{j}\right \vert <1$, so that there exists an integer $%
r=0,1,\ldots $, complex numbers $a_{i},$ $j=1,\ldots ,r$ and a $\eta >0$
such that 
\begin{equation*}
\min_{j}\left \vert \left \vert a_{j}\right \vert -1\right \vert \geq \eta
>0,
\end{equation*}%
where the case $r=0$ is interpreted as $A\left( z\right) =K_{0}$ being a
constant (in $z$) orthogonal matrix. Since we do not restrict the $a_{j}$ so
that $1/a_{j}^{\ast }$ matches a root of $\Psi \left( z\right) ,$
considering any BM $A\left( z\right) $ we can deal with both basic and
non-basic representations of VARMA models in the sense of Lippi and Reichlin
(1994).

To solve the problem that second order dynamics cannot identify the phase of 
$\Psi $, we resort to higher order moments as proposed by Lii and Rosenblatt
(1992). In Appendix~A we develop a compact representation of the spectral
density $f_{\mathbf{a},k}$ of $\left( Y_{t,\mathbf{a}\left( 1\right)
},\ldots ,Y_{t,\mathbf{a}\left( k\right) }\right) $ for any order $%
k=2,3,\ldots ,$ and $k$-tuple $\mathbf{a}=\left( \mathbf{a}\left( 1\right)
,\ldots ,\mathbf{a}\left( k\right) \right) ,$ where $Y_{t}$ follows a linear
model with IRF $\Psi (L)$ and innovations $\mathbf{\varepsilon }_{t}$ which
are \emph{iid} up to moments of order $k.$ The $k$-th order cumulants of the
vector $\mathbf{\varepsilon }_{t}$ can be characterized by the $d^{2}\times
d^{k-2}$ matrix v$\mathbf{\kappa }_{k}^{0}$, 
\begin{equation*}
\text{v}\mathbf{\kappa }_{k}^{0}:=\left[ \text{vec}\left( \mathbf{\kappa }%
_{\cdot \cdot 1\cdots 1}\right) \  \  \text{vec}\left( \mathbf{\kappa }_{\cdot
\cdot 2\cdots 1}\right) \  \  \cdots \  \  \text{vec}\left( \mathbf{\kappa }%
_{\cdot \cdot d\cdots d}\right) \right] ,
\end{equation*}%
where $\mathbf{\kappa }_{\cdot \cdot j\left( 3\right) \cdots j\left(
k\right) }$ is the $d\times d$ matrix with typical $\left( j\left( 1\right)
,j\left( 2\right) \right) $ element equal to the $k$-th order joint cumulant
cum$\left( \mathbf{\varepsilon }_{t,j\left( 1\right) },\mathbf{\varepsilon }%
_{t,j\left( 2\right) },\mathbf{\varepsilon }_{t,j\left( 3\right) },\ldots ,%
\mathbf{\varepsilon }_{t,j\left( k\right) }\right) ,$ $j\left( h\right) \in
\left \{ 1,\ldots ,d\right \} .$ Then we find that%
\begin{equation*}
f_{\mathbf{a},k}(\boldsymbol{\lambda })=\frac{1}{\left( 2\pi \right) ^{k-1}}%
\Psi _{\mathbf{a}}^{\otimes k}\left( \mathbf{\lambda }\right) \text{vec}%
\left( \text{v}\mathbf{\kappa }_{k}^{0}\right) ,
\end{equation*}%
where for $\boldsymbol{\lambda }=(\lambda _{1},\ldots ,\lambda _{k-1})$ we
define%
\begin{equation*}
\Psi _{\mathbf{a}}^{\otimes k}\left( \mathbf{\lambda }\right) :=\Psi _{%
\mathbf{a}\left( k\right) }\left( e^{i\left( \lambda _{1}+\cdots +\lambda
_{k-1}\right) }\right) \otimes \Psi _{\mathbf{a}\left( k-1\right) }\left(
e^{-i\lambda _{k-1}}\right) \otimes \cdots \otimes \Psi _{\mathbf{a}\left(
2\right) }\left( e^{-i\lambda _{2}}\right) \otimes \Psi _{\mathbf{a}\left(
1\right) }\left( e^{-i\lambda _{1}}\right)
\end{equation*}%
for the usual Kronecker product $\otimes $ on the rows $\Psi _{\mathbf{a}%
\left( j\right) }$ of $\Psi $. This representation produces the usual
spectral density for $k=2$ because 
\begin{equation*}
f_{(\mathbf{a}(1),\mathbf{a}(2)),2}(\lambda )=\frac{1}{2\pi }\Psi _{\mathbf{a%
}(2)}(e^{i\lambda })\otimes \Psi _{\mathbf{a}(1)}(e^{-i\lambda })\text{vec}%
\left( \text{v}\mathbf{\kappa }_{2}^{0}\right) =\frac{1}{2\pi }\Psi _{%
\mathbf{a}(1)}(e^{-i\lambda })\Psi _{\mathbf{a}(2)}^{\prime }(e^{i\lambda }),
\end{equation*}%
as in this case$\ $v$\mathbf{\kappa }_{2}^{0}=\ $vec$(E[\mathbf{\varepsilon }%
_{t}\mathbf{\varepsilon }_{t}^{\prime }])=\ $vec$(\mathbf{I}_{d})$, where $%
\mathbf{I}_{d}$ is the covariance matrix of $\mathbf{\varepsilon }_{t}$
under the imposed normalization.

We now discuss the intuition on why higher order spectral densities with $%
k\geq 3$ can achieve dynamics identification unlike for $k=2.$ Thus, 
\begin{equation*}
f_{\mathbf{a},3}(\boldsymbol{\lambda };A,\text{v}\mathbf{\kappa }_{3})=\frac{%
1}{\left( 2\pi \right) ^{2}}\Psi _{\mathbf{a}}^{\otimes 3}\left( \mathbf{%
\lambda }\right) A^{\otimes 3}\left( \mathbf{\lambda }\right) \text{vec}%
\left( \text{v}\mathbf{\kappa }_{3}\right)
\end{equation*}%
is the implied $k=3$ spectral density for any third order marginal cumulants
matrix v$\mathbf{\kappa }_{3}$ under the (wrong) assumption that $\mathbf{u}%
_{t}=A\left( L\right) ^{-1}\mathbf{\varepsilon }_{t}$ is an $iid\left( 0,%
\mathbf{I}_{d}\right) $ sequence in (\ref{u}) for any non-constant BM $%
A\left( z\right) $, as the true $\mathbf{\varepsilon }_{t}$, and not just
serially uncorrelated.

Under some identification (rank) assumptions $f_{\mathbf{a},3}(\boldsymbol{%
\lambda };A,\text{v}\mathbf{\kappa }_{3})$ does differ from the true density 
$f_{\mathbf{a},3}(\boldsymbol{\lambda })=f_{\mathbf{a},k}(\boldsymbol{%
\lambda };\mathbf{I}_{d},$v$\mathbf{\kappa }_{3}^{0})$ for all choices of v$%
\mathbf{\kappa }_{3}$ because $A^{\otimes 3}\left( \mathbf{\lambda }\right)
\,$vec$\, \left( \text{v}\mathbf{\kappa }_{3}\right) $ depends on $\mathbf{%
\lambda} $ in general, unlike $A^{\otimes 2}\left( \lambda \right) $vec$%
\left( \text{v}\mathbf{\kappa }_{2}\right) =\ $vec$\left( A\left(
e^{-i\lambda _{1}}\right) \mathbf{I}_{d}A^{\prime }\left( e^{i\lambda
_{1}}\right) \right) =$ vec$\left( \mathbf{I}_{d}\right) .$ In particular,
for $d=1,$ $A^{\otimes 3}\left( \mathbf{\lambda }\right) $ is the bispectral
density of an all-pass process, which is not constant unlike its second
order spectral density. Similar arguments apply for any higher spectral
density, so, following VL, we can set up an $L^{2}$ distance between $f_{%
\mathbf{a},k}(\boldsymbol{\lambda };A,$v$\mathbf{\kappa }_{k})\ $and $f_{%
\mathbf{a},k}(\boldsymbol{\lambda }),$ for all $k$-tuples $\mathbf{a}$ from $%
\{1,2,\ldots ,d\}$, 
\begin{equation*}
\mathcal{L}_{k}^{0}\left( A,\text{v}\mathbf{\kappa }_{k}\right) :=\sum_{%
\mathbf{a}}\int_{\Pi ^{k-1}}\left \vert f_{\mathbf{a},k}(\boldsymbol{\lambda 
};A,\text{v}\mathbf{\kappa }_{k})-f_{\mathbf{a},k}(\boldsymbol{\lambda }%
)\right \vert ^{2}d\boldsymbol{\lambda },
\end{equation*}%
with known $\Psi ,$ but unknown location of the roots of $\Theta \left(
z\right) $ (or $\Phi \left( z\right) $ if $q=0$) expressed by the factor $A$%
, possibly flipping inside or outside some of these roots, adding additional
all-pass dynamics, or simply rotating elements of $\mathbf{\varepsilon }%
_{t}. $ Here $\Pi =\left[ -\pi ,\pi \right] $ and $\Pi ^{k-1}$ is the ($k-1)$%
-th cartesian product of $\Pi $\textbf{.} Obviously $\mathcal{L}_{k}\left( 
\mathbf{I}_{d},\text{v}\mathbf{\kappa }_{k}^{0}\right) =0,$ but we have to
rule out the possibility that for some spectral factor $A\neq \mathbf{I}_{d}$
it is possible to choose some v$\mathbf{\kappa }_{k}$ such that $\mathcal{L}%
_{k}\left( A,\text{v}\mathbf{\kappa }_{k}\right) =0$ for a given $k>2$, as
indeed it is possible for $k=2$ just setting v$\mathbf{\kappa }_{2}=\mathbf{I%
}_{d}.$

To illustrate this problem notice that we can write $\mathcal{L}%
_{k}^{0}\left( A,\text{v}\mathbf{\kappa }_{k}\right) $ as%
\begin{equation*}
\int_{\Pi ^{k-1}}\left \{ \text{vec}\! \left( \text{v}\mathbf{\kappa }%
_{k}\right) ^{\prime }A^{\otimes k}\! \left( \mathbf{\lambda }\right) ^{\ast
}-\text{vec}\! \left( \text{v}\mathbf{\kappa }_{k}^{0}\right) ^{\prime
}\right \} \Upsilon _{k}^{0}\left( \mathbf{I}_{d},\boldsymbol{\lambda }%
\right) \left \{ A^{\otimes k}\left( \mathbf{\lambda }\right) \text{vec}%
\left( \text{v}\mathbf{\kappa }_{k}\right) -\text{vec}\left( \text{v}\mathbf{%
\kappa }_{k}^{0}\right) \right \} d\boldsymbol{\lambda },
\end{equation*}%
where for any $d\times d$ spectral factor $A$ and $\boldsymbol{\lambda }\in
\Pi ^{k-1}$ we define for $k=2,3,\ldots $ 
\begin{equation*}
\Upsilon _{k}^{0}\left( A,\boldsymbol{\lambda }\right) :=A^{\otimes k}\left( 
\mathbf{\lambda }\right) ^{\ast }\sum_{\mathbf{a}}\Psi _{\mathbf{a}%
}^{\otimes k}\left( \mathbf{\lambda }\right) ^{\ast }\Psi _{\mathbf{a}%
}^{\otimes k}\left( \mathbf{\lambda }\right) A^{\otimes k}\left( \mathbf{%
\lambda }\right) .
\end{equation*}%
Then, for identification of $\Psi $ under the assumption that $\Upsilon
_{k}^{0}\left( \mathbf{I}_{d},\boldsymbol{\lambda }\right) $ is full rank
for every $\boldsymbol{\lambda }$, we have to rule out the possibility that
for some BM $A\neq \mathbf{I}_{d}$, possibly constant, and some choice of v$%
\mathbf{\kappa }_{k}$ 
it holds that $A^{\otimes k}\left( \mathbf{\lambda }\right) $vec$\left( 
\text{v}\mathbf{\kappa }_{k}\right) -$vec$\left( \text{v}\mathbf{\kappa }%
_{k}^{0}\right) =0$ for all $\mathbf{\lambda }$ (except possibly in a set of
measure zero), implying that $\mathcal{L}_{k}^{0}\left( A,\text{v}\mathbf{%
\kappa }_{k}\right) =0.$

We now introduce rank conditions on $\text{v}\mathbf{\kappa }_{k}^{0}$ to
reduce to a minimum the range of situations where identification is lost for 
$k=3$ and $4$. For dynamic identification it is sufficient to use a rank
condition on $\text{v}\mathbf{\kappa }_{3}^{0}$, but when using $k=4$
moments, we need to impose ICA with nonzero marginal kurtosis coefficients
for all components of $\mathbf{\varepsilon }_{t}$, though no further
conditions on higher order moments are required in contrast to Chan et al.
(2006) and Gouri\'{e}roux et al. (2019). 
For static identification we also impose ICA of order $k$ among the
components of $\mathbf{\varepsilon }_{t}$ but we allow for at most one of
the components to have zero marginal cumulants, being possibly Gaussian,
resembling the result of Comon (1994). We also relax the serial \emph{idd}
assumption on $\mathbf{\varepsilon }_{t} $ to equal distribution and
independence up to moments of order $k$, which is sufficient for
stationarity of $Y_{t}$ of order $k$ and to specify the corresponding higher
order spectral densities without further conditions on probability
distributions or conditional moments. Denote by $\lambda _{\min }\left(
M\right) $ the minimum eigenvalue of a matrix $M$ and let $\mathbf{\alpha }%
_{k}^{0}=\left( \mathbf{\alpha }_{k1}^{0},\ldots ,\mathbf{\alpha }%
_{kd}^{0}\right) ^{\prime }$ be the true vector of marginal cumulants of
order $k$ of $\mathbf{\varepsilon }_{t}.$ \bigskip

\noindent \textbf{Assumption~1}$\left( k\right) $: The $\mathbf{\varepsilon }%
_{t}$ are stationary and serially independent up to $k$ moments$,$ $%
E\left
\Vert \mathbf{\varepsilon }_{t}\right \Vert ^{k}<\infty ,$ and are
standardized with zero mean and $\mathbf{I}_{d}$ covariance matrix. \bigskip

\noindent \textbf{Assumption~2}$\left( 3\right) $. Rank$\left( \text{v}%
\mathbf{\kappa }_{3}^{0}\right) =d.$\bigskip

\noindent \textbf{Assumption~3}$\left( k\right) $. The components of $%
\mathbf{\varepsilon }_{t}\ $are independent up to $k$ moments and

\begin{itemize}
\item For $k=3$ all marginal skewness coefficients are nonzero, $\mathbf{%
\alpha }_{3j}^{0}=\mathbf{\kappa }_{jjj}^{0}\neq 0,$ $j=1,\ldots ,d.$

\item For $k=4$ all marginal kurtosis coefficients are nonzero, $\mathbf{%
\alpha }_{4j}^{0}=\mathbf{\kappa }_{jjjj}^{0}\neq 0,$ $j=1,\ldots ,d.$%
\bigskip
\end{itemize}

\noindent \textbf{Assumption~4}: For some $\eta >0,$%
\begin{equation*}
\inf_{\left \vert z\right \vert =1}\lambda _{\min }\left( \Psi \left(
z\right) \right) \geq \eta >0.
\end{equation*}
\smallskip

Assumption~4 determines the full rank of the dynamic system excluding unit
roots on the AR and MA lag polynomials, so $\Upsilon _{k}^{0}\left( \mathbf{I%
}_{d},\mathbf{\lambda }\right) =\left( \Psi ^{\otimes k}\left( \mathbf{%
\lambda }\right) \right) ^{\ast }\Psi ^{\otimes k}\left( \mathbf{\lambda }%
\right) $ is positive definite as well as $\Upsilon _{k}^{0}\left( A,%
\boldsymbol{\lambda }\right) \ $for all $\mathbf{\lambda }$ and any BM $A,$
because $\Psi \left( z\right) $ has rank $d$ on the complex unit circle. 
However, note that despite that under Assumption~4 $\Upsilon _{k}^{0}\left(
A,\boldsymbol{\lambda }\right) >0$ for all $\mathbf{\lambda }$ and any BM $%
A, $ it is not sufficient for dynamics identification and we need to make
sure that v$\mathbf{\kappa }_{k}^{0}$ is rich enough through Assumptions~2$%
\left( 3\right) $ or~3$\left( k\right) $. 
In particular, Assumption~2$\left( 3\right) $ is equivalent to the linear
independence of the third order cumulant matrices $\left \{ \mathbf{\kappa }%
_{\cdot \cdot j}^{0}\right \} _{j=1}^{d}$ assumption of Chen, Choi and
Escanciano (2018) used to investigate the consistency of a fundamentalness
test by showing that the Wold (invertible) innovations of a nonfundamental
VARMA model cannot be a martingale difference sequence despite being white
noise. 

The ICA of order $k$ among the elements of $\mathbf{\  \varepsilon }_{t}$
contained in Assumption~3$(k)$ implies that all joint higher order cumulants
up to order $k$ are zero, i.e. $\mathbf{\kappa }_{abc}^{0}=0$ when $a,b,c$
are not all equal, so Assumption~3$\left( 3\right) $ implies Assumption~2$%
\left( 3\right) $ but imposes further structure on the multivariate skewness
of the vector $\mathbf{\varepsilon }_{t}$ given by 
\begin{equation}
\text{v}\mathbf{\kappa }_{3}^{0}=\text{v}\mathbf{\kappa }_{3}^{\text{IC}%
}\left( \mathbf{\alpha }_{3}^{0}\right) :=\left( \mathbf{\alpha } _{31}^{0} 
\mathbf{e}_{1}^{\otimes 2} ,\  \mathbf{\alpha } _{32}^{0} \mathbf{e}%
_{2}^{\otimes 2},\  \ldots ,\  \mathbf{\alpha } _{3d}^{0} \mathbf{e}%
_{d}^{\otimes 2} \right) =\sum_{j=1}^{d} \mathbf{\alpha} _{3j}^{0}\mathbf{e}%
_{j}^{\otimes 2}\mathbf{e}_{j}^{\prime } ,  \label{vk3}
\end{equation}%
where $\mathbf{\alpha }_{3}^{0}=\left( \mathbf{\alpha } _{31}^{0},\mathbf{%
\alpha } _{32}^{0},\ldots ,\mathbf{\alpha } _{3d}^{0}\right) ^{\prime }$ are
the marginal skewness coefficients of $\mathbf{\varepsilon }_{t}$ and $%
\mathbf{e}_{j}$ is the $j$-th column of $\mathbf{I}_{d}$. Therefore, under
Assumption~3(3), all $\mathbf{\kappa }_{\cdot \cdot j}^{0}$, $j=1,\ldots ,d,$
are $d\times d$ matrices of zeros with a unique nonzero element $\mathbf{%
\alpha }_{3j}^{0}=\mathbf{\kappa } _{jjj}^{0}$ in position $\left(
j,j\right) $. Note that orthogonal rotations $\mathbf{\eta }_{t}=K\mathbf{%
\varepsilon }_{t}$ have the same identity covariance matrix of $\mathbf{%
\varepsilon }_{t},$ but their components are not longer independent if $%
K\neq P_{d},$ a signed permutation matrix of dimension $d$, because e.g. v$%
\mathbf{\kappa }_{3}^{\eta }=K^{\otimes 2}\ $v$\mathbf{\kappa }_{3}^{0}\
K^{\prime }$ has not the same structure $\left( \ref{vk3}\right) $ of v$%
\mathbf{\kappa }_{3}^{0},$ despite it maintains its rank, see Appendix~A for
details.

By contrast, Assumption~3$\left( 4\right) $ implies that v$\mathbf{\kappa }%
_{4}^{0}$ satisfies 
\begin{equation}
\text{v}\mathbf{\kappa }_{4}^{0}=\text{v}\mathbf{\kappa }_{4}^{\text{IC}%
}\left( \mathbf{\alpha }_{4}^{0}\right) :=\sum_{j=1}^{d}\mathbf{\alpha }%
_{4j}^{0}\mathbf{e}_{j}^{\otimes 2}\mathbf{e}_{j}^{\otimes 2\prime }
\label{vk4}
\end{equation}%
for the kurtosis coefficients $\mathbf{\alpha }_4^0$ and has only rank $d$,
because, despite $\mathbf{\kappa }_{\cdot \cdot jj}^{0}$ are all matrices of
zeros with a unique nonzero element in position $\left( j,j\right) $ given
by the marginal kurtosis coefficient $\mathbf{\alpha }_{4j}^{0}=\mathbf{%
\kappa }_{jjjj}^{0}$, $j=1,\ldots ,d,$ we have that $\mathbf{\kappa }_{\cdot
\cdot hj}^{0}=\mathbf{0},$ $h\neq j,$ so v$\mathbf{\kappa }_{4}^{0}$ has at
most $d$ columns different from zero. Note that Assumption~A.6 in Gouri\'{e}%
roux et al. (2019) similarly needs that each component of $\varepsilon _{t}$
has a nonzero cumulant of order $k,$ $k\geq 3,$ but it further requires that
components of $\mathbf{\varepsilon }_{t}$ are full independent with a finite
moment of order $s,$ where $s$ is an even integer greater than $k.$

Finally, it is possible to obtain static identifying results when there is
at most one Gaussian structural shock as in Comon (1994) result, or a
non-Gaussian one with zero higher order cumulant of order $k$, as formalized
in the next weaker version of Assumption~$3\left( k\right) .$\bigskip

\noindent \textbf{Assumption~3}$^{\ast }\left( k\right) $. The components of 
$\mathbf{\varepsilon }_{t}\ $are independent up to $k$ moments and

\begin{itemize}
\item For $k=3$ the marginal skewness coefficients are nonzero, $\mathbf{%
\alpha }_{3j}^{0}=\mathbf{\kappa }_{jjj}^{0}\neq 0,$ for all $j=1,\ldots ,d$
but at most one index.

\item For $k=4$ the marginal kurtosis coefficients are nonzero, $\mathbf{%
\alpha }_{4j}^{0}=\mathbf{\kappa }_{jjjj}^{0}\neq 0,$ for all $j=1,\ldots ,d$
but at most one index.\bigskip
\end{itemize}

\section{3. Nonparametric Identification}

In this section we discuss general identification results for SVARMA models
based on the spectral loss functions $\mathcal{L}_{k}^{0}$ under the
assumption of known dynamics up to a Blaschke factor $A$ and the
corresponding cumulants of structural shocks. We consider first dynamic
identification using nonconstant BM, while we later move to the static
components identification using constant BM.\bigskip

\begin{theorem}
\label{1A} Under Assumptions 1$\left( k\right) ,$ 2$\left( 3\right) $ and 4,
for any nonconstant BM $A\left( z\right) $, there exists an $\epsilon >0$
such that, 
\begin{equation*}
\inf_{\text{v}\mathbf{\kappa }_{3}}\mathcal{L}_{3}^{0}\left( A,\text{v}%
\mathbf{\kappa }_{3}\right) \geq \epsilon >0.
\end{equation*}
\end{theorem}

\bigskip

All proofs are contained in Appendix~B, while auxiliary results are included
in Appendices~C and D. As when $d=1$ in VL, Theorem~\ref{1A} implies that
there is no way of choosing v$\mathbf{\kappa }_{3}$ such that third order
dynamics can be replicated after introducing a Blaschke factor inverting any
root of $\Psi $ (as can be done for $k=2).$ The conditions of Theorem~\ref%
{1A} allow for general nonconstant $A\left( z\right) $ which have an
infinite expansion in positive and/or negative powers of $z$ under a rank
condition on v$\mathbf{\kappa }_{3}^{0}.$ \bigskip

\begin{theorem}
\label{1ANEW} Under Assumptions 1$\left( k\right) ,$ 3$\left( k\right) $, $%
k=3$ or$\ 4$, and 4, for any nonconstant BM $A\left( z\right) $, there
exists an $\epsilon >0$ such that, 
\begin{equation*}
\inf_{\mathbf{\alpha }}\mathcal{L}_{k}^{0}\left( A,\text{v}\mathbf{\kappa }%
_{k}^{\text{IC}}\left( \mathbf{\alpha }\right) \right) \geq \epsilon >0.
\end{equation*}
\end{theorem}

\bigskip

Theorem~\ref{1ANEW} relies on the particular structure of v$\mathbf{\kappa }%
_{k}^{0}$ imposed by Assumption~$3\left( k\right) .$ In fact, for $k=3,$
this is just a particular case of Theorem~\ref{1A}, since rank$\left( \text{v%
}\mathbf{\kappa }_{k}^{0}\right) =d$ under Assumption~3$\left( k\right) ,$ $%
k=3,4$. However, the argument of Theorem~\ref{1ANEW} can not be extended
under a generic rank condition on v$\mathbf{\kappa }_{4}^{0}$ to cover $%
\mathcal{L}_{4}^{0}$ in Theorem~\ref{1A} without further structure, because
for any v$\mathbf{\kappa }_{4}^{0}$ it holds rank$\left( \text{v}\mathbf{%
\kappa }_{4}^{0}\right) \leq d\left( d+1\right) /2<d^{2}=\ $rank$\left(
A\right) ^{2}$ for $d>1.$

For a signed permutation matrix $P_{d}$ of dimension $d$ with all elements
equal to zero but a single term equal to $+1$ or $-1$ in each column and
row, let $P_{d}^{+}$ be equal to $P_{d}$ but taking the absolute value of
all its elements. Then $P_{d}\mathbf{\alpha }$ and $P_{d}^{+}\mathbf{\alpha }
$ are (signed) permutations of the elements of the vector $\mathbf{\alpha }.$
\bigskip

\begin{theorem}
\label{1B} Under Assumptions 1$\left( k\right) ,$ 3$\left( k\right) $, $k=3\ 
$or $4,$ and 4, for any constant BM $K$ different from a signed permutation,
i.e. $K \neq P_{d}$, there exists an $\epsilon >0$ such that%
\begin{equation*}
\inf_{\mathbf{\alpha }}\mathcal{L}_{k}^{0}\left( K,\text{v}\mathbf{\kappa }%
_{k}^{\text{IC}}\left( \mathbf{\alpha }\right) \right) \geq \epsilon >0.
\end{equation*}
\end{theorem}

\bigskip

Under marginal independence of order $k$, which implies the co-kurtosis
conditions of Lanne and Luoto (2019), Theorem~\ref{1B} shows that for any
orthogonal matrix $K$ different from $\mathbf{I}_{d}$ and any signed
permutation matrix $P_{d}$, it is not possible to find any $\mathbf{\alpha }$
so that $\mathcal{L}_{k}^{0}\left( K,\text{v}\mathbf{\kappa }_{k}^{\text{IC}%
}\left( \mathbf{\alpha }\right) \right) =0.$ This provides identification of
the components of $\mathbf{\varepsilon }_{t}$ up to signed permutations
because for any $P_{d}$ and any v$\mathbf{\kappa }_{k}^{0}=\ $v$\mathbf{%
\kappa }_{k}^{IC}\left( \mathbf{\alpha }_{k}^{0}\right) $ we could select $%
\mathbf{\alpha }_{3}=P_{d}^{\prime }\mathbf{\alpha }_{3}^{0}$ or $\mathbf{%
\alpha }_{4}=P_{d}^{+\prime }\mathbf{\alpha }_{4}^{0}$ to make $\mathcal{L}%
_{3}^{0}\left( P_{d},\text{v}\mathbf{\kappa }_{3}^{\text{IC}}\left(
P_{d}^{\prime }\mathbf{\alpha }_{3}^{0}\right) \right) =\mathcal{L}%
_{4}^{0}\left( P_{d},\text{v}\mathbf{\kappa }_{4}^{\text{IC}}\left(
P_{d}^{+\prime }\mathbf{\alpha }_{4}^{0}\right) \right) =0.$

The results of Theorems~\ref{1ANEW} and~\ref{1B} can be combined in the
following result that identifies SVARMA models under the assumption of $k$%
-order independence of innovation components. Denote by $\left \Vert
M\right
\Vert =trace\left( M^{\ast }M\right) ^{1/2}$ the Frobenious norm of
a matrix $M$ and by $\left \Vert M\right \Vert _{L^{2}}=\int_{\left \vert
z\right
\vert =1}\left \Vert M\left( z\right) \right \Vert ^{2}dz$ the $%
L^{2} $ norm of $\Vert M\left( z\right) \Vert $ over the unit circle.

\begin{theorem}
\label{1A1B} Under Assumptions 1$\left( k\right) ,$ 3$\left( k\right) ,$ $%
k=3\ $or $4,$ and 4, for all $\nu >0$ there exists an $\epsilon >0$ such that%
\begin{equation*}
\inf_{A,\mathbf{\alpha },P_{d}\mathbf{:}\left \Vert A-P_{d}\right \Vert
_{L^{2}}+\left \Vert \mathbf{\alpha }-P_{d}^{\prime }\mathbf{\alpha }%
_{k}^{0}\right \Vert \geq \nu >0}\mathcal{L}_{k}^{0}\left( A,\text{v}\mathbf{%
\kappa }_{k}^{\text{IC}}\left( \mathbf{\alpha }\right) \right) \geq \epsilon
>0.
\end{equation*}%
where $A$ is any BM and$\ P_{d}$ is any signed permutation matrix, and $%
P_{d}^{\prime }\mathbf{\alpha }_{k}^{0}$ has to be replaced by $%
P_{d}^{+\prime }\mathbf{\alpha }_{k}^{0}$ for $k=4.$\bigskip
\end{theorem}

Note that under Assumption~3$(k)$ we can deal simultaneously with both
dynamics and marginal identification for both $k=3$ and$~4$, providing
identification up to signed permutations of the specific components of $%
\mathbf{\varepsilon }_{t}$ under mutual independence of order $k$.

We now explore the possibility of relaxing Assumption$~3 $ by allowing for
some marginal cumulants to be zero in the static identification of Theorem~%
\ref{1B} and also for the dynamics identification of Theorems~\ref{1A},~\ref%
{1ANEW} or~\ref{1A1B}, but only under conditions which guarantee that every
single component of $\mathbf{\varepsilon }_{t}$ is non-Gaussian. Then, for
static identification robustness we explore, first, the situation when just
one of the marginal cumulants is zero for a given $k$ and, second, when some
further marginal skewness coefficients are zero, but the corresponding
kurtosis coefficients are not, or vice versa, so they mutually compensate
for the lack of identification of some particular component, both mechanisms
indicating that at most one Gaussian component in $\mathbf{\varepsilon }_{t}$
can be allowed as in Comon (1994).

\bigskip

\begin{corollary}
\label{L0} 
Under Assumptions 1$\left( k\right) ,$ 3*$\left( k\right) $, $k=3\ $or $4,$
and 4, the conclusions of Theorem~\ref{1B} hold. 
\bigskip
\end{corollary}

Then, when a single higher order cumulant of a given order $k$ is zero only
signed permutation matrices are not discarded, but this is not true if more
than one marginal element of $\mathbf{\alpha }_{k}^{0}$ is zero. The next
result investigates the case when possibly more than one marginal cumulant
of the same order $k=3$ or $k=4$ is zero, but the corresponding marginal
cumulants of the other order are nonzero, allowing for at most one component
to have simultaneously zero skewness and zero kurtosis for static
identification, while the others need to have at least one coefficient
different from zero. As before, for dynamics identification we need to
guarantee that all components are non-Gaussian up to order four. To avoid
the potential lack of identification provided by a single set of cumulants
of a given order, we need to consider a robustified loss function involving
both third and fourth moments simultaneously.\bigskip

\begin{corollary}
\label{L1} Under Assumptions 1$\left( k\right) ,$ 3$\left( k\right) $, $k=3\ 
$and $4,$ and 4, with index sets $\mathcal{I}_{3}$ and $\mathcal{I}_{4}$,
subsets of $\left \{ 1,2,\ldots ,d\right \} ,$ such that $\mathbf{\alpha }%
_{3j}^{0}=0$ for $j\in \mathcal{I}_{3}$, $\mathbf{\alpha }_{4j}^{0}=0$ for $%
j\in \mathcal{I}_{4},$ and \#$\left \{ \mathcal{I}_{3}\cap \mathcal{I}%
_{4}\right \} \leq 1$, for any constant BM $A$ different from a signed
permutation, i.e. $A\neq P_{d}$, there exists an $\epsilon >0$ such that 
\begin{equation}
\inf_{\mathbf{\alpha }_{3}}\mathcal{L}_{3}^{0}\left( A,\text{v}\mathbf{%
\kappa }_{3}^{\text{IC}}\left( \mathbf{\alpha }_{3}\right) \right) +\inf_{%
\mathbf{\alpha }_{4}}\mathcal{L}_{4}^{0}\left( A,\text{v}\mathbf{\kappa }%
_{4}^{\text{IC}}\left( \mathbf{\alpha }_{4}\right) \right) \geq \epsilon >0,
\label{L34}
\end{equation}%
while if \#$\left \{ \mathcal{I}_{3}\cap \mathcal{I}_{4}\right \} =0$, then (%
\ref{L34}) holds for all non constant BM $A.$
\end{corollary}

\bigskip

\section{4. Parameter Identification}

We assume that the observed $d$-dimensional SVARMA$\left( p,q\right) $
process $Y_t$ admits the following parameterization 
\begin{equation}
\Phi _{\mathbf{\theta }_{0}}\left( L\right) Y_{t}=\mu +\Theta _{\mathbf{%
\theta }_{0}}\left( L\right) \mathbf{\varepsilon }_{t},\  \  \  \mathbf{%
\varepsilon }_{t}\sim iid_{k}\left( \mathbf{0},\mathbf{I}_{d},\text{v}%
\mathbf{\kappa }_{k}^{\text{IC}}\left( \mathbf{\alpha }_{k}^{0}\right) ,k\in 
\mathcal{K}\right)  \label{Rep}
\end{equation}%
where the index set $\mathcal{K}\subseteq \left \{ 3,4\right \} $ is non
empty, the lag polynomials 
\begin{eqnarray*}
\Phi _{\mathbf{\theta }}\left( L\right) &:=&\mathbf{I}_{d}+\Phi _{1}\left( 
\mathbf{\theta }\right) L+\cdots +\Phi _{p}\left( \mathbf{\theta }\right)
L^{p} \\
\Theta _{\mathbf{\theta }}\left( L\right) &:=&\Theta _{0}\left( \mathbf{%
\theta }\right) +\Theta _{1}\left( \mathbf{\theta }\right) L+\cdots +\Theta
_{q}\left( \mathbf{\theta }\right) L^{q}
\end{eqnarray*}%
depend on a $m$-dimensional parameter $\mathbf{\theta }\in \mathcal{S}%
\subset \mathbb{R}^{m}$ and $iid_{k}$ means that Assumption~1$\left(
k\right) $ holds imposing serial $iid$-ness up to $k$ moments. The
parameterization v$\mathbf{\kappa }_{k}^{\text{IC}}\left( \mathbf{\alpha }%
_{k}\right) $ given in (\ref{vk3}) and (\ref{vk4}) impose the independence
component condition of Assumption~3$\left( k\right) $ on the arrays of $k=3$
and/or $k=4$ order cumulants of the standardized error sequence $\mathbf{%
\varepsilon }_{t}$ with vectors $\mathbf{\alpha }_{k}\  \mathbf{\in \ }%
\mathcal{D}_{k}\subset \mathbb{R}^{d}$ of marginal skewness $\left(
k=3\right) $ and kurtosis coefficients $\left( k=4\right) $. Then $\mathbf{%
\theta }_{0}$ and $\mathbf{\alpha }_{0}$ denote the true value of the
parameters and, while the level $\mu $ could be estimated by OLS or GLS
based on estimates of $\mathbf{\theta }$ as usual, since our methods are
invariant to $\mu $ we do not discuss this further.

The $k$-th order spectral density parametric model for each index $\mathbf{a}%
=\left( \mathbf{a}(1),\ldots ,\mathbf{a}(k)\right) $ of components of $Y_{t}$
with representation $\left( \ref{Rep}\right) $ is given for $k=3,4,$ by 
\begin{equation*}
f_{\mathbf{a},k}(\boldsymbol{\lambda };\mathbf{\theta },\mathbf{\alpha }) :
=\left( \Phi _{\mathbf{\theta }}^{-1}\Theta _{\mathbf{\theta }}\right) _{%
\mathbf{a}}^{\otimes k}\left( \mathbf{\lambda }\right) \text{vec}\left( 
\text{v}\mathbf{\kappa }_{k}^{\text{IC}}\left( \mathbf{\alpha }\right)
\right) =\mathbf{\Psi }\left( \lambda ;\mathbf{\theta }\right) _{\mathbf{a}%
}^{\otimes k}\mathbf{S}_{k}\mathbf{\alpha }
\end{equation*}%
where $\mathbf{S}_{k}:=\left( \mathbf{e}_{1}^{\otimes k},\mathbf{e}%
_{2}^{\otimes k},\ldots ,\mathbf{e}_{d}^{\otimes k}\right) $ is a rank $d$
selection matrix and $\mathbf{\Psi }\left( \lambda ;\mathbf{\theta }\right)
:=\Phi _{\mathbf{\theta }}^{-1}\left( e^{-i\lambda }\right) \Theta _{\mathbf{%
\theta }}\left( e^{-i\lambda }\right) $. For $k=2$ we replace $\text{v}%
\mathbf{\kappa }_{2}^{\text{IC}}\left( \mathbf{\alpha }\right) $ by $\mathbf{%
I}_{d}$ in $f_{\mathbf{a},2}(\lambda ;\mathbf{\theta })$, to impose
normalization and uncorrelation of the components of $\mathbf{\varepsilon }%
_{t}$.

We assume that the parameterization $\left( \ref{Rep}\right) $ satisfies the
following conditions. \bigskip

\noindent \textbf{Assumption~5} \label{Ass4}

\noindent \textbf{5.1.} For all $\mathbf{\theta }\in \mathcal{S},$ det$%
\left( \Phi _{\mathbf{\theta }}\left( z\right) \right) $det$\left( \Theta _{%
\mathbf{\theta }}\left( z\right) \right) \neq 0$ for $\left \vert
z\right
\vert =1$ and $\Theta _{0}\left( \mathbf{\theta }\right) $ is
nonsingular.

\noindent \textbf{5.2.} For all $\mathbf{\theta }\neq \mathbf{\theta }_{0},$ 
$\  \Phi _{\mathbf{\theta }}^{-1}\left( z\right) \Theta _{\mathbf{\theta }%
}\left( z\right) \neq \Phi _{\mathbf{\theta }_{0}}^{-1}\left( z\right)
\Theta _{\mathbf{\theta }_{0}}\left( z\right) $ in a subset of positive
measure of $\left \{ z\in \mathbb{C}:\left \vert z\right \vert =1\right \} .$

\noindent \textbf{5.3.} $\mathbf{\theta }_{0}\in \mathcal{S}\ $and $\mathcal{%
S}$ is compact.

\noindent \textbf{5.4.} $\Phi _{i}\left( \mathbf{\theta }\right) ,$ $%
i=0,\ldots ,p,\ $and $\Theta _{i}\left( \mathbf{\theta }\right) ,$ $%
i=0,\ldots ,q,$ are continuously differentiable for $\mathbf{\theta }\in 
\mathcal{S}$.\bigskip

Assumption$~5.1$ imposes Assumption~4 for each parameterized model, and with
the $iid$ condition of order $k$ on the sequence $\mathbf{\varepsilon }_{t}$
in Assumption~1$\left( k\right) $ guarantee that $Y_{t}$ with representation 
$\left( \ref{Rep}\right) $ is $k$-stationary. The identifiability conditions
in Assumption~5.2 are satisfied when the parameter space $\mathcal{S}$ is
sufficiently constrained (cf. Boubacar Mainassara and Francq, 2011) as for
restricted versions of causal and invertible VARMA models (e.g. echelon or
final equations forms) that guarantee that $\Phi _{\mathbf{\theta }}\left(
L\right) \ $and $\Theta _{\mathbf{\theta }}\left( L\right) $ are left
coprime and that the unique unimodular common left divisor of $\Phi _{%
\mathbf{\theta }}\left( L\right) $ and $\Theta _{\mathbf{\theta }}\left(
L\right) $ is the identity matrix, see e.g. Section~12.1 in L\"{u}tkepohl
(2005). Assumption~5.3 is a standard parameter space restriction and,
together with Assumption~5.4, allows for uniformity arguments.

Note that parameterizations covering both invertible and noninvertible (or
causal and noncausal) solutions are allowed by Assumption~5.2, which
identifies uniquely the parametric transfer and impulse response functions.
However, it is possible that $\Phi _{\mathbf{\theta }}^{-1}\left( z\right)
\Theta _{\mathbf{\theta }}\left( z\right) A(z)=\Phi _{\mathbf{\theta }%
_{0}}^{-1}\left( z\right) \Theta _{\mathbf{\theta }_{0}}\left( z\right) $
for almost all $z,\  \left \vert z\right \vert =1,\ $some $\mathbf{\theta }%
\neq \mathbf{\theta }_{0}$ and some BM $A\left( z\right) \neq \mathbf{I}%
_{d}, $ and therefore $f_{\mathbf{a},2}\left( \lambda ;\mathbf{\theta }%
_{0}\right) $ can not identify $\Phi _{\mathbf{\theta }_{0}}^{-1}\left(
z\right) \Theta _{\mathbf{\theta }_{0}}\left( z\right) .$ However,
Assumption~5.2 is not sufficient either to identify $\Phi _{\mathbf{\theta }%
}^{-1}\left( z\right) \Theta _{\mathbf{\theta }}\left( z\right) $ uniquely
from $f_{\mathbf{a},k}(\boldsymbol{\lambda };\mathbf{\theta }_{0},\mathbf{%
\alpha }_{k}^{0})$ for $k=3$ or 4, without further restrictions to discard
signed permutations $A\left( z\right) =P_{d}\neq \mathbf{I}_{d}$ since
Assumption~3$\left( k\right) $ already prevents for all nonconstant $A\left(
z\right) .$ With this aim, we introduce the following assumptions, which fix
the signs of the components of $\mathbf{\varepsilon }_{t}$ and impose an
ordering by either imposing a given structure on $\Theta _{0}\left( \mathbf{%
\theta }_{0}\right) $ (6A), or by restricting the allowed set of values for
the marginal third or fourth order cumulants $\mathbf{\alpha }_{k}^{0}$ (6B$%
(k)$), or by directly excluding any signed permutations among the columns of
the transfer function (6C). \bigskip

\noindent \textbf{Assumption~6A}. The diagonal elements of $\Theta
_{0}\left( \mathbf{\theta }_{0}\right) $ are all positive and the elements
of $Y_{t}$ are ordered so that there is no signed permutation $P_{d}$ such
that the absolute value of the product of the diagonal elements of $\Theta
_{0}\left( \mathbf{\theta }_{0}\right) P_{d}$ is equal or larger than that
of $\Theta _{0}\left( \mathbf{\theta }_{0}\right) .$\bigskip

The restriction on the diagonal elements of $\Theta _{0}\left( \mathbf{%
\theta }_{0}\right) ,$ cf. Pham and Garat (1997) and Lanne and Luoto (2019),
fixes the signs of $\mathbf{\varepsilon }_{t}$ so that a positive increment
in a component of $\mathbf{\varepsilon }_{t}$ corresponds to a positive
increment in the element of $Y_{t}$ associated to this shock. This
restriction can be described as a sign restriction (at lag $0)$ on the IRF
of each endogenous variable with respect to the corresponding error term,
giving a unique interpretation of the IRF (as it is automatically imposed
when using Cholesky identification). Further, it is imposed a unique
permutation that maximizes the absolute value of the product of the diagonal
elements of $\Theta _{0}\left( \mathbf{\theta }_{0}\right) P_{d}$ for $P_{d}=%
\mathbf{I}_d.$ Alternative ordering schemes based on the elements of $\Theta
_{0}$ are possible, see e.g. Lanne et al. (2017), as well as schemes that
use information for the ordering from the vectors $\mathbf{\alpha }_{3}^{0}$
and $\mathbf{\alpha }_{4}^{0},$ as the following novel condition.\bigskip

\noindent \textbf{Assumption~6B}$\left( k\right) $. The diagonal elements of 
$\Theta _{0}\left( \mathbf{\theta }_{0}\right) $ are all positive and it
holds for $k\in \mathcal{K}$,

\begin{itemize}
\item[\textbf{6B}$(3)$] \noindent For $k=3:$ $\mathbf{\alpha }_{3}^{0}\in 
\mathcal{D}_{3}$ compact where 
\begin{equation*}
\mathcal{D}_{3}\subseteq \left \{ \mathbf{\alpha }=\left( \alpha _{1},\ldots
,\alpha _{d}\right) ^{\prime }\in \mathbb{R}^{d}:-\infty <\alpha _{1}<\alpha
_{2}<\cdots <\alpha _{d}<\infty , \  \alpha _{j}\neq 0\right \} ,
\end{equation*}%
and there is no signed permutation $P_{d}\neq \mathbf{I}_{d}$ such that the
diagonal elements of $\Theta _{0}\left( \mathbf{\theta }_{0}\right) P_{d}$
are positive and $P_{d}^{\prime }\mathbf{\alpha }_{3}^{0}\in \mathcal{D}%
_{3}. $

\item[\textbf{6B}$(4)$] \noindent For $k=4:$ $\mathbf{\alpha }_{4}^{0}\in 
\mathcal{D}_{4}$ compact where 
\begin{equation*}
\mathcal{D}_{4}\subseteq \left \{ \mathbf{\alpha }=\left( \alpha _{1},\ldots
,\alpha _{d}\right) ^{\prime }\in \mathbb{R}^{d}:-2\leq \alpha _{1}<\alpha
_{2}<\cdots <\alpha _{d}<\infty ,\  \alpha _{j}\neq 0\right \} .
\end{equation*}
\end{itemize}

The sign of the diagonal elements of $\Theta _{0}\left( \mathbf{\theta }%
_{0}\right) $ also determines the sign of the elements of $\mathbf{\alpha }%
_{3}^{0}\in \mathcal{D}_{3},$ but not that of $\mathbf{\alpha }_{4}^{0}\in 
\mathcal{D}_{4},$ so any signed permutation $P_{d},$ $P_{d}^{+}\neq \mathbf{I%
}_{d},$ would alter the ordering of the elements of $\mathbf{\alpha }%
_{4}^{0} $, but it might not that of $\mathbf{\alpha }_{3}^{0}$ depending on
the values of $\Theta _{0}\left( \mathbf{\theta }_{0}\right) $ and the shift
of signs. Then, Assumption~6B($3$) imposes explicitly a unique sign-ordering
associated with an increasing sequence of asymmetry coefficients to remove
these situations. While Assumption~6B imposes the non-zero cumulant part of
Assumption~3$\left( k\right) $, it further requires some further knowledge
on the marginal distributions of the structural errors, e.g. related to the
asymmetry sign or the tail behaviour, that allows to order the skewness or
kurtosis coefficients of the components of $\mathbf{\varepsilon }_{t}$.

Assumptions~6A and 6B fix the ordering of the columns of $\Theta _{0}\left( 
\mathbf{\theta }_{0}\right) $ or of the cumulant coefficients, respectively,
to exclude permutations in the columns of the IRF (cf. Theorem~3) if $%
\mathcal{S}$ and the parameterization $\Phi _{\mathbf{\theta }}^{-1}\left(
z\right) \Theta _{\mathbf{\theta }}\left( z\right) $ allow for rotations
(i.e. constant Blaschke factors) different from $\mathbf{I}_{d}$, where
rotations different from $P_{d}$ are ruled out by $\mathcal{L}_{k}^{0},$ $%
k=3,4.$ If $\mathcal{S}$ restricts to only invertible and causal models,
then it would be possible to allow for a component of $\mathbf{\varepsilon }%
_{t}$ to be symmetric (and possibly Gaussian), i.e. $\mathbf{\alpha }%
_{3j}^{0}=0$ for a single $j\in \left \{ 1,\ldots ,d\right \} ,$ if second
order moments (to identify dynamics) were used in conjunction with third
order ones (to identify rotations), cf. Corollary~1. If in this case we
include information from $k=4$ simultaneously to $k=2$ and $k=3,$ we could
also allow for a component with zero skewness and kurtosis (and possibly
Gaussian), while the other ones should display at least one non-Gaussian
feature in their third and fourth order moments, i.e. $\left \vert \mathbf{%
\alpha }_{3j}^{0}\right \vert +\left \vert \mathbf{\alpha }%
_{4j}^{0}\right
\vert =0$ for at most one $j\in \left \{ 1,\ldots
,d\right
\} .$

In the next assumption, the parameterization and $\mathcal{S}$ are further
restricted not to allow that $\Phi _{\mathbf{\theta }}^{-1}\left( z\right)
\Theta _{\mathbf{\theta }}\left( z\right) =\Phi _{\mathbf{\theta }%
_{0}}^{-1}\left( z\right) \Theta _{\mathbf{\theta }_{0}}\left( z\right)
P_{d} $ for any $\mathbf{\theta }\neq \mathbf{\theta }_{0}$ and any constant
signed permutation matrix $P_{d}$ because some identifying restrictions have
been already imposed, such as a recursive system (Choleski) assumption with
a known ordering, so that Assumption~3$\left( k\right) $ becomes sufficient
for identification$,$ cf. Theorem~\ref{1A1B}. \bigskip

\noindent \textbf{Assumption~6C}. \label{Ass4NoRotC)}For all $\mathbf{\theta 
}\neq \mathbf{\theta }_{0},$ $\mathbf{\theta }\in \mathcal{S}$, $\  \Phi _{%
\mathbf{\theta }}^{-1}\left( z\right) \Theta _{\mathbf{\theta }}\left(
z\right) \neq \Phi _{\mathbf{\theta }_{0}}^{-1}\left( z\right) \Theta _{%
\mathbf{\theta }_{0}}\left( z\right) P_{d}$ for any signed permutation
matrix $P_{d}$ in a subset of positive measure of $\left \{ z\in \mathbb{C}%
:\left \vert z\right \vert =1\right \} $.\bigskip

Define the $L^{2}$ distance between $f_{\mathbf{a},k}(\boldsymbol{\lambda };%
\mathbf{\theta },\mathbf{\alpha })\ $and $f_{\mathbf{a},k}(\boldsymbol{%
\lambda })=f_{\mathbf{a},k}(\boldsymbol{\lambda };\mathbf{\theta }_{0},%
\mathbf{\alpha }_{k}^{0})$ for $k=3,4$ and all possible indices $\mathbf{a}$%
, 
\begin{equation*}
\mathcal{L}_{k}\left( \mathbf{\theta },\mathbf{\alpha }\right) :=\sum_{%
\mathbf{a}}\int_{\Pi ^{k-1}}\left \vert f_{\mathbf{a},k}(\boldsymbol{\lambda 
};\mathbf{\theta },\mathbf{\alpha })-f_{\mathbf{a},k}(\boldsymbol{\lambda }%
)\right \vert ^{2}d\boldsymbol{\lambda },
\end{equation*}%
in terms of the parameters $\mathbf{\theta }$ and the marginal cumulants $%
\mathbf{\alpha }$ of order $k,$ and define%
\begin{equation*}
\mathcal{L}_{2}\left( \mathbf{\theta }\right) :=\sum_{\mathbf{a}}\int_{\Pi
}\left \vert f_{\mathbf{a},2}(\lambda ;\mathbf{\theta })-f_{\mathbf{a}%
,2}(\lambda )\right \vert ^{2}d\lambda,
\end{equation*}%
where $f_{\mathbf{a},2}(\lambda ;\mathbf{\theta })$ is the parametric model
for the second order spectral density $f_{\mathbf{a},2}(\lambda )$ for all
pairs $\mathbf{a}=\left( \mathbf{a}(1),\mathbf{a}(2)\right) ,$ which only
depends on $\mathbf{\theta }$ because of the normalization $\mathbb{V}\left( 
\mathbf{\varepsilon }_{t}\right) =\mathbf{I}_{d}$ under $\left( \ref{Rep}%
\right) .$ The loss function $\mathcal{L}_{2}$ is only able to identify $%
\Phi _{\mathbf{\theta }}^{-1}\left( z\right) \Theta _{\mathbf{\theta }%
}\left( z\right) $ up to a Blaschke factor so it cannot identify $\mathbf{%
\theta }$ under Assumptions~3 and~5 in absence of more restrictions.

We now show that the conclusions of Theorem~\ref{1A1B} extend to these
parametric loss functions under any of the versions of Assumption~6 (A, B or
C) with the corresponding restrictions on the parameter space. Thus, let $%
\mathcal{S}^{+}\subseteq \mathcal{S}$ be a compact set where we impose the
restriction that the diagonal elements of $\Theta _{0}\left( \mathbf{\theta }%
\right) $ are positive, while in the compact $\mathcal{S}^{\max }\subseteq 
\mathcal{S}^{+}$ we further impose that $\Pi _{j=1}^{d}\Theta _{0,jj}\left( 
\mathbf{\theta }\right) >\max_{\mathbf{\theta }^{\prime }\in \mathcal{S}%
^{\max }}\Pi _{j=1}^{d}\Theta _{0,jj}\left( \mathbf{\theta }^{\prime
}\right) $ where $\mathbf{\theta }^{\prime }$ is any parameter vector which
describes a permutation of the columns of $\Theta _{\mathbf{\theta }}\left(
L\right) ,$ i.e. $\Theta _{\mathbf{\theta }^{\prime }}\left( L\right)
=\Theta _{\mathbf{\theta }}\left( L\right) P_{d}$ for some signed
permutation matrix $P_{d}.$ \bigskip

\begin{theorem}
\label{Th4} Under Assumptions 1$\left( k\right) ,$ 3$\left( k\right) ,$ 5
and 6$\left( k\right) $, $k\in \mathcal{K},$ for any $\nu >0,$ there exists
an $\epsilon >0$ such that,%
\begin{equation*}
\inf_{\mathbf{\theta }\in \mathcal{S},\mathbf{\alpha }:\left \Vert \mathbf{%
\theta }-\mathbf{\theta }_{0}\right \Vert +\left \Vert \mathbf{\alpha }-%
\mathbf{\alpha }_{k}^0\right \Vert \geq \nu >0}\mathcal{L}_{2}\left( \mathbf{%
\theta }\right) +\mathcal{L}_{k}\left( \mathbf{\theta },\mathbf{\alpha }%
\right) \geq \epsilon >0,
\end{equation*}%
where it is further imposed $\mathbf{\theta }\in \mathcal{S}^{\max }$ under
Assumption~6A and $\mathbf{\theta }\in \mathcal{S}^{+}\ $and $\mathbf{\alpha 
}\in \mathcal{D}_{k}$ under Assumption~6B$\left( k\right) $.
\end{theorem}

This result provides identification of the dynamics and scaling parameters $%
\mathbf{\theta }\ $and the marginal cumulant vector $\mathbf{\alpha}$ in
absence of knowledge on the possible noninvertibility (or noncausality) of $%
Y_{t}.$ Following the comments to Theorem~\ref{1A1B} and despite
Assumption~5.2 allows that for some $\mathbf{\theta \neq \theta }_{0}$ and
some BM $A,$ $\Phi _{\mathbf{\theta }}^{-1}\left( z\right) \Theta _{\mathbf{%
\theta }}\left( z\right) =\Phi _{\mathbf{\theta }_{0}}^{-1}\left( z\right)
\Theta _{\mathbf{\theta }_{0}}\left( z\right) A(z)$ for almost all $z$ in
the unit circle so that $\mathcal{L}_{2}\left( \mathbf{\theta }\right) =0,$ $%
\mathcal{L}_{k}$ is not minimized unless $A(z)=P_{d}$ and $\mathbf{\alpha }%
_{3}=P_{d}^{\prime }\mathbf{\alpha }_{3}^{0}\ $or $\mathbf{\alpha }%
_{4}=P_{d}^{+\prime }\mathbf{\alpha }_{4}^{0}$. However, Assumption~6
imposes a unique ordering (and sign) to discard any signed permutations of
the columns of $\Phi _{\mathbf{\theta }_{0}}^{-1}\left( z\right) \Theta _{%
\mathbf{\theta }_{0}}\left( z\right) $ that would generate the same (second
and $k$-order) dynamics of $Y_{t}$ based on a signed permuted version of $%
\mathbf{\varepsilon }_{t}.$ Further, similarly as in VL, $\mathcal{L}_{3}$
and $\mathcal{L}_{4}$ could be considered jointly. Thus, if the target
function includes contributions from $\mathcal{L}_{3}$ and $\mathcal{L}_{4},$
then both $\mathbf{\alpha }_{3}$ and $\mathbf{\alpha }_{4}$ are identified
by the ordering of the columns of $\Theta _{0}\left( \mathbf{\theta }%
_{0}\right) $ (Assumption~6A) or by the ordering information from a single
higher order cumulant (Assumption~6B$\left( k\right) )$.

By contrast with the univariate analysis of VL, in the multivariate case it
is not possible to identify jointly $\mathbf{\theta }$ and $\mathbf{\alpha }%
_{k}$ from only $\mathcal{L}_{k}$ without $\mathcal{L}_{2},$ because the MA
matrix polynomial $\Theta _{\mathbf{\theta }}\left( z\right) $ incorporates
the scaling in $\Theta _{0}\left( \mathbf{\theta }\right) $ as we set $\mathbb{V}
\left( \mathbf{\varepsilon }_{t}\right) =\mathbf{I}_{d},$ so that $\mathbf{%
\alpha }_{k}$ are skewness and kurtosis \textit{coefficients}. In
alternative parameterizations, e.g. when setting $\Theta _{0}\left( \mathbf{%
\theta }\right) =\mathbf{I}_{d}$ for all $\mathbf{\theta }$ and%
\begin{equation*}
\Phi _{\mathbf{\theta }}\left( L\right) Y_{t}=\Theta _{\mathbf{\theta }%
}\left( L\right) \Omega \mathbf{\varepsilon }_{t},\  \  \  \mathbf{\varepsilon }%
_{t}\sim iid_{k}\left( \mathbf{0},\mathbf{I}_{d},\text{v}\mathbf{\kappa }%
_{k}^{IC}\left( \mathbf{\alpha }_{k}^{0}\right) ,k\in \mathcal{K}\right),
\end{equation*}%
where $\Omega $ is a non-singular matrix parameterized independently of $%
\mathbf{\theta },$ the dynamics parameters $\mathbf{\theta }$ and the
marginal cumulants $\mathbf{\alpha }_{k}$ could be identified jointly by a
single $\mathcal{L}_{k},$ $k=3\ $or $4,$ for a given scaling rotation $%
\Omega ,$ which has still to be identified together with $\mathcal{L}_{2}.$%
\bigskip

\section{5. Parameter Minimum Distance Estimation}

Given a time series of $Y_{t}$, $t=1,\ldots ,T$, we define sample analogs of
the loss functions $\mathcal{L}_{k}\left( \mathbf{\theta },\mathbf{\alpha }%
\right) $ for $k=3,4$ and $\mathcal{L}_{2}\left( \mathbf{\theta }\right) $
as in VL (2018) and Brillinger (1985), 
\begin{equation*}
\mathcal{L}_{k,T}\left( \mathbf{\theta },\mathbf{\alpha }\right) :=\frac{%
\left( 2\pi \right) ^{2k-2}}{T^{k-1}} \sum_{\mathbf{a}}\sum_{\mathbf{\lambda 
}_{\mathbf{j}}}\left \vert f_{\mathbf{a},k}(\mathbf{\lambda }_{\mathbf{j}};%
\mathbf{\theta },\mathbf{\alpha })-I_{\mathbf{a},k}(\mathbf{\lambda }_{%
\mathbf{j}})\right \vert ^{2},
\end{equation*}%
replacing the true spectral densities by the sample higher order
periodograms $I_{\mathbf{a},k}$, 
\begin{equation*}
I_{\mathbf{a},k}\left( \lambda _{1},\ldots ,\lambda _{k-1}\right) :=\frac{1}{%
\left( 2\pi \right) ^{k-1}T}w_{T,\mathbf{a}(1)}\left( \lambda _{1}\right)
\cdots w_{T,\mathbf{a}(k-1)}\left( \lambda _{k-1}\right) w_{T,\mathbf{a}%
(k)}\left( -\lambda _{1}\cdots -\lambda _{k-1}\right),
\end{equation*}%
where $w_{T}\left( \lambda \right) =\sum_{t=1}^{T}e^{-i\lambda t}Y_{t}$ is
the discrete Fourier transform (DFT) of $Y_{t}.$ In $\mathcal{L}_{k,T},$ the
summation in $\mathbf{\lambda }_{\mathbf{j}}=(\lambda _{j_{1}},\ldots
,\lambda _{j_{k-1}}) $ for Fourier frequencies $\lambda _{j_n}=2\pi j_n/T$
runs for all $j_{n}=1,\ldots ,T-1,$ $n=1,\ldots ,k-1,$ excluding $%
j_{a}+j_{b}=0\func{mod}(T),$ $a\neq b, $ and $j_{a}+j_{b}+j_{c}=0\func{mod}%
(T),$ all $a$, $b$ and $c$ different, for sample mean correction. Similarly,%
\begin{equation*}
\mathcal{L}_{2,T}\left( \mathbf{\theta }\right) :=\frac{\left( 2\pi \right)
^{2}}{T}\sum_{\mathbf{a}}\sum_{j=1}^{T-1}\left \vert f_{\mathbf{a}%
,2}(\lambda _{j};\mathbf{\theta })-I_{\mathbf{a},2}(\lambda _{j})\right
\vert ^{2}
\end{equation*}%
for the usual periodogram $I_{\mathbf{a},2}$ and spectral density $f_{%
\mathbf{a},2}.$

The $d^{k}$-vector containing all $f_{\mathbf{a},k}$ spectral densities can
be written under Assumptions~1$\left( k\right) $ and~3$\left( k\right) $ as $%
\left( 2\pi \right) ^{1-k}\mathbf{\Psi }_{k}\left( \mathbf{\lambda ;\theta }%
\right) \mathbf{S}_{k}\mathbf{\alpha },\ $with$\  \mathbf{\Psi }_{k}\left( 
\mathbf{\lambda ;\theta }\right) =\mathbf{\Psi }^{\otimes k}\left( \mathbf{%
\lambda ;\theta }\right) ,$ so that for $k=3,4$%
\begin{equation*}
\mathcal{L}_{k,T}\left( \mathbf{\theta },\mathbf{\alpha }\right) =\frac{1}{%
T^{k-1}}\sum_{\mathbf{\lambda }_{\mathbf{j}}}\left( \mathbf{\Psi }_{k}\left( 
\mathbf{\lambda }_{\mathbf{j}}\mathbf{;\theta }\right) \mathbf{S}_{k}\mathbf{%
\alpha }-\mathbb{I}_{k}(\mathbf{\lambda }_{\mathbf{j}})\right) ^{\ast
}\left( \mathbf{\Psi }_{k}\left( \mathbf{\lambda ;\theta }\right) \mathbf{S}%
_{k}\mathbf{\alpha }-\mathbb{I}_{k}(\mathbf{\lambda }_{\mathbf{j}})\right) ,
\end{equation*}%
where $\mathbb{I}_{k}(\mathbf{\lambda })$ is the $d^{k}\times 1$ vector
stacking all ($k=3$ order) biperiodograms and ($k=4$ order) triperiodograms
of $Y_{t},$ normalized by $\left( 2\pi \right) ^{k-1},$ 
\begin{equation*}
\mathbb{I}_{k}(\mathbf{\lambda }) : =\frac{1}{T}w_{T}\left( -\lambda _{1}-
\cdots -\lambda _{k-1}\right)\otimes \cdots \otimes w_{T}\left( \lambda
_{2}\right) \otimes w_{T}\left( \lambda _{1}\right) .
\end{equation*}%
Correspondingly,%
\begin{equation*}
\mathcal{L}_{2,T}\left( \mathbf{\theta }\right) =\frac{1}{T}%
\sum_{j=1}^{T-1}\left( \mathbf{\Psi }_{2}\left( \lambda _{j}\mathbf{;\theta }%
\right) \text{vec}\left( \mathbf{I}_{d}\right) -\mathbb{I}_{2}(\lambda
_{j})\right) ^{\ast }\left( \mathbf{\Psi }_{2}\left( \lambda _{j}\mathbf{%
;\theta }\right) \text{vec}\left( \mathbf{I}_{d}\right) -\mathbb{I}%
_{2}(\lambda _{j})\right),
\end{equation*}%
where $\mathbb{I}_{2}\left( \lambda \right) =\frac{1}{T}w_{T}\left( -\lambda
\right) \otimes w_{T}\left( \lambda \right) =2\pi $vec$\left( I_{YY}\left(
\lambda \right) \right) $ and $I_{YY}\left( \lambda \right) =\left( 2\pi
T\right) ^{-1}w_{T}\left( \lambda \right) w_{T}^{\ast }\left( \lambda
\right) =\left \{ I_{(a,b),2}\right \} _{a,b=1,\ldots ,d}$ is the usual
periodogram matrix of $Y_{t}.$

We set the following minimum distance parameter estimates for weights $w_k$, 
\begin{equation*}
\left( \mathbf{\hat{\theta}}_{w,T},\  \mathbf{\hat{\alpha}}_{k,T},\ k\in 
\mathcal{K}\right) :=\arg \min_{\mathbf{\theta }\in \mathcal{S}, \mathbf{%
\alpha }_{k},k\in \mathcal{K}}\mathcal{L}_{2,T}\left( \mathbf{\theta }%
\right) +\sum_{k\in \mathcal{K}}w_{k}\mathcal{L}_{k,T}\left( \mathbf{\theta }%
,\mathbf{\alpha }_{k}\right) ,
\end{equation*}%
where min$_{k\in \mathcal{K}}w_{k}>0$ and we further restrict to $\mathbf{%
\theta }\in \mathcal{S}^{\max }$ under Assumption~6A and to $\mathbf{\theta }%
\in \mathcal{S}^{+}\ $and $\mathbf{\alpha }_{k} \in \mathcal{D}_{k}$ under
Assumption~6B$\left( k\right) ,$ $k\in \mathcal{K}$. The main purpose of
combining loss functions $\mathcal{L}_{k,T}$ involving cumulants of
different orders $k=3,4,$ is robustness to lack of identification due to
failure of the nonzero cumulant condition of Assumption~3$\left( k\right) $
for a single $k$, as efficiency gains are possible though difficult to
characterize, even in the univariate case, see Lobato and Velasco (2018).
Additionally, we always need to include $\mathcal{L}_{2,T}$ in our loss
function for scaling identification.

Solving the first order conditions for $\mathbf{\alpha }_{k},$ $k\in 
\mathcal{K}$, using 
\begin{equation*}
\frac{\partial }{\partial \mathbf{\alpha }}\mathcal{L}_{k,T}\left( \mathbf{%
\theta ,\alpha }\right) =\frac{2}{T^{k-1}}\sum_{\mathbf{\lambda }_{\mathbf{j}%
}}\func{Re}\left \{ \mathbf{S}_{k}^{\prime }\mathbf{\Psi }_{k}^{\ast }(%
\mathbf{\lambda }_{\mathbf{j}};\mathbf{\theta })\left( \mathbf{\Psi }_{k}(%
\mathbf{\lambda }_{\mathbf{j}};\mathbf{\theta })\mathbf{S}_{k}\mathbf{\alpha 
}-\mathbb{I}_{k}(\mathbf{\lambda }_{\mathbf{j}})\right) \right \} ,
\end{equation*}%
we can obtain the (unrestricted) estimate of $\mathbf{\alpha }_{k}$ for a
given $\mathbf{\theta }$ 
\begin{equation*}
\mathbf{\hat{\alpha}}_{k,T}\left( \mathbf{\theta }\right): =\left( \sum_{%
\mathbf{\lambda }_{\mathbf{j}}}\func{Re}\left \{ \mathbf{S}_{k}^{\prime }%
\mathbf{\Psi }_{k}^{\ast }(\mathbf{\lambda }_{\mathbf{j}};\mathbf{\theta })%
\mathbf{\Psi }_{k}(\mathbf{\lambda }_{\mathbf{j}};\mathbf{\theta })\mathbf{S}%
_{k}\right \} \right) ^{-1}\sum_{\mathbf{\lambda }_{\mathbf{j}}}\func{Re}%
\left \{ \mathbf{S}_{k}^{\prime }\mathbf{\Psi }_{k}^{\ast }(\mathbf{\lambda }%
_{\mathbf{j}};\mathbf{\theta })\mathbb{I}_{k}(\mathbf{\lambda }_{\mathbf{j}%
})\right \} ,
\end{equation*}%
and concentrate out $\mathbf{\alpha }_{k}$ in $\mathcal{L}_{k,T}\left( 
\mathbf{\theta },\mathbf{\alpha }_{k}\right) $ so that 
\begin{equation*}
\mathbf{\hat{\theta}}_{w,T}=\arg \min_{\mathbf{\theta }\in \mathcal{S}}%
\mathcal{L}_{w,T}\left( \mathbf{\theta }\right) ,
\end{equation*}%
or restricted to $\mathbf{\theta }\in \mathcal{S}^{\max }$ under
Assumption~6A, where%
\begin{equation*}
\mathcal{L}_{w,T}\left( \mathbf{\theta }\right) :=\mathcal{L}_{2,T}\left( 
\mathbf{\theta }\right) +\sum_{k\in \mathcal{K}}w_{k}\mathcal{\hat{L}}%
_{k,T}\left( \mathbf{\theta }\right) ,\  \  \mathcal{\hat{L}}_{k,T}\left( 
\mathbf{\theta }\right) :=\mathcal{L}_{k,T}\left( \mathbf{\theta },\mathbf{%
\hat{\alpha}}_{k,T}\left( \mathbf{\theta }\right) \right),
\end{equation*}%
and $\mathbf{\hat{\alpha}}_{k,T}=\mathbf{\hat{\alpha}}_{k,T}\left( \mathbf{%
\hat{\theta}}_{w,T}\right) .$ However, under Assumption~6B$\left( k\right) $
there is no guarantee that for $\mathbf{\theta }\in \mathcal{S}^{+}\ $it
holds that $\mathbf{\hat{\alpha}}_{k,T}\left( \mathbf{\theta }\right) \in 
\mathcal{D}_{k},$ $k\in \mathcal{K}$, and optimization should be done
simultaneously for $\left( \mathbf{\theta },\  \mathbf{\alpha }_{k}, \ k\in 
\mathcal{K}\right)$ in $\mathcal{S}^{+}\times \prod_{ \mathcal{K}}\mathcal{D}%
_{k}$. 
Then, consistency of estimates under our set of identifying conditions is
achieved exploiting that periodogram averages estimate consistently
integrals of the true spectral densities.\bigskip

\begin{theorem}
\label{Th5} Under Assumptions 1$\left( 2k\right) ,$ 3$\left( k\right) ,$ 5, 6%
$\left( k\right) $, $k\in \mathcal{K},$ as $T\rightarrow \infty ,$ 
\begin{equation*}
\left( \mathbf{\hat{\theta}}_{w,T},\  \mathbf{\hat{\alpha}}_{k,T},\ k\in 
\mathcal{K}\right) \rightarrow _{p}\left( \mathbf{\theta }_{0}, \  \mathbf{\
\alpha }_{k}^{0},\ k\in \mathcal{K}\right) .
\end{equation*}
\end{theorem}

Note that the result for $\mathbf{\  \hat{\alpha}}_{k,T}$ also holds for $k
\notin \mathcal{K} $ from the consistency of $\mathbf{\hat{\theta}}_{w,T}$
when Assumptions~1$\left( 2k \right) $ and 3$\left( k \right)$ hold for both 
$k=3,4$, but only one set of cumulants needs to satisfy the nonzero
conditions to guarantee identification. Further, independence and equal
distribution of order $2k$ in Assumption~1$\left( k\right) $ are used to
facilitate the asymptotic analysis despite are not necessary for
identification.\bigskip

\textbf{Asymptotic Distribution}\bigskip

Now we consider optimal weighting of higher order periodograms replacing $%
\mathcal{L}_{k,T}\left( \mathbf{\theta },\mathbf{\alpha }\right) $ by a
weighted loss function 
\begin{equation*}
\mathcal{L}_{k,T}^{EFF}\left( \mathbf{\theta },\mathbf{\alpha }\right) :=%
\frac{1}{T^{k-1}}\sum_{\mathbf{\lambda }_{\mathbf{j}}}\left( \mathbf{\Psi }%
_{k}(\mathbf{\lambda }_{\mathbf{j}};\mathbf{\theta })\mathbf{S}_{k}\mathbf{%
\alpha }-\mathbb{I}_{k}(\mathbf{\lambda }_{\mathbf{j}})\right) ^{\ast }%
\mathbf{W}_{k}(\mathbf{\lambda }_{\mathbf{j}};\mathbf{\tilde{\theta}}%
_{T})\left( \mathbf{\Psi }_{k}(\mathbf{\lambda }_{\mathbf{j}};\mathbf{\theta 
})\mathbf{S}_{k}\mathbf{\alpha }-\mathbb{I}_{k}(\mathbf{\lambda }_{\mathbf{j}%
})\right) ,
\end{equation*}%
where $\mathbf{W}_{k}(\mathbf{\lambda };\mathbf{\tilde{\theta}}_{T}):=\left( 
\mathbf{\Psi }_{k}(\mathbf{\lambda };\mathbf{\tilde{\theta}}_{T})\mathbf{%
\Psi }_{k}^{\ast }(\mathbf{\lambda };\mathbf{\tilde{\theta}}_{T})\right)
^{-1}$ gives weights inversely proportional to the modulus of the higher
order transfer function $\mathbf{\Psi }_{k}(\mathbf{\lambda };\mathbf{\theta 
})$ evaluated at some preliminary estimate $\mathbf{\tilde{\theta}}%
_{T}\rightarrow _{p}\mathbf{\theta }_{0}.$ This weighting does not require
structural identification since $\mathbf{W}_{k}$ produces the same weighting
when $\mathbf{\Psi }(e^{-i\lambda };\mathbf{\tilde{\theta}}_{T})$ is
replaced by $\mathbf{\Psi }(e^{-i\lambda };\mathbf{\tilde{\theta}}%
_{T})A\left( e^{-i\lambda }\right) $ for a BM $A\left( z\right) $, because%
\begin{equation*}
\left( \mathbf{\Psi }(e^{-i\lambda };\mathbf{\theta })A(e^{-i\lambda
})\right) ^{\otimes k}=\mathbf{\Psi }(e^{-i\lambda };\mathbf{\theta }%
)^{\otimes k}A(e^{-i \lambda })^{\otimes k}=\mathbf{\Psi }_{k}(\mathbf{%
\lambda };\mathbf{\theta })A^{\otimes k}(\mathbf{\lambda }).
\end{equation*}%
Further, instead of using the solution from the first order condition for $%
\mathbf{\alpha }$, we propose to use instead the simpler estimate%
\begin{equation*}
\mathbf{\hat{\alpha}}_{k,T}^{\dag }\left( \mathbf{\theta }\right) :=\frac{1}{%
T^{k-1}}\sum_{\mathbf{\lambda }_{\mathbf{j}}}\func{Re}\left \{ \mathbf{S}%
_{k}^{\prime }\mathbf{\Psi }_{k}^{-1}(\mathbf{\lambda }_{\mathbf{j}};\mathbf{%
\theta })\mathbb{I}_{k}(\mathbf{\lambda }_{\mathbf{j}})\right \} ,
\end{equation*}%
exploiting that $\mathbf{S}_{k}^{\prime }\mathbf{\Psi }_{k}^{\ast }(\mathbf{%
\lambda };\mathbf{\theta }_{0})\mathbf{W}_{k}(\mathbf{\lambda };\mathbf{%
\theta }_{0})\mathbf{\Psi }_{k}(\mathbf{\lambda };\mathbf{\theta }_{0})%
\mathbf{S}_{k}=\mathbf{S}_{k}^{\prime }\mathbf{S}_{k}=\mathbf{I}_{d},$
similar to the proposal in VL (2018), and replace $\mathcal{L}_{k,T}^{EFF}$
by the pseudo-profile loss function $\mathcal{\hat{L}}_{k,T}^{\dag }\left( 
\mathbf{\theta }\right) :=\mathcal{L}_{k,T}^{EFF}\left( \mathbf{\theta ,\hat{%
\alpha}}_{k,T}^{\dag }\left( \mathbf{\theta }\right) \right) .$ 

The consistency of parameter estimates%
\begin{equation*}
\mathbf{\hat{\theta}}_{w,T}^{\dag }\ := \  \arg \min_{\mathbf{\theta }\in 
\mathcal{S}} \mathcal{\hat{L}}_{w,T}^{\dag } \left( \mathbf{\theta } \right)
, \  \  \mathcal{\hat{L}}_{w,T}^{\dag } \left( \mathbf{\theta } \right) := 
\mathcal{L}_{2,T}\left( \mathbf{\theta }\right) +\sum_{k\in \mathcal{K}}w_{k}%
\mathcal{\hat{L}}_{k,T}^{\dag }\left( \mathbf{\theta }\right) ,\ 
\end{equation*}%
restricted to $\mathbf{\theta }\in \mathcal{S}^{\max }$ under Assumption~6A,
can be deduced by the same arguments as for $\mathbf{\hat{\theta}}_{w,T}$.
We focus on estimates $\mathbf{\hat{\theta}}_{w,T}^{\dag }$ based on $%
\mathcal{\hat{L}}_{w,T}^{\dag }$ because its analysis is simpler due to the
efficient weighting scheme and a more straightforward estimation effect from 
$\mathbf{\hat{\alpha}}_{k,T}^{\dag }\left( \mathbf{\theta }\right) ,$ but it
is immediate to show the asymptotic equivalence with estimates based on
minimizing $\mathcal{L}_{k,T}^{EFF}\left( \mathbf{\theta },\mathbf{\alpha }%
\right) $ for $\mathbf{\theta }\in \mathcal{S}^{+}\ $and $\mathbf{\alpha }%
\in \mathcal{D}_{k},$ $k\in \mathcal{K}$, under Assumption~6B$\left(
k\right) $, as well as that of $\mathbf{\hat{\alpha}}_{k,T}^{\dag }\left( 
\mathbf{\theta }\right) $ with $\mathbf{\hat{\alpha}}_{k,T}^{EFF}\left( 
\mathbf{\theta }\right) .$

To investigate the asymptotic distribution of parameter estimates we need
further restrictions on the parameterization and a local identification
condition. Define 
\begin{equation*}
\mathbf{H}_{k}\left( \mathbf{\theta }\right) :=\left( 2\pi \right)
^{1-k}\int_{\Pi ^{k-1}}\func{Re}\left \{ \mathbf{B}_{k}^{\ast }\left( 
\mathbf{\lambda };\mathbf{\theta }\right) \mathbf{B}_{k}\left( \mathbf{%
\lambda };\mathbf{\theta }\right) \right \} d\mathbf{\lambda }
\end{equation*}%
where for $k=3,4$ we set%
\begin{equation*}
\mathbf{B}_{k}\left( \mathbf{\lambda };\mathbf{\theta }\right) :=\mathbf{%
\Psi }_{k}^{-1}\left( \mathbf{\lambda };\mathbf{\theta }\right) \mathbf{\dot{%
\Psi}}_{k}(\mathbf{\lambda };\mathbf{\theta })-\mathbf{S}_{k}\mathbf{S}%
_{k}^{\prime }\mathbf{\bar{\Lambda}}_{k}\left( \mathbf{\theta }\right)
=\sum_{j=1}^{k}\mathbf{B}_{k,j}\left( \lambda _{j};\mathbf{\theta }\right)
,\  \  \ 
\end{equation*}%
with $\mathbf{\bar{\Lambda}}_{k}\left( \mathbf{\theta }\right) :=\left( 2\pi
\right) ^{1-k}\int_{\Pi ^{k-1}}\func{Re}\left \{ \mathbf{\Psi }%
_{k}^{-1}\left( \mathbf{\lambda };\mathbf{\theta }\right) \mathbf{\dot{\Psi}}%
_{k}(\mathbf{\lambda };\mathbf{\theta })\right \} d\mathbf{\lambda }$, and 
$\mathbf{B}_{2}\left( \lambda ;\mathbf{\theta }\right) :=\mathbf{\Psi }%
_{2}^{-1}\left( \lambda ;\mathbf{\theta }\right) \mathbf{\dot{\Psi}}%
_{2}(\lambda ;\mathbf{\theta })$\\ $=\mathbf{B}_{2,1}\left( \lambda ;\mathbf{%
\theta }\right) +\mathbf{B}_{2,2}\left( -\lambda ;\mathbf{\theta }\right) ,$
with%
\begin{equation*}
\mathbf{\dot{\Psi}}_{k}(\mathbf{\lambda };\mathbf{\theta }):=\left( \mathbf{%
\dot{\Psi}}_{k}^{\left( 1\right) }(\mathbf{\lambda };\mathbf{\theta }%
),\ldots ,\mathbf{\dot{\Psi}}_{k}^{\left( m\right) }(\mathbf{\lambda };%
\mathbf{\theta })\right) ,\  \  \  \mathbf{\dot{\Psi}}_{k}^{\left( \ell \right)
}(\mathbf{\lambda };\mathbf{\theta }):=\frac{\partial }{\partial \mathbf{%
\theta }_{\ell }}\mathbf{\Psi }_{k}(\mathbf{\lambda };\mathbf{\theta }),
\end{equation*}%
where $\mathbf{B}_{k,j}\left( \lambda _{j};\mathbf{\theta }\right) ,$ $%
j=1,\ldots ,k,$ are obtained at once by the $k$-fold multiplicative
structure of $\mathbf{B}_{k}.$ Introduce for $\mathbf{\alpha }=\left \{ 
\mathbf{\alpha }_{k}\right \} _{k\in \mathcal{K}}$%
\begin{equation*}
\mathbf{\Sigma }\left( \mathbf{\theta ,\alpha }\right) :=\left( \mathbf{I}%
_{m}\otimes \text{vec}\left( \mathbf{I}_{d}\right) \right) ^{\prime }\mathbf{%
H}_{2}\left( \mathbf{\theta }\right) \left( \mathbf{I}_{m}\otimes \text{vec}%
\left( \mathbf{I}_{d}\right) \right) +\sum_{k\in \mathcal{K}} w_{k}\left( 
\mathbf{I}_{m}\otimes \mathbf{S}_{k}\mathbf{\alpha }_{k}\right) ^{\prime }%
\mathbf{H}_{k}\left( \mathbf{\theta }\right) \left( \mathbf{I}_{m}\otimes 
\mathbf{S}_{k}\mathbf{\alpha }_{k}\right)
\end{equation*}%
and the following assumption which imposes a rank condition on $\mathbf{%
\Sigma }\left( \mathbf{\theta }_{0},\mathbf{\alpha }_{0}\right) $ and
reinforces the smoothness conditions of the parameterization.\bigskip

\noindent \textbf{Assumption~7}.

\noindent \textbf{7.1.}\emph{\ Let }$\Phi _{i}\left( \mathbf{\theta }\right)
,$\emph{\ }$i=0,\ldots ,p,$ \emph{and\ }$\Theta _{i}\left( \mathbf{\theta }%
\right) ,$\emph{\ }$i=0,\ldots ,q,$\emph{\ have continuous third order
derivatives for all }$\theta \in S.$

\noindent \textbf{7.2.} $\mathbf{\theta }_{0}\in Int\left( \mathcal{S}^{\max
}\right) .$

\noindent \textbf{7.3. }$\mathbf{\Sigma }\left( \mathbf{\theta }_{0},\mathbf{%
\alpha }_{0}\right) >0.$\bigskip

Assumptions 7.1 and 7.2 are standard for the analysis of asymptotic
properties of extremum estimates, while sufficient conditions for the local
identification Assumption~7.3 are the full rank of individual Hessian
matrices, $\mathbf{H}_{k}\left( \mathbf{\theta }_{0}\right) >0$, for at
least one $k=2,3,4$ with $w_{k}>0$ and Assumption~3($k$) for $k\in \mathcal{K%
}.$ $\mathbf{H}_{2}$ is similar to the Hessian of PMLE estimates under
Gaussianity (and causality and invertibility), like Whittle estimates, which
only use second order information, noting that here $\mathbf{H}_{2}$
includes also the scaling parameters in $\mathbf{\Theta }_{0}$. This
indicates that usual methods are sufficient for local identification and
local asymptotic inference, but non-Gaussian information is key to achieve
global identification and potential efficiency improvements. Further, the
centering terms in $\mathbf{\bar{\Lambda}}_{k}\left( \mathbf{\theta }\right) 
$ reflect the higher order cumulant estimation effect, with the ones
corresponding to non scaling parameters being identically zero for causal
and invertible processes.

To simplify the asymptotic variance of $\mathbf{\hat{\theta}}_{w,T}^{\dag }$
we could strengthen Assumption~3$(k),$ $k\in \mathcal{K}$, to ICA of order $%
k=8$, as we do in Appendix~D to obtain explicit formulae for the variance of
the following vectors of powers of $\mathbf{\varepsilon }_{t}$ and $\mathbf{%
\varepsilon }_{r}$ for $t\neq r,$ 
\begin{eqnarray*}
\mathbf{\varepsilon }_{t,r}^{\left[ 2\right] } &:=&\left[ 
\begin{array}{c}
\mathbf{\varepsilon }_{t}\otimes \mathbf{\varepsilon }_{r} \\ 
\mathbf{\varepsilon }_{r}\otimes \mathbf{\varepsilon }_{t}%
\end{array}%
\right] \\
\mathbf{\varepsilon }_{t,r}^{\left[ 3\right] } &:=&\left[ 
\begin{array}{c}
\mathbf{\varepsilon }_{t}\otimes \mathbf{\varepsilon }_{t}\otimes \mathbf{%
\varepsilon }_{r} \\ 
\mathbf{\varepsilon }_{t}\otimes \mathbf{\varepsilon }_{r}\otimes \mathbf{%
\varepsilon }_{t} \\ 
\mathbf{\varepsilon }_{r}\otimes \mathbf{\varepsilon }_{t}\otimes \mathbf{%
\varepsilon }_{t}%
\end{array}%
\right] -\sum_{a=1}^{d}\left[ 
\begin{array}{c}
\mathbf{e}_{a}\otimes \mathbf{e}_{a}\otimes \mathbf{\varepsilon }_{r} \\ 
\mathbf{e}_{a}\otimes \mathbf{\varepsilon }_{r}\otimes \mathbf{e}_{a} \\ 
\mathbf{\varepsilon }_{r}\otimes \mathbf{e}_{a}\otimes \mathbf{e}_{a}%
\end{array}%
\right] ,\  \  \\
\mathbf{\varepsilon }_{t,r}^{\left[ 4\right] } &:=&\left[ 
\begin{array}{c}
\mathbf{\varepsilon }_{t}\otimes \mathbf{\varepsilon }_{t}\otimes \mathbf{%
\varepsilon }_{t}\otimes \mathbf{\varepsilon }_{r} \\ 
\mathbf{\varepsilon }_{t}\otimes \mathbf{\varepsilon }_{t}\otimes \mathbf{%
\varepsilon }_{r}\otimes \mathbf{\varepsilon }_{t} \\ 
\mathbf{\varepsilon }_{t}\otimes \mathbf{\varepsilon }_{r}\otimes \mathbf{%
\varepsilon }_{t}\otimes \mathbf{\varepsilon }_{t} \\ 
\mathbf{\varepsilon }_{r}\otimes \mathbf{\varepsilon }_{t}\otimes \mathbf{%
\varepsilon }_{t}\otimes \mathbf{\varepsilon }_{t}%
\end{array}%
\right] -\sum_{a=1}^{d}\mathbf{\alpha }_{3a}^{0}\left[ 
\begin{array}{c}
\mathbf{e}_{a}\otimes \mathbf{e}_{a}\otimes \mathbf{e}_{a}\otimes \mathbf{%
\varepsilon }_{r} \\ 
\mathbf{e}_{a}\otimes \mathbf{e}_{a}\otimes \mathbf{\varepsilon }_{r}\otimes 
\mathbf{e}_{a} \\ 
\mathbf{e}_{a}\otimes \mathbf{\varepsilon }_{r}\otimes \mathbf{e}_{a}\otimes 
\mathbf{e}_{a} \\ 
\mathbf{\varepsilon }_{r}\otimes \mathbf{e}_{a}\otimes \mathbf{e}_{a}\otimes 
\mathbf{e}_{a}%
\end{array}%
\right] .
\end{eqnarray*}%
However, there is no need of nonzero assumptions on marginal cumulants for
order larger than max$_{\mathcal{K}}k$, apart of existence of $2k$ moments
for $k\in \mathcal{K}$, as, e.g., it is straightforward to justify estimates
with $w_{4}=0$ based on up most $k=3$ information with only six bounded
moments.

Define $\mathbf{C}_{k}\left( 0\right) :=\left( 2\pi \right) ^{1-k}\int_{\Pi
^{k-1}}\mathbf{B}_{k}^{\ast }\left( \mathbf{\lambda };\mathbf{\theta }%
_{0}\right) d\mathbf{\lambda }$ for $k=2,3,4$ and the row block matrices for 
$j=\pm 1,\pm 2,\ldots ,$ 
\begin{equation*}
\mathbf{\mathbf{C}}_{k}\left( j\right) :=\left \{ \left( 2\pi \right)
^{1-k}\int_{\Pi ^{k-1}}\mathbf{B}_{k,a}^{\ast }\left( \mathbf{\lambda };%
\mathbf{\theta }_{0}\right) e^{-ij\lambda _{a}}d\mathbf{\lambda }\right \}
_{a=1,\ldots ,k},
\end{equation*}%
and%
\begin{equation*}
\mathbf{\Omega }\left( \mathbf{\theta }_{0};\mathbf{C}\right) :=\left \{ 
\QATOP{{}}{{}}\mathbf{\Omega }_{ab}\left( \mathbf{\theta }_{0};\mathbf{C}%
\right) \QATOP{{}}{{}}\right \} _{a,b=2,3,4},
\end{equation*}%
where $\mathbf{\Omega }_{ab}\left( \mathbf{\theta }_{0};\mathbf{C}\right) :=%
\mathbf{\Phi }_{ab}^{0}\left( \mathbf{\theta }_{0};\mathbf{C}\right) +%
\mathbf{\Phi }_{ab}\left( \mathbf{\theta }_{0};\mathbf{C}\right) +\mathbf{%
\Phi }_{ab}^{\dag }\left( \mathbf{\theta }_{0};\mathbf{C}\right) $ with
\begin{eqnarray*}
\mathbf{\Phi }_{ab}^{0}\left( \mathbf{\theta }_{0};\mathbf{C}\right) &:= &%
\mathbf{C}_{a}\left( 0\right) \mathbb{C}\left[ \mathbf{\varepsilon }%
_{t}^{\otimes a},\mathbf{\varepsilon }_{t}^{\otimes b}\right] \mathbf{C}%
_{b}^{\prime }\left( 0\right) \\
\mathbf{\Phi }_{ab}\left( \mathbf{\theta }_{0};\mathbf{C}\right) &:=
&\sum_{j=-\infty ,\neq 0}^{\infty }\mathbf{C}_{a}\left( j\right) \mathbb{C}%
\left[ \mathbf{\varepsilon }_{t,r}^{\left[ a\right] },\mathbf{\varepsilon }%
_{r,t}^{\left[ b\right] }\right] \mathbf{C}_{b}^{\prime }\left( j\right) \\
\mathbf{\Phi }_{ab}^{\dag }\left( \mathbf{\theta }_{0};\mathbf{C}\right) &:=
&\sum_{j=-\infty ,\neq 0}^{\infty }\mathbf{C}_{a}\left( -j\right) \mathbb{C}%
\left[ \mathbf{\varepsilon }_{t,r}^{\left[ a\right] },\mathbf{\varepsilon }%
_{r,t}^{\left[ b\right] }\right] \mathbf{C}_{b}^{\prime }\left( j\right)
\end{eqnarray*}%
and variance-covariance matrices not depending on $t$ or $r$ by
stationarity of order $2k,$ $k\in \mathcal{K}$. These definitions provide
multivariate generalizations of the score variance expressions developed in
VL for averages of higher order periodograms accounting for possibly
nonfundamental solutions and allow a compact presentation of the asymptotic
distribution of $\mathbf{\hat{\theta}}_{w,T}^{\dag }$ in next theorem. See
Appendix~D for explicit expressions for each $k$.

\begin{theorem}
\label{Th6} Under Assumptions 1$\left( 2k\right) ,$ 3$\left( k\right) ,$ 5, 6%
$\left( k\right) ,$ 7, $k\in \mathcal{K},$ $\min_{\mathcal{K}}w_{k}>0,$ as $%
T\rightarrow \infty ,$ 
\begin{equation*}
\sqrt{T}\left( \mathbf{\hat{\theta}}_{w,T}^{\dag }-\mathbf{\theta }%
_{0}\right) \rightarrow _{d}N_{m}\left( 0\mathbf{,\Sigma }^{-1}\left( 
\mathbf{\theta }_{0},\mathbf{\alpha }_{0}\right) \mathbf{\delta }\left( 
\mathbf{\alpha }_{0}\right) \mathbf{\Omega }\left( \mathbf{\theta }_{0};%
\mathbf{B}\right) \mathbf{\delta }^{\prime }\left( \mathbf{\alpha }%
_{0}\right) \mathbf{\Sigma }^{-1}\left( \mathbf{\theta }_{0},\mathbf{\alpha }%
_{0}\right) \right) ,
\end{equation*}%
where%
\begin{equation*}
\mathbf{\delta }\left( \mathbf{\alpha }_{0}\right) :=\left[ \  \  \left( 
\mathbf{I}_{m}\otimes \text{vec}\left( \mathbf{I}_{d}\right) \right)
^{\prime }\  \  \left \vert \  \ w_{3}\left( \mathbf{I}_{m}\otimes \mathbf{S}%
_{3}\mathbf{\alpha }_{3}^{0}\right) ^{\prime }\  \  \right \vert \  \
w_{4}\left( \mathbf{I}_{m}\otimes \mathbf{S}_{4}\mathbf{\alpha }%
_{4}^{0}\right) ^{\prime }\  \  \right] .
\end{equation*}
\end{theorem}

As in VL, the terms $\mathbf{\Phi }_{ab}^{\dag }$ are only different from
zero for noncausal or noninvertible models for which $\mathbf{B}_{k}\left(
\lambda ;\mathbf{\theta }\right) $ have a representation with terms in $%
e^{-ij\lambda }$ for $j=-1,-2,\ldots .$ Note also that $\mathbf{\Phi }%
_{22}^{0}$ incorporates the scaling estimation effect, which is treated
separately in other parameterizations, as in the univariate model in VL and,
partially, in the one discussed in Appendix~E, because, e.g. for causal and
invertible models, only the contributions from elements of $\mathbf{\theta }$
affecting $\Theta _{0}\left( \mathbf{\theta }\right) $ are different from
zero in $\mathbf{C}_{2}\left( 0\right) $. In the same line, the form of $\mathbf{\Phi }_{ab}^{0}$, $a,b>2$, accounts for the estimation of higher order cumulants under the ICA restriction of Assumption~3$(k)$, see Appendix~D for details.

We next present the asymptotic distribution of $\mathbf{\hat{\alpha}}%
_{k,T}^{\dag }\left( \mathbf{\hat{\theta}}_{w,T}^{\dag }\right) $ for a
particular $k\in \left \{ 3,4\right \} $ without assumptions on $\mathbf{%
\alpha }_{k}^{0},$ which could contain many zeros, obtaining an easy test of
overidentification. The only requisite is that $\mathbf{\hat{\theta}}%
_{w,T}^{\dag }$ is consistent and asymptotic normal with identification
provided by, possibly, a different set of cumulants. We introduce the
following vector of fourth powers of errors which include all permutations
of pairs, where $E_{t}$ stands for expectation conditional on $\mathbf{%
\varepsilon }_{t}$, $t\neq r,$ which is relevant for the variance of the
sample kurtosis vector coefficients $\mathbf{\hat{\alpha}}_{4,T}^{\dag
}\left( \mathbf{\theta }_{0}\right) ,$%
\begin{equation*}
\mathbf{\varepsilon }_{t}^{\left[ 4\right] }:=E_{t}\left[ 
\begin{array}{c}
\mathbf{\varepsilon }_{t}\otimes \mathbf{\varepsilon }_{t}\otimes \mathbf{%
\varepsilon }_{r}\otimes \mathbf{\varepsilon }_{r} \\ 
\mathbf{\varepsilon }_{r}\otimes \mathbf{\varepsilon }_{r}\otimes \mathbf{%
\varepsilon }_{t}\otimes \mathbf{\varepsilon }_{t} \\ 
\mathbf{\varepsilon }_{t}\otimes \mathbf{\varepsilon }_{r}\otimes \mathbf{%
\varepsilon }_{t}\otimes \mathbf{\varepsilon }_{r} \\ 
\mathbf{\varepsilon }_{r}\otimes \mathbf{\varepsilon }_{t}\otimes \mathbf{%
\varepsilon }_{r}\otimes \mathbf{\varepsilon }_{t} \\ 
\mathbf{\varepsilon }_{t}\otimes \mathbf{\varepsilon }_{r}\otimes \mathbf{%
\varepsilon }_{r}\otimes \mathbf{\varepsilon }_{t} \\ 
\mathbf{\varepsilon }_{r}\otimes \mathbf{\varepsilon }_{t}\otimes \mathbf{%
\varepsilon }_{t}\otimes \mathbf{\varepsilon }_{r}%
\end{array}%
\right] ,
\end{equation*}%
where, e.g., $\mathbf{\varepsilon }_{t,1}^{\left[ 4\right] }=\sum_{a,b=1}^{d}%
\mathbf{\varepsilon }_{t,a}^{2}\left( \mathbf{e}_{a}\otimes \mathbf{e}%
_{a}\otimes \mathbf{e}_{b}\otimes \mathbf{e}_{b}\right) $ given that $E\left[
\mathbf{\varepsilon }_{r}\otimes \mathbf{\varepsilon }_{r}\right]
=\sum_{b=1}^{d}E\left[ \mathbf{\varepsilon }_{r,b}^{2}\right] \left( \mathbf{%
e}_{b}\otimes \mathbf{e}_{b}\right) $\\ $=\sum_{b=1}^{d}\left( \mathbf{e}%
_{b}\otimes \mathbf{e}_{b}\right) .$ Define $\mathbf{\Delta }_{k}:=\mathbf{1}%
_{3}^{\prime }\otimes \mathbf{S}_{k}^{\prime }$ and $\mathbf{\Omega }\left( 
\mathbf{\theta }_{0};\mathbf{D}_{k}\right) $ as $\mathbf{\Omega }\left( 
\mathbf{\theta }_{0};\mathbf{C}\right) $ with $\mathbf{C}_{h}$ replaced by $%
\mathbf{D}_{k,h},$ $h=2,3,4,$ where for $j=\pm 1,\pm 2,\ldots ,$%
\begin{equation*}
\mathbf{D}_{k,h}\left( j\right) :=-w_{h}\mathbf{\bar{\Lambda}}_{k}\left( 
\mathbf{\theta }_{0}\right) \left( \mathbf{I}_{m}\otimes \mathbf{S}_{k}%
\mathbf{\alpha }_{k}^{0}\right) \mathbf{\Sigma }^{-1}\left( \mathbf{\theta }%
_{0},\mathbf{\alpha }_{0}\right) \left( \mathbf{I}_{m}\otimes \mathbf{S}_{h}%
\mathbf{\alpha }_{h}^{0}\right) ^{\prime }\mathbf{\mathbf{C}}_{h}\left(
j\right)
\end{equation*}%
with $\mathbf{S}_{2}\mathbf{\alpha }_{2}^{0}=\ $vec$\left( \mathbf{I}%
_{d}\right) $ and $w_{2}=1,$ and for $j=0,$%
\begin{equation*}
\mathbf{D}_{k,h}\left( 0\right) :=\mathbf{I}_{d^{k}}1_{\left \{ k=h\right \}
}-w_{h}\mathbf{\bar{\Lambda}}_{k}\left( \mathbf{\theta }_{0}\right) \left( 
\mathbf{I}_{m}\otimes \mathbf{S}_{k}\mathbf{\alpha }_{k}^{0}\right) \mathbf{%
\Sigma }^{-1}\left( \mathbf{\theta }_{0},\mathbf{\alpha }_{0}\right) \left( 
\mathbf{I}_{m}\otimes \mathbf{S}_{h}\mathbf{\alpha }_{h}^{0}\right) ^{\prime
}\mathbf{\mathbf{C}}_{h}\left( 0\right) .
\end{equation*}%
Define also $\mathbf{\bar{\Delta}}_{4}:=\mathbf{1}_{4}^{\prime }\otimes 
\mathbf{S}_{4}^{\prime }\ $and the block matrix $\mathbf{\bar{\Phi}}%
^{0}\left( \mathbf{\theta }_{0};\mathbf{D}_{4}\right) :=\left \{ \mathbf{%
\Phi }_{ab}^{0}\left( \mathbf{\theta }_{0};\mathbf{D}_{4}\right) \right \}
_{a,b=2,3,4,5}$ being equal to $\mathbf{\Phi }^{0}$ adding an extra row and
column$\ $given by, $b=2,3,4,$%
\begin{eqnarray*}
\mathbf{\Phi }_{5b}^{0}\left( \mathbf{\theta }_{0};\mathbf{D}_{4}\right) &:=
&\left( \mathbf{1}_{6}^{\prime }\otimes \mathbf{I}_{d^{4}}\right) \mathbb{C}%
\left[ \mathbf{\varepsilon }_{t}^{\left[ 4\right] },\mathbf{\varepsilon }%
_{t}^{\otimes b}\right] \mathbf{D}_{k,b}^{\prime }\left( 0\right) \ = \ 
\mathbf{\Phi }_{b5}^{0}\left( \mathbf{\theta }_{0};\mathbf{D}_{4}\right)
^{\prime } \\
\mathbf{\Phi }_{55}^{0}\left( \mathbf{\theta }_{0};\mathbf{D}_{4}\right) &:=
&\left( \mathbf{1}_{6}^{\prime }\otimes \mathbf{I}_{d^{4}}\right) \mathbb{V}%
\left[ \mathbf{\varepsilon }_{t}^{\left[ 4\right] }\right] \left( \mathbf{1}%
_{6}^{\prime }\otimes \mathbf{I}_{d^{4}}\right) ^{\prime }.
\end{eqnarray*}

\begin{theorem}
\label{Th7} Under Assumptions 1$\left( 2h\right) ,$ $3\left( h\right) ,$ 5, 6%
$\left( h\right) ,$ 7, $h\in \mathcal{K}\cup \{k\},$ $\min_{\mathcal{K}%
}w_{h}>0,$ as $T\rightarrow \infty ,$ for $k=3$, 
\begin{equation*}
T^{1/2}\left( \mathbf{\hat{\alpha}}_{3,T}^{\dag }\left( \mathbf{\hat{\theta}}%
_{w,T}^{\dag }\right) -\mathbf{\alpha }_{3}^{0}\right) \rightarrow
_{d}N_{d}\left( 0,\mathbf{\Delta }_{3}\mathbf{\Omega }\left( \mathbf{\theta }%
_{0};\mathbf{D}_{3}\right) \mathbf{\Delta }_{3}^{\prime }\right)
\end{equation*}%
and for $k=4$, 
\begin{equation*}
T^{1/2}\left( \mathbf{\hat{\alpha}}_{4,T}^{\dag }\left( \mathbf{\hat{\theta}}%
_{w,T}^{\dag }\right) -\mathbf{\alpha }_{4}^{0}\right) \rightarrow
_{d}N_{d}\left( 0,\mathbf{\bar{\Delta}}_{4}\mathbf{\bar{\Phi}}^{0}\left( 
\mathbf{\theta }_{0};\mathbf{D}_{4}\right) \mathbf{\bar{\Delta}}_{4}^{\prime
}+\mathbf{\Delta }_{4}\left( \mathbf{\Phi }\left( \mathbf{\theta }_{0};%
\mathbf{D}_{4}\right) +\mathbf{\Phi }^{\dag }\left( \mathbf{\theta }_{0};%
\mathbf{D}_{4}\right) \right) \mathbf{\Delta }_{4}^{\prime }\right) .
\end{equation*}
\end{theorem}

\section{6. Parameter GMM Estimation and Bootstrap Approximations}

In this section we propose estimates that exploit efficiently all the
information used by the minimum distance estimates of the previous section
by minimizing simultaneously the score functions of $\mathcal{L}_{k,T}\left( 
\mathbf{\theta },\mathbf{\alpha }\right) $ for $k=2,3,4$ without need to
specify $w.$ Denote as in VL the gradient vector of the concentrated loss
functions for all $k=2,3,4$ by%
\begin{equation*}
\mathbb{S}_{T}\left( \mathbf{\theta }\right) : =\left( 
\begin{array}{c}
\frac{\partial }{\partial \theta }\mathcal{L}_{2,T}(\mathbf{\theta }) \\ 
\frac{\partial }{\partial \theta }\mathcal{\hat{L}}_{3,T}^{\dag }(\mathbf{%
\theta }) \\ 
\frac{\partial }{\partial \theta }\mathcal{\hat{L}}_{4,T}^{\dag }(\mathbf{%
\theta })%
\end{array}%
\right) ,
\end{equation*}%
and denote the asymptotic variance of $T^{1/2}\mathbb{S}_{T}\left( \mathbf{%
\theta }_{0}\right) $ by $\mathcal{V}$ depending on $\mathbf{\delta }$ and $%
\mathbf{\Omega }$ from Theorem~7, so we can consider the objective function 
\begin{equation}
\mathcal{Q}_{T}\left( \mathbf{\theta }\right) =\mathbb{S}_{T}\left( \mathbf{%
\theta }\right) ^{\prime }\mathcal{\hat{V}}_{T}^{-}\, \mathbb{S}_{T}\left( 
\mathbf{\theta }\right) ,  \label{gmmof}
\end{equation}%
to optimally weight the information on $\mathbf{\theta }_{0}$ contained in
the score vector $\mathbb{S}_{T}\left( \mathbf{\theta }\right) .$ Here $%
\mathcal{\hat{V}}_{T}^{-}$ is a consistent estimator of the matrix $\mathcal{%
V}^{-},$ which is a reflexive generalized inverse of $\mathcal{V}$, and
hence, satisfies $\mathcal{VV}^{-}\mathcal{V}=\mathcal{V}$ and $\mathcal{V}%
^{-}\mathcal{VV}^{-}=\mathcal{V}^{-}$ as the Moore-Penrose inverse. We
employ generalized inverses to take account of cases where the asymptotic
variance of $\mathbb{S}_{T}\left( \mathbf{\theta }\right) $ is default rank
when identification rank conditions in Assumption~3$\left( k\right) $ fail
for some $k=3,4$. Hence, the proposed efficient estimator of $\mathbf{\theta 
}_{0}$ is a Newton-Raphson step using (\ref{gmmof}), 
\begin{equation}
\mathbf{\hat{\theta}}_{GMM,T}:=\mathbf{\tilde{\theta}}_{T}-\left( \mathbb{H}%
_{T}\left( \mathbf{\tilde{\theta}}_{T}\right) ^{\prime }\mathcal{\hat{V}}%
_{T}^{-}\, \mathbb{H}_{T}\left( \mathbf{\tilde{\theta}}_{T}\right) \right)
^{-1}\mathbb{H}_{T}\left( \mathbf{\tilde{\theta}}_{T}\right) ^{\prime }%
\mathcal{\hat{V}}_{T}^{-}\, \mathbb{S}_{T}\left( \mathbf{\tilde{\theta}}%
_{T}\right) ,  \label{theestimator}
\end{equation}%
where $\mathbb{H}_{T}\left( \mathbf{\theta }\right) :=\left( \partial
/\partial \mathbf{\theta }^{\prime }\right) \mathbb{S}_{T}\left( \mathbf{%
\theta }\right) ,$ and the initial estimate $\tilde{\mathbf{\theta }}_{T}$
satisfying 
\begin{equation}
\mathbf{\tilde{\theta}}_{T}-\mathbf{\theta }_{0}=O_{p}\left( T^{-1/2}\right)
,  \label{condis}
\end{equation}%
could be $\mathbf{\hat{\theta}}_{w,T}^{\dag }$ or any other PML or GMM
estimate obtained under appropriate identifying conditions. Given (\ref%
{condis}), the consistency of $\mathbf{\hat{\theta}}_{GMM,T}$ is trivial and
the next theorem states its asymptotic distribution defining the column
block matrix with with $\mathbf{S}_{2}\mathbf{\alpha }_{2}^{0}=\ $vec$\left( 
\mathbf{I}_{d}\right) .$ 
\begin{eqnarray*}
\mathbb{H} &:=&p\lim_{T\rightarrow \infty }\mathbb{H}_{T}\left( \mathbf{%
\tilde{\theta}}_{T}\right) =p\lim_{T\rightarrow \infty }\frac{\partial }{%
\partial \mathbf{\theta }^{\prime }}\mathbb{S}_{T}\left( \mathbf{\theta }%
_{0}\right) \\
&=&\left \{ \left( \mathbf{I}_{m}\otimes \mathbf{S}_{k}\mathbf{\alpha }%
_{k}^{0}\right) ^{\prime }\mathbf{H}_{k}\left( \mathbf{\theta }_{0}\right)
\left( \mathbf{I}_{m}\otimes \mathbf{S}_{k}\mathbf{\alpha }_{k}^{0}\right)
\right \} _{k=2,3,4}
\end{eqnarray*}

\begin{theorem}
\label{THGMM} Under Assumptions 1$\left( 2k\right) ,$ 3$(k),$ $k\in \left \{
3,4\right \} $, 5, 7, and $\mathcal{\hat{V}}_{T}^{-}\rightarrow _{p}\mathcal{%
V}^{-}$ as $T\rightarrow \infty ,$%
\begin{equation*}
\sqrt{T}(\mathbf{\hat{\theta}}_{GMM,T}-\mathbf{\theta }_{0})\rightarrow
_{d}N\left( 0,(\mathbb{H}^{\prime }\mathcal{V}^{-}\mathbb{H})^{-1}\right) .
\end{equation*}
\end{theorem}

Note that Theorem~\ref{THGMM} holds irrespective of zero values in $\mathbf{%
\alpha }_{k}^{0}$ for $k=3$ and $4,$ so no need of this part of Assumption~3(%
$k$) for both $k=3,4$ as far as (\ref{condis}) holds, and though estimation
procedures to obtain candidates for $\mathbf{\tilde{\theta}}_{T}$ can rely
on some form of non-Gaussianity, they might use only second-order
information complemented by economic identification restrictions. For the
same reason, Assumption~6 is not needed under (\ref{condis}).

To perform inference avoiding possibly imprecise estimation of $\mathbf{%
\Omega }$ or $\mathcal{V}$ it is possible to use a simple parametric
bootstrap based on resampling from the empirical distribution of model
residuals. However this procedure would require model simulation and
re-estimation in each resample, which could be costly for large high
dimensional models. Alternatively, we could resample the following
linearization of the estimates in terms of the higher order periodograms of
true errors $\mathbf{\varepsilon }_{t},$ 
\begin{equation*}
\mathbf{\hat{\theta}}_{w,T}^{\dag }-\mathbf{\theta }_{0}=\mathbf{\Sigma }%
^{-1}\left( \mathbf{\theta }_{0},\mathbf{\alpha }_{0}\right) \sum _{k\in 
\mathcal{K} \cup 2}\frac{w_{k}}{T^{k-1}}\left( \mathbf{I}_{m}\otimes \mathbf{%
S}_{k}\mathbf{\alpha }_{k}^{0}\right) ^{\prime }\sum_{\mathbf{\lambda }_{%
\mathbf{j}}}\func{Re}\left \{ \mathbf{B}_{k}^{\ast }\left( \mathbf{\lambda }%
_{\mathbf{j}};\mathbf{\theta }_{0}\right) \mathbb{I}_{k}^{\mathbf{%
\varepsilon }}(\mathbf{\lambda }_{\mathbf{j}})\right \} +o_{p}\left(
T^{-1/2}\right) ,
\end{equation*}%
where for $k=2$ we set $\mathbf{S}_{2}\mathbf{\alpha }_{2}^{0}=$ vec$\left( 
\mathbf{I}_{d}\right) $, $w_{2}=1,$ and replace $\mathbb{I}_{2}^{\mathbf{%
\varepsilon }}(\mathbf{\lambda }_{\mathbf{j}})$ by $\mathbb{I}_{2}^{\mathbf{%
\varepsilon }}(\mathbf{\lambda }_{\mathbf{j}})-$vec$\left( \mathbf{I}%
_{2}\right) $. Resampled versions of the estimates are obtained replacing by 
$\left( \mathbf{\theta }_{0},\mathbf{\alpha }_{k}^{0}, k\in \mathcal{K}
\right) $ by $\left( \mathbf{\hat{\theta}}_{w,T}^{\dag },\mathbf{\hat{\alpha}%
}_{k,T}^{\dag }, k\in \mathcal{K} \right) $ and $\mathbb{I}_{k}^{\mathbf{%
\varepsilon }}$ by the periodograms $\mathbb{I}_{k}^{\mathbf{\hat{\varepsilon%
}}^{\star }}$ of resampled residuals $\mathbf{\hat{\varepsilon}}_{t}^{\star
} $ from the empirical distribution of residuals $\mathbf{\hat{\varepsilon}}%
_{t}=\mathbf{\varepsilon }_{t}\left( \mathbf{\hat{\theta}}_{w,T}^{\dag
}\right) ,$ $t=1,\ldots ,T,$ properly standardized, 
\begin{equation*}
\mathbf{\hat{\theta}}_{w,T}^{\star }:=\mathbf{\hat{\theta}}_{w,T}^{\dag }+%
\mathbf{\Sigma }^{-1}\left( \mathbf{\hat{\theta}}_{w,T}^{\dag },\mathbf{\hat{%
\alpha}}_{T}^{\dag }\right) \sum _{k\in \mathcal{K} \cup 2}\frac{w_{k}}{%
T^{k-1}}\left( \mathbf{I}_{m}\otimes \mathbf{S}_{k}\mathbf{\hat{\alpha}}%
_{k,T}^{\dag }\right) ^{\prime }\sum_{\mathbf{\lambda }_{\mathbf{j}}}\func{Re%
}\left \{ \mathbf{B}_{k}^{\ast }\left( \mathbf{\lambda }_{\mathbf{j}};%
\mathbf{\hat{\theta}}_{w,T}^{\dag }\right) \mathbb{I}_{k}^{\mathbf{\hat{%
\varepsilon}}^{\star }}(\mathbf{\lambda }_{\mathbf{j}})\right \} ,
\end{equation*}%
which requires only one computation of $\mathbf{\Sigma }\left( \mathbf{\hat{%
\theta}}_{w,T}^{\dag },\mathbf{\hat{\alpha}}_{T}^{\dag }\right) $ and $%
\mathbf{B}_{k}\left( \mathbf{\lambda },\mathbf{\hat{\theta}}_{w,T}^{\dag
}\right) ,$ but not simulation or parameter re-estimation in each resample.
Interestingly, residuals $\mathbf{\hat{\varepsilon}}_{t}$ are obtained
directly in the frequency domain after inversion of the residual DFT $%
\mathbf{\Psi }(\lambda ;\mathbf{\hat{\theta}}_{w,T}^{\dag
})^{-1}w_{T}(\lambda ),$ without need to care about noninvertible or
noncausal roots or imposing ad-hoc factorizations. Resampling of residuals $%
\mathbf{\hat{\varepsilon}}_{t}$ and the FFT can be implemented efficiently
to obtain periodograms of resampled residuals $\mathbb{I}_{k}^{\mathbf{\hat{%
\varepsilon}}^{\star }}$ for any $k,$ as well as set $\mathbf{\hat{\alpha}}%
_{kT}^{\star }=\mathbf{\hat{\alpha}}_{kT}^{\dag }\left( \mathbf{\hat{\theta}}%
_{w,T}^{\star }\right) $ evaluated at the periodogram of data simulated from
resampled $\mathbf{\hat{\varepsilon}}_{t}^{\star }$ to approximate the
distribution of cumulant estimates.

Similar resampling methods to construct $\mathbf{\hat{\theta}}_{GMM,T}$ (and
approximate its finite sample distribution) can be implemented from (\ref%
{theestimator}) in the same way as for $\mathbf{\hat{\theta}}_{w,T}^{\dag }$%
, since the elements of $\mathbb{S}_{T}$ and $\mathbb{H}_{T}$ are the same
as those involved in the linearization of $\mathbf{\hat{\theta}}_{w,T}^{\dag
},$ but pooled instead of aggregated. In particular, it is straightforward
to estimate $\mathcal{V}$ by the sample covariance of bootstrap versions of $%
\mathbb{S}_{T}$ evaluated at a preliminary inefficient $\mathbf{\hat{\theta}}%
_{w,T}^{\dag }.$ Then, it is possible to construct bootstrap standard
errors, significance tests on parameter values and overidentification tests
based on redundant cumulants or restrictions on functions of the IRF.\bigskip

\section{7. Simulations}

In this section we consider several Monte Carlo experiments to check the
finite sample performance of our identification and estimation procedures.
We simulate bivariate $\left( d=2\right) $ SVARMA$\left( p,q\right) $
systems with $\left( p,q\right) =\left( 1,0\right) $ and $\left( 1,1\right)
, $ $\Theta \left( L\right) =\left( I-B_{1}L\right) \Omega $ and 
\begin{equation*}
\Phi _{1}=\left( 
\begin{array}{cc}
0.9 & 0 \\ 
-0.4 & 0.7%
\end{array}%
\right) ,\ B_{1}=\left( 
\begin{array}{cc}
-\zeta _{1}^{\rho _{1}} & 0 \\ 
0 & -\zeta _{2}^{\rho _{2}}%
\end{array}%
\right) ,\  \  \Omega =\left( 
\begin{array}{cc}
10 & 4 \\ 
-2 & 5%
\end{array}%
\right).
\end{equation*}%
Then, $\Theta _{0}=\Omega $ satisfies Assumption~6A and $\Theta
_{1}=-B_{1}\Omega $ with $\left \vert \zeta _{j}\right \vert <1$ and $\rho
_{j}=\pm 1,$ so the two roots of $\det \left( \Theta \left( z\right) \right) 
$ are $-\zeta _{j}^{-\rho _{j}},$ $j=1,2$. We consider all configurations of
invertible and non-invertible roots (including mixed cases) when $p=0.$ We
also simulate the SVARMA$\left( 0,1\right) $ model considered by Gouri\'{e}%
roux et al. (2019) with mixed roots as for our VARMA$\left( 1,1\right) $
model. We set $\mathbf{\theta }_{0}=\left( \text{vec}\left( \Phi _{1}\right)
^{\prime },\text{vec}\left( \Omega \right) ^{\prime },\text{vec}\left(
B_{1}\right) ^{\prime }\right) ^{\prime }$ for $p=1$ and $\mathbf{\theta }%
_{0}=\left( \text{vec}\left( \Omega \right) ^{\prime },\text{vec}\left(
B_{1}\right) ^{\prime }\right) ^{\prime }$ for $p=0$ and simulate two sample
sizes $T=100,200$ with 1000 and 500 replications, respectively.

We consider three sets of innovations. The first type are mutually
independent standardized $\chi _{n_{j}}^{2}$ variates, $(n_{1},n_{2})=(6,1)$, 
designed to satisfy\ Assumption~6B$\left( k\right) $ for both $k=3$ or $k=4$
together with the value of $\Theta _{0}$ with $\mathbf{\alpha }%
_{3}^{0}=(1.155,2.828)^{\prime }$ and $\mathbf{\alpha }_{4}^{0}=(2,12)^{%
\prime }$. We also consider two further shock vectors, one composed of two
standardized $t_{n_{j}}$ variates, $(n_{1},n_{2})=(6,5)$ with $\mathbf{%
\alpha }_{3}^{0}=(0,0)^{\prime }$ and $\mathbf{\alpha }_{4}^{0}=(3,6)^{%
\prime }$ and one with a mixed normal and a $t_{6}$ variable as in Gouri\'{e}%
roux et al. (2019) with $\mathbf{\alpha }_{3}^{0}=(1,0)^{\prime }$ and $%
\mathbf{\alpha }_{4}^{0}=(6,3)^{\prime }.$ Note that these last two
distributions do not satisfy Assumption~3$\left( k\right) $ for $k=3$
because at least one component is symmetric, so (dynamic) identification
only relies on information from $k=4$ order cumulants. Also note that
Assumption~1$(2k)$ does not hold for these $t$ distributions, as they have 5
finite moments at most.

To identify both the correct rotation of the innovations and the location of
the MA lag polynomial roots we implement the following algorithm which
obtains different estimates in each step:

\begin{enumerate}
\item $\mathbf{\hat{\theta}}_{2,T}$: Causal and invertible reduced form
VARMA estimation using $\mathcal{L}_{2,T}$ loss function with Whittle
initial estimates to define the spectral weighting imposing uncorrelation
among components of $\mathcal{\varepsilon} _{t}$. If $p>0,$ preliminary IV
estimation of VAR parameters is performed as in Gouri\'{e}roux et al. (2019)
under the (correct) assumption of causality.

\item $\mathbf{\hat{\theta}}_{w,T}^\dag$: Minimum distance higher order
spectral estimation with unrestricted lag polynomial root location.

\begin{enumerate}
\item[2.1] Computation of all $2^{dq}=4$ MA basic representations of $%
\mathbf{\hat{\theta}}_{2,T}$ obtained by combinations of possibly inverted
MA roots using the procedure of Baggio and Ferrante (2019).

\item[2.2] Initial approximation of the rotation of the residuals for each
representation closest to component independence using a reconstruction ICA
algorithm as implemented in MATLAB RICA function.

\item[2.3] Higher order spectral estimation imposing independence of order $%
k=3$ or $k=3$ and $4$ in the components of $\mathbf{\varepsilon }_{t}$ by
minimizing $\mathcal{L}^\dag_{w,T}\left( \mathbf{\theta }\right) $ with $%
w\in \left \{ \left( 1,0\right) ,\left( 1,1\right) \right \} $ using the
rotated parameters from 2.2 of each specific root-model configurations as
initial estimates.

\item[2.4] Global minimum: choose the root configuration that minimizes $%
\mathcal{L}^\dag_{w,T}\left( \mathbf{\theta }\right) .$

\item[2.5] Permutation of components: choose the signed permutation of error
components to match Assumption~6A among all $2^{d}d!=8$ possible ones.
\end{enumerate}

\item $\mathbf{\hat{\theta}}_{GMM,T}$: Local GMM estimation based on
preliminary estimates $\mathbf{\hat{\theta}}_{w,T}^\dag $ from Step 2 and
bootstrap estimated variance of the score.\bigskip
\end{enumerate}

For Step 2.1 we also tried to minimize a penalized version of $\mathcal{L}%
_{2,T}$ that forces each specific root location configuration, with similar
results in most situations in terms of root identification, but less
efficient estimates. Alternatively, we also imposed Assumption~6B$\left(
k\right) $ in Step~2.5, but since it was not possible to get estimates of $%
\Theta _{0}$ with the right configuration of the implied cumulant estimates
in some replications, we do not report its results here. For comparison
purposes, we also report for invertible models the outcomes for an
infeasible version of the estimates $\mathbf{\hat{\theta}}_{2,T}$ obtained
in Step 1 (and after a RICA rotation to enforce component independence as in
Step 2.2) by finding the signed permuted version that minimizes the $%
\mathcal{L}^{2}$ distance to the true parameters. All $\mathcal{L}_{k,T}
\left( \mathbf{\theta }\right) $ objective functions minimized in the
different steps use the same causal and invertible weighting based on
Whittle estimates obtained in Step $1$ (as this is invariant to rotations of
innovations and flipping of polynomial roots). For all methods investigated
we also obtain estimates of $\mathbf{\alpha }_{k}$ and we report bias and
Root MSE across simulations of all estimates (and average absolute bias and
MSE for all elements of the vector $\mathbf{\theta }$).

\setlength{\tabcolsep}{3.3pt} 
\renewcommand*{\arraystretch}{.8}

\begin{table}[h]
\centering
\noindent \textbf{Table 1. }Percentage of correct
identification of the location of the MA roots by $\mathcal{L}_{k,T}.$ 
\begin{tabular}{lllllcccccc}
\hline \hline
&  &  &  &  & \multicolumn{4}{c}{SVARMA$(0,1)$} &  & SVARMA$(1,1)$ \\ 
\cline{6-9}\cline{11-11}
& Innovation & $T$ & $k$ &  & Invert. & Mixed & Non-Inv. & Mix-GMR &  & Mixed
\\ \hline \hline
& $\left( \chi _{6}^{2},\chi _{1}^{2}\right) $ & $100$ & $\mathcal{L}_{3}$ & 
& \  \  \ 92.3\  \  \  \  & 83.3 & 92.0 & 66.7 &  & 61.7 \\ 
&  &  & $\mathcal{L}_{3}+\mathcal{L}_{4}$ &  & 91.4 & 79.5 & 92.1 & 67.3 & 
& 60.6 \\ 
&  & $200$ & $\mathcal{L}_{3}$ &  & 99.0 & 99.0 & 99.6 & 76.6 &  & 70.0 \\ 
&  &  & $\mathcal{L}_{3}+\mathcal{L}_{4}$ &  & 89.2 & 86.8 & 88.4 & 76.4 & 
& 71.4 \\ \hline
& $(MN,t_{6})$ & $100$ & $\mathcal{L}_{3}$ &  & 78.0 & 75.3 & 77.3 & 36.8 & 
& 61.7 \\ 
&  &  & $\mathcal{L}_{3}+\mathcal{L}_{4}$ &  & 78.3 & 76.1 & 76.3 & 30.0 & 
& 59.8 \\ 
&  & $200$ & $\mathcal{L}_{3}$ &  & 89.4 & 69.8 & 86.6 & 35.8 &  & 72.0 \\ 
&  &  & $\mathcal{L}_{3}+\mathcal{L}_{4}$ &  & 89.4 & 76.6 & 87.4 & 37.6 & 
& 73.5 \\ \hline
& $(t_{6},t_{5})$ & $100$ & $\mathcal{L}_{3}$ &  & 52.2 & 51.0 & 48.3 & 43.7
&  & 39.3 \\ 
&  &  & $\mathcal{L}_{3}+\mathcal{L}_{4}$ &  & 54.0 & 52.0 & 50.4 & 42.9 & 
& 39.3 \\ 
&  & $200$ & $\mathcal{L}_{3}$ &  & 55.4 & 40.6 & 56.2 & 46.6 &  & 42.0 \\ 
&  &  & $\mathcal{L}_{3}+\mathcal{L}_{4}$ &  & 62.8 & 54.0 & 59.4 & 47.4 & 
& 41.2 \\ \hline \hline
\end{tabular}
\caption{ \footnotesize 
\noindent 
Identification of SVARMA$\left( p, 1 \right)$ models
 by minimizing $w_3 \mathcal{L}_{3,T}  +w_4 \mathcal{L}_{4,T}$, $p\in \{0,1\}$, $d= 2 ,$
 $\zeta _{1}= \zeta _{2}= 0.5 ,$ $\rho_{1}= 1 ,  \rho _{2}= 1 $. Mix-GMR is the bivariate SVARMA$(0,1)$ model simulated in Gouri\'{e}roux et al. (2019) with
$\mathbf{\theta} =(0,1,1,0.5,-0.5,1,0,-2)'$ and
 mixed MA roots $(0.5,2)$. Innovations are $\left( \chi _{6}^{2},\chi
_{1}^{2}\right), (MN(2.12,1.41^2;-0.24,0.58^2;0.1,0.9),t_6),$ and $(t_6,t_5)$.\  
$\mathcal{L}_{3}$ and $\mathcal{L}_{3}+\mathcal{L}_{4}$ 
identification use $w=(1,0)$ and $w=(1,1)$, respectively.}
\end{table}
\bigskip

In Table 1 we report the percentage of simulations that identified the right
location of roots for all the parameter configurations for the $p=0$ case
(just mixed case for Gouri\'{e}roux et al. (2019) and SVARMA$\left(
1,1\right) $ models), innovations, sample sizes and weights $w$ in the loss
functions. From the results for chi squared innovations we find that
skewness provides a very precise information for identification in
multivariate systems as was found for univariate models in VL, while the
additional use of kurtosis information does not help much for the sample
sizes considered and, even for the smallest ones, can introduce further
noise. For mixed normal and $t$ shocks our identification procedure is still
able to gain substantial information on the location of roots when
Assumption~3$\left( 3\right) $ fails for one component, but obtains quite
little when it fails for all components as for bivariate $t$ distributions
(about 50\% of right locations for each of the three configurations). When
adding information from kurtosis, i.e. $k=4$ cumulants, in these $k=3$
identification-failure situations, the results get better in almost all
cases, but with modest improvements in general given the small sample sizes
considered. Across all setups simulated, we found that mixed models with MA
roots both inside and outside the unit circle are more difficult to identify
than pure invertible or non-invertible systems, though only marginally in
same situations.

The results for the SVARMA$\left( 0,1\right) $ model simulated in Gouri\'{e}%
roux et al. (2019) also report significant differences for each set of
innovations, with kurtosis being relatively helpful for improving the
performance both for symmetric and asymmetric distributions. We observe a
similar pattern for SVARMA$\left( 1,1\right) $ models, but in general
identification of the MA polynomial roots is more complicated in the
presence of autorregressive dynamics. Still, kurtosis information seems more
valuable when skewness is default rank, but, when all series are perfectly
symmetric, identification results are poorer and only improve very slowly
with sample size.

\begin{table}[h!]
\renewcommand*{\arraystretch}{.75}
\centering
\noindent \textbf{Table 2.} Estimation of invertible SVARMA$\left( 0,1\right) .$ $T=100$ 
\begin{tabular}{rrrrrrrrr}
\hline \hline
& \multicolumn{8}{c}{Bias} \\ \cline{2-9}
& \multicolumn{2}{c}{$\left( \chi _{6}^{2},\chi _{1}^{2}\right) $} &  & 
\multicolumn{2}{c}{$(MN,t_{6})$} &  & \multicolumn{2}{c}{$(t_{6},t_{5})$} \\ 
\cline{2-3}\cline{5-6}\cline{8-9}
$\mathbf{\theta }_{0}$ & $\mathbf{\hat{\theta}}_{2,T}$ & $\mathbf{\hat{\theta%
}}_{w,T}$ &  & $\mathbf{\hat{\theta}}_{2,T}$ & $\mathbf{\hat{\theta}}_{w,T}$
&  & $\mathbf{\hat{\theta}}_{2,T}$ & $\mathbf{\hat{\theta}}_{w,T}$ \\ \hline
10 & -0.612 & -0.299 &  & -0.197 & -0.087 &  & -0.304 & -0.648 \\ 
-2 & 0.075 & 0.046 &  & -0.012 & 0.055 &  & -0.075 & 0.542 \\ 
4 & -0.162 & -0.083 &  & -0.078 & -0.180 &  & -0.061 & -1.464 \\ 
5 & -0.193 & -0.049 &  & -0.151 & -0.099 &  & -0.171 & -0.284 \\ 
-0.5 & -0.016 & -0.005 &  & -0.013 & -0.005 &  & -0.017 & 0.005 \\ 
0 & -0.003 & -0.001 &  & -0.002 & 0.001 &  & 0.002 & 0.010 \\ 
0 & -0.006 & -0.002 &  & -0.003 & 0.004 &  & -0.005 & -0.016 \\ 
-0.5 & -0.018 & -0.001 &  & -0.019 & -0.012 &  & -0.011 & 0.016 \\ 
AVE & 0.134 & 0.061 &  & 0.060 & 0.055 &  & 0.080 & 0.373 \\ 
$\mathbf{\alpha }_{3}^{0}$ & $\mathbf{\hat{\alpha}}_{32,T}$ & $\mathbf{\hat{%
\alpha}}_{3w,T}$ & $\mathbf{\alpha }_{3}^{0}$ & $\mathbf{\hat{\alpha}}%
_{32,T} $ & $\mathbf{\hat{\alpha}}_{3w,T}$ & $\mathbf{\alpha }_{3}^{0}$ & $%
\mathbf{\hat{\alpha}}_{32,T}$ & $\mathbf{\hat{\alpha}}_{3w,T}$ \\ \hline
1.155 & -0.141 & -0.128 & 2 & -0.129 & -0.174 & 0 & -0.050 & -0.008 \\ 
2.828 & -0.503 & -0.420 & 0 & -0.011 & -0.036 & 0 & -0.008 & -0.038 \\ 
$\mathbf{\alpha }_{4}^{0}$ & $\mathbf{\hat{\alpha}}_{42,T}$ & $\mathbf{\hat{%
\alpha}}_{4w,T}$ & $\mathbf{\alpha }_{4}^{0}$ & $\mathbf{\hat{\alpha}}%
_{42,T} $ & $\mathbf{\hat{\alpha}}_{4w,T}$ & $\mathbf{\alpha }_{4}^{0}$ & $%
\mathbf{\hat{\alpha}}_{42,T}$ & $\mathbf{\hat{\alpha}}_{4w,T}$ \\ \hline
2 & -0.481 & -0.568 & 6 & -0.413 & -0.600 & 3 & -0.411 & 0.449 \\ 
12 & -4.286 & -4.073 & 3 & -1.233 & -1.144 & 6 & -2.704 & -1.554 \\ 
\hline \hline
& \multicolumn{8}{c}{RMSE} \\ \cline{2-9}
& \multicolumn{2}{c}{$\left( \chi _{6}^{2},\chi _{1}^{2}\right) $} &  & 
\multicolumn{2}{c}{$(MN,t_{6})$} &  & \multicolumn{2}{c}{$(t_{6},t_{5})$} \\ 
\cline{2-3}\cline{5-6}\cline{8-9}
$\mathbf{\theta }_{0}$ & $\mathbf{\hat{\theta}}_{2,T}$ & $\mathbf{\hat{\theta%
}}_{w,T}$ &  & $\mathbf{\hat{\theta}}_{2,T}$ & $\mathbf{\hat{\theta}}_{w,T}$
&  & $\mathbf{\hat{\theta}}_{2,T}$ & $\mathbf{\hat{\theta}}_{w,T}$ \\ \hline
10 & 1.498 & 1.046 &  & 1.720 & 1.607 &  & 1.598 & 2.027 \\ 
-2 & 0.983 & 0.470 &  & 0.785 & 0.921 &  & 1.030 & 1.651 \\ 
4 & 2.233 & 1.340 &  & 1.402 & 1.694 &  & 2.140 & 3.740 \\ 
5 & 1.146 & 1.009 &  & 0.647 & 0.593 &  & 0.865 & 1.178 \\ 
-0.5 & 0.146 & 0.095 &  & 0.124 & 0.089 &  & 0.127 & 0.200 \\ 
0 & 0.059 & 0.051 &  & 0.052 & 0.047 &  & 0.057 & 0.088 \\ 
0 & 0.227 & 0.188 &  & 0.225 & 0.199 &  & 0.247 & 0.442 \\ 
-0.5 & 0.167 & 0.100 &  & 0.125 & 0.101 &  & 0.120 & 0.184 \\ 
AVE & 0.808 & 0.537 &  & 0.635 & 0.656 &  & 0.773 & 1.189 \\ 
$\mathbf{\alpha }_{3}^{0}$ & $\mathbf{\hat{\alpha}}_{32,T}$ & $\mathbf{\hat{%
\alpha}}_{3w,T}$ & $\mathbf{\alpha }_{3}^{0}$ & $\mathbf{\hat{\alpha}}%
_{32,T} $ & $\mathbf{\hat{\alpha}}_{3w,T}$ & $\mathbf{\alpha }_{3}^{0}$ & $%
\mathbf{\hat{\alpha}}_{32,T}$ & $\mathbf{\hat{\alpha}}_{3w,T}$ \\ \hline
1.155 & 0.419 & 0.394 & 2 & 0.483 & 0.510 & 0 & 0.797 & 0.887 \\ 
2.828 & 0.941 & 0.846 & 0 & 0.612 & 0.636 & 0 & 0.887 & 0.893 \\ 
$\mathbf{\alpha }_{4}^{0}$ & $\mathbf{\hat{\alpha}}_{42,T}$ & $\mathbf{\hat{%
\alpha}}_{4w,T}$ & $\mathbf{\alpha }_{4}^{0}$ & $\mathbf{\hat{\alpha}}%
_{42,T} $ & $\mathbf{\hat{\alpha}}_{4w,T}$ & $\mathbf{\alpha }_{4}^{0}$ & $%
\mathbf{\hat{\alpha}}_{42,T}$ & $\mathbf{\hat{\alpha}}_{4w,T}$ \\ \hline
2 & 2.137 & 2.116 & 6 & 2.851 & 2.815 & 3 & 3.556 & 3.219 \\ 
12 & 7.502 & 7.288 & 3 & 2.866 & 2.892 & 6 & 5.653 & 3.399 \\ \hline \hline
\end{tabular}
\caption{ \footnotesize 
\noindent Estimation of SVARMA$\left( 0, 1 \right)$ models by GMM, $d= 2 ,$ $\zeta _{1}= \zeta _{2}= 0.5 ,$ $%
\rho_{1}= 1 , \rho _{2}= 1 $. Innovations $\left( \chi _{6}^{2},\chi
_{1}^{2}\right), (MN(2.12,1.41^2;-0.24,0.58^2;0.1,0.9),t_6),$ and $(t_6,t_5)$.\  $\mathcal{L}_{2}+\mathcal{L}_{3}+\mathcal{L}_{4}$ identification and estimation, $w=(1,1)$. AVE is the average of the absolute value of the column.}
\end{table}

In Tables 2-5 we report the results on the performance of parameter
estimates only for the replications which correctly identified the MA root
location to make the comparisons meaningful. We can check from the output in
Table 2 for SVARMA$\left( 0,1\right) $ invertible models that the behaviour
of parameter estimates $\mathbf{\hat{\theta}}_{GMM,T}$ compares well in
terms of bias and RMSE with the unfeasible estimate $\mathbf{\hat{\theta}}%
_{2,T},$ which uses information on the true value of the parameters. Among
the different distributions, the case of the symmetric $t$ distribution
seems to be the most complicated, both in terms of bias and variability,
given that its heavy tails violate the moment condition in Theorem~6, while
the $\chi ^{2}$ distribution is the most informative given its strong
skewness. The results for non-invertible models (mixed and non-invertible
cases) in Table~3 are similar but in general estimates are less precise
comparing, for instance, the estimation of $\Omega $ contained in the first
four components of $\mathbf{\theta }.$

\begin{table}[h!]
\renewcommand*{\arraystretch}{.75} 
\centering
\noindent \textbf{Table 3.} Estimation of non-invertible SVARMA$\left( 0,1\right) .$ $T=100$ 
\begin{tabular}{rrrrrrrrrrrr}
\hline \hline
& \multicolumn{11}{c}{Bias} \\ \cline{2-12}
& \multicolumn{3}{c}{$\left( \chi _{6}^{2},\chi _{1}^{2}\right) $} &  & 
\multicolumn{3}{c}{$(MN,t_{6})$} &  & \multicolumn{3}{c}{$(t_{6},t_{5})$} \\ 
\cline{2-4}\cline{6-8}\cline{10-12}
$\mathbf{\theta }_{0}$ &  & Mixed & P.inv &  &  & Mixed & P.inv &  &  & Mixed
& P.inv \\ \hline
10 &  & -0.407 & -0.145 &  &  & -0.832 & 0.098 &  &  & -1.686 & 0.094 \\ 
-2 &  & 0.072 & 0.115 &  &  & 0.553 & -0.018 &  &  & 1.105 & 0.972 \\ 
4 &  & -0.208 & -0.124 &  &  & -1.242 & -0.046 &  &  & -2.576 & -1.869 \\ 
5 &  & 0.141 & 0.031 &  &  & 0.021 & 0.054 &  &  & 1.175 & 0.049 \\ 
-0.5/-2 &  & -0.046 & -0.034 &  &  & -0.049 & -0.039 &  &  & -0.413 & -0.003
\\ 
0 &  & 0.017 & 0.018 &  &  & 0.056 & -0.014 &  &  & -0.008 & -0.004 \\ 
0 &  & 0.012 & 0.039 &  &  & 0.061 & 0.017 &  &  & -0.016 & -0.009 \\ 
-2 &  & 0.017 & -0.053 &  &  & 0.050 & -0.008 &  &  & 0.409 & -0.021 \\ 
AVE &  & 0.115 & 0.070 &  &  & 0.358 & 0.037 &  &  & 0.923 & 0.378 \\ 
& $\mathbf{\alpha }_{3}^{0}$ &  &  &  & $\mathbf{\alpha }_{3}^{0}$ &  &  & 
& $\mathbf{\alpha }_{3}^{0}$ &  &  \\ \hline
& 1.155 & -0.120 & -0.103 &  & 2 & -0.095 & -0.175 &  & 0 & -0.047 & -0.030
\\ 
& 2.828 & -0.394 & -0.435 &  & 0 & -0.118 & -0.047 &  & 0 & 0.031 & 0.007 \\ 
& $\mathbf{\alpha }_{4}^{0}$ &  &  &  & $\mathbf{\alpha }_{4}^{0}$ &  &  & 
& $\mathbf{\alpha }_{4}^{0}$ &  &  \\ \hline
& 2 & -0.491 & -0.340 &  & 6 & 0.480 & -0.741 &  & 3 & -1.392 & -1.080 \\ 
& 12 & -3.722 & -4.060 &  & 3 & 0.436 & -1.157 &  & 6 & -3.934 & -3.296 \\ 
\hline \hline
& \multicolumn{11}{c}{RMSE} \\ \cline{2-12}
& \multicolumn{3}{c}{$\left( \chi _{6}^{2},\chi _{1}^{2}\right) $} &  & 
\multicolumn{3}{c}{$(MN,t_{6})$} &  & \multicolumn{3}{c}{$(t_{6},t_{5})$} \\ 
\cline{2-4}\cline{6-8}\cline{10-12}
$\mathbf{\theta }_{0}$ &  & Mixed & P.inv &  &  & Mixed & P.inv &  &  & Mixed
& P.inv \\ \hline
10 &  & 1.204 & 2.095 &  &  & 2.132 & 2.501 &  &  & 2.678 & 2.225 \\ 
-2 &  & 1.559 & 1.218 &  &  & 2.034 & 1.434 &  &  & 3.589 & 2.576 \\ 
4 &  & 1.675 & 2.752 &  &  & 3.829 & 2.882 &  &  & 4.952 & 5.020 \\ 
5 &  & 1.504 & 1.449 &  &  & 1.840 & 1.161 &  &  & 2.331 & 1.411 \\ 
-0.5/-2 &  & 0.138 & 0.430 &  &  & 0.257 & 0.406 &  &  & 0.657 & 0.389 \\ 
0 &  & 0.285 & 0.224 &  &  & 1.013 & 0.188 &  &  & 0.439 & 0.223 \\ 
0 &  & 0.281 & 0.959 &  &  & 0.580 & 0.980 &  &  & 0.473 & 0.759 \\ 
-2 &  & 0.522 & 0.498 &  &  & 4.381 & 0.452 &  &  & 0.789 & 0.417 \\ 
AVE &  & 0.896 & 1.203 &  &  & 2.008 & 1.250 &  &  & 1.989 & 1.628 \\ 
& $\mathbf{\alpha }_{3}^{0}$ &  &  &  & $\mathbf{\alpha }_{3}^{0}$ &  &  & 
& $\mathbf{\alpha }_{3}^{0}$ &  &  \\ \hline
& 1.155 & 0.414 & 0.415 &  & 2 & 0.764 & 1.178 &  & 0 & 0.718 & 0.711 \\ 
& 2.828 & 0.948 & 0.908 &  & 0 & 0.937 & 0.396 &  & 0 & 0.731 & 0.753 \\ 
& $\mathbf{\alpha }_{4}^{0}$ &  &  &  & $\mathbf{\alpha }_{4}^{0}$ &  &  & 
& $\mathbf{\alpha }_{4}^{0}$ &  &  \\ \hline
& 2 & 2.189 & 2.318 &  & 6 & 2.059 & 6.607 &  & 3 & 2.711 & 2.488 \\ 
& 12 & 7.606 & 7.409 &  & 3 & 3.186 & 3.365 &  & 6 & 4.657 & 4.051 \\ 
\hline \hline
\end{tabular}
\caption{\footnotesize 
\noindent  Estimation of SVARMA$\left( 0, 1 \right)$ models by GMM, $d= 2 ,$ $\zeta _{1}= \zeta _{2}= 0.5 ,$ $\rho_{1}= 1(\text{Mixed})/-1(\text{P.inv})$, $ \rho _{2}= -1 $. Innovations $\left( \chi _{6}^{2},\chi
_{1}^{2}\right), (MN(2.12,1.41^2;-0.24,0.58^2;0.1,0.9),t_6),$ and $(t_6,t_5)$.\    $\mathcal{L}_{2}+\mathcal{L}_{3}+\mathcal{L}_{4}$ identification and estimation, $w=(1,1)$. AVE is the average of the absolute value of the column.}
\end{table}
\bigskip

Summary results on estimation of higher order cumulants are severely
affected by some extreme replications, and, in general, confirm the usual
intuition that estimation of higher order moments is quite difficult for
small and moderate sample sizes and heavy tail distributions, even more
complicated for the kurtosis than for the skewness. Cumulant estimates are
biased towards zero, i.e. underestimate the non-Gaussianity, the bias
growing with the magnitude of the cumulant, and work similarly when based on
the unfeasible estimates $\mathbf{\hat{\theta}}_{2,T}.$

Table 4 reports the results for sample size $T=200$ and the same set-up, but
only for chi square innovations, confirming that models with non-invertible
roots tend to report more imprecise estimates and that estimation exploiting
higher order moments information in the invertible case can outperform
estimates using only second order moments.

The simulations on estimation of the SVARMA$\left( 1,1\right) $ mixed-root
model described in Table~5 for both $T=100$ and $200$ and $\chi ^{2}$
innovations using either $w=\left(1,0\right) $ and $\left( 1,1\right) $ to select different sets
of moments, confirm that kurtosis does not contribute much on top of
skewness in terms of bias and variability of estimates, in parallel with its
reduced identifying information. Cumulant estimation also becomes more
difficult with model complexity, though it improves substantially with the
sample size. \bigskip

\begin{table}[h!]
\centering
\renewcommand*{\arraystretch}{.75} 
\noindent \noindent \textbf{Table 4.} Estimation of $\left( \chi _{6}^{2},\chi _{1}^{2}\right) $ SVARMA$\left( 0,1\right) .$ $T=200$ 
\begin{tabular}{rrrrrrrrr}
\hline \hline
\multicolumn{9}{c}{Bias} \\ \hline
\multicolumn{3}{c}{Invertible} &  & \multicolumn{2}{c}{Mixed} &  & 
\multicolumn{2}{c}{Non-invert.} \\ \cline{1-3}\cline{5-6}\cline{8-9}
$\mathbf{\theta }_{0}$ & $\mathbf{\hat{\theta}}_{2,T}$ & $\mathbf{\hat{\theta%
}}_{w,T}$ &  & $\mathbf{\theta }_{0}$ & $\mathbf{\hat{\theta}}_{w,T}$ &  & $%
\mathbf{\theta }_{0}$ & $\mathbf{\hat{\theta}}_{w,T}$ \\ \hline
10 & -0.107 & -0.041 &  & 10 & -0.237 &  & 10 & 0.087 \\ 
-2 & 0.153 & 0.047 &  & -2 & -0.080 &  & -2 & -0.037 \\ 
4 & -0.324 & -0.076 &  & 4 & 0.007 &  & 4 & 0.095 \\ 
5 & -0.069 & -0.043 &  & 5 & -0.073 &  & 5 & -0.069 \\ 
-0.5 & -0.003 & -0.004 &  & -0.5 & -0.009 &  & -2 & 0.006 \\ 
0 & -0.001 & -0.001 &  & 0 & -0.021 &  & 0 & -0.007 \\ 
0 & -0.027 & -0.027 &  & 0 & 0.001 &  & 0 & -0.027 \\ 
-0.5 & 0.005 & 0.004 &  & -2 & -0.019 &  & -2 & -0.035 \\ 
AVE & 0.086 & 0.030 &  &  & 0.056 &  &  & 0.045 \\ 
$\mathbf{\alpha }_{3}^{0}$ & $\mathbf{\hat{\alpha}}_{32,T}$ & $\mathbf{\hat{%
\alpha}}_{3w,T}$ &  &  & $\mathbf{\hat{\alpha}}_{3w,T}$ &  &  & $\mathbf{%
\hat{\alpha}}_{3w,T}$ \\ \hline
1.155 & -0.054 & -0.051 &  &  & -0.010 &  &  & -0.035 \\ 
2.828 & -0.258 & -0.215 &  &  & -0.060 &  &  & -0.115 \\ 
$\mathbf{\alpha }_{4}^{0}$ & $\mathbf{\hat{\alpha}}_{42,T}$ & $\mathbf{\hat{%
\alpha}}_{4w,T}$ &  &  & $\mathbf{\hat{\alpha}}_{4w,T}$ &  &  & $\mathbf{%
\hat{\alpha}}_{4w,T}$ \\ \hline
2 & -0.210 & -0.275 &  &  & -0.213 &  &  & -0.267 \\ 
12 & -2.501 & -2.308 &  &  & -1.171 &  &  & -1.454 \\ \hline \hline
\multicolumn{9}{c}{RMSE} \\ \hline
\multicolumn{3}{c}{Invertible} &  & \multicolumn{2}{c}{Mixed} &  & 
\multicolumn{2}{c}{Non-invert.} \\ \cline{1-3}\cline{5-6}\cline{8-9}
$\mathbf{\theta }_{0}$ & $\mathbf{\hat{\theta}}_{2,T}$ & $\mathbf{\hat{\theta%
}}_{w,T}$ &  & $\mathbf{\theta }_{0}$ & $\mathbf{\hat{\theta}}_{w,T}$ &  & $%
\mathbf{\theta }_{0}$ & $\mathbf{\hat{\theta}}_{w,T}$ \\ \hline
10 & 0.984 & 0.741 &  & 10 & 0.737 &  & 10 & 1.700 \\ 
-2 & 0.789 & 0.314 &  & -2 & 0.773 &  & -2 & 0.895 \\ 
4 & 1.780 & 1.000 &  & 4 & 1.092 &  & 4 & 1.752 \\ 
5 & 0.742 & 0.704 &  & 5 & 1.067 &  & 5 & 1.078 \\ 
-0.5 & 0.064 & 0.065 &  & -0.5 & 0.080 &  & -2 & 0.312 \\ 
0 & 0.028 & 0.029 &  & 0 & 0.136 &  & 0 & 0.162 \\ 
0 & 0.121 & 0.121 &  & 0 & 0.157 &  & 0 & 0.591 \\ 
-0.5 & 0.063 & 0.064 &  & -2 & 0.297 &  & -2 & 0.304 \\ 
AVE & 0.571 & 0.380 &  &  & 0.542 &  &  & 0.849 \\ 
$\mathbf{\alpha }_{3}^{0}$ & $\mathbf{\hat{\alpha}}_{32,T}$ & $\mathbf{\hat{%
\alpha}}_{3w,T}$ &  &  & $\mathbf{\hat{\alpha}}_{3w,T}$ &  &  & $\mathbf{%
\hat{\alpha}}_{3w,T}$ \\ \hline
1.155 & 0.335 & 0.308 &  &  & 0.284 &  &  & 0.273 \\ 
2.828 & 0.729 & 0.679 &  &  & 0.827 &  &  & 0.814 \\ 
$\mathbf{\alpha }_{4}^{0}$ & $\mathbf{\hat{\alpha}}_{42,T}$ & $\mathbf{\hat{%
\alpha}}_{4w,T}$ &  &  & $\mathbf{\hat{\alpha}}_{4w,T}$ &  &  & $\mathbf{%
\hat{\alpha}}_{4w,T}$ \\ \hline
2 & 2.019 & 1.933 &  &  & 1.750 &  &  & 1.681 \\ 
12 & 6.961 & 6.978 &  &  & 9.409 &  &  & 9.148 \\ \hline \hline
\end{tabular}
\caption{\footnotesize 
\noindent   Estimation of SVARMA$\left( 0, 1 \right)$ models by GMM, $d= 2 ,$ $\zeta _{1}= \zeta _{2}= 0.5 ,$ $%
\rho_{1}= \pm 1 , \rho _{2}= \pm 1 $. Innovations $\left( \chi _{6}^{2},\chi
_{1}^{2}\right).$\   $\mathcal{L}_{2}+\mathcal{L}_{3}+\mathcal{L}_{4}$ identification and estimation, $w=(1,1)$. AVE is the average of the absolute value of the column.}
\end{table}

\begin{table}[h!]
\centering
\renewcommand*{\arraystretch}{.58} 
\noindent \textbf{Table 5.} Estimation of $\left( \chi
_{6}^{2},\chi _{1}^{2}\right) $-SVARMA$\left( 1,1\right) .$ Mixed Roots. 
\begin{tabular}{rrcrrrrcrr}
\hline \hline
& \multicolumn{9}{c}{Bias} \\ \hline
& \multicolumn{4}{c}{$T=100$} &  & \multicolumn{4}{c}{$T=200$} \\ 
\cline{2-5}\cline{7-10}
 $\mathbf{\theta }_{0}$&   & $\mathcal{L}_{3}$ &  & \multicolumn{1}{c}{$\mathcal{L}_{3}+\mathcal{L}%
_{4}$} &  &  & $\mathcal{L}_{3}$ &  & \multicolumn{1}{c}{$\mathcal{L}_{3}+%
\mathcal{L}_{4}$} \\ \hline
0.9 &  & \multicolumn{1}{r}{-0.045} &  & -0.002 &  &  & \multicolumn{1}{r}{
-0.026} &  & -0.027 \\ 
-0.4 &  & \multicolumn{1}{r}{0.155} &  & 0.162 &  &  & \multicolumn{1}{r}{
0.170} &  & 0.156 \\ 
0 &  & \multicolumn{1}{r}{0.028} &  & 0.039 &  &  & \multicolumn{1}{r}{-0.007
} &  & -0.007 \\ 
0.7 &  & \multicolumn{1}{r}{-0.041} &  & -0.052 &  &  & \multicolumn{1}{r}{
-0.011} &  & -0.017 \\ 
10 &  & \multicolumn{1}{r}{-1.100} &  & -1.178 &  &  & \multicolumn{1}{r}{
-0.384} &  & -0.700 \\ 
-2 &  & \multicolumn{1}{r}{0.200} &  & 0.159 &  &  & \multicolumn{1}{r}{0.051
} &  & 0.283 \\ 
4 &  & \multicolumn{1}{r}{-1.773} &  & -1.853 &  &  & \multicolumn{1}{r}{
-0.704} &  & -1.010 \\ 
5 &  & \multicolumn{1}{r}{0.401} &  & 0.823 &  &  & \multicolumn{1}{r}{-0.037
} &  & 0.224 \\ 
-0.5 &  & \multicolumn{1}{r}{-0.159} &  & -0.264 &  &  & \multicolumn{1}{r}{
-0.104} &  & -0.149 \\ 
0 &  & \multicolumn{1}{r}{-0.087} &  & -0.035 &  &  & \multicolumn{1}{r}{
-0.055} &  & -0.029 \\ 
0 &  & \multicolumn{1}{r}{-0.008} &  & 0.032 &  &  & \multicolumn{1}{r}{0.008
} &  & 0.044 \\ 
-2 &  & \multicolumn{1}{r}{-0.079} &  & 0.203 &  &  & \multicolumn{1}{r}{
-0.038} &  & -0.010 \\ 
AVE &  & \multicolumn{1}{r}{0.340} &  & 0.400 &  &  & \multicolumn{1}{r}{
0.133} &  & 0.221 \\ 
$\mathbf{\alpha }_{3}^{0}$ &  & \multicolumn{1}{r}{$\mathbf{\hat{\alpha}}%
_{3w,T}$} &  & $\mathbf{\hat{\alpha}}_{3w,T}$ &  &  & \multicolumn{1}{r}{$%
\mathbf{\hat{\alpha}}_{3w,T}$} &  & $\mathbf{\hat{\alpha}}_{3w,T}$ \\ \hline
1.155 &  & \multicolumn{1}{r}{-0.192} &  & -0.265 &  &  & \multicolumn{1}{r}{
-0.014} &  & -0.162 \\ 
2.828 &  & \multicolumn{1}{r}{-1.005} &  & -1.227 &  &  & \multicolumn{1}{r}{
-0.511} &  & -0.934 \\ 
$\mathbf{\alpha }_{4}^{0}$ &  & \multicolumn{1}{r}{$\mathbf{\hat{\alpha}}%
_{4w,T}$} &  & $\mathbf{\hat{\alpha}}_{4w,T}$ &  &  & \multicolumn{1}{r}{$%
\mathbf{\hat{\alpha}}_{4w,T}$} &  & $\mathbf{\hat{\alpha}}_{4w,T}$ \\ \hline
2 &  & \multicolumn{1}{r}{0.277} &  & 0.060 &  &  & \multicolumn{1}{r}{0.788}
&  & -0.239 \\ 
12 &  & \multicolumn{1}{r}{-10.415} &  & -7.540 &  &  & \multicolumn{1}{r}{
-2.112} &  & -6.346 \\ \hline \hline
& \multicolumn{9}{c}{RMSE} \\ \hline
& \multicolumn{4}{c}{$T=100$} &  & \multicolumn{4}{c}{$T=200$} \\ 
\cline{2-5}\cline{7-10}
 $\mathbf{\theta }_{0}$ & & $\mathcal{L}_{3}$ &  & \multicolumn{1}{c}{$\mathcal{L}_{3}+\mathcal{L}%
_{4}$} &  &  & $\mathcal{L}_{3}$ &  & \multicolumn{1}{c}{$\mathcal{L}_{3}+%
\mathcal{L}_{4}$} \\ \hline
0.9 &  & \multicolumn{1}{r}{0.117} &  & 0.420 &  &  & \multicolumn{1}{r}{
0.071} &  & 0.118 \\ 
-0.4 &  & \multicolumn{1}{r}{0.172} &  & 0.213 &  &  & \multicolumn{1}{r}{
0.174} &  & 0.165 \\ 
0 &  & \multicolumn{1}{r}{0.179} &  & 0.217 &  &  & \multicolumn{1}{r}{0.104}
&  & 0.161 \\ 
0.7 &  & \multicolumn{1}{r}{0.113} &  & 0.128 &  &  & \multicolumn{1}{r}{
0.054} &  & 0.109 \\ 
10 &  & \multicolumn{1}{r}{2.378} &  & 2.790 &  &  & \multicolumn{1}{r}{1.355
} &  & 1.644 \\ 
-2 &  & \multicolumn{1}{r}{2.942} &  & 2.986 &  &  & \multicolumn{1}{r}{2.398
} &  & 2.457 \\ 
4 &  & \multicolumn{1}{r}{4.549} &  & 4.816 &  &  & \multicolumn{1}{r}{2.744}
&  & 3.204 \\ 
5 &  & \multicolumn{1}{r}{2.210} &  & 2.462 &  &  & \multicolumn{1}{r}{1.534}
&  & 1.857 \\ 
-0.5 &  & \multicolumn{1}{r}{0.856} &  & 0.837 &  &  & \multicolumn{1}{r}{
0.259} &  & 0.361 \\ 
0 &  & \multicolumn{1}{r}{0.951} &  & 0.977 &  &  & \multicolumn{1}{r}{0.439}
&  & 0.412 \\ 
0 &  & \multicolumn{1}{r}{0.853} &  & 0.977 &  &  & \multicolumn{1}{r}{0.471}
&  & 0.382 \\ 
-2 &  & \multicolumn{1}{r}{1.071} &  & 1.475 &  &  & \multicolumn{1}{r}{0.542
} &  & 0.635 \\ 
AVE &  & \multicolumn{1}{r}{1.366} &  & 1.525 &  &  & \multicolumn{1}{r}{
0.845} &  & 0.959 \\ 
$\mathbf{\alpha }_{3}^{0}$ &  & \multicolumn{1}{r}{$\mathbf{\hat{\alpha}}%
_{3w,T}$} &  & $\mathbf{\hat{\alpha}}_{3w,T}$ &  &  & \multicolumn{1}{r}{$%
\mathbf{\hat{\alpha}}_{3w,T}$} &  & $\mathbf{\hat{\alpha}}_{3w,T}$ \\ \hline
1.155 &  & \multicolumn{1}{r}{0.809} &  & 0.936 &  &  & \multicolumn{1}{r}{
0.747} &  & 0.422 \\ 
2.828 &  & \multicolumn{1}{r}{1.838} &  & 1.799 &  &  & \multicolumn{1}{r}{
1.392} &  & 1.338 \\ 
$\mathbf{\alpha }_{4}^{0}$ &  & \multicolumn{1}{r}{$\mathbf{\hat{\alpha}}%
_{4w,T}$} &  & $\mathbf{\hat{\alpha}}_{4w,T}$ &  &  & \multicolumn{1}{r}{$%
\mathbf{\hat{\alpha}}_{4w,T}$} &  & $\mathbf{\hat{\alpha}}_{4w,T}$ \\ \hline
2 &  & \multicolumn{1}{r}{7.205} &  & 2.674 &  &  & \multicolumn{1}{r}{6.950}
&  & 2.007 \\ 
12 &  & \multicolumn{1}{r}{66.500} &  & 8.112 &  &  & \multicolumn{1}{r}{
9.998} &  & 6.859 \\ \hline
\% root loc. &  & \multicolumn{1}{r}{61.7} &  & 60.7 &  &  & 
\multicolumn{1}{r}{70.2} &  & 67.0 \\ \hline \hline
\end{tabular}
\caption{\footnotesize 
\noindent  Estimation of SVARMA$\left( 1, 1 \right)$ models by GMM, $d= 2 ,$ $\zeta _{1}= \zeta _{2}= 0.5 ,$ $%
\rho_{1}=  1 , \rho _{2}= -1 $. Innovations $\left( \chi _{6}^{2},\chi
_{1}^{2}\right).$\  $\mathcal{L}_{2}+\mathcal{L}_{3}$  and $\mathcal{L}_{2}+\mathcal{L}_{3}+\mathcal{L}_{4}$ identification and estimation, $w=(1,0)$ and $w=(1,1)$, respectively. Last line provides the percentage of simulated paths with correct MA roots location. AVE is the average of the absolute value of the column. $T=100,200$.}
\end{table}

\section{8. Empirical analysis}

We apply our identification and estimation methods to Blanchard and Quah
(1989) bivariate system for the US real GNP growth and unemployment rate
after linear detrending. They fit a SVAR with 8 lags (1948Q2-1987Q4, $T=159)$
and use for identification a long-run restriction by which the demand shock
has no long-run effect on real GNP in the same way as both supply and demand
shocks have no long-run effect on unemployment. Lippi and Reichlin (1994)
alternatively propose that these long VAR dynamics could be better
approximated by a VARMA model and explore the properties of the different
versions of the IRF obtained by inverting the MA roots of the fundamental
VARMA$\left( 1,1\right) $ representation deduced from the fitted VAR$\left(
8\right) $ parameterization. Gouri\'{e}roux et al. (2019) fit SVARMA$\left(
p,1\right) $ models to the same dataset for $p\in \left \{ 1,\ldots
,6\right
\} $ by PMLE using GMM initial estimates based on IV estimation of
the AR parameters. The shocks are assumed mixed Gaussian distributed and $%
p=4 $ is the order chosen by a combination of model selection criteria and
correlation diagnostics, finding a representation with a mixed
invertible/noninvertible MA roots.

We fit with our methods a series of simple SVARMA$\left( p,q\right) $ models
to the original dataset to investigate the possible nonfundamentalness of
the dynamics of the system using our higher order cumulant identification.
We follow the same MA polynomial parameterization and the same procedure as
in the Monte Carlo simulations and base the preliminary choice of the MA
roots configuration on minimization of $\mathcal{L}_{3,T}^\dag+\mathcal{L}%
_{4,T}^\dag$, $w=(1,1)$, among all the versions of the model obtained by MA
root flipping of Whittle estimates using Baggio and Ferrante (2019)
algorithm imposing causality. Then, a local of minimization of $\mathcal{L}%
_{w,T}^\dag=\mathcal{L}_{2,T}+\mathcal{L}_{3,T}^\dag+\mathcal{L}_{4,T}^\dag$
is performed to find $\mathbf{\hat{\theta}}_{w,T}^\dag$ and finally GMM
estimates $\mathbf{\hat{\theta}}_{GMM,T}$ are obtained with $\mathcal{V}$
estimated by bootstrap using the linear representation of the score. Each $%
\mathcal{L}^\dag _{k,T}$ is normalized by the number of spectral densities
of order $k$ in the system, $d^{k},$ and by $T^{k-2},$ so that their values
are close to one and can be compared easily across $k$ and model orders,
though the larger $k, $ the larger the relative contribution from the sample
variation of the $k $-order periodogram and the smaller the contribution of
the particular parameter value at which the loss function is evaluated.

We report in Table 6 the values of $\mathcal{L}_{2,T}$, $\mathcal{L}%
_{3,T}^\dag$, $\mathcal{L}_{4,T}^\dag$ and $\mathcal{L}_{w,T}^\dag$
evaluated at the final estimates with weighting defined by initial Whittle
estimates (and the value $\mathcal{L}_{2,T}^{0}$ obtained by minimizing only 
$\mathcal{L}_{2,T}$ as a benchmark for best linear fit), the modulus of the
MA roots of the estimated parameterization and the estimates of $\mathbf{%
\alpha} _{3}$ and $\mathbf{\alpha} _{4},$ using the scheme of Assumption~6A
for component identification, which imposes positive diagonal values of $%
\Theta _{0}$ and maximizes their product. We also report the estimates of $%
\Omega $ for the instantaneous impact of the shocks on the endogenous
variables. We use a similar identification strategy as Gouri\'{e}roux et al.
(2019) to label error components as transitory or demand and permanent or
supply shocks, facilitating an easy comparison to previous analysis.

The model which appears to best fit the data is the SVARMA$\left( 1,1\right)
, $ reporting the smallest value for the overall loss function $\mathcal{L}%
_{w,T}^\dag$ after joint optimization and an efficient GMM Newton-Raphson
step, just third best for $\mathcal{L}_{2,T},$ outperformed by the much
larger SVARMA$\left( 4,1\right) $ and VAR$\left( 8\right) $ models, which
did also a better job when fitting only second order dynamics attending to $%
\mathcal{L}_{2,T}^{0}$. Note that this comparison is made against the
initial Whittle estimates found only using second moments for the same order 
$\left( p,q\right) $, so additional local optimization of $\mathcal{L}%
_{3,T}^\dag$ and $\mathcal{L}_{4,T}^\dag$ to match higher order dynamics is
made by penalizing the second order goodness-of-fit at the same that imposes
ICA and the correct location of MA roots. For the same reason, and despite
the model nesting, larger models do not necessarily provide better fit in
finite samples attending to initial-estimates weighted $\mathcal{L}%
_{k,T}^\dag$ loss functions.

The SVARMA$\left( 1,1\right) $ model provides an invertible solution, but
some more complex models present non-invertible dynamics. Attending to
higher cumulants, invertible solutions indicate that both shocks have
moderate skewness, possibly of different signs. The SVARMA$\left( 1,1\right) 
$ identifies the first shock as the transitory with negative skewness, while
the second shock would be the permanent one with much larger positive
asymmetry, so in both cases negative news (decreasing GNP growth and
increasing unemployment) tend to be more extreme than positive ones.
Further, typically one shock displays large positive kurtosis (the permanent
one for the SVARMA$\left( 1,1\right) )$, and there is no conclusive evidence
about the kurtosis of the other one, given the large bootstrap standard
errors which make most estimates not significatively different from zero. At
least for the SVARMA$\left( 1,1\right) $ model, these estimation results
would confirm model dynamics identification using jointly third and fourth
cumulants by Corollary~2, while none of them in isolation would be
sufficient.

\begin{table}[t]
\renewcommand*{\arraystretch}{.8}  
\setlength{\tabcolsep}{2.8pt} 
\textbf{Table 6.} US real GNP growth
and unemployment rate. GMM SVARMA$(p,q)$ estimates. 
\begin{tabular}{lrrrrrrrrr}
\hline \hline
$(p,q)$ & $(1,0)\ $ & $(4,0)\ $ & $(8,0)\ $ & $(0,1)\ $ & $(1,1)\ $ & $%
(2,1)\ $ & $(4,1)\ $ & $\left( 1,2\right) \ $ & $\left( 2,2\right) \ $ \\ 
\hline \hline
$\mathcal{L}_{2,T}^{0}$ & $0.9608$ & $0.9991$ & $0.7694$ & $2.6118$ & $%
0.9136 $ & $0.9342$ & $0.8928$ & $0.9395$ & $0.9371$ \\ 
$\mathcal{L}_{2,T}$ & $0.9608$ & $1.0393$ & $0.8325$ & $2.6118$ & $0.9148$ & 
$0.9329$ & $0.8885$ & $0.9395$ & $0.9364$ \\ 
$\mathcal{\hat{L}}_{3,T}^\dag$ & $1.0258$ & $1.0345$ & $1.0343$ & $1.3832$ & $%
1.0284$ & $1.0349$ & $1.0335$ & $1.0285$ & $1.0297$ \\ 
$\mathcal{\hat{L}}_{4,T}^\dag$ & $1.1066$ & $1.0784$ & $1.2122$ & $1.8540$ & $%
1.0116$ & $1.0693$ & $1.1070$ & $1.0578$ & $1.0622$ \\ 
$\mathcal{\hat{L}}_{w,T}^\dag$ & $3.0932$ & $3.1522$ & $3.0790$ & $5.8490$ & $%
2.9548$ & $3.0370$ & $3.0290$ & $3.0258$ & $3.0284$ \\ \hline
$\left \vert \text{MA\ roots}\right \vert $ & $-$ & $-$ & $-$ & $3.1072$ & $%
8.7673$ & $0.7873$ & $1.0002$ & $0.25,0.51$ & $1.13,1.32$ \\ 
& $-$ & $-$ & $-$ & $3.1072$ & $8.7673$ & $5.2772$ & $1.0002$ & $0.51,2.81$
& $1.32,2.99$ \\ 
& $-$ & $-$ & $-$ & $(Inv.)$ & $(Inv.)$ & $(Mix.)$ & $(Inv.)$ & $(Mix.)$ & $%
(Inv.)$ \\ \hline
$\hat{\mathbf{\alpha }}_{3}$ & $-0.4474$ & $-0.5422$ & $-1.3776$ & $%
-0.4611 $ & $-0.5260 $ & $-0.3023^{{}}$ & $-0.2722$ & $%
-0.8550 $ & $-0.2959$ \\ 
& $1.4965$ & $1.5745$ & $-0.3435$ & $0.5114 $ & $1.2839$ & $1.4849$ & 
$1.9261$ & $-1.2217$ & $1.4230$ \\ 
$\hat{\mathbf{\alpha }}_{4}$ & $-1.6447$ & $0.1856$ & $6.3477 $ & $%
0.4471 $ & $-0.2134$ & $-1.0077$ & $-2.6953$ & $-0.0563$ & $-0.9281$
\\ 
& $5.4016$ & $8.1305 $ & $-1.8098$ & $1.1451 $ & $6.8220$ & $7.1990 $ & $-1.8966$ & $7.9476 $ & $7.0659 $ \\ 
\hline
$\hat{\Omega}_{1,1}$ & $\underset{\left( .103\right) }{0.8648}$ & $\underset{%
\left( .190\right) }{0.7465}$ & $\underset{\left( .067\right) }{0.4192}$ & $%
\underset{\left( .072\right) }{1.4697}$ & $\underset{\left( .088\right) }{%
0.8077}$ & $\underset{\left( .091\right) }{0.9184}$ & $\underset{\left(
.026\right) }{0.7116}$ & $\underset{\left( .305\right) }{0.2338}$ & $%
\underset{\left( .287\right) }{0.8178}$ \\ 
$\hat{\Omega}_{2,1}$ & $\underset{\left( .068\right) }{-0.2131}$ & $\underset%
{\left( .088\right) }{-0.1993}$ & $\underset{\left( .070\right) }{-0.0356}$
& $\underset{\left( .040\right) }{-0.6311}$ & $\underset{\left( .136\right) }%
{-0.2504}$ & $\underset{\left( .122\right) }{-0.2046}$ & $\underset{\left(
.027\right) }{-0.2568}$ & $\underset{\left( .194\right) }{-0.0364}$ & $%
\underset{\left( .173\right) }{-0.2042}$ \\ 
$\hat{\Omega}_{1,2}$ & $\underset{\left( .308\right) }{0.0886}$ & $\underset{%
\left( .202\right) }{0.1438}$ & $\underset{\left( .074\right) }{0.0379}$ & $%
\underset{\left( .194\right) }{0.3036}$ & $\underset{\left( .619\right) }{%
0.2445}$ & $\underset{\left( .187\right) }{-0.0575}$ & $\underset{\left(
.022\right) }{0.2805}$ & $\underset{\left( .394\right) }{-0.0704}$ & $%
\underset{\left( .261\right) }{0.1052}$ \\ 
$\hat{\Omega}_{2,2}$ & $\underset{\left( .111\right) }{0.1910}$ & $\underset{%
\left( .097\right) }{0.1688}$ & $\underset{\left( .073\right) }{0.0565}$ & $%
\underset{\left( .208\right) }{0.1961}$ & $\underset{\left( .140\right) }{%
0.1468}$ & $\underset{\left( .156\right) }{0.2301}$ & $\underset{\left(
.057\right) }{0.1094}$ & $\underset{\left( .170\right) }{0.0701}$ & $%
\underset{\left( .192\right) }{0.1825}$ \\ \hline \hline
\end{tabular}
\caption{\footnotesize  SVARMA$(p,q)$ model fitting for US real GNP growth and unemployment rate after linear detrending, 1948Q2-1987Q4, $T=159$. GMM estimates $\hat{\mathbf{\theta}}_{GMM,T}$ in (\ref{theestimator}) obtained with $\tilde{\mathbf{\theta}}_T= \hat{\mathbf{\theta}}_{w,T}$ with $w=(1,1)$ and $\hat{\mathbb{V}}_T^-$ obtained by inverting the sample variance of the joint score $\mathbb{S}_T$ evaluated at $\hat{\mathbf{\theta}}_{w,T}$ for 400 draws on the empirical distribution of residuals $\hat{\mathbf{\varepsilon}}_t$. Component identification uses Assumption 6A.}
\end{table}

The estimates of the lag zero impact matrix $\Omega $ of our preferred
model, as for other models with only invertible dynamics in Table 6,
indicate a much higher impact of the first (transitory) shock on output
compared to the impacts of the second (permanent) shock on both output and
employment, while the transitory shock has a negative short run effect on
unemployment. This last feature is preserved in models which identify
non-invertible components, but the sign of the instantaneous impact of the
permanent shock on output is reversed.

Finally, we provide in Figure 1 the plots of the IRF's identified by our
estimation methods for the SVARMA $\left( 1,1\right) $ model together with
those reported by Blanchard and Quah (1989) and Gouri\'{e}roux et al.
(2019). We can observe that the effect of the supply (or permanent) shock on
both endogenous variables is very close to the results of Blanchard and Quah
for the SVAR(8) model, though there is a slight delay in the maximum effect
and a quite more persistent effect on output. For the demand (transitory)
shock, the shape of both IRF's and the timing of maximal effects are almost
the same as for the SVAR model, but now the long run effect on output is not
restricted to zero and the IRF displays a very slow rate of decay but close
to the horizontal axis. The bootstrap tests of the long-run identification
restriction used by Blanchard and Quah (1989), based on the distance at lags
40 and 100 of the IRF of GNP growth on this transitory shock, can not reject
the hypothesis of zero long run effect, confirming the plausibility of this
identification strategy. On the other hand, the IRF obtained by Gouri\'{e}%
roux et al. (2019) with noninvertible dynamics retain some of the previous
properties for the unemployment response, but GNP growth behaves quite
differently in the long run after either type of structural shocks hitting
the system. \bigskip

\begin{center}
{\normalsize \includegraphics[width=1.05%
\textwidth]{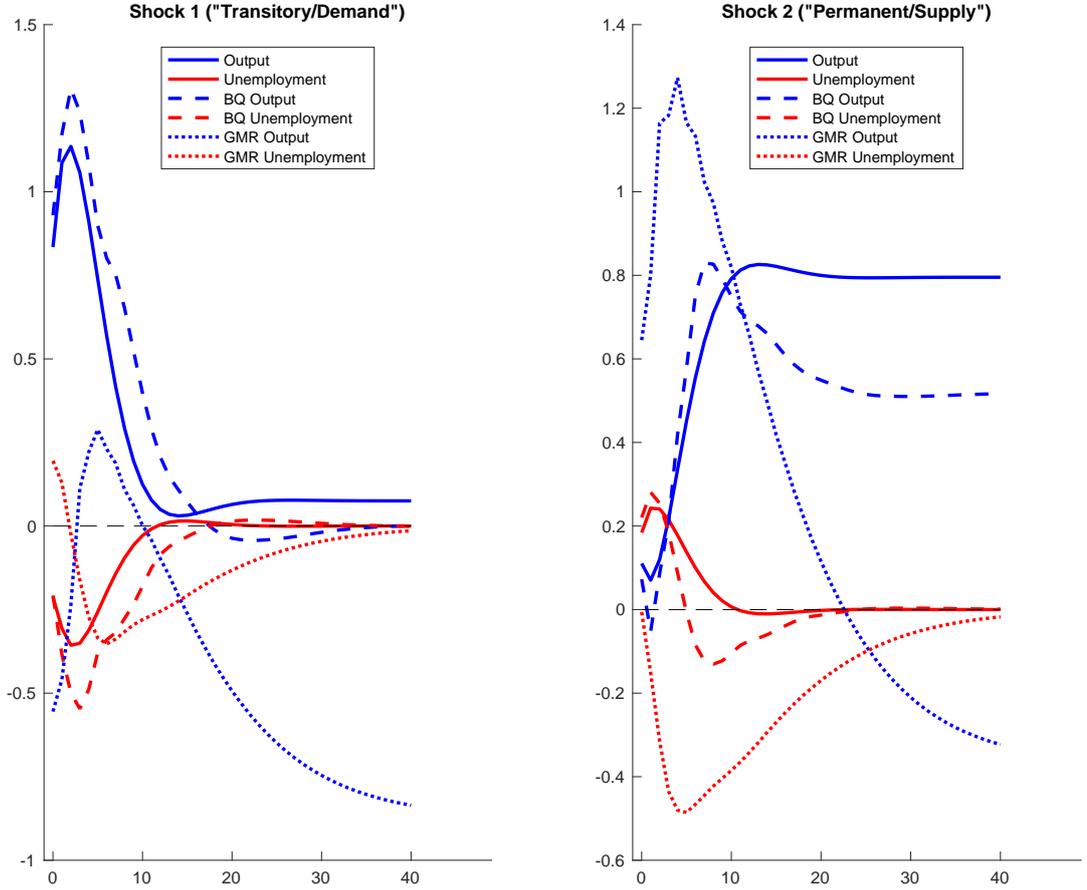} 
\captionof{figure}{IRF for US GNP
growth and unemployment based on our GMM estimates (\ref{theestimator}) for the SVARMA$\left(1,1\right) $ model, BQ: Blanchard and Quah (1989) SVAR$\left( 8\right) $ and
GMR: Gouri\'{e}roux et al. (2019) SVARMA$\left( 4,1\right) $ model.} }
\end{center}

\section{9. Conclusions}

In this paper we have showed how to achieve identification of non-Gaussian
SVARMA models using basic restrictions on higher order moments under serial
and component independence conditions of finite order $k=3$ or 4 on the
structural shocks sequence. We use an identification criterion in the
frequency domain that leads to easy to check global and local identification
conditions and permits the design of consistent and asymptotically normal
parameter estimates which exploit all dynamic and static information in
second, third and fourth order moments. These results provide consistent
estimation of IRFs without need to specify the fundamentalness of the system
and can be combined with different sources of information to proper label
the structural shocks of the model or to test relevant hypothesis and
overidentification conditions.\bigskip

\section{\protect \small Appendix A: Higher order cumulants and spectral
densities}

{\small \noindent \textbf{Cumulants of linear combinations. }The third order
cumulant matrices $\mathbf{\kappa }_{\cdot \cdot j}$,\ $j=1,\ldots ,d,$ of a
random vector $\mathbf{\varepsilon }_{t}$ can be updated easily under linear
transformations for a $d\times d$ matrix $K,$ where $K$ is orthogonal, $%
KK^{\prime }=\mathbf{I}_{d}, $ to maintain the covariance structure of $%
\mathbf{\varepsilon }_{t},$%
\begin{equation*}
\eta _{t}=K \mathbf{\varepsilon} _{t}
\end{equation*}%
so that it holds%
\begin{equation*}
\text{vec}\left( \text{v}\mathbf{\kappa }_{3}^{\eta }\right) =K_{1}^{\otimes
3}\text{vec}\left( \text{v}\mathbf{\kappa }_{3}^{0}\right)
\end{equation*}%
or, alternatively, 
\begin{equation*}
\text{v}\mathbf{\kappa }_{3}^{\eta }=K^{\otimes 2}\  \text{v}\mathbf{\kappa }%
_{3}^{0}\ K^{\prime },
\end{equation*}%
and in general, for any $k\geq 2$, see Jammalamadaka, Rao and Terdik (2006), 
\begin{equation*}
\text{vec}\left( \text{v}\mathbf{\kappa }_{k}^{\eta }\right) =K^{\otimes k}%
\text{vec}\left( \text{v}\mathbf{\kappa }_{k}^{0}\right) =\text{vec}\left(
K^{\otimes 2}\  \text{v}\mathbf{\kappa }_{k}^{0}\  \left( K^{\otimes
(k-2)}\right) ^{\prime }\right) ,
\end{equation*}%
with $K^{\otimes 0}=\mathbf{I}_{d}$ using that vec$\left( ABC\right) =\left(
C^{\prime }\otimes A\right) $vec$\left( B\right) $ for conformable matrices $%
A,B,C.$ Note however that, for $k\geq 4,$ v$\mathbf{\kappa }_{k}^{0}$
contains many repeated columns, since for instance, for $k=4,$ vec$\left( 
\mathbf{\kappa }_{\cdot \cdot jh}\right) =\ $vec$\left( \mathbf{\kappa }%
_{\cdot \cdot hj}\right) ,$ so in this particular case rank$\left( \text{v}%
\mathbf{\kappa }_{k}^{0}\right) \leq d(d+1)/2.\ $\newline
}

{\small \noindent \textbf{Higher order spectral densities.} Consider a $d$%
-dimensional stationary stochastic process $\{Y_{t}\}_{t\in Z}$ with $%
E[\left \Vert Y_{t}\right \Vert ^{k}]<\infty ,$ some $k\geq 3,$ and call $%
\mu =E[Y_{t}].$ Define the autocovariance of order $j$ as 
\begin{equation*}
\Gamma \left( j\right) =Cov[Y_{t},Y_{t-j}]=E[(Y_{t}-\mu )(Y_{t-j}-\mu
)^{\prime }],\qquad \text{for }j=0,\pm 1,\ldots ,
\end{equation*}%
and the spectral density matrix, $f(\lambda ),$ is defined implicitly as%
\begin{equation*}
\Gamma \left( j\right) =\int_{-\pi }^{\pi }f(\lambda )\exp (-ij\lambda
)d\lambda .
\end{equation*}%
The autocovariance sequence and the spectral density are measures of the
dependence of the stochastic process based on second moments, hence they are
the objects of interest of usual time series analysis. The dependence
contained in higher order moments can also be described by the cumulants
which are defined in terms of higher order moments as%
\begin{equation*}
\text{cum}\left( Y_{t\left( 1\right) ,\mathbf{a}\left( 1\right) },\ldots
,Y_{t\left( k\right) ,\mathbf{a}\left( k\right) }\right) =\sum
(-1)^{p-1}(p-1)!E(\Pi _{j\in v_{1}}Y_{t\left( j\right) ,\mathbf{a}\left(
j\right) })\cdots E(\Pi _{j\in v_{p}}Y_{t\left( j\right) ,\mathbf{a}\left(
j\right) }),\  \ k=1,2,\ldots ,
\end{equation*}%
where $v_{1},\ldots ,v_{p}$ is a partition of $(1,2,\ldots ,k),$ and the sum
runs over all these partitions, $\mathbf{a}\left( j\right) \in \left \{
1,\ldots ,d\right \} ,$ $t\left( j\right) =0,\pm 1,\ldots ,$ for $j=1.\ldots
,k,$ see Brillinger (1975) or Rosenblatt (1985, p.~34). Hence, the first and
second cumulants are the mean and the variance, respectively. }

{\small We also define the $k$-th order cumulant spectral density $%
k=2,3,\ldots ,$ for \textbf{$a$}$=\left( \mathbf{a}\left( 1\right) ,\ldots ,%
\mathbf{a}\left( k\right) \right) $ which is the Fourier transform of the $k$%
-th order cumulants for elements $\left( \mathbf{a}\left( 1\right) ,\ldots ,%
\mathbf{a}\left( k\right) \right) $ of the vector $Y_{t},$ $f_{\mathbf{a},k}(%
\boldsymbol{\lambda })=f_{\mathbf{a}\left( 1\right) ,\ldots ,\mathbf{a}%
\left( k\right) }(\lambda _{1},\ldots ,\lambda _{k-1}),$ as 
\begin{equation}
f_{\mathbf{a},k}(\boldsymbol{\lambda })=\left( 2\pi \right)
^{1-k}\sum_{j_{1},\ldots ,j_{k-1}=-\infty }^{\infty }\text{cum}(Y_{\mathbf{a}%
(1),t},Y_{\mathbf{a}(2),t+j_{1}},\ldots ,Y_{\mathbf{a}(k),t+j_{k-1}})\exp
\left( -\sum_{s=1}^{k-1}ij_{s}\lambda _{s}\right) ,  \notag
\end{equation}%
introducing for simplicity, when there is no ambiguity, the notation $%
\boldsymbol{\lambda }=(\lambda _{1},\ldots ,\lambda _{k-1})$. Existence of $%
f_{\mathbf{a},k}(\boldsymbol{\lambda })$ can be guaranteed by summability
conditions on cumulants, as those implied by a linear process condition.
Note that the elements of the usual spectral density (matrix) are recovered
for $k=2$, while, in general, $f_{\mathbf{a},k}$ can be complex valued for
any $k\geq 2$, except the diagonal elements of the usual spectral density
matrix $f=f_{2}$, i.e. $f_{a,a}\left( \lambda \right) ,$ $a=1,\ldots ,d,$
which are always real valued.$\ $\newline
}

{\small \noindent \textbf{VARMA higher order spectral densities}. The second
order spectral density matrix of $Y_{t}$ is given by%
\begin{eqnarray*}
f\left( \lambda \right) &=&\frac{1}{2\pi }\Phi ^{-1}\left( e^{-i\lambda
}\right) \Theta \left( e^{-i\lambda }\right) \left( \Phi ^{-1}\left(
e^{-i\lambda }\right) \Theta \left( e^{-i\lambda }\right) \right) ^{\ast } \\
&=&\frac{1}{2\pi }\Psi \left( e^{-i\lambda }\right) \Psi ^{\ast }\left(
e^{i\lambda }\right) ,
\end{eqnarray*}%
where $^{\ast }$ denotes complex conjugation and transposition, and in
particular for any two components \textbf{$a$}$\mathbf{=}\left( \mathbf{a}%
\left( 1\right) ,\mathbf{a}\left( 2\right) \right) $ of $Y_{t},$ its
spectral density $f_{\mathbf{a}}\left( \lambda \right) =f_{\left( a\left(
1\right) ,a\left( 2\right) \right) }\left( \lambda \right) $ satisfies 
\begin{eqnarray*}
f_{\mathbf{a}}\left( \lambda \right) &=&\frac{1}{2\pi }\Psi _{\mathbf{a}%
\left( 1\right) }\left( e^{-i\lambda }\right) \Psi _{\mathbf{a}\left(
2\right) }^{\prime }\left( e^{i\lambda }\right) \\
&=&\frac{1}{2\pi }\sum_{h,j=1}^{d}\Psi _{\mathbf{a}\left( 1\right) ,h}\left(
e^{-i\lambda }\right) \Psi _{\mathbf{a}\left( 2\right) ,j}\left( e^{i\lambda
}\right)
\end{eqnarray*}%
where $\Psi _{j}\left( z\right) =\left( \Psi _{j,1}\left( z\right) ,\ldots
,\Psi _{j,d}\left( z\right) \right) $ is the $j$-th row of $\Psi \left(
z\right) ,$ because the $j$-th element of the vector $Y_{t}$ is obtained as $%
Y_{t,j}=\Psi _{j}\left( L\right) \mathbf{\varepsilon }_{t}.$ }

{\small For any triplet \textbf{$a$}$\mathbf{=}\left( \mathbf{a}\left(
1\right) ,\mathbf{a}\left( 2\right) ,\mathbf{a}\left( 3\right) \right) $,
the third order spectral density of the 3-dimensional vector $\left( Y_{t,%
\mathbf{a}(1)},Y_{t,\mathbf{a}(2)},Y_{t,\mathbf{a}(3)}\right) ,$ $f_{\mathbf{%
a},3}(\boldsymbol{\lambda })=f_{\left( \mathbf{a}\left( 1\right) ,\mathbf{a}%
\left( 2\right) ,\mathbf{a}\left( 3\right) \right) }(\lambda _{1},\lambda
_{2}),$ is given by%
\begin{equation*}
f_{\mathbf{a},3}(\boldsymbol{\lambda })=\left( 2\pi \right) ^{-2}\sum_{%
\mathbf{j}=1}^{d}\Psi _{\mathbf{a}\left( 1\right) ,j(1)}\left( e^{-i\lambda
_{1}}\right) \Psi _{\mathbf{a}\left( 2\right) ,j(2)}\left( e^{-i\lambda
_{2}}\right) \Psi _{\mathbf{a}\left( 3\right) ,j(3)}\left( e^{i\left(
\lambda _{1}+\lambda _{2}\right) }\right) \mathbf{\kappa}
_{j(1),j(2),j(3)}^{0}
\end{equation*}%
where $\mathbf{\kappa }_{\mathbf{j}}^{0}=\ $cum$\left( \mathbf{\varepsilon }%
_{t,j(1)},\mathbf{\varepsilon }_{t,j(2)},\mathbf{\varepsilon }%
_{t,j(3)}\right) $, $\mathbf{j}=\left( j(1),j(2),j(3)\right) ,$ is the joint
third order cumulants of the innovations indexed by \textbf{$a$}. This third
order spectral density can be we written more compactly as%
\begin{equation*}
f_{\mathbf{a},3}(\boldsymbol{\lambda })=\left( 2\pi \right) ^{-2}\Psi _{%
\mathbf{a}}^{\otimes 3}\left( \mathbf{\lambda }\right) \text{vec}\left( 
\text{v}\mathbf{\kappa }_{3}^{0}\right) ,
\end{equation*}%
where, using Kronecker product $\otimes ,$%
\begin{equation*}
\Psi _{\mathbf{a}}^{\otimes 3}\left( \mathbf{\lambda }\right) =\Psi _{%
\mathbf{a}\left( 3\right) }\left( e^{i\left( \lambda _{1}+\lambda
_{2}\right) }\right) \otimes \Psi _{\mathbf{a}\left( 2\right) }\left(
e^{-i\lambda _{2}}\right) \otimes \Psi _{\mathbf{a}\left( 1\right) }\left(
e^{-i\lambda _{1}}\right)
\end{equation*}%
and, denoting $\mathbf{\kappa }_{\cdot \cdot j}=E\left[ \mathbf{\varepsilon }%
_{t}\mathbf{\varepsilon }_{t}^{\prime }\mathbf{\varepsilon }_{t,j}\right] ,\
j=1,\ldots ,d,$%
\begin{equation*}
\text{v}\mathbf{\kappa }_{3}^{0}=\left[ \text{vec}\left( \mathbf{\kappa }%
_{\cdot \cdot 1}\right) \  \  \cdots \  \  \text{vec}\left( \mathbf{\kappa }%
_{\cdot \cdot d}\right) \right]
\end{equation*}%
is a $d^{2}\times d$ matrix, where vec is the usual operator stacking all
the columns of a matrix in a single column vector, and the indexes in the $%
d^{3}$-dimensional vector vec$\left( \text{v}\mathbf{\kappa }_{3}^{0}\right)
=\left \{ \mathbf{\kappa }_{hj\ell }\right \} $ run first from left to right
from 1 to $d.$ }

{\small This notation extends readily to any $k$-order spectral density, $%
k=2,3,\ldots ,$%
\begin{equation*}
f_{\mathbf{a},k}(\boldsymbol{\lambda })=\frac{1}{\left( 2\pi \right) ^{k-1}}%
\Psi _{\mathbf{a}}^{\otimes k}\left( \mathbf{\lambda }\right) \text{vec}%
\left( \text{v}\mathbf{\kappa }_{k}^{0}\right) ,
\end{equation*}%
where for $\boldsymbol{\lambda }=(\lambda _{1},\ldots ,\lambda _{k-1})$ we
define%
\begin{equation*}
\Psi _{\mathbf{a}}^{\otimes k}\left( \mathbf{\lambda }\right) =\Psi _{%
\mathbf{a}\left( k\right) }\left( e^{i\left( \lambda _{1}+\cdots +\lambda
_{k-1}\right) }\right) \otimes \Psi _{\mathbf{a}\left( k-1\right) }\left(
e^{-i\lambda _{k-1}}\right) \otimes \cdots \otimes \Psi _{\mathbf{a}\left(
2\right) }\left( e^{-i\lambda _{2}}\right) \otimes \Psi _{\mathbf{a}\left(
1\right) }\left( e^{-i\lambda _{1}}\right)
\end{equation*}%
and the $d^{2}\times d^{k-2}$ matrix v$\mathbf{\kappa }_{k}^{0}$ satisfies%
\begin{equation*}
\text{v}\mathbf{\kappa }_{k}^{0}=\left[ \text{vec}\left( \mathbf{\kappa }%
_{\cdot \cdot 1\cdots 1}\right) \  \  \text{vec}\left( \mathbf{\kappa }_{\cdot
\cdot 2\cdots 1}\right) \  \  \cdots \  \  \text{vec}\left( \mathbf{\kappa }%
_{\cdot \cdot d\cdots d}\right) \right]
\end{equation*}%
and $\mathbf{\kappa }_{\cdot \cdot j\left( 3\right) \cdots j\left( k\right)
} $ is a $d\times d$ matrix with typical $\left( j\left( 1\right) ,j\left(
2\right) \right) $ element equal to the $k$-th order joint cumulant cum$%
\left( \mathbf{\varepsilon }_{t,j\left( 1\right) },\mathbf{\varepsilon }%
_{t,j\left( 2\right) },\mathbf{\varepsilon }_{t,j\left( 3\right) },\ldots ,%
\mathbf{\varepsilon }_{t,j\left( k\right) }\right) ,$ $j\left( h\right) \in
\left \{ 1,\ldots ,d\right \} ,$ $h=1,\ldots ,k.$ \bigskip }

\section{\protect \small Appendix B: Proofs of Results}

{\small \noindent \textbf{Proof of Theorem~\ref{1A}. } }

{\small \noindent \textbf{Proof for }$k=3.$\textbf{\ }The integrand of $%
\mathcal{L}_{3}^{0}\left( A,\text{v}\mathbf{\kappa }_{3}\right) ,$
satisfies, uniformly for $\mathbf{\lambda }\in \Pi ^{2},$%
\begin{eqnarray*}
&&\left \{ \text{vec}\! \left( \text{v}\mathbf{\kappa }_{3}\right) ^{\prime
}A^{\otimes 3}\! \left( \mathbf{\lambda }\right) ^{\ast }-\text{vec}\!
\left( \text{v}\mathbf{\kappa }_{3}^{0}\right) ^{\prime }\right \} \Upsilon
_{3}^{0}\left( \mathbf{I}_{d},\mathbf{\lambda }\right) \left \{ A^{\otimes
3}\left( \mathbf{\lambda }\right) \text{vec}\left( \text{v}\mathbf{\kappa }%
_{3}\right) -\text{vec}\left( \text{v}\mathbf{\kappa }_{3}^{0}\right) \right
\} \\
&\geq &\inf_{\mathbf{\lambda }}\lambda _{\min }\left( \Upsilon
_{3}^{0}\left( \mathbf{I}_{d},\boldsymbol{\lambda }\right) \right)
\inf_{\lambda }\left \vert \lambda _{\min }\left( A\left( \lambda \right)
\otimes \mathbf{I}_{d^{2}}\right) \right \vert ^{2}\left \Vert \text{vec}%
\left( \left( A^{\otimes 2}\left( \lambda _{2},\lambda _{1}\right) \text{v}%
\mathbf{\kappa }_{3}-\text{v}\mathbf{\kappa }_{3}^{0}A\left( \lambda
_{1}+\lambda _{2}\right) \right) \right) \right \Vert ^{2} \\
&\geq &\eta ^{6}\left \Vert \text{vec}\left( A^{\otimes 2}\left( \lambda
_{2},\lambda _{1}\right) \text{v}\mathbf{\kappa }_{3}-\text{v}\mathbf{\kappa 
}_{3}^{0}A\left( \lambda _{1}+\lambda _{2}\right) \right) \right \Vert ^{2}
\end{eqnarray*}%
because of Assumption~4, $A\left( z\right) $ is a BM (with unitary
eigenvalues for $\left \vert z\right \vert =1)$ and 
\begin{eqnarray*}
A^{\otimes 3}\left( \mathbf{\lambda }\right) \text{vec}\left( \text{v}%
\mathbf{\kappa }_{3}\right) -\text{vec}\left( \text{v}\mathbf{\kappa }%
_{3}^{0}\right) &=&\text{vec}\left( A^{\otimes 2}\left( \lambda _{2},\lambda
_{1}\right) \text{v}\mathbf{\kappa }_{3}A\left( -\lambda _{1}-\lambda
_{2}\right) ^{\prime }-\text{v}\mathbf{\kappa }_{3}^{0}\right) \\
&=&\text{vec}\left( \left( A^{\otimes 2}\left( \lambda _{2},\lambda
_{1}\right) \text{v}\mathbf{\kappa }_{3}-\text{v}\mathbf{\kappa }%
_{3}^{0}A^{\ast }\left( \lambda _{1}+\lambda _{2}\right) \right) A\left(
-\lambda _{1}-\lambda _{2}\right) ^{\prime }\right) \\
&=&\left( A\left( -\lambda _{1}-\lambda _{2}\right) \otimes \mathbf{I}%
_{d^{2}}\right) \text{vec}\left( A^{\otimes 2}\left( \lambda _{2},\lambda
_{1}\right) \text{v}\mathbf{\kappa }_{3}-\text{v}\mathbf{\kappa }%
_{3}^{0}A^{\ast }\left( \lambda _{1}+\lambda _{2}\right) \right)
\end{eqnarray*}%
denoting $A^{\otimes 2}\left( \lambda _{2},\lambda _{1}\right) =A\left(
e^{-i\lambda _{2}}\right) \otimes A\left( e^{-i\lambda _{1}}\right) .$ }

{\small Then, under the assumption that $A\left( z\right) $ is a
non-constant BM (even if $a_{n}={a}_{m}^{\star-1}$ for some $n\neq m$ and
therefore $g_{a_{n}}\left( z\right) g_{a_{m}}\left( z\right) =1$), and
because v$\mathbf{\kappa }_{3}^{0}$ is full rank $d,$ there is at least one
pair $\left( j,h\right) ,$ $j\in \left \{ 1,\ldots ,d^{2}\right \} ,$ $h\in
\left
\{ 1,\ldots ,d\right \} ,$ such that $\xi _{3}\left( \mathbf{\lambda }%
\right) :=\left \{ \text{v}\mathbf{\kappa }_{3}^{0}A^{\ast }\left( \lambda
_{1}+\lambda _{2}\right) \right \} _{j,h}$ depends on $\mathbf{\lambda }$
through $\lambda _{1}+\lambda _{2}$ in a linear combination of products of $%
g_{a_{i}}\left( \lambda _{1}+\lambda _{2}\right) $ functions, $i=1,\ldots
,r. $\footnote{%
Note that in the representation $\left( \ref{BM}\right) $ of a GBM we can
allow for $a_{n}={a}_{m}^{\star-1}$ for some $n\neq m$, despite this implies
that $g_{a}\left( z\right) g_{{a}^{\star-1}}\left( z\right) =1$ for a real $%
a. $ Therefore, if $A\left( z\right) $ is assumed not constant in $z$, a
full cancellation of all roots is not allowed and the effect of $A$ is
shifting these flipped roots in different components of $\Psi $ when $d>1$,
avoiding cases like $A\left( z\right) =K_{0}R\left( {a}_{n}^{\star-1},z%
\right) R\left( a_{n},z\right) K_{0}^{\prime }=\mathbf{I}_{d},$ which leave
unchanged the dynamics as any constant $A$.} 

Then, for such pair $\left( j,h\right) $ depending on the form of $A$
and v$\mathbf{\kappa }_{3}^{0},$ 
\begin{equation*}
\left \Vert \text{vec}\left( A^{\otimes 2}\left( \lambda _{2},\lambda
_{1}\right) \text{v}\mathbf{\kappa }_{3}-\text{v}\mathbf{\kappa }%
_{3}^{0}A^{\ast }\left( \lambda _{1}+\lambda _{2}\right) \right) \right
\Vert \geq \left \vert \zeta _{3}\left( \mathbf{\lambda }\right) -\xi
_{3}\left( \mathbf{\lambda }\right) \right \vert ,
\end{equation*}%
where $\zeta _{3}\left( \mathbf{\lambda }\right) :=\left \{ A^{\otimes
2}\left( \lambda _{2},\lambda _{1}\right) \text{v}\mathbf{\kappa }%
_{3}\right
\} _{j,h}$ is a linear combination of products of 
\begin{equation*}
1,g_{a_{i}}\left( \lambda _{1}\right) ,g_{a_{i}}\left( \lambda _{2}\right)
,\  \  \ i=1,\ldots ,r,
\end{equation*}%
possibly zero or constant even if v$\mathbf{\kappa }_{3}$ is full rank.
Therefore%
\begin{equation*}
\mathcal{L}_{3}^{0}\left( A,\text{v}\mathbf{\kappa }_{3}\right) \geq \eta
^{6}\int_{\Pi ^{2}}\left \vert \zeta _{3}\left( \mathbf{\lambda }\right)
-\xi _{3}\left( \mathbf{\lambda }\right) \right \vert ^{2}d\mathbf{\lambda }%
\geq \epsilon >0
\end{equation*}%
for some $\epsilon $ not depending on v$\mathbf{\kappa }_{3},$ but depending
on $\eta ,$ v$\mathbf{\kappa }_{3}^{0}$ and $A$, because the functions $%
\zeta _{3}\left( \mathbf{\lambda }\right) \ $and $\xi _{3}\left( \mathbf{%
\lambda }\right) $ differ a.e. for any choice of v$\mathbf{\kappa }_{3}$
because $\xi _{3}\left( \mathbf{\lambda }\right) $ has an infinite expansion
on powers of $\exp \left( i\left( \lambda _{1}+\lambda _{2}\right) \right) $
which depends on $\mathbf{\lambda }$ only through $\lambda _{1}+\lambda _{2}$
(and can not be factorized in separated functions of $\lambda _{1}\ $and $%
\lambda _{2}$) while $\zeta _{3}\left( \mathbf{\lambda }\right) $ depends on 
$\mathbf{\lambda }$ only through products of functions with infinite
representation on exponential functions of a single $\lambda _{i},\ i=1,2.$

Note that for $r=1$ it holds that $A\left( z\right) =K_{0}R\left(
a_{1},z\right) K_{1},$ $\left \vert a_{1}\right \vert \neq 1,$ $a_{1}$ real,
and we can chose $\left( j,h\right) $ such that for $\widetilde{\text{v}%
\mathbf{\kappa }}_{3}^{0}=\ $v$\mathbf{\kappa }_{3}^{0}K_{0}$ and for some $%
c_{i}=c_{i}\left( \text{v}\mathbf{\kappa }_{3}^{0},K_{0},K_{1}\right) \in 
\mathbb{R},$ $i=1,2,$ with $c_{1}\neq 0,$%
\begin{equation*}
\xi _{3}\left( \mathbf{\lambda }\right) =\left \{ \text{v}\mathbf{\kappa }%
_{3}^{0}A\left( \lambda _{1}+\lambda _{2}\right) \right \} _{j,h}=\left \{ 
\widetilde{\text{v}\mathbf{\kappa }}_{3}^{0}R\left( a_{1},z\right)
K_{1}\right \} _{j,h}=c_{0}+c_{1}g_{a_{1}}\left( \lambda _{1}+\lambda
_{2}\right)
\end{equation*}%
and for some constants $d_{i,\ell }=d_{i,\ell }\left( \text{v}\mathbf{\kappa 
}_{3},K_{0},K_{1}\right) \in \mathbb{R},$%
\begin{equation*}
\zeta _{3}\left( \mathbf{\lambda }\right) =\left(
d_{1,0}+d_{1,1}g_{a_{1}}\left( \lambda _{1}\right) \right) \left(
d_{2,0}+d_{2,1}g_{a_{1}}\left( \lambda _{2}\right) \right) ,
\end{equation*}%
so that for all $c_{0},c_{1}\neq 0$ and $d_{i,\ell },$ the function $\zeta
_{3}\left( \mathbf{\lambda }\right) -\xi _{3}\left( \mathbf{\lambda }\right)
\neq 0$ a.e. for any $a_{1},$ $\left \vert a_{1}\right \vert \neq 1.$%

In particular, if $r=d=1,$ then $\left( j,h\right) =\left(
1,1\right) $ and $\xi _{3}\left( \mathbf{\lambda }\right) =\kappa
_{3}^{0}g_{a_{1}}\left( \lambda _{1}+\lambda _{2}\right) ,$ $\zeta
_{3}\left( \mathbf{\lambda }\right) =\kappa _{3}g_{a_{1}}\left( \lambda
_{1}\right) g_{a_{1}}\left( \lambda _{2}\right) ,$ with $\kappa _{3}^{0}\neq
0,$ so that 
\begin{eqnarray*}
\int_{\Pi ^{2}}\left \vert \zeta _{3}\left( \mathbf{\lambda }\right) -\xi
_{3}\left( \mathbf{\lambda }\right) \right \vert ^{2}d\mathbf{\lambda } &%
\mathbf{=}&\int_{\Pi ^{2}}\left \vert g_{a_{1}}\left( \lambda _{1}+\lambda
_{2}\right) \right \vert ^{2}\left \vert \kappa _{3}g_{a_{1}}\left( \lambda
_{1}\right) g_{a_{1}}\left( \lambda _{2}\right) g_{a_{1}}\left( -\lambda
_{1}-\lambda _{2}\right) -\kappa _{3}^{0}\right \vert ^{2}d\mathbf{\lambda }
\\
&=&\left( 2\pi \right) ^{2}\int_{\Pi ^{2}}\left \vert f_{3}\left( \mathbf{%
\lambda };\kappa _{3},a_{1}\right) -\left( 2\pi \right) ^{-2}\kappa
_{3}^{0}\right \vert ^{2}d\mathbf{\lambda >}0
\end{eqnarray*}%
because $\left \vert g_{a_{1}}\left( \lambda _{1}+\lambda _{2}\right)
\right
\vert ^{2}=1$ and where $f_{3}\left( \mathbf{\lambda };\kappa
_{3},a_{1}\right) $ is the bispectrum of an all-pass ARMA$\left( 1,1\right) $
model where $a_{1}$ is the root of the MA$\left( 1\right) $ polynomial (and $%
{a}_{1}^{\star-1}$ that of the AR$\left( 1\right) $ one$)$ and the third
order cumulant is equal to $\kappa _{3},$ whose $\mathcal{L}^{2}$ distance
to the nonzero constant function $\left( 2\pi \right) ^{-2}\kappa _{3}^{0}$
(the bispectrum of an independent series) is positive for all $\kappa
_{3}\in \mathbb{R}$, as showed by VL for any $r.$ $\Box $\bigskip }

{\small \noindent \textbf{Proof of Theorem~\ref{1ANEW}. } }

{\small \noindent \textbf{Proof for }$k=3.$\textbf{\ }It follows from
Theorem~\ref{1A}, because Assumption~3 implies Assumption~2 for $k=3.$ Note
in particular, that under Assumption~$3\left( k=3\right) $ it is easy to
check that v$\mathbf{\kappa }_{3}^{0}A^{\ast }\left( \lambda _{1}+\lambda
_{2}\right) $ contains only $d$ non-zero rows, $1,$ $d+2,\ldots ,d^{2},$
equal to the rows of $A^{\ast }\left( \lambda _{1}+\lambda _{2}\right) $ in
the same order multiplied by the corresponding marginal cumulant of $\mathbf{%
\kappa }_{3}^{0}$, so there is always at least one element $\left \{ \text{v}%
\mathbf{\kappa }_{3}^{0}A^{\ast }\left( \lambda _{1}+\lambda _{2}\right)
\right \} _{j,h}$ depending on $\lambda _{1}+\lambda _{2}$ through linear
combinations of products of $g_{a_{i}^{\ast }}\left( \lambda _{1}+\lambda
_{2}\right) $ functions, $i=1,\ldots ,r$.\bigskip }

{\small \noindent \textbf{Proof for }$k=4.$ Arguing in a similar way as in
the proof of Theorem~\ref{1A}, denoting $\lambda _{4}:=\lambda _{1}+\lambda
_{2}+\lambda _{3}$ and recalling $A^{\otimes 2}\left( \lambda _{a},\lambda
_{b}\right) =A\left( e^{-i\lambda _{a}}\right) \otimes A\left( e^{-i\lambda
_{b}}\right) ,$ we find that 
\begin{eqnarray*}
\left \Vert A^{\otimes 4}\left( \mathbf{\lambda }\right) \text{vec}\left( 
\text{v}\mathbf{\kappa }_{4}\right) -\text{vec}\left( \text{v}\mathbf{\kappa 
}_{4}^{0}\right) \right \Vert &=&\left \Vert \text{vec}\left( A^{\otimes
2}\left( \lambda _{2},\lambda _{1}\right) \text{v}\mathbf{\kappa }%
_{4}A^{\otimes 2}\left( -\lambda _{4},\lambda _{3}\right) ^{\prime }-\text{v}%
\mathbf{\kappa }_{4}^{0}\right) \right \Vert \\
&=&\left \Vert \text{vec}\left( \left( A^{\otimes 2}\left( \lambda
_{2},\lambda _{1}\right) \text{v}\mathbf{\kappa }_{4}-\text{v}\mathbf{\kappa 
}_{4}^{0}A^{\otimes 2\ast }\left( \lambda _{4},-\lambda _{3}\right) \right)
A^{\otimes 2}\left( -\lambda _{4},\lambda _{3}\right) ^{\prime }\right)
\right \Vert \\
&\geq &\left \Vert \text{vec}\left( A^{\otimes 2}\left( \lambda _{2},\lambda
_{1}\right) \text{v}\mathbf{\kappa }_{4}-\text{v}\mathbf{\kappa }%
_{4}^{0}A^{\otimes 2\ast }\left( \lambda _{4},-\lambda _{3}\right) \right)
\right \Vert .
\end{eqnarray*}
}

{\small Under Assumption~$3\left( k=4\right) $ and with the Kronecker
structure of $A^{\otimes 2\ast }\left( \lambda _{4},-\lambda _{3}\right) ,$
it is easy to check that v$\mathbf{\kappa }_{4}^{0}A^{\otimes 2\ast }\left(
\lambda _{4},-\lambda _{3}\right) $ contains only $d$ non-zero rows, $1,$ $%
d+2,\ldots ,d^{2},$ each containing all cross-products between the elements
of the same row of $A^{\ast }\left( \lambda _{4}\right) $ and of $A^{\ast
}\left( -\lambda _{3}\right) ,$ $A_{j}^{\ast }\left( e^{-i\lambda
_{4}}\right) \otimes A_{j}^{\ast }\left( e^{i\lambda _{3}}\right) ,$ $%
j=1,\ldots ,d$ (multiplied by the corresponding $j$-th element of $\mathbf{%
\kappa }_{4}^{0}$) so there is always at least one element $\left \{ \text{v}%
\mathbf{\kappa }_{4}^{0}A^{\otimes 2\ast }\left( \lambda _{4},-\lambda
_{3}\right) \right \} _{j,h}:=\xi _{4}\left( \mathbf{\lambda }\right) ,$
say, depending on $\left( \lambda _{3},\lambda _{4}\right) $ in the form of
a product of two functions of $\lambda _{3}$ and $\lambda _{4}$,
respectively, each with infinite series expansion in powers of $\lambda _{i}$%
. Then, following the same argument as for $k=3$ in Theorem~\ref{1A},
defining $\zeta _{4}\left( \mathbf{\lambda }\right) :=\left \{ A^{\otimes
2}\left( \lambda _{2},\lambda _{1}\right) \text{v}\mathbf{\kappa }%
_{4}\right
\} _{j,h}, $ we can show that 
\begin{equation*}
\mathcal{L}_{4}^{0}\left( A,\text{v}\mathbf{\kappa }_{4}\right) \geq \eta
^{8}\int_{\Pi ^{3}}\left \vert \zeta _{4}\left( \mathbf{\lambda }\right)
-\xi _{4}\left( \mathbf{\lambda }\right) \right \vert ^{2}d\mathbf{\lambda }%
\geq \epsilon >0
\end{equation*}%
and the theorem follows.\ $\Box $\bigskip }

{\small \noindent \textbf{Proof of Theorem~\ref{1B}.} }

{\small \noindent \textbf{Proof for} $k=3.$ For $A=K$ orthogonal$,$ $%
KK^{\prime }=\mathbf{I}_{d},$ we need to show that $K=P_{d}$ and v$\mathbf{%
\kappa }_{3}=\ $v$\mathbf{\kappa }_{3}^{\text{IC}}\left( P_{d}^{\prime }%
\mathbf{\alpha }_{3}^{0}\right) $ for signed permutation matrices $P_{d}$
are the only solutions that make true vec$\left( \text{v}\mathbf{\kappa }%
_{3}^{0}\right) =K^{\otimes 3}\ $vec$\left( \text{v}\mathbf{\kappa }%
_{3}\right) $ or equivalently v$\mathbf{\kappa }_{3}^{0}=K^{\otimes 2}\ $v$%
\mathbf{\kappa }_{3}\ K^{\prime },$ or 
\begin{equation}
\text{v}\mathbf{\kappa }_{3}^{0}\ K=K^{\otimes 2}\  \text{v}\mathbf{\kappa }%
_{3}.  \label{star}
\end{equation}%
Take any $d\geq 2$ and%
\begin{equation*}
\text{v}\mathbf{\kappa }_{3}^{0}:=\ v\mathbf{\kappa }_{3}^{\text{IC}}\left( 
\mathbf{\alpha }_{3}^{0}\right) =\left( \alpha _{1}^{0}\mathbf{e}%
_{1}^{\otimes 2},\  \  \alpha _{2}^{0}\mathbf{e}_{2}^{\otimes 2},\ldots
,\alpha _{d}^{0}\mathbf{e}_{d}^{\otimes 2}\right) ,\  \ K=\left \{
K_{ab}\right \} _{a,b=1}^{d}
\end{equation*}%
where the constants $\alpha _{j}^{0}$ satisfy $\alpha _{1}^{0}\alpha
_{2}^{0}\cdots \alpha _{d}^{0}\neq 0$ and $\mathbf{e}_{1}=\left( 1,0,\ldots
,0\right) ^{\prime }$ and so on are the unitary vectors of dimension $d$,
and we impose the same restriction to v$\mathbf{\kappa }_{3},$ i.e. 
\begin{equation*}
\text{v}\mathbf{\kappa }_{3}=\left( \alpha _{1}\mathbf{e}_{1}^{\otimes 2},\
\  \alpha _{2}\mathbf{e}_{2}^{\otimes 2},\ldots ,\alpha _{d}\mathbf{e}%
_{d}^{\otimes 2}\right) ,
\end{equation*}%
with $\alpha _{j}$ satisfying $\alpha _{1}\alpha _{2}\cdots \alpha _{d}\neq
0 $ and by orthogonality%
\begin{eqnarray}
K_{a1}^{2}+K_{a2}^{2}+\cdots +K_{ad}^{2} &=&1,\  \ a=1,\ldots ,d  \label{O1.a}
\\
K_{a1}K_{b1}+K_{a2}K_{b2}+\cdots +K_{ad}K_{bd} &=&0,\  \text{\  \ }%
b<a=2,\ldots ,d  \label{O2.ab}
\end{eqnarray}%
so that also $K_{ab}^{2}\leq 1.$ }

{\small Then we have that (\ref{star}) is equivalent to%
\begin{equation*}
\left( 
\begin{array}{cccc}
\alpha _{1}^{0}K_{11} & \alpha _{1}^{0}K_{12} & \cdots & \alpha
_{1}^{0}K_{1d} \\ 
0 & 0 &  & 0 \\ 
\vdots & \vdots &  & \vdots \\ 
0 & 0 &  & 0 \\ 
\alpha _{2}^{0}K_{21} & \alpha _{2}^{0}K_{22} & \cdots & \alpha
_{2}^{0}K_{2d} \\ 
0 & 0 &  & 0 \\ 
\vdots & \vdots &  & \vdots \\ 
&  &  &  \\ 
0 & 0 &  & 0 \\ 
\alpha _{d}^{0}K_{d1} & \alpha _{d}^{0}K_{d2} & \cdots & \alpha
_{d}^{0}K_{dd}%
\end{array}%
\right) =\left( 
\begin{array}{cccc}
\alpha _{1}K_{11}^{2} & \alpha _{2}K_{12}^{2} & \cdots & \alpha
_{d}K_{1d}^{2} \\ 
\vdots & \vdots & \cdots & \vdots \\ 
\alpha _{1}K_{11}K_{d1} & \alpha _{2}K_{12}K_{d2} & \cdots & \alpha
_{d}K_{1d}K_{dd} \\ 
\alpha _{1}K_{21}K_{11} & \alpha _{2}K_{22}K_{12} &  & \alpha
_{d}K_{2d}K_{1d} \\ 
\alpha _{1}K_{21}^{2} & \alpha _{2}K_{22}^{2} &  & \alpha _{d}K_{2d}^{2} \\ 
\vdots & \vdots &  & \vdots \\ 
\alpha _{1}K_{21}K_{d1} & \alpha _{2}K_{22}K_{d2} &  & \alpha
_{d}K_{2d}K_{dd} \\ 
\alpha _{1}K_{31}K_{11} & \alpha _{2}K_{32}K_{12} &  & \alpha
_{d}K_{3d}K_{1d} \\ 
\vdots & \vdots &  & \vdots \\ 
\alpha _{1}K_{d1}^{2} & \alpha _{2}K_{d2}^{2} & \cdots & \alpha
_{d}K_{dd}^{2}%
\end{array}%
\right)
\end{equation*}%
so that the restrictions imposed for all $a,b=1,\ldots ,d$ are%
\begin{equation}
\alpha _{a}^{0}K_{ab}=\alpha _{b}K_{ab}^{2}\  \  \Rightarrow \  \ K_{ab}=0\  \ 
\text{or\  \ }=\alpha _{a}^{0}/\alpha _{b}\ast 1\left \{ 0<\left \vert \alpha
_{a}^{0}/\alpha _{b}\right \vert \leq 1\right \} .  \label{V1.ab}
\end{equation}%
Note also that it must hold that $\alpha _{b}\neq 0$ for all $b,$ as
otherwise the $b$-th column of $K$ would be exactly zero.

We note the following consequences of these restrictions on $K$: }

\begin{enumerate}
\item {\small For $a=b,$ $K_{aa}=0$ or $=\alpha _{a}^{0}/\alpha _{a}$ with $%
0<\left
\vert \alpha _{a}^{0}/\alpha _{a}\right \vert \leq 1$ by (\ref{V1.ab}%
). }

\item {\small If $K_{aa}=\pm 1$ $\Rightarrow \ $ $\alpha _{a}=\pm \alpha
_{a}^{0},$ with$\  \ K_{ab}=0,b\neq a,$ by (\ref{O1.a}) ($\alpha _{b}$ no
restricted), i.e. the $a$-th row of $K$ is, up to sign, the unitary vector $%
\mathbf{e}_{a}.$ }

\item {\small If $K_{ab}=\pm 1\  \Rightarrow \  \  \alpha _{a}^{0}/\alpha
_{b}=\pm 1\ (a\neq b),\ $with$\  \ K_{aa}=0$ and $K_{aj}=0$ for all $j\neq b$
by (\ref{O1.a}). }

\item {\small If $K_{aa}\neq \pm 1$ $\  \  \Rightarrow \  \  \exists b\neq a,\  \ 
$s.t. $K_{ab}\neq 0$ by (\ref{O1.a}).\bigskip }
\end{enumerate}

{\small Then, to show that the only solutions to the equations (\ref{star})
are matrices $K$ which are signed permutations between row $b$ and row $a$
with signs given by $K_{ab}=\alpha _{a}^{0}/\alpha _{b}=\pm 1,$ $K_{aj}=0,$ $%
j\neq b,$ we have to show that it is not possible to select for any $%
a=1,\ldots ,d$ a set of indexes $\mathcal{I}_{a}=\left \{ j_{a,1},\ldots
,j_{a,p\left( a\right) }\right \} ,$ $\# \mathcal{I}_{a}=p\left( a\right)
\geq 2,$ $j_{a,\ell }\in \left \{ 1,\ldots ,d\right \} $ and values $%
K_{a\ell }\neq 0$ for $\ell \in \mathcal{I}_{a}$ that satisfy%
\begin{eqnarray}
\sum_{\ell \in \mathcal{I}_{a}}K_{a\ell }^{2} &=&\sum_{\ell \in \mathcal{I}%
_{a}}\left( \frac{\alpha _{a}^{0}}{\alpha _{\ell }}\right) ^{2}=1,\  \  \ 
\label{C1} \\
\sum_{\ell \in \mathcal{I}_{a}\cap \mathcal{I}_{b}}K_{a\ell }K_{b\ell }
&=&\sum_{\ell \in \mathcal{I}_{a}\cap \mathcal{I}_{b}}\frac{\alpha _{a}^{0}}{%
\alpha _{\ell }}\frac{\alpha _{b}^{0}}{\alpha _{\ell }}=0,\  \  \ b\neq a,
\label{C2}
\end{eqnarray}%
i.e. the normalization and orthogonalization conditions of $K,$ (\ref{O1.a}%
)-(\ref{O2.ab}). }

{\small Then note that condition (\ref{C1}) excludes permutation matrices
because \#$\mathcal{I}_{a}=p\left( a\right) \geq 2$ ($p\left( b\right) \geq
1,$ $b\neq a $) and that condition (\ref{C2}) implies that for all $a\neq b$%
\begin{equation*}
0=\sum_{\ell \in \mathcal{I}_{a}\cap \mathcal{I}_{b}}\frac{\alpha _{a}^{0}}{%
\alpha _{\ell }}\frac{\alpha _{b}^{0}}{\alpha _{\ell }}=\alpha
_{a}^{0}\alpha _{b}^{0}\sum_{\ell \in \mathcal{I}_{a}\cap \mathcal{I}_{b}}%
\frac{1}{\alpha _{\ell }^{2}}
\end{equation*}%
which, given that $\alpha _{\ell }\neq 0$ for all $\ell \in \mathcal{I}%
_{a}\cap \mathcal{I}_{b},$ is only feasible if $\mathcal{I}_{a}\cap \mathcal{%
I}_{b}=\emptyset $ for all $a\neq b,$ but since $p\left( a\right) \geq 2$,
even if $p\left( b\right) =1$, $b\neq a$, there must be some $b$ for which $%
\mathcal{I}_{a}\cap \mathcal{I}_{b}\not=\emptyset $ and therefore we
conclude that it is not possible to make such selection of indexes to
construct $K$ with elements different from $\pm 1.$\bigskip }

{\small \noindent \textbf{Proof for }$k=4.\ $Consider $d\geq 2,$ denote $%
\beta ^{0}=\mathbf{\alpha }_{4}^{0}$ and%
\begin{equation*}
\text{v}\mathbf{\kappa }_{4}^{0}:=\ v\mathbf{\kappa }_{4}^{\text{IC}}\left( 
\mathbf{\alpha }_{4}^{0}\right) =\left( \beta _{1}^{0}\mathbf{e}_{1}\ 0\
\cdots \ 0\  \beta _{2}^{0}\mathbf{e}_{d+2}\ 0\  \cdots \ 0\  \beta _{d}^{0}%
\mathbf{e}_{d^{2}}\right) ,\  \ 
\end{equation*}%
where the constants $\beta _{j}^{0}$ satisfy $\beta _{1}^{0}\beta
_{2}^{0}\cdots \beta _{d}^{0}\neq 0$, v$\mathbf{\kappa }_{4}=\ $v$\mathbf{%
\kappa }_{4}^{\text{IC}}\left( \mathbf{\beta }\right) $ has the same
structure for some constants $\beta _{j}$ satisfying $\beta _{1}\beta
_{2}\cdots \beta _{d}\neq 0$ and $K$ satisfies the normalization and
orthogonality conditions (\ref{O1.a}) and (\ref{O2.ab}). }

{\small Then we have that vec$\left( \text{v}\mathbf{\kappa }_{4}^{0}\right)
=K^{\otimes 4}\ $vec$\left( \text{v}\mathbf{\kappa }_{4}\right) $ is
equivalent to v$\mathbf{\kappa }_{4}^{0}\ K^{\otimes 2}=K^{\otimes 2}\ $v$%
\mathbf{\kappa }_{4},$ or to%
\begin{equation*}
\left( 
\begin{array}{cccccccc}
\beta _{1}^{0}K_{11}^{2} & \beta _{1}^{0}K_{11}K_{12} & \cdots & \beta
_{1}^{0}K_{11}K_{1d} & \beta _{1}^{0}K_{12}K_{11} & \beta _{1}^{0}K_{12}^{2}
& \cdots & \beta _{1}^{0}K_{1d}^{2} \\ 
0 & 0 & \cdots & 0 & 0 & 0 & \cdots & 0 \\ 
\vdots & \vdots &  & \vdots & \vdots & \vdots &  & \vdots \\ 
0 & 0 & \cdots & 0 & 0 & 0 & \cdots & 0 \\ 
\beta _{2}^{0}K_{21}^{2} & \beta _{2}^{0}K_{21}K_{22} & \cdots & \beta
_{2}^{0}K_{21}K_{1d} & \beta _{2}^{0}K_{22}K_{21} & \beta _{2}^{0}K_{22}^{2}
& \cdots & \beta _{2}^{0}K_{2d}^{2} \\ 
0 & 0 & \cdots & 0 & 0 & 0 & \cdots & 0 \\ 
\vdots & \vdots &  & \vdots & \vdots & \vdots &  & \vdots \\ 
0 & 0 & \cdots & 0 & 0 & 0 & \cdots & 0 \\ 
\beta _{d}^{0}K_{d1}^{2} & \beta _{d}^{0}K_{d1}K_{d2} & \cdots & \beta
_{d}^{0}K_{d1}K_{dd} & \beta _{d}^{0}K_{d2}K_{d1} & \beta _{d}^{0}K_{d2}^{2}
& \cdots & \beta _{d}^{0}K_{dd}^{2}%
\end{array}%
\right)
\end{equation*}%
\begin{equation*}
=\left( 
\begin{array}{ccccccccc}
\beta _{1}K_{11}^{2} & 0 & \cdots & 0 & \beta _{2}K_{12}^{2} & 0 & \cdots & 0
& \beta _{d}K_{1d}^{2} \\ 
\beta _{1}K_{11}K_{21} & 0 & \cdots & 0 & \beta _{2}K_{12}K_{22} & 0 & \cdots
& 0 & \beta _{d}K_{1d}K_{2d} \\ 
\vdots & \vdots &  & \vdots & \vdots & \vdots &  & \vdots & \vdots \\ 
\beta _{1}K_{11}K_{d1} & 0 & \cdots & 0 & \beta _{2}K_{12}K_{d2} & 0 & \cdots
& 0 & \beta _{d}K_{1d}K_{dd} \\ 
\beta _{1}K_{21}K_{11} & 0 & \cdots & 0 & \beta _{2}K_{22}K_{12} & 0 & \cdots
& 0 & \beta _{d}K_{2d}K_{1d} \\ 
\beta _{1}K_{21}^{2} & 0 & \cdots & 0 & \beta _{2}K_{22}^{2} & 0 & \cdots & 0
& \beta _{d}K_{2d}^{2} \\ 
\vdots & \vdots &  & \vdots & \vdots & \vdots &  & \vdots & \vdots \\ 
\beta _{1}K_{21}K_{d1} & 0 & \cdots & 0 & \beta _{2}K_{22}K_{d2} & 0 & \cdots
& 0 & \beta _{d}K_{2d}K_{dd} \\ 
\beta _{1}K_{31}K_{11} & 0 & \cdots & 0 & \beta _{2}K_{32}K_{12} & 0 & \cdots
& 0 & \beta _{d}K_{3d}K_{1d} \\ 
\vdots & \vdots &  & \vdots & \vdots & \vdots &  & \vdots & \vdots \\ 
\beta _{1}K_{d1}^{2} & 0 & \cdots & 0 & \beta _{2}K_{d2}^{2} & 0 & \cdots & 0
& \beta _{d}K_{dd}^{2}%
\end{array}%
\right)
\end{equation*}%
so for each pair $\left( a,b\right) $%
\begin{equation*}
K_{ab}=0\  \  \text{or\  \ }\beta _{b}=\beta _{a}^{0}\neq 0,
\end{equation*}%
and%
\begin{equation*}
K_{ab}K_{ac}=0,\ b\neq c,
\end{equation*}%
i.e. in each row $a$ there must be one single non-zero element, say $%
K_{a\ell }$, equal to $\pm 1$ to fulfill the orthogonality and normalization
restrictions, so all the restrictions can only hold when $K$ is a
permutation matrix $P_{d}$ and the assignment $\beta _{\ell }=\beta _{a}^{0}$
is done attending to the location of these nonzero elements $K_{a\ell }=\pm
1,$ i.e. $\mathbf{\beta }=P_{d}^{+\prime }\mathbf{\beta .}$ \ }

{\small This provides identification of the components of $\mathbf{%
\varepsilon }_{t}$ up to signed permutations because for any $A\left(
z\right) =P_{d}$ and any v$\mathbf{\kappa }_{k}^{0}=\ $v$\mathbf{\kappa }%
_{k}^{IC}\left( \mathbf{\alpha }_{k}^{0}\right) $ we could select v$\mathbf{%
\kappa }_{3}=\ $v$\mathbf{\kappa }_{3}^{IC}\left( P_{d}^{\prime }\mathbf{%
\alpha }_{3}^{0}\right) $ or v$\mathbf{\kappa }_{4}=\ $v$\mathbf{\kappa }%
_{4}^{IC}\left( P_{d}^{+\prime }\mathbf{\alpha }_{4}^{0}\right) $ so that $%
\mathcal{L}_{3}^{0}\left( P_{d},\text{v}\mathbf{\kappa }_{3}^{\text{IC}%
}\left( P_{d}^{\prime }\mathbf{\alpha }_{3}^{0}\right) \right) =\mathcal{L}%
_{4}^{0}\left( P_{d},\text{v}\mathbf{\kappa }_{4}^{\text{IC}}\left(
P_{d}^{+\prime }\mathbf{\alpha }_{4}^{0}\right) \right) =0\ $because that
for these choices of v$\mathbf{\kappa }_{k},$ 
\begin{equation*}
P_{d}^{\otimes 3}\text{vec}\left( \text{v}\mathbf{\kappa }_{3}\right) =\text{%
vec}\left( P_{d}^{\otimes 2}\text{v}\mathbf{\kappa }_{3}^{IC}\left(
P_{d}^{\prime }\mathbf{\alpha }_{3}^{0}\right) P_{d}^{\prime }\right) =\text{%
vec}\left( \text{v}\mathbf{\kappa }_{3}^{IC}\left( \mathbf{\alpha }%
_{3}^{0}\right) \right) =\text{vec}\left( \text{v}\mathbf{\kappa }%
_{3}^{0}\right)
\end{equation*}%
and%
\begin{equation*}
P_{d}^{\otimes 4}\text{vec}\left( \text{v}\mathbf{\kappa }_{4}\right) =\text{%
vec}\left( P_{d}^{\otimes 2}\text{v}\mathbf{\kappa }_{4}^{IC}\left(
P_{d}^{+\prime }\mathbf{\alpha }_{4}^{0}\right) P_{d}^{\otimes 2\prime
}\right) =\text{vec}\left( \text{v}\mathbf{\kappa }_{4}^{IC}\left( \mathbf{%
\alpha }_{4}^{0}\right) \right) =\text{vec}\left( \text{v}\mathbf{\kappa }%
_{4}^{0}\right) ,
\end{equation*}%
reflecting that kurtosis is identified independently of the sign of $\mathbf{%
\varepsilon }_{t},$ but skewness is not, because $\alpha _{3,j}^{0}=\mathbf{e%
}_{j}^{\prime }\mathbf{\alpha }_{3}^{0}$ and from (\ref{vk3}) it follows that%
\begin{equation*}
P_{d}^{\otimes 2}\text{v}\mathbf{\kappa }_{3}^{IC}\left( P_{d}^{\prime }%
\mathbf{\alpha }_{3}^{0}\right) P_{d}^{\prime }=\sum_{j=1}^{d}\mathbf{e}%
_{j}^{\prime }\left( P_{d}^{\prime }\mathbf{\alpha }_{3}^{0}\right) \left(
P_{d}\mathbf{e}_{j}\right) ^{\otimes 2}\left( P_{d}\mathbf{e}_{j}\right)
^{\prime }=\sum_{j=1}^{d}\alpha _{3j}^{0}\mathbf{e}_{j}^{\otimes 2}\mathbf{e}%
_{j}^{\prime }=\text{v}\mathbf{\kappa }_{3}^{IC}\left( \mathbf{\alpha }%
_{3}^{0}\right) .
\end{equation*}%
$\Box $\bigskip }

{\small \noindent \textbf{Proof of Theorem~\ref{1A1B}}. The results follows
from combining Theorems~\ref{1A} and~\ref{1ANEW} for non-constant and
constant $A\left( z\right) ,$ respectively. $\Box $\bigskip }

{\small \noindent \textbf{Proof of Corollary \ref{L0}. }For $k=3,$
proceeding as in the proof of Theorem~\ref{1B}, and assuming w.l.o.g. $%
\alpha _{1}^{0}=0,$ the set of assumptions imposed by (\ref{V1.ab}) only
affect now to rows $a=2,\ldots ,d,$ of \ $K,$ while for $a=1\ $the
implication is 
\begin{equation}
0=\alpha _{b}K_{1b}^{2}\  \  \Rightarrow \  \ K_{1b}=0\  \  \text{and/or}\  \
\alpha _{b}=0\  \  \text{(any }K_{1b}\text{).}  \label{V1.ab*}
\end{equation}%
Therefore, for $a=2,\ldots ,d,$ it is only possible to set $K_{ab}=\alpha
_{a}^{0}/\alpha _{b}\in \left[ -1,1\right] /\left \{ 0\right \} $ when $%
\alpha _{b}\neq 0$ (which implies $K_{1b}=0$ by (\ref{V1.ab*})), otherwise,
if $\alpha _{b}=0,$ then $K_{ab}=0$ for all $a>1.$ Hence, at most there can
be one $b$ such that $\alpha _{b}=0,$ otherwise there would be more than one
row in $K$ with a single nonzero element in the first column (and $K$ would
be no full rank). On the other hand, there must be at least one $\alpha
_{b}=0$ because otherwise the first row of $K$ would be zero. Then, for this
unique $b,$ $K_{1b}=\pm 1,$ $K_{1\ell }=0,$ $\ell \neq b,$ $K_{ab}=0,$ $%
a\neq 1,$ while for the rest of rows and columns of $K$ we can apply the
same argument as in the proof of Theorem~\ref{1B} to show that it must be
also a sign-permutation. }

{\small When there are two or more values of $a$ such that $\alpha
_{a}^{0}=0,$ e.g. $\alpha _{1}^{0}=\alpha _{2}^{0}=0,$ then it is possible
to set up to two values $\alpha _{b}=0,$ for $b=1,2,$ say, where the top
left corner of $K$ can be chosen freely as any orthogonal $2\times 2$
matrix, with the off diagonal blocks of $K$ being zero, while the bottom
right corner should remain of permutation type. }

{\small For $k=4$ and $\beta _{1}^{0}=0,$ we find again that $K_{1b}=0$
and/or $\beta _{b}=0$ (any $K_{1b}$), so the same reasoning as for $k=3,$
applies, since we can only set up one $\beta _{b}=0,$ and the first row (and
column) of $K$ will only contain a single nonzero element equal to $\pm 1,$
the rest of $K$ being also a permutation matrix orthogonal to this first
row.\bigskip }

{\small \noindent \textbf{Proof of Corollary \ref{L1}. }We only consider the
case $\mathcal{I}_{3}\cap \mathcal{I}_{4}=\varnothing ,$ i.e. \#$\left \{ 
\mathcal{I}_{3}\cap \mathcal{I}_{4}\right \} =0,$ the case with $\alpha
_{j}^{0}=\beta _{j}^{0}=0$ for a single $j$ can be dealt with using the
arguments in Corollary~\ref{L0} to show static identification. }

{\small Proceeding as in the proof of Theorem~\ref{1B}, the restrictions
given by v$\mathbf{\kappa }_{3}^{0}\ K=K^{\otimes 2}\ $v$\mathbf{\kappa }%
_{3}\ $for $k=3\ $in the typical case when all $\alpha _{j}^{0}\neq 0,$ $%
j=1,\ldots ,m,$ but $\alpha _{m+1}^{0}=\cdots =\alpha _{d}=0,$ so that 
\begin{equation*}
\text{v}\mathbf{\kappa }_{3}^{0}\ K=\left( 
\begin{array}{cccc}
\alpha _{1}^{0}K_{11} & \alpha _{1}^{0}K_{12} & \cdots & \alpha
_{1}^{0}K_{1d} \\ 
0 & 0 &  & 0 \\ 
\vdots & \vdots &  & \vdots \\ 
0 & 0 &  & 0 \\ 
\alpha _{2}^{0}K_{21} & \alpha _{2}^{0}K_{22} & \cdots & \alpha
_{2}^{0}K_{2d} \\ 
0 & 0 &  & 0 \\ 
\vdots & \vdots &  & \vdots \\ 
0 & 0 & \cdots & 0 \\ 
\alpha _{m}^{0}K_{m1} & \alpha _{m}^{0}K_{m1} &  & \alpha _{m}^{0}K_{21} \\ 
0 & 0 & \cdots & 0 \\ 
\vdots &  &  & \vdots \\ 
0 & 0 & \cdots & 0%
\end{array}%
\right)
\end{equation*}%
are 
\begin{equation*}
K_{ab}=0\  \  \text{or\ }=\alpha _{a}^{0}/\alpha _{b}\  \  \text{for\  \ }%
a=1,\ldots ,m,
\end{equation*}%
while the rows $j=m+1,\ldots ,d$ of $K$ are essentially unrestricted when $%
d-m\geq 2,$ cf. proof of Corollary~\ref{L0}. This implies by the same
argument as in Theorem~\ref{1B} that these restrictions can only be
satisfied for $K$ whose top $m$ rows are of permutation type if we show that
the remaining rows have to be also of permutation type. }

{\small The restrictions v$\mathbf{\kappa }_{4}^{0}\ K^{\otimes
2}=K^{\otimes 2}\ $v$\mathbf{\kappa }_{4}$ for $k=4$, under independence and 
$\beta _{j}^{0}\neq 0,$ $j=m+1,\ldots ,d$ and $\beta _{1}^{0}=\cdots =\beta
_{m}^{0}=0,$ where now%
\begin{equation*}
v\mathbf{\kappa }_{4}^{0}\ K^{\otimes 2}=\left( 
\begin{array}{cccc}
0 & 0 & \cdots & 0 \\ 
0 & 0 & \cdots & 0 \\ 
\vdots & \vdots &  & \vdots \\ 
0 & 0 & \cdots & 0 \\ 
\beta _{m+1}^{0}K_{m+11}^{2} & \beta _{m+1}^{0}K_{m+11}K_{m+12} & \cdots & 
\beta _{m+1}^{0}K_{m+1d}^{2} \\ 
\vdots & \vdots &  & \vdots \\ 
\beta _{d}^{0}K_{d1}^{2} & \beta _{d}^{0}K_{d1}K_{d2} & \cdots & \beta
_{d}^{0}K_{dd}^{2}%
\end{array}%
\right) ,
\end{equation*}%
are, following as in Theorem~\ref{1B},%
\begin{eqnarray*}
K_{ab} &=&0\  \text{or otherwise }\beta _{b}=\beta _{a}^{0}=0,\  \  \text{for }%
a=1,\ldots ,m \\
K_{ab} &=&0\  \text{or otherwise }\beta _{b}=\beta _{a}^{0}\neq 0,\  \text{for 
}a=m+1,\ldots ,d,
\end{eqnarray*}%
and still%
\begin{equation*}
K_{ab}K_{ac}=0,\  \ b\neq c,\  \  \ a=m+1,\ldots ,d,
\end{equation*}%
i.e. each of the last $d-m$ rows must contain one single non-zero element,
equal to $\pm 1$ to fulfill the orthogonality and normalization
restrictions. Then the restrictions for $k=3$ imply that the top $m$ rows
have to be also of signed permutation type, and therefore, considering both
types of restrictions together, we obtain that they can only hold
simultaneously for $K$ equal to a signed permutation $P_{d},$ where the
non-zero and zero values of $\alpha _{j}^{0}$ and $\beta _{j}^{0}$ in v$%
\mathbf{\kappa }_{3}^{0}$ and v$\mathbf{\kappa }_{4}^{0}$ (and signs for $%
\alpha _{j}^{0}$) are permuted in the same order by $\mathbf{\alpha }%
=P_{d}^{\prime }\mathbf{\alpha }^{0}$ and $\mathbf{\beta }=P_{d}^{+\prime }%
\mathbf{\beta }^{0}.$ }

{\small To show dynamic identification when $\mathcal{I}_{3}\cap \mathcal{I}%
_{4}=\varnothing ,$ we take in the proof of Theorem~\ref{1ANEW} the simplest
case in which $\alpha _{41}^{0}=0,$ but all order marginal cumulants of
order $k=4$ are nonzero. In this case, the argument for identification fails
because there are BM $A(z)$ for which the only row of $A^{\otimes 2\ast
}\left( \lambda _{4},-\lambda _{3}\right) $ which depends on $\left( \lambda
_{3},\lambda _{4}\right) $ is the first one, so that $\xi _{4}(\mathbf{%
\lambda })$ can not be constructed and it would be possible to choose $\text{%
v}\mathbf{\kappa }_{4}$ to achieve $\left \Vert \text{vec}\left( A^{\otimes
2}\left( \lambda _{2},\lambda _{1}\right) \text{v}\mathbf{\kappa }_{4}-\text{%
v}\mathbf{\kappa }_{4}^{0}A^{\otimes 2\ast }\left( \lambda _{4},-\lambda
_{3}\right) \right) \right \Vert =0$ a.e.. However, if at the same time $%
\alpha _{31}^{0}\neq 0$, we can take $j=1$ in the argument of the proof of
Theorem~\ref{1ANEW} to define $\xi _{3}(\lambda _{1}+\lambda _{2})$ for an
appropriate $h$ so that it has an infinite expansion in $\lambda
_{1}+\lambda _{2}$ that can not be matched by any function $\zeta
_{3}(\lambda _{1},\lambda _{2})$ for any choice of $\text{v}\mathbf{\kappa }%
_{3}$ and such BM $A\left( z\right) .$ The same reasoning applies in the
other direction and when more than one marginal cumulant is zero for any
order $k$ as far as the corresponding cumulants of the other order are
different from zero. $\square \bigskip $ }

{\small \noindent \textbf{Proof of Theorem~\ref{Th4}}. By Theorem~3 of Lippi
and Reichlin (1994), all second-order equivalent representations of $%
Y_{t}=\Phi _{\mathbf{\theta }_{0}}^{-1}\left( z\right) \Theta _{\mathbf{%
\theta }_{0}}\left( z\right) \mathbf{\varepsilon }_{t},$ $\mathbf{%
\varepsilon }_{t}=iid\left( \mathbf{0},\mathbf{I}_{d},\text{v}\mathbf{\kappa 
}_{k}^{IC}\left( \mathbf{\alpha }_{k}^{0}\right) \right) ,$ giving $\mathcal{%
L}_{2}\left( \mathbf{\theta }\right) =0$ for some $\mathbf{\theta \neq
\theta }_{0}\mathbf{,\  \theta }\in \mathcal{S}$, involve a matrix polynomial 
$\Phi _{\mathbf{\theta }}^{-1}\left( z\right) \Theta _{\mathbf{\theta }%
}\left( z\right) =\Phi _{\mathbf{\theta }_{0}}^{-1}\left( z\right) \Theta _{%
\mathbf{\theta }_{0}}\left( z\right) K$ up to an orthogonal matrix $K$ when $%
\Theta _{\mathbf{\theta }}\left( z\right) $ has the same roots as $\Theta _{%
\mathbf{\theta }_{0}}\left( z\right) ,$ while by Theorem~2 of Lippi and
Reichlin (1994) the transfer functions differ by a (non constant) BM $%
A\left( z\right) $ when some of the roots of $\Theta _{\mathbf{\theta }%
_{0}}\left( z\right) $ are flipped, $\Phi _{\mathbf{\theta }}^{-1}\left(
z\right) \Theta _{\mathbf{\theta }}\left( z\right) =\Phi _{\mathbf{\theta }%
_{0}}^{-1}\left( z\right) \Theta _{\mathbf{\theta }_{0}}\left( z\right)
A\left( z\right) .$ }

{\small Then the proof follows from our Theorem~\ref{1A1B} which shows that $%
\mathcal{L}_{k}\left( \mathbf{\theta ,\alpha }_{k}\right) =0$ could only
hold if such $A$ is equal to a signed permutation $P_{d}$ for $k=3$ or $k=4$
and $\mathbf{\alpha }_{3}=P_{d}^{\prime }\mathbf{\alpha }_{3}^{0}$ or $%
\mathbf{\alpha }_{4}=P_{d}^{+\prime }\mathbf{\alpha }_{4}^{0}$ (with $%
P_{d}^{+}$ being equal to $P_{d}$ with positive entries), while Assumption~6
discards all such signed permutations except the identity one by fixing an
ordering and a sign pattern, and Assumption~5.2 excludes all values $\mathbf{%
\theta }\neq \mathbf{\theta }_{0}$ such that $\Phi _{\mathbf{\theta }%
}^{-1}\left( z\right) \Theta _{\mathbf{\theta }}\left( z\right) =\Phi _{%
\mathbf{\theta }_{0}}^{-1}\left( z\right) \Theta _{\mathbf{\theta }%
_{0}}\left( z\right) A\left( z\right) $ a.e. for $A\left( z\right) =\mathbf{I%
}_{d}.$ $\Box $\bigskip }

{\small \noindent \textbf{Proof of Theorem~\ref{Th5}}. As in VL (2018) we
can show uniformly for $\mathbf{\theta \in \ }\mathcal{S}$ and for $k=3\ $or 
$4$ that%
\begin{equation*}
\mathbf{\hat{\alpha}}_{k,T}\left( \mathbf{\theta }\right) \rightarrow _{p}%
\mathbf{\alpha }_{k}\left( \mathbf{\theta }\right)
\end{equation*}%
where 
\begin{equation*}
\mathbf{\alpha }_{k}\left( \mathbf{\theta }\right) :=\left( \int_{\Pi ^{k-1}}%
\func{Re}\left \{ \mathbf{S}_{k}^{\prime }\mathbf{\Psi }_{k}^{\ast }(\mathbf{%
\lambda };\mathbf{\theta })\mathbf{\Psi }_{k}(\mathbf{\lambda };\mathbf{%
\theta })\mathbf{S}_{k}\right \} d\mathbf{\lambda }\right) ^{-1}\int_{\Pi
^{k-1}}\func{Re}\left \{ \mathbf{S}_{k}^{\prime }\mathbf{\Psi }_{k}^{\ast }(%
\mathbf{\lambda };\mathbf{\theta })\mathbf{\Psi }_{k}(\mathbf{\lambda };%
\mathbf{\theta }_{0})\mathbf{S}_{k}\right \} d\mathbf{\lambda \  \mathbf{%
\alpha }}_{k}^{0}
\end{equation*}%
satisfies $\mathbf{\alpha }_{k}\left( \mathbf{\theta }_{0}\right) =\mathbf{%
\alpha }_{k}^{0},$ $\mathcal{\hat{L}}_{k,T}\left( \mathbf{\theta }\right) $
converges uniformly to 
\begin{equation*}
\mathcal{\hat{L}}_{k}\left( \mathbf{\theta }\right) =\left( 2\pi \right)
^{1-k}\int_{\Pi ^{k-1}}\left( \mathbf{\Psi }_{k}\left( \mathbf{\lambda
;\theta }\right) \mathbf{S}_{k}\mathbf{\alpha }_{k}\left( \mathbf{\theta }%
\right) -\mathbf{\Psi }_{k}\left( \mathbf{\lambda ;\theta }_{0}\right) 
\mathbf{S}_{k}\mathbf{\alpha }_{k}^{0}\right) ^{\ast }\left( \mathbf{\Psi }%
_{k}\left( \mathbf{\lambda ;\theta }\right) \mathbf{S}_{k}\mathbf{\alpha }%
_{k}\left( \mathbf{\theta }\right) -\mathbf{\Psi }_{k}\left( \mathbf{\lambda
;\theta }_{0}\right) \mathbf{S}_{k}\mathbf{\alpha }_{k}^{0}\right) d\mathbf{%
\lambda }
\end{equation*}%
while $\mathcal{L}_{2,T}\left( \mathbf{\theta }\right) $ converges uniformly
to%
\begin{equation*}
\mathcal{L}_{2}\left( \mathbf{\theta }\right) =\left( 2\pi \right)
^{-1}\int_{\Pi }\text{vec}\left( \mathbf{I}_{d}\right) ^{\prime }\left( 
\mathbf{\Psi }_{2}\left( \lambda \mathbf{;\theta }\right) -\mathbf{\Psi }%
_{2}\left( \lambda \mathbf{;\theta }_{0}\right) \right) ^{\ast }\left( 
\mathbf{\Psi }_{2}\left( \lambda \mathbf{;\theta }\right) -\mathbf{\Psi }%
_{2}\left( \lambda \mathbf{;\theta }_{0}\right) \right) \text{vec}\left( 
\mathbf{I}_{d}\right) d\lambda
\end{equation*}%
which is minimized for $\mathbf{\theta =\theta }_{0}$ and for all $\mathbf{%
\theta \neq \theta }_{0}$ such that $\mathbf{\Psi }\left( e^{-i\lambda }%
\mathbf{;\theta }\right) =\mathbf{\Psi }\left( e^{-i\lambda }\mathbf{;\theta 
}_{0}\right) A(e^{-i\lambda })\ $where the BM factor satisfies $%
A(e^{-i\lambda })\neq \mathbf{I}_{d}$ in a set of positive measure by
Assumption~5.2, so that $\mathbf{\Psi }_{2}\left( \lambda \mathbf{;\theta }%
\right) =\mathbf{\Psi }_{2}\left( \lambda \mathbf{;\theta }_{0}\right) $
a.e. and $\mathcal{L}_{2}\left( \mathbf{\theta }\right) =0.$ }

{\small Then for those $\mathbf{\theta \neq \theta }_{0}$ for which $%
\mathcal{L}_{2}\left( \mathbf{\theta }\right) =0$, \textbf{$\Psi $}$_{k}(%
\mathbf{\lambda };\mathbf{\theta })\mathbf{S}_{k}\mathbf{\alpha }_{k}\left( 
\mathbf{\theta }\right) =\mathbf{\Psi }_{k}(\mathbf{\lambda };\mathbf{\theta 
}_{0})A^{\otimes k}(\mathbf{\lambda })\mathbf{S}_{k}\mathbf{\alpha }%
_{k}\left( \mathbf{\theta }\right) ,$ and therefore for $\mathcal{\hat{L}}%
_{k}\left( \mathbf{\theta }\right) =0$ to be true given Assumption~5.2, it
must hold that $A^{\otimes k}(\mathbf{\lambda })\mathbf{S}_{k}\mathbf{\alpha 
}_{k}\left( \mathbf{\theta }\right) =\mathbf{S}_{k}\mathbf{\alpha }_{k}^{0}$
a.e., which by Theorem~3 can only hold when $A(e^{-i\lambda })$ is a signed
permutation matrix $P_{d}$ a.e.. }

{\small Therefore, under Assumption~6A, all $\mathbf{\theta \neq \theta }%
_{0} $ such that $\mathbf{\Psi }\left( e^{-i\lambda }\mathbf{;\theta }%
\right) =\mathbf{\Psi }\left( e^{-i\lambda }\mathbf{;\theta }_{0}\right)
P_{d}$ with $P_{d}\neq \mathbf{I}_{d}$ are discarded as $\mathbf{\theta \not
\in }\mathcal{S}^{\max }$ because the product of the absolute value of the
diagonal elements of $\Theta _{0}\left( \mathbf{\theta }_{0}\right) $ is a
unique maximum up to permutations, and the consistency of $\mathbf{\hat{%
\theta}}_{k,T}$ follows by the standard argument. Alternatively,
Assumption~6C directly discard that for any $\mathbf{\theta \neq \theta }%
_{0} $ in $\mathcal{S}$ there exists a signed permutation $P_{d}$ so that $%
\mathbf{\Psi }\left( e^{-i\lambda }\mathbf{;\theta }\right) =\mathbf{\Psi }%
\left( e^{-i\lambda }\mathbf{;\theta }_{0}\right) P_{d}$ a.e. exists, so it
must hold that $\mathcal{\hat{L}}_{k}\left( \mathbf{\theta }\right) >0$. }

{\small Further, Lemma~\ref{LemmaA} shows that if $\mathbf{\Psi }\left(
e^{-i\lambda }\mathbf{;\theta }\right) =\mathbf{\Psi }\left( e^{-i\lambda }%
\mathbf{;\theta }_{0}\right) P_{d}$ a.e., then $\mathbf{\alpha }_{3}\left( 
\mathbf{\theta }\right) =P_{d}^{\prime }\mathbf{\alpha }_{3}^{0}\ $and $%
\mathbf{\alpha }_{4}\left( \mathbf{\theta }\right) =P_{d}^{+\prime }\mathbf{%
\alpha }_{4}^{0}$ (where $P_{d}^{+}$ is equal to $P_{d}$ with positive
entries)$,$ and therefore $\mathcal{L}_{k}\left( \mathbf{\theta ,\alpha }%
_{k}\left( \mathbf{\theta }\right) \right) =0$.\ But for such $\mathbf{%
\theta }\neq \mathbf{\theta }_{0},$ even if $\mathbf{\theta \in }\mathcal{S}%
^{+},$ it holds $\mathbf{\alpha }_{k}\left( \mathbf{\theta }\right) \neq 
\mathbf{\alpha }_{k}^{0}$ under Assumption~6B$\left( k\right) $ because $%
P_{d}\neq I_{d}$ and $P_{d}^{+}\neq I_{d}$ for $k=3$ and $4,$ respectively,
by Assumption~5.2, so that for such $\mathbf{\theta ,}$ $\mathbf{\alpha }%
_{k}\left( \mathbf{\theta }\right) \not \in \mathcal{D}_{k}$ and min$_{%
\mathbf{\alpha }\in \mathcal{D}_{k}}\mathcal{L}_{k}\left( \mathbf{\theta
,\alpha }\right) >0\ $because of the compactness of $\mathcal{D}_{k},$ and
therefore $\mathcal{L}_{2}\left( \mathbf{\theta }\right) +\mathcal{L}%
_{k}\left( \mathbf{\theta ,\alpha }\right) $ is uniquely miminized at $%
\left( \mathbf{\theta ,\alpha }\right) =\left( \mathbf{\theta }_{0},\mathbf{%
\alpha }_{0}\right) \ $in $\mathcal{S}^{+}\times \mathcal{D}_{k}.$ }

{\small Notice that the case where $P_{d}\neq \mathbf{I}_{d}$ is a pure
sign-flipping diagonal matrix so that $P_{d}^{+}=\mathbf{I}_{d}$ and $%
\mathbf{\Psi }_{4}(\mathbf{\lambda };\mathbf{\theta })=\mathbf{\Psi }_{4}(%
\mathbf{\lambda };\mathbf{\theta }_{0})$ is excluded by Assumption~6B$\left(
4\right) $ by imposing all the diagonal elements of $\mathbf{\Theta }%
_{0}\left( \mathbf{\theta }_{0}\right) $ being strictly positive, which is
not true for $\mathbf{\Psi }(e^{-i\lambda };\mathbf{\theta }_{0})P_{d}$ and
such diagonal $P_{d}.$ However, Assumption~6B$\left( 3\right) $ needs to
explicit prevent alternative orderings in $\mathbf{\alpha }_{3}$ due to sign
changes in $P_{d}^{\prime }\mathbf{\alpha }_{3}^{0}.$ }

{\small Finally, the consistency of $\mathbf{\hat{\alpha}}_{k,T}\left( 
\mathbf{\hat{\theta}}_{k,T}\right) \ $follows from the consistency of $%
\mathbf{\hat{\theta}}_{k,T}$ using similar methods as in VL and
Assumption~5.4.$\  \  \Box $\bigskip }

{\small \noindent \textbf{Proof of Theorem~\ref{Th6}}. The score of the
concentrated loss function wrt to each component of $\mathbf{\theta }$ is
given by, $k=3,4,$%
\begin{eqnarray*}
\frac{\partial }{\partial \mathbf{\theta }_{\ell }}\mathcal{\hat{L}}%
_{k,T}^{\dag }\left( \mathbf{\theta }\right) &=&\frac{2}{T^{k-1}}\sum_{%
\mathbf{\lambda }_{\mathbf{j}}}\func{Re}\left \{ \left( \mathbf{\Psi }%
_{k}\left( \mathbf{\lambda _{\mathbf{j}};\theta }\right) \mathbf{S}_{k}%
\mathbf{\hat{\alpha}}_{k,T}^{\dag }\left( \mathbf{\theta }\right) -\mathbb{I}%
_{k}(\mathbf{\lambda }_{\mathbf{j}})\right) ^{\ast }\mathbf{W}_{k}(\mathbf{%
\lambda }_{\mathbf{j}};\mathbf{\tilde{\theta}}_{T})\mathbf{\dot{\Psi}}%
_{k}^{\left( \ell \right) }(\mathbf{\lambda }_{\mathbf{j}};\mathbf{\theta })%
\mathbf{S}_{k}\mathbf{\hat{\alpha}}_{k,T}^{\dag }\left( \mathbf{\theta }%
\right) \right \} \\
&+&\frac{2}{T^{k-1}}\sum_{\mathbf{\lambda }_{\mathbf{j}}}\func{Re}\left \{
\left( \mathbf{\Psi }_{k}(\mathbf{\lambda }_{\mathbf{j}};\mathbf{\theta })%
\mathbf{S}_{k}\mathbf{\hat{\alpha}}_{k,T}^{\dag }\left( \mathbf{\theta }%
\right) -\mathbb{I}_{k}(\mathbf{\lambda }_{\mathbf{j}})\right) ^{\ast }%
\mathbf{W}_{k}(\mathbf{\lambda }_{\mathbf{j}};\mathbf{\tilde{\theta}}_{T})%
\mathbf{\Psi }_{k}(\mathbf{\lambda }_{\mathbf{j}};\mathbf{\theta })\mathbf{S}%
_{k}\frac{\partial }{\partial \mathbf{\theta }_{\ell }}\mathbf{\hat{\alpha}}%
_{k,T}^{\dag }\left( \mathbf{\theta }\right) \right \}
\end{eqnarray*}%
where the estimation effect of $\mathbf{\alpha }_{k}$ is similar to VL,%
\begin{equation*}
\frac{\partial }{\partial \mathbf{\theta }_{\ell }}\mathbf{\hat{\alpha}}%
_{k,T}^{\dag }\left( \mathbf{\theta }\right) =-\frac{1}{T^{k-1}}\sum_{%
\mathbf{\lambda }_{\mathbf{j}}}\func{Re}\left \{ \mathbf{S}_{k}^{\prime }%
\mathbf{\Psi }_{k}^{-1}(\mathbf{\lambda }_{\mathbf{j}};\mathbf{\theta })%
\mathbf{\dot{\Psi}}_{k}^{\left( \ell \right) }(\mathbf{\lambda }_{\mathbf{j}%
};\mathbf{\theta })\mathbf{\Psi }_{k}^{-1}(\mathbf{\lambda }_{\mathbf{j}};%
\mathbf{\theta })\mathbb{I}_{k}(\mathbf{\lambda }_{\mathbf{j}})\right \} ,
\end{equation*}%
which shows that estimation of $\mathbf{\theta }$ is not independent of
estimation of $\mathbf{\alpha }_{k}$ unlike with second order methods as we
can check that 
\begin{equation*}
\frac{\partial }{\partial \mathbf{\theta }_{\ell }}\mathbf{\hat{\alpha}}%
_{k,T}^{\dag }\left( \mathbf{\theta }_{0}\right) \rightarrow _{p}-\mathbf{S}%
_{k}^{\prime }\mathbf{\bar{\Lambda}}_{k}^{\left( \ell \right) }\left( 
\mathbf{\theta }_{0}\right) \mathbf{S}_{k}\mathbf{\alpha }_{0},
\end{equation*}%
which is also the limit of $\frac{\partial }{\partial \mathbf{\theta }_{\ell
}}\mathbf{\hat{\alpha}}_{k,T}^{EFF}\left( \mathbf{\theta }_{0}\right) ,$ so $%
\mathbf{\hat{\alpha}}_{k,T}^{\dag }\left( \mathbf{\theta }_{0}\right) $ and $%
\mathbf{\hat{\alpha}}_{k,T}^{EFF}\left( \mathbf{\theta }_{0}\right) $ share
the same asymptotic distribution because $\frac{\partial }{\partial \mathbf{%
\tilde{\theta}}_{T,\ell }}\mathbf{\hat{\alpha}}_{k,T}^{EFF}\left( \mathbf{%
\theta }_{0}\right) \rightarrow _{p}0.$ }

{\small Similarly for $k=2,$ we have%
\begin{equation*}
\frac{\partial }{\partial \mathbf{\theta }_{\ell }}\mathcal{L}_{2,T}\left( 
\mathbf{\theta }\right) =\frac{2}{T}\sum_{\lambda _{j}}\func{Re}\left \{
\left( \mathbf{\Psi }_{2}\left( \lambda _{j}\mathbf{;\theta }\right) \text{%
vec}\left( \mathbf{I}_{d}\right) -\mathbb{I}_{2}(\lambda _{j})\right) ^{\ast
}\mathbf{W}_{2}(\lambda _{j};\mathbf{\tilde{\theta}}_{T})\mathbf{\dot{\Psi}}%
_{2}^{\left( \ell \right) }(\lambda _{j};\mathbf{\theta )}\text{vec}\left( 
\mathbf{I}_{d}\right) \right \} .
\end{equation*}
}

{\small Then, using that $\mathbf{S}_{k}\mathbf{S}_{k}^{\prime }T^{1-k}\sum_{%
\mathbf{\lambda }_{\mathbf{j}}}\func{Re}\left \{ \mathbf{B}_{k}\left( 
\mathbf{\lambda }_{\mathbf{j}};\mathbf{\theta }_{0}\right) \right \} =%
\mathbf{S}_{k}\mathbf{S}_{k}^{\prime }\left( 2\pi \right) ^{1-k}\int \func{Re%
}\left \{ \mathbf{B}_{k}\left( \mathbf{\lambda };\mathbf{\theta }_{0}\right)
\right \} d\mathbf{\lambda }+O\left( T^{-1}\right) =O\left( T^{-1}\right) 
\mathbf{,\ }$because $\mathbf{S}_{k}\mathbf{S}_{k}^{\prime }\left( 2\pi
\right) ^{1-k}\int_{\Pi ^{k-1}}\func{Re}\left \{ \mathbf{B}_{k}\left( 
\mathbf{\lambda };\mathbf{\theta }_{0}\right) \right \} d\mathbf{\lambda }=0$
as $\mathbf{S}_{k}^{\prime }\mathbf{S}_{k}=\mathbf{I}_{d},$ up to $%
o_{p}\left( 1\right) $ terms, 
\begin{eqnarray*}
&&T^{1/2}\mathbf{\Sigma }\left( \mathbf{\theta }_{0},\mathbf{\alpha }%
_{0}\right) \left( \mathbf{\hat{\theta}}_{w,T}^{\dag }-\mathbf{\theta }%
_{0}\right) \\
&=&-\sum_{k=3}^{4}w_{k}\frac{T^{1/2}}{T^{k-1}}\left( \mathbf{I}_{m}\otimes 
\mathbf{\alpha }_{k}^{0\prime }\mathbf{S}_{k}^{\prime }\right) \sum_{\mathbf{%
\lambda }_{\mathbf{j}}}\func{Re}\left \{ \mathbf{B}_{k}^{\ast }\left( 
\mathbf{\lambda }_{\mathbf{j}};\mathbf{\theta }_{0}\right) \mathbf{\Psi }%
_{k}^{-1}(\mathbf{\lambda }_{\mathbf{j}};\mathbf{\theta }_{0})\left( \mathbf{%
\Psi }_{k}(\mathbf{\lambda }_{\mathbf{j}};\mathbf{\theta }_{0})\mathbf{S}_{k}%
\mathbf{\hat{\alpha}}_{k,T}^{\dag }\left( \mathbf{\theta }_{0}\right) -%
\mathbb{I}_{k}(\mathbf{\lambda }_{\mathbf{j}})\right) \right \} \\
&&-\frac{1}{T^{1/2}}\left( \mathbf{I}_{m}\otimes \text{vec}\left( \mathbf{I}%
_{d}\right) ^{\prime }\right) \sum_{\lambda _{j}}\func{Re}\left \{ \mathbf{B}%
_{2}^{\ast }\left( \lambda _{j}\mathbf{;\theta }_{0}\right) \mathbf{\Psi }%
_{2}^{-1}(\lambda _{j};\mathbf{\theta }_{0})\left( \mathbf{\Psi }_{2}\left(
\lambda _{j}\mathbf{;\theta }_{0}\right) \text{vec}\left( \mathbf{I}%
_{d}\right) -\mathbb{I}_{2}(\lambda _{j})\right) \right \} \\
&=&\sum_{k=3}^{4}w_{k}\frac{T^{1/2}}{T^{k-1}}\left( \mathbf{I}_{m}\otimes 
\mathbf{\alpha }_{k}^{0\prime }\mathbf{S}_{k}^{\prime }\right) \sum_{\mathbf{%
\lambda }_{\mathbf{j}}}\func{Re}\left \{ \mathbf{B}_{k}^{\ast }\left( 
\mathbf{\lambda }_{\mathbf{j}};\mathbf{\theta }_{0}\right) \left \{ \mathbf{%
\Psi }_{k}^{-1}(\mathbf{\lambda }_{\mathbf{j}};\mathbf{\theta }_{0})\mathbb{I%
}_{k}(\mathbf{\lambda }_{\mathbf{j}})-\mathbf{S}_{k}\mathbf{\hat{\alpha}}%
_{k,T}^{\dag }\left( \mathbf{\theta }_{0}\right) \right \} \right \} \\
&&+\frac{1}{T^{1/2}}\left( \mathbf{I}_{m}\otimes \text{vec}\left( \mathbf{I}%
_{d}\right) ^{\prime }\right) \sum_{\lambda _{j}}\func{Re}\left \{ \mathbf{B}%
_{2}^{\ast }\left( \lambda _{j}\mathbf{;\theta }_{0}\right) \left( \mathbf{%
\Psi }_{2}^{-1}(\lambda _{j};\mathbf{\theta }_{0})\mathbb{I}_{2}(\lambda
_{j})-\text{vec}\left( \mathbf{I}_{d}\right) \right) \right \}
\end{eqnarray*}
\begin{eqnarray*}
&=&\sum_{k=3}^{4}w_{k}\frac{T^{1/2}}{T^{k-1}}\left( \mathbf{I}_{m}\otimes 
\mathbf{\alpha }_{k}^{0\prime }\mathbf{S}_{k}^{\prime }\right) \sum_{\mathbf{%
\lambda }_{\mathbf{j}}}\func{Re}\left \{ \mathbf{B}_{k}^{\ast }\left( 
\mathbf{\lambda }_{\mathbf{j}};\mathbf{\theta }_{0}\right) \left \{ \mathbb{I%
}_{k}^{\mathbf{\varepsilon }}(\mathbf{\lambda }_{\mathbf{j}})-\mathbf{S}_{k}%
\mathbf{S}_{k}^{\prime }\frac{1}{T^{k-1}}\sum_{\mathbf{\lambda }_{\mathbf{j}%
}}\mathbb{I}_{k}^{\mathbf{\varepsilon }}(\mathbf{\lambda }_{\mathbf{j}%
})\right \} \right \} \\
&&+\frac{1}{T^{1/2}}\left( \mathbf{I}_{m}\otimes \text{vec}\left( \mathbf{I}%
_{d}\right) ^{\prime }\right) \sum_{\lambda _{j}}\func{Re}\left \{ \mathbf{B}%
_{2}^{\ast }\left( \lambda _{j}\mathbf{;\theta }_{0}\right) \left( \mathbb{I}%
_{2}^{\mathbf{\varepsilon }}(\lambda _{j})-\text{vec}\left( \mathbf{I}%
_{d}\right) \right) \right \} \\
&=&\sum_{k=3}^{4}w_{k}\frac{T^{1/2}}{T^{k-1}}\left( \mathbf{I}_{m}\otimes 
\mathbf{\alpha }_{k}^{0\prime }\mathbf{S}_{k}^{\prime }\right) \sum_{\mathbf{%
\lambda }_{\mathbf{j}}}\func{Re}\left \{ \mathbf{B}_{k}^{\ast }\left( 
\mathbf{\lambda }_{\mathbf{j}};\mathbf{\theta }_{0}\right) \left \{ \mathbb{I%
}_{k}^{\mathbf{\varepsilon }}(\mathbf{\lambda }_{\mathbf{j}})-E\left[ 
\mathbb{I}_{k}^{\mathbf{\varepsilon }}(\mathbf{\lambda }_{\mathbf{j}})\right]
\right \} \right \} \\
&&+\frac{1}{T^{1/2}}\left( \mathbf{I}_{m}\otimes \text{vec}\left( \mathbf{I}%
_{d}\right) ^{\prime }\right) \sum_{\lambda _{j}}\func{Re}\left \{ \mathbf{B}%
_{2}^{\ast }\left( \lambda _{j}\mathbf{;\theta }_{0}\right) \left( \mathbb{I}%
_{2}^{\mathbf{\varepsilon }}(\lambda _{j})-E\left[ \mathbb{I}_{2}^{\mathbf{%
\varepsilon }}(\lambda _{j})\right] \right) \right \} \\
& \rightarrow _{d} &N_{m}\left( 0,\mathbf{\delta }\left( \mathbf{\alpha }%
_{0}\right) \mathbf{\Omega }\left( \mathbf{\theta }_{0}\right) \mathbf{%
\delta }\left( \mathbf{\alpha }_{0}\right) ^{\prime }\right) ,
\end{eqnarray*}%
applying the CLT in Appendix~D to weighted sums of $\mathbb{I}_{2}^{\mathbf{%
\varepsilon }}(\lambda _{j})$ and $\mathbb{I}_{k}^{\mathbf{\varepsilon }}(%
\mathbf{\lambda }_{\mathbf{j}})$ as in VL, because for $k=3,4,$ $%
T^{1-k}\sum_{\mathbf{\lambda }_{\mathbf{j}}}\func{Re}\left \{ \mathbf{B}%
_{k}^{\ast }\left( \mathbf{\lambda }_{\mathbf{j}};\mathbf{\theta }%
_{0}\right) E\left[ \mathbb{I}_{k}^{\mathbf{\varepsilon }}(\mathbf{\lambda }%
_{\mathbf{j}})\right] \right \} =T^{1-k}\sum_{\mathbf{\lambda }_{\mathbf{j}}}%
\func{Re}\left \{ \mathbf{B}_{k}^{\ast }\left( \mathbf{\lambda }_{\mathbf{j}%
};\mathbf{\theta }_{0}\right) \right \} \mathbf{S}_{k}\mathbf{\alpha }_{0}=%
\mathbf{\bar{\Lambda}}_{k}^{\prime }\left( \mathbf{\theta }\right) \mathbf{S}%
_{k}\mathbf{\alpha }_{0}-\mathbf{\bar{\Lambda}}_{k}^{\prime }\left( \mathbf{%
\theta }\right) \mathbf{S}_{k}\mathbf{S}_{k}^{\prime }\mathbf{S}_{k}\mathbf{%
\alpha }_{0}+O\left( T^{-1}\right) =O\left( T^{-1}\right) $ as $\mathbf{S}%
_{k}\mathbf{S}_{k}^{\prime }\mathbf{S}_{k}=\mathbf{S}_{k}.$ The proof is
completed using Lemma~\ref{Hessian} for the convergence of the Hessian. $%
\mathbf{\  \  \Box }$\bigskip }

{\small \noindent \textbf{Proof of Theorem~\ref{Th7}}. We can write%
\begin{equation*}
\mathbf{\hat{\alpha}}_{k,T}^{\dag }\left( \mathbf{\hat{\theta}}_{k,T}^{\dag
}\right) -\mathbf{\alpha }_{k}^{0}=\mathbf{\hat{\alpha}}_{k,T}^{\dag }\left( 
\mathbf{\hat{\theta}}_{w,T}^{\dag }\right) -\mathbf{\hat{\alpha}}%
_{k,T}^{\dag }\left( \mathbf{\theta }_{0}\right) +\mathbf{\hat{\alpha}}%
_{k,T}^{\dag }\left( \mathbf{\theta }_{0}\right) -\mathbf{\alpha }_{k}^{0}
\end{equation*}%
where%
\begin{equation*}
\mathbf{\hat{\alpha}}_{k,T}^{\dag }\left( \mathbf{\hat{\theta}}_{k,T}^{\dag
}\right) -\mathbf{\hat{\alpha}}_{k,T}^{\dag }\left( \mathbf{\theta }%
_{0}\right) =\frac{\partial }{\partial \mathbf{\theta }^{\prime }}\mathbf{%
\hat{\alpha}}_{k,T}^{\dag }\left( \mathbf{\theta }_{T}\right) \left( \mathbf{%
\hat{\theta}}_{w,T}^{\dag }-\mathbf{\theta }_{0}\right)
\end{equation*}%
for some $\mathbf{\theta }_{T}\rightarrow _{p}\mathbf{\theta }_{0}$ and%
\begin{eqnarray*}
T^{1/2}\left( \mathbf{\hat{\alpha}}_{k,T}^{\dag }\left( \mathbf{\theta }%
_{0}\right) -\mathbf{\alpha }_{k}^{0}\right) &=&\frac{T^{1/2}}{T^{k-1}}\sum_{%
\mathbf{\lambda }_{\mathbf{j}}}\func{Re}\left \{ \mathbf{S}_{k}^{\prime }%
\mathbf{\Psi }_{k}^{-1}(\mathbf{\lambda }_{\mathbf{j}};\mathbf{\theta }_{0})%
\mathbb{I}_{k}(\mathbf{\lambda }_{\mathbf{j}})-\mathbf{\alpha }%
_{k}^{0}\right \} \\
&=&\frac{T^{1/2}}{T^{k-1}}\sum_{\mathbf{\lambda }_{\mathbf{j}}}\mathbf{S}%
_{k}^{\prime }\func{Re}\left \{ \mathbb{I}_{k}^{\mathbf{\varepsilon }}(%
\mathbf{\lambda }_{\mathbf{j}})-E\left[ \mathbb{I}_{k}^{\mathbf{\varepsilon }%
}(\mathbf{\lambda }_{\mathbf{j}})\right] \right \} +o_{p}\left( 1\right) .
\end{eqnarray*}%
Then%
\begin{eqnarray*}
\frac{\partial }{\partial \mathbf{\theta }_{\ell }}\mathbf{\hat{\alpha}}%
_{k,T}^{\dag }\left( \mathbf{\theta }_{T}\right) &=&-\frac{1}{T^{k-1}}\sum_{%
\mathbf{\lambda }_{\mathbf{j}}}\func{Re}\left \{ \mathbf{S}_{k}^{\prime }%
\mathbf{\Psi }_{k}^{-1}(\mathbf{\lambda }_{\mathbf{j}};\mathbf{\theta }_{T})%
\mathbf{\dot{\Psi}}_{k}^{\left( \ell \right) }(\mathbf{\lambda }_{\mathbf{j}%
};\mathbf{\theta }_{T})\mathbf{\Psi }_{k}^{-1}(\mathbf{\lambda }_{\mathbf{j}%
};\mathbf{\theta }_{T})\mathbb{I}_{k}(\mathbf{\lambda }_{\mathbf{j}})\right
\} \\
& \rightarrow _{p} &-\frac{1}{T^{k-1}}\sum_{\mathbf{\lambda }_{\mathbf{j}}}%
\func{Re}\left \{ \mathbf{S}_{k}^{\prime }\mathbf{\Psi }_{k}^{-1}(\mathbf{%
\lambda }_{\mathbf{j}};\mathbf{\theta }_{0})\mathbf{\dot{\Psi}}_{k}^{\left(
\ell \right) }(\mathbf{\lambda }_{\mathbf{j}};\mathbf{\theta }_{0})\mathbf{%
\Psi }_{k}^{-1}(\mathbf{\lambda }_{\mathbf{j}};\mathbf{\theta }_{0})\mathbb{I%
}_{k}(\mathbf{\lambda }_{\mathbf{j}})\right \} \\
& \rightarrow _{p} &-\frac{1}{T^{k-1}}\sum_{\mathbf{\lambda }_{\mathbf{j}}}%
\func{Re}\left \{ \mathbf{S}_{k}^{\prime }\mathbf{\Psi }_{k}^{-1}(\mathbf{%
\lambda }_{\mathbf{j}};\mathbf{\theta }_{0})\mathbf{\dot{\Psi}}_{k}^{\left(
\ell \right) }(\mathbf{\lambda }_{\mathbf{j}};\mathbf{\theta }_{0})\mathbb{I}%
_{k}^{\mathbf{\varepsilon }}(\mathbf{\lambda }_{\mathbf{j}})\right \} \\
& \rightarrow _{p} &-\mathbf{S}_{k}^{\prime }\mathbf{\bar{\Lambda}}%
_{k}^{\left( \ell \right) }\left( \mathbf{\theta }_{0}\right) \mathbf{S}_{k}%
\mathbf{\alpha }_{k}^{0}
\end{eqnarray*}%
so that%
\begin{equation*}
\frac{\partial }{\partial \mathbf{\theta }^{\prime }}\mathbf{\hat{\alpha}}%
_{k,T}^{\dag }\left( \mathbf{\theta }_{T}\right) \rightarrow _{p}-\mathbf{S}%
_{k}^{\prime }\mathbf{\bar{\Lambda}}_{k}\left( \mathbf{\theta }_{0}\right)
\left( \mathbf{I}_{m}\otimes \mathbf{S}_{k}\mathbf{\alpha }_{k}^{0}\right) .
\end{equation*}%
Then, pooling all results, with $ \mathbf{\Sigma }_0 := \mathbf{\Sigma }\left( \mathbf{\theta }_{0},%
\mathbf{\alpha }_{k}^0 \right), $
\begin{eqnarray*}
&&T^{1/2}\left( \mathbf{\hat{\alpha}}_{k,T}^{\dag }\left( \mathbf{\hat{\theta%
}}_{w,T}^{\dag }\right) -\mathbf{\alpha }_{k}^{0}\right) \\
&=&-T^{1/2}\mathbf{S}_{k}^{\prime }\mathbf{\bar{\Lambda}}_{k}\left( \mathbf{%
\theta }_{0}\right) \left( \mathbf{I}_{m}\otimes \mathbf{S}_{k}\mathbf{%
\alpha }_{k}^{0}\right) \mathbf{\Sigma }^{-1}_0 \sum_{h=3}^{4}\frac{w_{h}}{T^{h-1}}\left( 
\mathbf{I}_{m}\otimes \mathbf{\alpha }_{h}^{0\prime }\mathbf{S}_{hk}^{\prime
}\right) \sum_{\mathbf{\lambda }_{\mathbf{j}}}\func{Re}\left \{ \mathbf{B}%
_{h}^{\ast }\left( \mathbf{\lambda }_{\mathbf{j}};\mathbf{\theta }%
_{0}\right) \left \{ \mathbb{I}_{h}^{\mathbf{\varepsilon }}(\mathbf{\lambda }%
_{\mathbf{j}})-E\left[ \mathbb{I}_{h}^{\mathbf{\varepsilon }}(\mathbf{%
\lambda }_{\mathbf{j}})\right] \right \} \right \} \\
&&-\frac{1}{T^{1/2}}\mathbf{S}_{k}^{\prime }\mathbf{\bar{\Lambda}}_{k}\left( 
\mathbf{\theta }_{0}\right) \left( \mathbf{I}_{m}\otimes \mathbf{S}_{k}%
\mathbf{\alpha }_{k}^{0}\right) \mathbf{\Sigma }^{-1}_0 \left( \mathbf{I}_{m}\otimes \text{vec}%
\left( \mathbf{I}_{d}\right) \right) ^{\prime }\sum_{\lambda _{j}}\func{Re}%
\left \{ \mathbf{B}_{2}^{\ast }\left( \lambda _{j}\mathbf{;\theta }%
_{0}\right) \left( \mathbb{I}_{2}^{\mathbf{\varepsilon }}(\lambda _{j})-E%
\left[ \mathbb{I}_{2}^{\mathbf{\varepsilon }}(\lambda _{j})\right] \right)
\right \} \\
&&+\frac{T^{1/2}}{T^{k-1}}\sum_{\mathbf{\lambda }_{\mathbf{j}}}\mathbf{S}%
_{k}^{\prime }\func{Re}\left \{ \mathbb{I}_{k}^{\mathbf{\varepsilon }}(%
\mathbf{\lambda }_{\mathbf{j}})-E\left[ \mathbb{I}_{k}^{\mathbf{\varepsilon }%
}(\mathbf{\lambda }_{\mathbf{j}})\right] \right \} +o_{p}\left( 1\right) \\
&=&T^{1/2}\sum_{h=3}^{4}\frac{1}{T^{h-1}}\sum_{\mathbf{\lambda }_{\mathbf{j}%
}}\func{Re}\left \{ \mathbf{S}_{k}^{\prime }\mathbf{D}_{k,h}^{\ast }\left( 
\mathbf{\lambda }_{\mathbf{j}};\mathbf{\theta }_{0}\right) \left \{ \mathbb{I%
}_{h}^{\mathbf{\varepsilon }}(\mathbf{\lambda }_{\mathbf{j}})-E\left[ 
\mathbb{I}_{h}^{\mathbf{\varepsilon }}(\mathbf{\lambda }_{\mathbf{j}})\right]
\right \} \right \} \\
&&+\frac{1}{T^{1/2}}\sum_{\lambda _{j}}\func{Re}\left \{ \mathbf{S}%
_{k}^{\prime }\mathbf{D}_{k,2}^{\ast }\left( \lambda _{j}\mathbf{;\theta }%
_{0}\right) \left( \mathbb{I}_{2}^{\mathbf{\varepsilon }}(\lambda _{j})-E%
\left[ \mathbb{I}_{2}^{\mathbf{\varepsilon }}(\lambda _{j})\right] \right)
\right \} +o_{p}\left( 1\right)
\end{eqnarray*}%
and the result follows as Theorem~\ref{Th6} and Appendix~D for $k=2,3,$
while for $k=4$ we have to consider the extra term $\mathbf{\eta }_{t}$
coming for the decomposition of%
\begin{eqnarray*}
\frac{1}{T^{3}}\sum_{\mathbf{\lambda }_{\mathbf{j}}}\mathbb{I}_{4}^{\mathbf{%
\varepsilon }}(\mathbf{\lambda }_{\mathbf{j}}) &=&\frac{1}{T}\sum_{t=1}^{T}%
\mathbf{\varepsilon }_{t}^{\otimes 4}-\frac{1}{T^{2}}\sum_{t=1}^{T}%
\sum_{r=1}^{T}\left[ 
\begin{array}{c}
\mathbf{\varepsilon }_{t}\otimes \mathbf{\varepsilon }_{t}\otimes \mathbf{%
\varepsilon }_{r}\otimes \mathbf{\varepsilon }_{r} \\ 
+\mathbf{\varepsilon }_{t}\otimes \mathbf{\varepsilon }_{r}\otimes \mathbf{%
\varepsilon }_{t}\otimes \mathbf{\varepsilon }_{r} \\ 
+\mathbf{\varepsilon }_{t}\otimes \mathbf{\varepsilon }_{r}\otimes \mathbf{%
\varepsilon }_{r}\otimes \mathbf{\varepsilon }_{t}%
\end{array}%
\right] +O_{p}\left( T^{-1}\right) \\
&=&\frac{1}{T}\sum_{t=1}^{T}\mathbf{\varepsilon }_{t}^{\otimes 4}-\frac{1}{T}%
\sum_{t=1}^{T}\left( \mathbf{1}_{6}^{\prime }\otimes \mathbf{I}%
_{d^{4}}\right) \mathbf{\eta }_{t}+E_{4}^{\left( a\right) }+E_{4}^{\left(
b\right) }+E_{4}^{\left( c\right) }+O_{p}\left( T^{-1}\right)
\end{eqnarray*}%
because the sum $\mathbf{\lambda }_{\mathbf{j}}$ does not include the terms $%
j_{a}=0\func{mod}T$, e.g., for $E_{4}^{\left( a\right) }:=E\left[ \mathbf{%
\varepsilon }_{t}\otimes \mathbf{\varepsilon }_{t}\otimes \mathbf{%
\varepsilon }_{r}\otimes \mathbf{\varepsilon }_{r}\right] =E\left[ \mathbf{%
\varepsilon }_{t}\otimes \mathbf{\varepsilon }_{t}\right] \otimes E\left[ 
\mathbf{\varepsilon }_{r}\otimes \mathbf{\varepsilon }_{r}\right]
=\sum_{ab}\left( \mathbf{e}_{a}\otimes \mathbf{e}_{a}\otimes \mathbf{e}%
_{b}\otimes \mathbf{e}_{b}\right) :=E_{2}\otimes E_{2},$ say, for $t\neq r,$
we find that 
\begin{eqnarray*}
&&\frac{1}{T^{2}}\sum_{t=1}^{T}\sum_{r=1}^{T}\left( \left( \mathbf{%
\varepsilon }_{t}\otimes \mathbf{\varepsilon }_{t}\pm E_{2}\right) \otimes
\left( \mathbf{\varepsilon }_{r}\otimes \mathbf{\varepsilon }_{r}\pm
E_{2}\right) \right) \\
&=&\frac{1}{T}\sum_{t=1}^{T}\left( \mathbf{\varepsilon }_{t}\otimes \mathbf{%
\varepsilon }_{t}-E_{2}\right) \otimes \frac{1}{T}\sum_{r=1}^{T}\left( 
\mathbf{\varepsilon }_{r}\otimes \mathbf{\varepsilon }_{r}-E_{2}\right)
+E_{4}^{\left( a\right) } \\
&&+\frac{1}{T^{2}}\sum_{t=1}^{T}\sum_{r=1}^{T}\left( \left( \mathbf{%
\varepsilon }_{t}\otimes \mathbf{\varepsilon }_{t}-E_{2}\right) \otimes
E_{2}\right) +\frac{1}{T^{2}}\sum_{t=1}^{T}\sum_{r=1}^{T}\left( E_{2}\otimes
\left( \mathbf{\varepsilon }_{r}\otimes \mathbf{\varepsilon }%
_{r}-E_{2}\right) \right) \\
&=&\frac{1}{T^{2}}\sum_{t=1}^{T}\sum_{r=1}^{T}\left( \mathbf{\varepsilon }%
_{t}\otimes \mathbf{\varepsilon }_{t}\otimes E_{2}\right) +\frac{1}{T^{2}}%
\sum_{t=1}^{T}\sum_{r=1}^{T}\left( E_{2}\otimes \mathbf{\varepsilon }%
_{r}\otimes \mathbf{\varepsilon }_{r}\right) -E_{4}^{\left( a\right)
}+O_{p}\left( T^{-1}\right)
\end{eqnarray*}%
where $\left( \mathbf{\varepsilon }_{t}\otimes \mathbf{\varepsilon }%
_{t}\otimes E_{2}\right) =E_{r}\left[ \mathbf{\varepsilon }_{t}\otimes 
\mathbf{\varepsilon }_{t}\otimes \mathbf{\varepsilon }_{r}\otimes \mathbf{%
\varepsilon }_{r}\right] ,$ and we can proceed similarly for $E_{4}^{\left(
b\right) }:=E\left[ \mathbf{\varepsilon }_{t}\otimes \mathbf{\varepsilon }%
_{r}\otimes \mathbf{\varepsilon }_{t}\otimes \mathbf{\varepsilon }_{r}\right]
$ and $E_{4}^{\left( c\right) }:=E\left[ \mathbf{\varepsilon }_{t}\otimes 
\mathbf{\varepsilon }_{r}\otimes \mathbf{\varepsilon }_{r}\otimes \mathbf{%
\varepsilon }_{t}\right] ,$ where these constant terms cancel with the
expectation of $\mathbb{I}_{h}^{\mathbf{\varepsilon }}(\mathbf{\lambda }_{%
\mathbf{j}})\ $and do not contribute to the variance of $\  \mathbf{\hat{%
\alpha}}_{4,T}^{\dag }.\  \  \Box $\bigskip }

\section{\protect \small Appendix C: Auxiliary results}

\begin{lemma}
{\small \label{LemmaA} Under the conditions of Theorem~\ref{Th5}, for any
signed permutation matrix $P_{d}$, if $\mathbf{\Psi }\left( e^{-i\lambda }%
\mathbf{;\theta }\right) =\mathbf{\Psi }\left( e^{-i\lambda }\mathbf{;\theta 
}_{0}\right) P_{d}$ a.e., then $\mathbf{\alpha }_{3}\left( \mathbf{\theta }%
\right) =P_{d}^{\prime }\mathbf{\alpha }_{3}^{0}\ $and $\mathbf{\alpha }%
_{4}\left( \mathbf{\theta }\right) =P_{d}^{+\prime }\mathbf{\alpha }_{4}^{0}$%
.\bigskip }
\end{lemma}

{\small \noindent \textbf{Proof of Lemma~\ref{LemmaA}}. We can show using $%
P_{d}P_{d}^{\prime }=\mathbf{I}_{d}$ and $P_{d}^{+}P_{d}^{+\prime }=\mathbf{I%
}_{d}$ 
\begin{eqnarray*}
\mathbf{\alpha }_{k}\left( \mathbf{\theta }\right) &=&\left( \int_{\Pi
^{k-1}}\func{Re}\left \{ \mathbf{S}_{k}^{\prime }P_{d}^{\otimes k\prime }%
\mathbf{\Psi }_{k}^{\ast }(\mathbf{\lambda };\mathbf{\theta }_{0})\mathbf{%
\Psi }_{k}(\mathbf{\lambda };\mathbf{\theta }_{0})P_{d}^{\otimes k}\mathbf{S}%
_{k}\right \} d\mathbf{\lambda }\right) ^{-1} \\
&& \times \int_{\Pi ^{k-1}}\func{Re}\left \{ \mathbf{S}_{k}^{\prime
}P_{d}^{\otimes k\prime }\mathbf{\Psi }_{k}^{\ast }(\mathbf{\lambda };%
\mathbf{\theta }_{0})\mathbf{\Psi }_{k}(\mathbf{\lambda };\mathbf{\theta }%
_{0})\mathbf{S}_{k}\right \} d\mathbf{\lambda \ }P_{d}P_{d}^{\prime }\mathbf{%
\alpha }_{k}^{0}
\end{eqnarray*}
is equal to $P_{d}^{\prime }\mathbf{\alpha }_{k}^{0}$ or $P_{d}^{+\prime }%
\mathbf{\alpha }_{k}^{0}$ iff%
\begin{equation*}
P_{d}^{\otimes 3}\mathbf{S}_{3}=\mathbf{S}_{3}P_{d}\  \  \text{or\  \ }P_{4}%
\mathbf{S}_{4}=\mathbf{S}_{4}P_{d}^{+},\  \  \text{respectively.}
\end{equation*}%
Then writing $P_{d}=\left( s_{j_{1}}\mathbf{e}_{j_{1}},\ldots ,s_{j_{d}}%
\mathbf{e}_{j_{d}}\right) $ and $P_{d}=\left( \mathbf{e}_{j_{1}},\ldots ,%
\mathbf{e}_{j_{d}}\right) $ for $j_{a}\in \left \{ 1,\ldots ,d\right \} ,$ $%
j_{a}\neq j_{b}$ for $a\neq b,$ and $s_{j_{a}}=\pm 1,$ we notice that%
\begin{eqnarray*}
\mathbf{S}_{k}P_{d} &=&\left( \mathbf{e}_{1}^{\otimes k},\ldots ,\mathbf{e}%
_{d}^{\otimes k}\right) \left( s_{j_{1}}\mathbf{e}_{j_{1}},\ldots ,s_{j_{d}}%
\mathbf{e}_{j_{d}}\right) =\left( s_{j_{1}}\mathbf{e}_{j_{1}}^{\otimes
k},\ldots ,s_{j_{d}}\mathbf{e}_{j_{d}}^{\otimes k}\right) \\
\mathbf{S}_{k}P_{d}^{+} &=&\left( \mathbf{e}_{1}^{\otimes k},\ldots ,\mathbf{%
e}_{d}^{\otimes k}\right) \left( \mathbf{e}_{j_{1}},\ldots ,\mathbf{e}%
_{j_{d}}\right) =\left( \mathbf{e}_{j_{1}}^{\otimes k},\ldots ,\mathbf{e}%
_{j_{d}}^{\otimes k}\right)
\end{eqnarray*}%
is a reordering of the columns of $\mathbf{S}_{k}$ with the appropriate
sign, while%
\begin{equation*}
P_{d}^{\otimes k}\mathbf{S}_{k}=P_{d}^{\otimes k}\left( \mathbf{e}%
_{1}^{\otimes k},\ldots ,\mathbf{e}_{d}^{\otimes k}\right) =\left[ \left(
P_{d}\mathbf{e}_{1}\right) ^{\otimes k},\ldots ,\left( P_{d}\mathbf{e}%
_{d}\right) ^{\otimes k}\right] =\left[ s_{j_{1}}^{k}\mathbf{e}%
_{j_{1}}^{\otimes k},\ldots ,s_{j_{d}}^{k}\mathbf{e}_{j_{d}}^{\otimes k}%
\right] ,
\end{equation*}%
so that $P_{d}^{\otimes 3}\mathbf{S}_{3}=\left[ s_{j_{1}}\mathbf{e}%
_{j_{1}}^{\otimes 3},\ldots ,s_{j_{d}}\mathbf{e}_{j_{d}}^{\otimes 3}\right] =%
\mathbf{S}_{3}P_{d}$ and $P_{d}^{\otimes 4}\mathbf{S}_{4}=\left[ \mathbf{e}%
_{j_{1}}^{\otimes 4},\ldots ,\mathbf{e}_{j_{d}}^{\otimes 4}\right] =\mathbf{S%
}_{4}P_{d}^{+}.\  \  \Box $\bigskip }

\begin{lemma}
{\small \label{Hessian} Under the Assumptions of Theorem~\ref{Th6}, for $%
\tilde{\mathbf{\theta}}_{T}\rightarrow _{p}\mathbf{\theta }_{0},$ $k=3,4,$ 
\begin{eqnarray*}
\frac{\partial ^{2}}{\partial \mathbf{\theta }\partial \mathbf{\theta }%
^{\prime }}\mathcal{\hat{L}}_{k,T}^{\dag }\left( \mathbf{\theta }_{T}\right)
& \rightarrow _{p} &\left( \mathbf{I}_{m}\otimes \mathbf{S}_{k}\mathbf{%
\alpha }_{k}^{0}\right) ^{\prime }\mathbf{H}_{k}\left( \mathbf{\theta }%
_{0}\right) \left( \mathbf{I}_{m}\otimes \mathbf{S}_{k}\mathbf{\alpha }%
_{k}^{0}\right) \\
\frac{\partial ^{2}}{\partial \mathbf{\theta }\partial \mathbf{\theta }%
^{\prime }}\mathcal{L}_{2,T}\left( \mathbf{\theta }_{T}\right)
&\rightarrow_{p} &\left( \mathbf{I}_{m}\otimes \text{vec}\left( \mathbf{I}%
_{d}\right) \right) ^{\prime }\mathbf{H}_{2}\left( \mathbf{\theta }%
_{0}\right) \left( \mathbf{I}_{m}\otimes \text{vec}\left( \mathbf{I}%
_{d}\right) \right) .
\end{eqnarray*}
}
\end{lemma}

{\small \noindent \textbf{Proof of Lemma~\ref{Hessian}}. We give the proof
for $k=3,4,$ the case for $k=2$ is much simpler,%
\begin{eqnarray*}
&&\frac{\partial ^{2}}{\partial \mathbf{\theta }_{\ell }\partial \mathbf{%
\theta }_{p}}\mathcal{\hat{L}}_{k,T}^{\dag }\left( \mathbf{\theta }\right) \\
&=&\frac{2}{T^{k-1}}\sum_{\mathbf{\lambda }_{\mathbf{j}}}\func{Re}\left \{
\left( \mathbf{\Psi }_{k}(\mathbf{\lambda }_{\mathbf{j}};\mathbf{\theta })%
\mathbf{S}_{k}\mathbf{\hat{\alpha}}_{k,T}^{\dag }\left( \mathbf{\theta }%
\right) -\mathbb{I}_{k}(\mathbf{\lambda }_{\mathbf{j}})\right) ^{\ast }%
\mathbf{W}_{k}(\mathbf{\lambda }_{\mathbf{j}};\mathbf{\tilde{\theta}}_{T})%
\mathbf{\dot{\Psi}}_{k}^{\left( \ell ,p\right) }(\mathbf{\lambda }_{\mathbf{j%
}};\mathbf{\theta })\mathbf{S}_{k}\mathbf{\hat{\alpha}}_{k,T}^{\dag }\left( 
\mathbf{\theta }\right) \right \} \\
&&+\frac{2}{T^{k-1}}\sum_{\mathbf{\lambda }_{\mathbf{j}}}\func{Re}\left \{
\left( \mathbf{\dot{\Psi}}_{k}^{\left( p\right) }(\mathbf{\lambda }_{\mathbf{%
j}};\mathbf{\theta })\mathbf{S}_{k}\mathbf{\hat{\alpha}}_{k,T}^{\dag }\left( 
\mathbf{\theta }\right) \right) ^{\ast }\mathbf{W}_{k}(\mathbf{\lambda }_{%
\mathbf{j}};\mathbf{\tilde{\theta}}_{T})\mathbf{\dot{\Psi}}_{k}^{\left( \ell
\right) }(\mathbf{\lambda }_{\mathbf{j}};\mathbf{\theta })\mathbf{S}_{k}%
\mathbf{\hat{\alpha}}_{k,T}^{\dag }\left( \mathbf{\theta }\right) \right \}
\\
&&+\frac{2}{T^{k-1}}\sum_{\mathbf{\lambda }_{\mathbf{j}}}\func{Re}\left \{
\left( \mathbf{\Psi }_{k}(\mathbf{\lambda }_{\mathbf{j}};\mathbf{\theta })%
\mathbf{S}_{k}\frac{\partial }{\partial \mathbf{\theta }_{p}}\mathbf{\hat{%
\alpha}}_{k,T}^{\dag }\left( \mathbf{\theta }\right) \right) ^{\ast }\mathbf{%
W}_{k}(\mathbf{\lambda }_{\mathbf{j}};\mathbf{\tilde{\theta}}_{T})\mathbf{%
\dot{\Psi}}_{k}^{\left( \ell \right) }(\mathbf{\lambda }_{\mathbf{j}};%
\mathbf{\theta })\mathbf{S}_{k}\mathbf{\hat{\alpha}}_{k,T}^{\dag }\left( 
\mathbf{\theta }\right) \right \} \\
&&+\frac{2}{T^{k-1}}\sum_{\mathbf{\lambda }_{\mathbf{j}}}\func{Re}\left \{
\left( \mathbf{\Psi }_{k}(\mathbf{\lambda }_{\mathbf{j}};\mathbf{\theta })%
\mathbf{S}_{k}\mathbf{\hat{\alpha}}_{k,T}^{\dag }\left( \mathbf{\theta }%
\right) -\mathbb{I}_{k}(\mathbf{\lambda }_{\mathbf{j}})\right) ^{\ast }%
\mathbf{W}_{k}(\mathbf{\lambda }_{\mathbf{j}};\mathbf{\tilde{\theta}}_{T})%
\mathbf{\dot{\Psi}}_{k}^{\left( \ell \right) }(\mathbf{\lambda }_{\mathbf{j}%
};\mathbf{\theta })\mathbf{S}_{k}\frac{\partial }{\partial \mathbf{\theta }%
_{p}}\mathbf{\hat{\alpha}}_{k,T}^{\dag }\left( \mathbf{\theta }\right)
\right \} \\
&&+\frac{2}{T^{k-1}}\sum_{\mathbf{\lambda }_{\mathbf{j}}}\func{Re}\left \{
\left( \mathbf{\Psi }_{k}(\mathbf{\lambda }_{\mathbf{j}};\mathbf{\theta })%
\mathbf{S}_{k}\mathbf{\hat{\alpha}}_{k,T}^{\dag }\left( \mathbf{\theta }%
\right) -\mathbb{I}_{k}(\mathbf{\lambda }_{\mathbf{j}})\right) ^{\ast }%
\mathbf{W}_{k}(\mathbf{\lambda }_{\mathbf{j}};\mathbf{\tilde{\theta}}_{T})%
\mathbf{\Psi }_{k}(\mathbf{\lambda }_{\mathbf{j}};\mathbf{\theta })\mathbf{S}%
_{k}\frac{\partial ^{2}}{\partial \mathbf{\theta }_{\ell }\partial \mathbf{%
\theta }_{p}}\mathbf{\hat{\alpha}}_{k,T}^{\dag }\left( \mathbf{\theta }%
\right) \right \} \\
&&+\frac{2}{T^{k-1}}\sum_{\mathbf{\lambda }_{\mathbf{j}}}\func{Re}\left \{
\left( \mathbf{\Psi }_{k}(\mathbf{\lambda }_{\mathbf{j}};\mathbf{\theta })%
\mathbf{S}_{k}\mathbf{\hat{\alpha}}_{k,T}^{\dag }\left( \mathbf{\theta }%
\right) -\mathbb{I}_{k}(\mathbf{\lambda }_{\mathbf{j}})\right) ^{\ast }%
\mathbf{W}_{k}(\mathbf{\lambda }_{\mathbf{j}};\mathbf{\tilde{\theta}}_{T})%
\mathbf{\dot{\Psi}}_{k}^{(p)}(\mathbf{\lambda }_{\mathbf{j}};\mathbf{\theta }%
)\mathbf{S}_{k}\frac{\partial }{\partial \mathbf{\theta }_{\ell }}\mathbf{%
\hat{\alpha}}_{k,T}^{\dag }\left( \mathbf{\theta }\right) \right \} \\
&&+\frac{2}{T^{k-1}}\sum_{\mathbf{\lambda }_{\mathbf{j}}}\func{Re}\left \{
\left( \mathbf{\Psi }_{k}(\mathbf{\lambda }_{\mathbf{j}};\mathbf{\theta })%
\mathbf{S}_{k}\frac{\partial }{\partial \mathbf{\theta }_{p}}\mathbf{\hat{%
\alpha}}_{k,T}^{\dag }\left( \mathbf{\theta }\right) \right) ^{\ast }\mathbf{%
W}_{k}(\mathbf{\lambda }_{\mathbf{j}};\mathbf{\tilde{\theta}}_{T})\mathbf{%
\Psi }_{k}(\mathbf{\lambda }_{\mathbf{j}};\mathbf{\theta })\mathbf{S}_{k}%
\frac{\partial }{\partial \mathbf{\theta }_{\ell }}\mathbf{\hat{\alpha}}%
_{k,T}^{\dag }\left( \mathbf{\theta }\right) \right \} \\
&&+\frac{2}{T^{k-1}}\sum_{\mathbf{\lambda }_{\mathbf{j}}}\func{Re}\left \{
\left( \mathbf{\dot{\Psi}}_{k}^{\left( p\right) }(\mathbf{\lambda }_{\mathbf{%
j}};\mathbf{\theta })\mathbf{S}_{k}\mathbf{\hat{\alpha}}_{k,T}^{\dag }\left( 
\mathbf{\theta }\right) \right) ^{\ast }\mathbf{W}_{k}(\mathbf{\lambda }_{%
\mathbf{j}};\mathbf{\tilde{\theta}}_{T})\mathbf{\Psi }_{k}(\mathbf{\lambda }%
_{\mathbf{j}};\mathbf{\theta })\mathbf{S}_{k}\frac{\partial }{\partial 
\mathbf{\theta }_{\ell }}\mathbf{\hat{\alpha}}_{k,T}^{\dag }\left( \mathbf{%
\theta }\right) \right \}
\end{eqnarray*}%
where the limits of lines 1, 4, 5 and 6 of the rhs do not contribute
asymptotically when evaluated at $\tilde{\mathbf{\theta }}_{T} \rightarrow_{p} \mathbf{\mathbf{\theta }}_{0}$ and the other ones converge to the
probability limit of%
\begin{eqnarray*}
&&\frac{\partial ^{2}}{\partial \mathbf{\theta }_{\ell }\partial \mathbf{%
\theta }_{p}}\mathcal{\hat{L}}_{k,T}^{\dag }\left( \mathbf{\theta }%
_{0}\right) \\
&=&\frac{2}{T^{k-1}}\sum_{\mathbf{\lambda }_{\mathbf{j}}}\func{Re}\left \{
\left( \mathbf{\dot{\Psi}}_{k}^{(p)}(\mathbf{\lambda }_{\mathbf{j}};\mathbf{%
\theta }_{0})\mathbf{S}_{k}\mathbf{\alpha }_{0}\right) ^{\ast }\mathbf{W}%
_{T}(\mathbf{\lambda }_{\mathbf{j}};\mathbf{\theta }_{0})\mathbf{\dot{\Psi}}%
_{k}^{\left( \ell \right) }(\mathbf{\lambda }_{\mathbf{j}};\mathbf{\theta }%
_{0})\mathbf{S}_{k}\mathbf{\alpha }_{0}\right \} \\
&&+\frac{2}{T^{k-1}}\sum_{\mathbf{\lambda }_{\mathbf{j}}}\func{Re}\left \{
\left( \mathbf{\Psi }_{k}(\mathbf{\lambda }_{\mathbf{j}};\mathbf{\theta }%
_{0})\mathbf{S}_{k}\left( -\mathbf{S}_{k}^{\prime }\mathbf{\bar{\Lambda}}%
_{k}^{\left( p\right) }\left( \mathbf{\theta }_{0}\right) \mathbf{S}_{k}%
\mathbf{\alpha }_{0}\right) \right) ^{\ast }\mathbf{W}_{k}(\mathbf{\lambda }%
_{\mathbf{j}};\mathbf{\theta }_{0})\mathbf{\dot{\Psi}}_{k}^{\left( \ell
\right) }(\mathbf{\lambda }_{\mathbf{j}};\mathbf{\theta }_{0})\mathbf{S}_{k}%
\mathbf{\alpha }_{0}\right \} \\
&&+\frac{2}{T^{k-1}}\sum_{\mathbf{\lambda }_{\mathbf{j}}}\func{Re}\left \{
\left( \mathbf{\Psi }_{k}(\mathbf{\lambda }_{\mathbf{j}};\mathbf{\theta }%
_{0})\mathbf{S}_{k}\left( -\mathbf{S}_{k}^{\prime }\mathbf{\bar{\Lambda}}%
_{k}^{\left( p\right) }\left( \mathbf{\theta }_{0}\right) \mathbf{S}_{k}%
\mathbf{\alpha }_{0}\right) \right) ^{\ast }\mathbf{W}_{k}(\mathbf{\lambda }%
_{\mathbf{j}};\mathbf{\theta }_{0})\mathbf{\Psi }_{k}(\mathbf{\lambda }_{%
\mathbf{j}};\mathbf{\theta }_{0})\mathbf{S}_{k}\left( -\mathbf{S}%
_{k}^{\prime }\mathbf{\bar{\Lambda}}_{k}^{\left( \ell \right) }\left( 
\mathbf{\theta }_{0}\right) \mathbf{S}_{k}\mathbf{\alpha }_{0}\right) \right
\} \\
&&+\frac{2}{T^{k-1}}\sum_{\mathbf{\lambda }_{\mathbf{j}}}\func{Re}\left \{
\left( \mathbf{\dot{\Psi}}_{k}^{(p)}(\mathbf{\lambda }_{\mathbf{j}};\mathbf{%
\theta }_{0})\mathbf{S}_{k}\mathbf{\alpha }_{0}\right) ^{\ast }\mathbf{W}%
_{k}(\mathbf{\lambda }_{\mathbf{j}};\mathbf{\theta }_{0})\mathbf{\Psi }_{k}(%
\mathbf{\lambda }_{\mathbf{j}};\mathbf{\theta }_{0})\mathbf{S}_{k}\left( -%
\mathbf{S}_{k}^{\prime }\mathbf{\bar{\Lambda}}_{k}^{\left( \ell \right)
}\left( \mathbf{\theta }_{0}\right) \mathbf{S}_{k}\mathbf{\alpha }%
_{0}\right) \right \} +o_{p}\left( 1\right) ,
\end{eqnarray*}
which converges to $2\mathbf{\alpha }_{k}^{0\prime }\mathbf{S}_{k}^{\prime
}\left \{ \mathbf{H}_{k}\left( \mathbf{\theta }_{0} \right) \right \}
_{\left( p,\ell \right) }\mathbf{S}_{k}\mathbf{\alpha }_{k}^0$ because%
\begin{eqnarray*}
&&\left( 2\pi \right) ^{k-1}\left \{ \mathbf{H}_{k}\left( \mathbf{\theta }%
\right) \right \} _{\left( p,\ell \right) } \\
&=&\int \func{Re}\left \{ \mathbf{\dot{\Psi}}_{k}^{(p)}(\mathbf{\lambda };%
\mathbf{\theta })^{\ast }\mathbf{W}_{k}(\mathbf{\lambda };\mathbf{\theta })%
\mathbf{\dot{\Psi}}_{k}^{\left( \ell \right) }(\mathbf{\lambda };\mathbf{%
\theta })\right \} d\mathbf{\lambda } \\
&&-\int \func{Re}\left \{ \mathbf{\bar{\Lambda}}_{k}^{\left( p\right)
}\left( \mathbf{\theta }\right) ^{\ast }\mathbf{S}_{k}\mathbf{S}_{k}^{\prime
}\mathbf{\Psi }_{k}^{\ast }(\mathbf{\lambda };\mathbf{\theta })\mathbf{W}%
_{k}(\mathbf{\lambda };\mathbf{\theta })\mathbf{\dot{\Psi}}_{k}^{\left( \ell
\right) }(\mathbf{\lambda };\mathbf{\theta })\right \} d\mathbf{\lambda } \\
&&+\int \func{Re}\left \{ \mathbf{\bar{\Lambda}}_{k}^{\left( p\right)
}\left( \mathbf{\theta }\right) ^{\ast }\mathbf{S}_{k}\mathbf{S}_{k}^{\prime
}\mathbf{\Psi }_{k}^{\ast }(\mathbf{\lambda };\mathbf{\theta })\mathbf{W}%
_{k}(\mathbf{\lambda };\mathbf{\theta })\mathbf{\Psi }_{k}(\mathbf{\lambda };%
\mathbf{\theta })\mathbf{S}_{k}\mathbf{S}_{k}^{\prime }\mathbf{\bar{\Lambda}}%
_{k}^{\left( \ell \right) }\left( \mathbf{\theta }\right) \right \} d\mathbf{%
\lambda } \\
&&-\int \func{Re}\left \{ \mathbf{\dot{\Psi}}_{k}^{(p)}(\mathbf{\lambda };%
\mathbf{\theta })^{\ast }\mathbf{W}_{k}(\mathbf{\lambda };\mathbf{\theta })%
\mathbf{\Psi }_{k}(\mathbf{\lambda };\mathbf{\theta })\mathbf{S}_{k}\mathbf{S%
}_{k}^{\prime }\mathbf{\bar{\Lambda}}_{k}^{\left( \ell \right) }\left( 
\mathbf{\theta }\right) \right \} d\mathbf{\lambda } \\
&=&\int \func{Re}\left \{ \left( \mathbf{\dot{\Psi}}_{k}^{(p)}(\mathbf{%
\lambda };\mathbf{\theta })-\mathbf{\Psi }_{k}(\mathbf{\lambda };\mathbf{%
\theta })\mathbf{S}_{k}\mathbf{S}_{k}^{\prime }\mathbf{\bar{\Lambda}}%
_{k}^{\left( p\right) }\left( \mathbf{\theta }\right) \right) ^{\ast }%
\mathbf{W}_{k}(\mathbf{\lambda };\mathbf{\theta })\left( \mathbf{\dot{\Psi}}%
_{k}^{(\ell )}(\mathbf{\lambda };\mathbf{\theta })-\mathbf{\Psi }_{k}(%
\mathbf{\lambda };\mathbf{\theta })\mathbf{S}_{k}\mathbf{S}_{k}^{\prime }%
\mathbf{\bar{\Lambda}}_{k}^{\left( \ell \right) }\left( \mathbf{\theta }%
\right) \right) \right \} d\mathbf{\lambda } \\
&=&\int \func{Re}\left \{ \left( \mathbf{\Psi }_{k}^{-1}(\mathbf{\lambda };%
\mathbf{\theta })\mathbf{\dot{\Psi}}_{k}^{(p)}(\mathbf{\lambda };\mathbf{%
\theta })-\mathbf{S}_{k}\mathbf{S}_{k}^{\prime }\mathbf{\bar{\Lambda}}%
_{k}^{\left( p\right) }\left( \mathbf{\theta }\right) \right) ^{\ast }\left( 
\mathbf{\Psi }_{k}^{-1}(\mathbf{\lambda };\mathbf{\theta })\mathbf{\dot{\Psi}%
}_{k}^{(\ell )}(\mathbf{\lambda };\mathbf{\theta })-\mathbf{S}_{k}\mathbf{S}%
_{k}^{\prime }\mathbf{\bar{\Lambda}}_{k}^{\left( \ell \right) }\left( 
\mathbf{\theta }\right) \right) \right \} d\mathbf{\lambda .}
\end{eqnarray*}
 $\Box $\bigskip}

\section{\protect \small Appendix D: Asymptotics of averages of higher order
periodograms}

{\small Define for $k=2,3,4,$ the average of periodograms of order }$k,$%
{\small 
\begin{equation*}
\mathbf{Z}_{k,T}:=\frac{T^{1/2}}{T^{k-1}}\sum_{\mathbf{\lambda }_{\mathbf{j}%
}}\func{Re}\left \{ \mathbf{B}_{k}^{\ast }\left( \mathbf{\lambda }_{\mathbf{j%
}}\right) \left( \mathbb{I}_{k}^{\varepsilon }(\mathbf{\lambda }_{\mathbf{j}%
})-E\left[ \mathbb{I}_{k}^{\varepsilon }(\mathbf{\lambda }_{\mathbf{j}})%
\right] \right) \right \}
\end{equation*}%
where 
\begin{equation*}
\mathbb{I}_{k}^{\mathbf{\varepsilon }}(\mathbf{\lambda })=\frac{1}{T}w_{T}^{%
\mathbf{\varepsilon }}\left( -\lambda _{1}-\cdots -\lambda _{k-1}\right)
\otimes w_{T}^{\mathbf{\varepsilon }}\left( \lambda _{k-1}\right) \otimes
\cdots \otimes w_{T}^{\mathbf{\varepsilon }}\left( \lambda _{1}\right) ,\  \
w_{T}^{\mathbf{\varepsilon }}\left( \lambda \right)
=\sum_{t=1}^{T}e^{-i\lambda t}\mathbf{\varepsilon }_{t}.
\end{equation*}%
is the $k$-periodogram of shocks satisfying Assumptions 1$(2k)$ and 3$(k)$.

Denote $\mathbf{\dot{\Lambda}}^{(\ell )}\left( \lambda ;\mathbf{%
\theta }\right) :=\mathbf{\Psi }^{-1}\left( e^{-i\lambda }\mathbf{;\theta }%
\right) \mathbf{\dot{\Psi}}^{(\ell )}\left( e^{-i\lambda }\mathbf{;\theta }%
\right) $ and $\mathbf{\dot{\Lambda}}\left( \lambda ;\mathbf{\theta }\right)
:=\left( \mathbf{\dot{\Lambda}}^{(1)}\left( \lambda ;\mathbf{\theta }\right)
,\ldots ,\mathbf{\dot{\Lambda}}^{(m)}\left( \lambda ;\mathbf{\theta }\right)
\right) ,$ and set for $j=0,\pm 1,\ldots ,$ 
\begin{equation*}
\mathbf{\mathbf{C}}\left( j\right) :=\left( 2\pi \right) ^{-1}\int_{-\pi
}^{\pi }\mathbf{\dot{\Lambda}}^{\ast }\left( \lambda ;\mathbf{\theta }%
_{0}\right) \exp \left( -ij\lambda \right) d\lambda
\end{equation*}%
so for $\mathbf{C}_{k}\left( 0\right) :=\left( 2\pi \right) ^{1-k}\int_{\Pi
^{k-1}}\mathbf{B}_{k}^{\ast }\left( \mathbf{\lambda };\mathbf{\theta }%
_{0}\right) d\mathbf{\lambda }$ for $k=2,3,4,$ we find that $\mathbf{C}%
_{2}\left( 0\right) =\mathbf{I}_{d}\otimes \mathbf{\mathbf{C}}\left(
0\right) +\mathbf{\mathbf{C}}\left( 0\right) \otimes \mathbf{I}_{d}$ while%
\begin{eqnarray*}
\mathbf{C}_{3}\left( 0\right) &=&\left[ \mathbf{I}_{d}\otimes \mathbf{I}%
_{d}\otimes \mathbf{\mathbf{C}}\left( 0\right) +\mathbf{I}_{d}\otimes 
\mathbf{\mathbf{C}}\left( 0\right) \otimes \mathbf{I}_{d}+\mathbf{\mathbf{C}}%
\left( 0\right) \otimes \mathbf{I}_{d}\otimes \mathbf{I}_{d}\right] \left( 
\mathbf{I}_{d^{3}}-\mathbf{S}_{3}\mathbf{S}_{3}^{\prime }\right) \\
\mathbf{C}_{4}\left( 0\right) &=&\left[ \mathbf{I}_{d}\otimes \mathbf{I}%
_{d}\otimes \mathbf{I}_{d}\otimes \mathbf{\mathbf{C}}\left( 0\right) +%
\mathbf{I}_{d}\otimes \mathbf{I}_{d}\otimes \mathbf{\mathbf{C}}\left(
0\right) \otimes \mathbf{I}_{d}+\cdots +\mathbf{\mathbf{C}}\left( 0\right)
\otimes \mathbf{I}_{d}\otimes \mathbf{I}_{d}\otimes \mathbf{I}_{d}\right]
\left( \mathbf{I}_{d^{4}}-\mathbf{S}_{4}\mathbf{S}_{4}^{\prime }\right) ,\  \
\  \ 
\end{eqnarray*}%
and for the block row matrices, $j\neq 0,$ 
\begin{equation*}
\mathbf{\mathbf{C}}_{k}\left( j\right) :=\left \{ \left( 2\pi \right)
^{1-k}\int_{\Pi ^{k-1}}\mathbf{B}_{k,a}^{\ast }\left( \mathbf{\lambda };%
\mathbf{\theta }_{0}\right) \exp \left( -ij\lambda _{a}\right) d\mathbf{%
\lambda }\right \} _{a=1,\ldots ,k}\mathbf{,}
\end{equation*}%
we find that%
\begin{eqnarray*}
\mathbf{C}_{2}\left( j\right) &=&\left[ \mathbf{I}_{d}\otimes \mathbf{%
\mathbf{C}}\left( j\right) \  \QATOP{{}}{{}}\  \mathbf{\mathbf{C}}\left(
j\right) \otimes \mathbf{I}_{d}\right] \\
\mathbf{C}_{3}\left( j\right) &=&\left[ \mathbf{I}_{d}\otimes \mathbf{I}%
_{d}\otimes \mathbf{\mathbf{C}}\left( j\right) \  \QATOP{{}}{{}}\  \mathbf{I}%
_{d}\otimes \mathbf{\mathbf{C}}\left( j\right) \otimes \mathbf{I}_{d}\ 
\QATOP{{}}{{}}\  \mathbf{\mathbf{C}}\left( j\right) \otimes \mathbf{I}%
_{d}\otimes \mathbf{I}_{d}\right] \\
\mathbf{C}_{4}\left( j\right) &=&\left[ \mathbf{I}_{d}\otimes \mathbf{I}%
_{d}\otimes \mathbf{I}_{d}\otimes \mathbf{\mathbf{C}}\left( j\right) \ 
\QATOP{{}}{{}}\  \mathbf{I}_{d}\otimes \mathbf{I}_{d}\otimes \mathbf{\mathbf{C%
}}\left( j\right) \otimes \mathbf{I}_{d}\QATOP{{}}{{}}\  \cdots \  \QATOP{{}}{%
{}}\mathbf{\mathbf{C}}\left( j\right) \otimes \mathbf{I}_{d}\otimes \mathbf{I%
}_{d}\otimes \mathbf{I}_{d}\right] .
\end{eqnarray*}%
The factors $\left( \mathbf{I}_{d^{k}}-\mathbf{S}_{k}\mathbf{S}_{k}^{\prime }\right) $that appear in $\mathbf{C}_{3}\left( 0\right) $, $\mathbf{C}_{4}\left(
0\right) $ and $\mathbf{\Phi }_{ab}^{0}$, $a$ or $b>2$, indicate that the ICA
restriction is imposed when estimating higher order cumulants, so not all
elements in the cumulant array v$\mathbf{\kappa}_k$ are estimated. In an
unrestricted estimation, based e.g. on Assumption~2(3), $\mathbf{S}_{k}$
would be replaced by $\mathbf{I}_{d^{k}}$ so the contributions from $\Phi
_{ab}^{0}\left( \mathbf{\theta }_{0};\mathbf{C}\right) $ for $a>2\ $or $b>2$
cancel out.

For $k=2,$ 
\begin{eqnarray*}
\mathbf{B}_{2,1}^{\left( \ell \right) }\left( \lambda ;\mathbf{\theta }%
\right) &=&\mathbf{B}_{2,1}^{\left( \ell \right) }\left( \lambda ;\mathbf{%
\theta }\right) =\mathbf{I}_{d}\otimes \left \{ \mathbf{\dot{\Lambda}}%
^{(\ell )}\left( \lambda \mathbf{;\theta }\right) \right \} \\
\mathbf{B}_{2,2}^{\left( \ell \right) }\left( \lambda ;\mathbf{\theta }%
\right) &=&\mathbf{B}_{2,2}^{\left( \ell \right) }\left( \lambda ;\mathbf{%
\theta }\right) =\left \{ \mathbf{\dot{\Lambda}}^{(\ell )}\left( \lambda 
\mathbf{;\theta }\right) \right \} \otimes \mathbf{I}_{d}
\end{eqnarray*}%
because 
$
\mathbf{\Psi }_{2}^{-1}(\lambda ;\mathbf{\theta })\mathbf{\dot{\Psi}}%
_{2}^{\left( \ell \right) }(\lambda ;\mathbf{\theta })=\left( \mathbf{I}%
_{d}\otimes \mathbf{\dot{\Lambda}}^{(\ell )}\left( \lambda \mathbf{;\theta }%
\right) \right) +\left( \mathbf{\dot{\Lambda}}^{(\ell )}\left( \lambda 
\mathbf{;\theta }\right) \otimes \mathbf{I}_{d}\right) ,\  
$
and 
\begin{eqnarray*}
\mathbb{I}_{2}^{\mathbf{\varepsilon }}(\lambda ) &=&\frac{1}{T}w_{T}\left(
-\lambda \right) \otimes w_{T}\left( \lambda \right) =\frac{1}{T}%
\sum_{t,r=1}^{T}\exp \left( -i\left( t-r\right) \lambda \right) \left( 
\mathbf{\varepsilon }_{r}\otimes \mathbf{\varepsilon }_{t}\right) \\
&=&\frac{1}{T}\sum_{t=1}^{T}\left( \mathbf{\varepsilon }_{t}\otimes \mathbf{%
\varepsilon }_{t}\right) +\frac{1}{T}\sum_{t,r=1}^{T}%
\sum_{r=1}^{t-1}A_{t,r}^{\left[ 2\right] }\left( \lambda _{j}\right) \left[ 
\begin{array}{c}
\mathbf{\varepsilon }_{r}\otimes \mathbf{\varepsilon }_{t} \\ 
\mathbf{\varepsilon }_{t}\otimes \mathbf{\varepsilon }_{r}%
\end{array}%
\right]
\end{eqnarray*}%
where 
\begin{equation*}
A_{t,r}^{\left[ 2\right] }\left( \lambda _{j}\right) :=\left[ \exp \left \{
-i\left( t-r\right) \lambda _{j}\right \} \  \QATOP{{}}{{}}\  \  \exp \left \{
-i\left( r-t\right) \lambda _{j}\right \} \right] .
\end{equation*}%

Then 
\begin{equation*}
\mathbf{Z}_{2,T}=\frac{T^{1/2}}{T}\sum_{\lambda _{j}}\func{Re}\left \{ 
\mathbf{B}_{2}^{\ast }(\lambda _{j})\left( \mathbb{I}_{2}^{\varepsilon
}(\lambda _{j})-E\left[ \mathbb{I}_{2}^{\varepsilon }(\lambda _{j})\right]
\right) \right \} =\sum_{t=1}^{T}Z_{2,t}
\end{equation*}%
where $Z_{2,t}$ is a martingale difference sequence (MDS) 
\begin{eqnarray*}
Z_{2,t} &:= &\frac{1}{T^{1/2}}\mathbf{C}_{T}^{\left[ 2,0\right] }\left \{
\left( \mathbf{\varepsilon }_{t}\otimes \mathbf{\varepsilon }_{t}\right) -%
\text{vec}\left( \mathbf{I}_{d}\right) \right \} \\
&&+\frac{1}{T^{1/2}}\sum_{r=1}^{t-1}\left[ \mathbf{C}_{T}^{\left[ 2,1\right]
}\left( r-t\right) \  \  \  \mathbf{C}_{T}^{\left[ 2,2\right] }\left(
r-t\right) \right] \left[ \QATOP{\mathbf{\varepsilon }_{t}\otimes \mathbf{%
\varepsilon }_{r}}{\mathbf{\varepsilon }_{r}\otimes \mathbf{\varepsilon }_{t}%
}\right] \\
&&+\frac{1}{T^{1/2}}\sum_{r=1}^{t-1}\left[ \mathbf{C}_{T}^{\left[ 2,1\right]
}\left( t-r\right) \  \  \  \  \mathbf{C}_{T}^{\left[ 2,2\right] }\left(
t-r\right) \right] \left[ \QATOP{\mathbf{\varepsilon }_{r}\otimes \mathbf{%
\varepsilon }_{t}}{\mathbf{\varepsilon }_{t}\otimes \mathbf{\varepsilon }_{r}%
}\right] \\
&=&\frac{1}{T^{1/2}}\mathbf{C}_{T}^{\left[ 2,0\right] }\left \{ \left( 
\mathbf{\varepsilon }_{t}\otimes \mathbf{\varepsilon }_{t}\right) -\text{vec}%
\left( \mathbf{I}_{d}\right) \right \} +\frac{1}{T^{1/2}}\sum_{r=1}^{t-1}%
\mathbf{C}_{T}^{\left[ 2\right] }\left( r-t\right) \mathbf{\varepsilon }%
_{t,r}^{\left[ 2\right] }+\frac{1}{T^{1/2}}\sum_{r=1}^{t-1}\mathbf{C}_{T}^{%
\left[ 2\right] }\left( t-r\right) \mathbf{\varepsilon }_{r,t}^{\left[ 2%
\right] }
\end{eqnarray*}%
with%
\begin{equation*}
\mathbf{C}_{T}^{\left[ 2,0\right] }:=\frac{1}{T}\sum_{\lambda _{j}}\func{Re}%
\mathbf{B}_{2}^{\ast }\left( \lambda _{j};\mathbf{\theta }\right)
\rightarrow \mathbf{C}_{2}\left( 0\right) :=\frac{1}{2\pi }\int_{-\pi }^{\pi
}\mathbf{B}_{2}^{\ast }\left( \lambda ;\mathbf{\theta }\right) d\lambda ,
\end{equation*}%
which in general is different from zero since $\mathbf{B}_{2}$ is not
centered as the scaling is incorporated in $\mathbf{\Psi }$ and second order
cumulants (or covariances) are not estimated separately due to the
normalization assumption, and 
\begin{eqnarray*}
\mathbf{C}_{T}^{\left[ 2\right] }\left( r-t\right)  &:=&\left[ \mathbf{C}%
_{T}^{\left[ 2,1\right] }\left( r-t\right) \  \  \  \  \mathbf{C}_{T}^{\left[ 2,2%
\right] }\left( r-t\right) \right] \\
&:=&\left[ \frac{1}{T}\sum_{\lambda _{j}}\func{Re}\left \{ \mathbf{B}%
_{2,1}^{\ast }\left( \lambda _{j};\mathbf{\theta }\right) \exp \left(
i\left( t-r\right) \lambda _{j}\right) \right \} \  \  \  \frac{1}{T}%
\sum_{\lambda _{j}}\func{Re}\left \{ \mathbf{B}_{2,2}^{\ast }\left( \lambda
_{j};\mathbf{\theta }\right) \exp \left( -i\left( t-r\right) \lambda
_{j}\right) \right \} \right] \\
&=&\mathbf{\mathbf{C}}_{2}\left( r-t\right) +O\left( T^{-1}\right) ,
\end{eqnarray*}%
as $T\rightarrow \infty $ where%
\begin{eqnarray*}
\mathbf{C}_{2}\left( j\right) &:=&\left[ \mathbf{C}_{\left[ 2,1\right]
}\left( j\right) \  \  \  \mathbf{C}_{\left[ 2,2\right] }\left( j\right) \right]
:=\left[ \left( \mathbf{I}_{d}\otimes \mathbf{\mathbf{C}}\left( j\right)
\right) \  \  \  \left( \mathbf{\mathbf{C}}\left( j\right) \otimes \mathbf{I}%
_{d}\right) \right] \\
\mathbf{\mathbf{C}}\left( j\right)  &:=&\frac{1}{2\pi }\int_{-\pi }^{\pi }%
\mathbf{\dot{\Lambda}}^{\ast }\left( \lambda \mathbf{;\theta }\right) \exp
\left( -ij\lambda \right) d\lambda .
\end{eqnarray*}%

Then we can write $\mathbf{Z}_{2,T}=\sum_{t=1}^{T}Z_{2,t}^{0}+o_{p}%
\left( 1\right) ,$ where%
\begin{equation*}
Z_{2,t}^{0}:=\frac{1}{T^{1/2}}\sum_{r=1}^{t-1}\left[ \mathbf{C}_{2}\left(
0\right) \  \  \  \mathbf{C}_{2}\left( r-t\right) \  \  \  \mathbf{C}_{2}\left(
t-r\right) \right] \left[ 
\begin{array}{c}
\frac{1}{t-1}\left( \mathbf{\varepsilon }_{t}^{\otimes 2}-\text{vec}\left( 
\mathbf{I}_{d}\right) \right) \\ 
\mathbf{\varepsilon }_{t,r}^{\left[ 2\right] } \\ 
\mathbf{\varepsilon }_{r,t}^{\left[ 2\right] }%
\end{array}%
\right]
\end{equation*}%
where%
\begin{equation*}
\mathbf{\varepsilon }_{t,r}:=\left[ 
\begin{array}{c}
\frac{1}{t-1}\left( \mathbf{\varepsilon }_{t}^{\otimes 2}-\text{vec}\left( 
\mathbf{I}_{d}\right) \right) \\ 
\mathbf{\varepsilon }_{t,r}^{\left[ 2\right] } \\ 
\mathbf{\varepsilon }_{r,t}^{\left[ 2\right] }%
\end{array}%
\right] :=\left[ 
\begin{array}{c}
\frac{1}{t-1}\left( \mathbf{\varepsilon }_{t}\otimes \mathbf{\varepsilon }%
_{t}-\text{vec}\left( \mathbf{I}_{d}\right) \right) \\ 
\mathbf{\varepsilon }_{t}\otimes \mathbf{\varepsilon }_{r} \\ 
\mathbf{\varepsilon }_{r}\otimes \mathbf{\varepsilon }_{t}%
\end{array}%
\right]
\end{equation*}%
so that exploiting symmetries, $\mathbb{C}\left[ \mathbf{\varepsilon }%
_{t}^{\otimes 2},\mathbf{\varepsilon }_{r,t}^{\left[ 2\right] }\right] =%
\mathbb{C}\left[ \mathbf{\varepsilon }_{t}^{\otimes 2},\mathbf{\varepsilon }%
_{t,r}^{\left[ 2\right] }\right] =0,$ $\mathbb{V}\left[ \mathbf{Z}_{2,T}%
\right] =$}$\sum_{t=1}^{T}\mathbb{V}\left[ Z_{2,t}^{0}\right] ${\small \ can
be approximated by 
\begin{eqnarray*}
&&\mathbf{C}_{2}\left( 0\right) V_{0}\mathbf{C}_{2}^{\prime }\left( 0\right)
\\
&&+\frac{1}{T}\sum_{t=1}^{T}\sum_{r=1}^{t-1}\left[ \mathbf{C}^{\left[ 2%
\right] }\left( r-t\right) \  \  \  \mathbf{C}^{\left[ 2\right] }\left(
t-r\right) \right] \left[ 
\begin{array}{cc}
\mathbb{V}\left[ \mathbf{\varepsilon }_{t,r}^{\left[ 2\right] }\right] & 
\mathbb{C}\left[ \mathbf{\varepsilon }_{t,r}^{\left[ 2\right] },\mathbf{%
\varepsilon }_{r,t}^{\left[ 2\right] }\right] \\ 
\mathbb{C}\left[ \mathbf{\varepsilon }_{r,t}^{\left[ 2\right] },\mathbf{%
\varepsilon }_{t,r}^{\left[ 2\right] }\right] & \mathbb{V}\left[ \mathbf{%
\varepsilon }_{r,t}^{\left[ 2\right] }\right]%
\end{array}%
\right] \left[ 
\begin{array}{c}
\mathbf{C}^{\left[ 2\right] }\left( r-t\right) ^{\prime } \\ 
\mathbf{C}^{\left[ 2\right] }\left( t-r\right) ^{\prime }%
\end{array}%
\right] \\
&=&\mathbf{C}_{2}\left( 0\right) V_{0}\mathbf{C}_{2}^{\prime }\left(
0\right) +\frac{1}{T}\sum_{t=1}^{T}\sum_{r=1}^{t-1}\mathbf{C}^{\left[ 2%
\right] }\left( r-t\right) \mathbb{V}\left[ \mathbf{\varepsilon }_{t,r}^{%
\left[ 2\right] }\right] \mathbf{C}^{\left[ 2\right] }\left( r-t\right)
^{\prime }+\mathbf{C}^{\left[ 2\right] }\left( t-r\right) \mathbb{V}\left[ 
\mathbf{\varepsilon }_{r,t}^{\left[ 2\right] }\right] \mathbf{C}^{\left[ 2%
\right] }\left( t-r\right) ^{\prime } \\
&&+ \mathbf{C}^{\left[ 2\right] }\left( r-t\right) \mathbb{C}\left[ \mathbf{%
\varepsilon }_{t,r}^{\left[ 2\right] },\mathbf{\varepsilon }_{r,t}^{\left[ 2%
\right] }\right] \mathbf{C}^{\left[ 2\right] }\left( t-r\right) ^{\prime }+%
\mathbf{C}^{\left[ 2\right] }\left( t-r\right) \mathbb{C}\left[ \mathbf{%
\varepsilon }_{r,t}^{\left[ 2\right] },\mathbf{\varepsilon }_{t,r}^{\left[ 2%
\right] }\right] \mathbf{C}^{\left[ 2\right] }\left( r-t\right) ^{\prime } \\
&\rightarrow &\mathbf{C}_{2}\left( 0\right) V_{0}\mathbf{C}_{2}^{\prime
}\left( 0\right) +\sum_{j=-\infty }^{\infty }\mathbf{C}^{\left[ 2\right]
}\left( j\right) \mathbb{V}\left[ \mathbf{\varepsilon }_{t,r}^{\left[ 2%
\right] }\right] \mathbf{C}^{\left[ 2\right] }\left( j\right) ^{\prime
}+\sum_{j=-\infty }^{\infty }\mathbf{C}^{\left[ 2\right] }\left( -j\right) 
\mathbb{C}\left[ \mathbf{\varepsilon }_{t,r}^{\left[ 2\right] },\mathbf{%
\varepsilon }_{r,t}^{\left[ 2\right] }\right] \mathbf{C}^{\left[ 2\right]
}\left( j\right) ^{\prime } \\
&=&\mathbf{\Phi }_{22}^{0}+\mathbf{\Phi }_{22}+\mathbf{\Phi }_{22}^{\dag }
\end{eqnarray*}%
where, under Assumption 3$(4)$, $V_{0}=V_{0}\left( \mathbf{\kappa }%
_{4}^{0}\right) =\mathbb{V}\left[ \mathbf{\varepsilon }_{t}^{\otimes 2}%
\right] =\sum_{a,b}\left( \mathbf{e}_{a}\mathbf{e}_{a}^{\prime }\otimes 
\mathbf{e}_{b}\mathbf{e}_{b}^{\prime }+\mathbf{e}_{a}\mathbf{e}_{b}^{\prime
}\otimes \mathbf{e}_{b}\mathbf{e}_{a}^{\prime }\right) +\sum_{a}\mathbf{%
\kappa }_{4a}^{0}\left( \mathbf{e}_{a}\mathbf{e}_{a}^{\prime }\otimes 
\mathbf{e}_{a}\mathbf{e}_{a}^{\prime }\right) ,$ 
\begin{equation*}
\mathbb{V}\left[ \mathbf{\varepsilon }_{t,r}^{\left[ 2\right] }\right] =%
\mathbb{V}\left[ \mathbf{\varepsilon }_{r,t}^{\left[ 2\right] }\right]
=\left( \mathbf{I}_{d}\otimes \mathbf{I}_{d}\right) =\mathbf{I}_{d^{2}}
\end{equation*}%
and for $t\neq r,$%
\begin{equation*}
\mathbb{C}\left[ \mathbf{\varepsilon }_{t,r}^{\left[ 2\right] },\mathbf{%
\varepsilon }_{r,t}^{\left[ 2\right] }\right] =\sum_{a,b}\left( \mathbf{e}%
_{a}\mathbf{e}_{b}^{\prime }\otimes \mathbf{e}_{b}\mathbf{e}_{a}^{\prime
}\right) .
\end{equation*}%

For $k=3$ and $\mathbf{\bar{\Lambda}}^{(\ell )}\left( \mathbf{\theta 
}\right) =\left( 2\pi \right) ^{-1}\int_{\Pi }\mathbf{\dot{\Lambda}}^{(\ell
)}\left( \lambda \mathbf{;\theta }\right) d\lambda $ we obtain%
\begin{eqnarray*}
\mathbf{B}_{3,1}^{\left( \ell \right) }\left( \mathbf{\lambda };\mathbf{%
\theta }\right) &=&\mathbf{B}_{3,1}^{\left( \ell \right) }\left( \lambda
_{1};\mathbf{\theta }\right) =\mathbf{I}_{d}\otimes \mathbf{I}_{d}\otimes
\left \{ \mathbf{\dot{\Lambda}}^{(\ell )}\left( \lambda _{1}\mathbf{;\theta }%
\right) -\mathbf{S}_{3}\mathbf{S}_{3}^{\prime }\mathbf{\bar{\Lambda}}^{(\ell
)}\left( \mathbf{\theta }\right) \right \} \\
\mathbf{B}_{3,2}^{\left( \ell \right) }\left( \mathbf{\lambda };\mathbf{%
\theta }\right) &=&\mathbf{B}_{3,2}^{\left( \ell \right) }\left( \lambda
_{2};\mathbf{\theta }\right) =\mathbf{I}_{d}\otimes \left \{ \mathbf{\dot{%
\Lambda}}^{(\ell )}\left( \lambda _{2}\mathbf{;\theta }\right) -\mathbf{S}%
_{3}\mathbf{S}_{3}^{\prime }\mathbf{\bar{\Lambda}}^{(\ell )}\left( \mathbf{%
\theta }\right) \right \} \otimes \mathbf{I}_{d} \\
\mathbf{B}_{3,2}^{\left( \ell \right) }\left( \mathbf{\lambda };\mathbf{%
\theta }\right) &=&\mathbf{B}_{3,1}^{\left( \ell \right) }\left( \lambda
_{3};\mathbf{\theta }\right) =\left \{ \mathbf{\dot{\Lambda}}^{(\ell
)}\left( -\lambda _{1}-\lambda _{2}\mathbf{;\theta }\right) -\mathbf{S}_{3}%
\mathbf{S}_{3}^{\prime }\mathbf{\bar{\Lambda}}^{(\ell )}\left( \mathbf{%
\theta }\right) \right \} \otimes \mathbf{I}_{d}\otimes \mathbf{I}_{d},
\end{eqnarray*}%
because 
\begin{equation*}
\mathbf{\Psi }_{3}^{-1}(\mathbf{\lambda };\mathbf{\theta })\mathbf{\dot{\Psi}%
}_{3}^{\left( \ell \right) }(\mathbf{\lambda };\mathbf{\theta })=\left( 
\mathbf{I}_{d}\otimes \mathbf{I}_{d}\otimes \mathbf{\dot{\Lambda}}^{(\ell
)}\left( \lambda _{1}\mathbf{;\theta }\right) \right) +\left( \mathbf{I}%
_{d}\otimes \mathbf{\dot{\Lambda}}^{(\ell )}\left( \lambda _{2}\mathbf{%
;\theta }\right) \otimes \mathbf{I}_{d}\right) +\left( \mathbf{\dot{\Lambda}}%
^{(\ell )}\left( -\lambda _{1}-\lambda _{2}\mathbf{;\theta }\right) \otimes 
\mathbf{I}_{d}\otimes \mathbf{I}_{d}\right) ,\  \  \ 
\end{equation*}%
and the third order periodogram }$\mathbb{I}_{3}^{\mathbf{\varepsilon }}(%
\mathbf{\lambda })=T^{-1}w_{T}^{\mathbf{\varepsilon }}\left( -\lambda
_{1}-\lambda _{2}\right) \otimes w_{T}^{\mathbf{\varepsilon }}\left( \lambda
_{2}\right) \otimes w_{T}^{\mathbf{\varepsilon }}\left( \lambda _{1}\right) $
can be written as 
\begin{eqnarray*}
&&\frac{1}{T}\sum_{t,r,s=1}^{T}\exp \left( is\left( \lambda _{1}+\lambda
_{2}\right) -ir\lambda _{2}-it\lambda _{1}\right) \left( \mathbf{\varepsilon 
}_{s}\otimes \mathbf{\varepsilon }_{r}\otimes \mathbf{\varepsilon }%
_{t}\right) \\ \! \! \! \! \! \!
&=&\frac{1}{T}\sum_{t=1}^{T}\left( \mathbf{\varepsilon }_{t}\otimes \mathbf{%
\varepsilon }_{t}\otimes \mathbf{\varepsilon }_{t}\right) +\frac{1}{T}%
\sum_{t=1}^{T}\sum_{r=1}^{t-1}A_{t,r}^{\left[ 3\right] }\left( \mathbf{%
\lambda }\right) \! \left[ \!
\begin{array}{c}
\mathbf{\varepsilon }_{t}\otimes \mathbf{\varepsilon }_{t}\otimes \mathbf{%
\varepsilon }_{r} \\ 
\mathbf{\varepsilon }_{t}\otimes \mathbf{\varepsilon }_{r}\otimes \mathbf{%
\varepsilon }_{t} \\ 
\mathbf{\varepsilon }_{r}\otimes \mathbf{\varepsilon }_{t}\otimes \mathbf{%
\varepsilon }_{t}%
\end{array} \!
\right] \!+\frac{1}{T}\sum_{t=1}^{T}\sum_{r,s=1}^{t-1}G_{t,r,s}^{\left[ 3%
\right] }\left( \mathbf{\lambda }\right) \! \left[ \!
\begin{array}{c}
\mathbf{\varepsilon }_{s}\otimes \mathbf{\varepsilon }_{r}\otimes \mathbf{%
\varepsilon }_{t} \\ 
\mathbf{\varepsilon }_{s}\otimes \mathbf{\varepsilon }_{t}\otimes \mathbf{%
\varepsilon }_{r} \\ 
\mathbf{\varepsilon }_{t}\otimes \mathbf{\varepsilon }_{s}\otimes \mathbf{%
\varepsilon }_{r}%
\end{array}%
\! \right] \ \ \ \
\end{eqnarray*}%
where 
\begin{eqnarray*}
A_{t,r}^{\left[ 3\right] }\left( \mathbf{\lambda }_{\mathbf{j}}\right) &:= &%
\left[ \exp \left \{ -i\left( \left( r-t\right) \lambda _{j_{1}}\right)
\right \} \QATOP{{}}{{}}\exp \left \{ -i\left( \left( r-t\right) \lambda
_{j_{2}}\right) \right \} \QATOP{{}}{{}}\exp \left \{ -i\left( \left(
t-r\right) \left( \lambda _{j_{1}}+\lambda _{j_{2}}\right) \right) \right \} %
\right] , \\
G_{t,r,s}^{\left[ 3\right] }\left( \mathbf{\lambda }_{\mathbf{j}}\right) &:=&
\left[ 
\begin{array}{c}
\exp \left \{ -i\left( \left( t-s\right) \lambda _{j_{1}}+\left( r-s\right)
\lambda _{j_{2}}\right) \right \} \\ 
\exp \left \{ -i\left( \left( r-s\right) \lambda _{j_{1}}+\left( t-s\right)
\lambda _{j_{2}}\right) \right \} \\ 
\exp \left \{ -i\left( \left( r-t\right) \lambda _{j_{1}}+\left( s-t\right)
\lambda _{j_{2}}\right) \right \}%
\end{array}%
\right] ^{\prime }.
\end{eqnarray*}%

Then 
\begin{equation*}
\mathbf{Z}_{3,T}=\frac{T^{1/2}}{T^{2}}\sum_{\mathbf{\lambda }_{\mathbf{j}}}%
\func{Re}\left \{ \mathbf{B}_{3}^{\ast }\left( \mathbf{\lambda }_{\mathbf{j}%
};\mathbf{\theta }_{0}\right) \left( \mathbb{I}_{3}^{\varepsilon }(\mathbf{%
\lambda }_{\mathbf{j}})-E\left[ \mathbb{I}_{3}^{\varepsilon }(\mathbf{%
\lambda }_{\mathbf{j}})\right] \right) \right \} +o_{p}\left( 1\right)
=\sum_{t=1}^{T}Z_{3,t}+o_{p}\left( 1\right)
\end{equation*}%
where, vec$\left( \mathbf{I}_{d^{3}}\right) =\sum_{a=1}^{d}\left( \mathbf{e}%
_{a}\otimes \mathbf{e}_{a}\otimes \mathbf{e}_{a}\right) ,$ $\mathbf{e}_{a}$
equal to the $a$-th column of $\mathbf{I}_{d},$ 
\begin{eqnarray*}
Z_{3,t}&:= &\frac{1}{T^{1/2}}\mathbf{C}_{T}^{\left[ 3,0\right] }\left \{
\left( \mathbf{\varepsilon }_{t}\otimes \mathbf{\varepsilon }_{t}\otimes 
\mathbf{\varepsilon }_{t}\right) -\sum_{a=1}^{d}\mathbf{\kappa }_{3,a}\left( 
\mathbf{e}_{a}\otimes \mathbf{e}_{a}\otimes \mathbf{e}_{a}\right) \right \}
\\
&&+\frac{1}{T^{1/2}}\sum_{r=1}^{t-1}\left \{ \frac{1}{T^{2}}\sum_{\mathbf{%
\lambda }_{\mathbf{j}}}\mathbf{B}_{3}^{\ast }\left( \mathbf{\lambda }_{%
\mathbf{j}};\mathbf{\theta }_{0}\right) A_{t,r}^{\left[ 3\right] }\left( 
\mathbf{\lambda }_{\mathbf{j}}\right) \right \} \left[ 
\begin{array}{c}
\mathbf{\varepsilon }_{t}\otimes \mathbf{\varepsilon }_{t}\otimes \mathbf{%
\varepsilon }_{r} \\ 
\mathbf{\varepsilon }_{t}\otimes \mathbf{\varepsilon }_{r}\otimes \mathbf{%
\varepsilon }_{t} \\ 
\mathbf{\varepsilon }_{r}\otimes \mathbf{\varepsilon }_{t}\otimes \mathbf{%
\varepsilon }_{t}%
\end{array}%
\right] \\
&&+\frac{1}{T^{1/2}}\sum_{r,s=1}^{t-1}\left \{ \frac{1}{T^{2}}\sum_{\mathbf{%
\lambda }_{\mathbf{j}}}\mathbf{B}_{3}^{\ast }\left( \mathbf{\lambda }_{%
\mathbf{j}};\mathbf{\theta }_{0}\right) G_{t,r,s}^{\left[ 3\right] }\left( 
\mathbf{\lambda }_{\mathbf{j}}\right) \right \} \left[ 
\begin{array}{c}
\mathbf{\varepsilon }_{s}\otimes \mathbf{\varepsilon }_{r}\otimes \mathbf{%
\varepsilon }_{t} \\ 
\mathbf{\varepsilon }_{s}\otimes \mathbf{\varepsilon }_{t}\otimes \mathbf{%
\varepsilon }_{r} \\ 
\mathbf{\varepsilon }_{t}\otimes \mathbf{\varepsilon }_{s}\otimes \mathbf{%
\varepsilon }_{r}%
\end{array}%
\right]
\end{eqnarray*}%
where%
\begin{equation*}
\mathbf{C}_{T}^{\left[ 3,0\right] }:=\frac{1}{T^{2}}\sum_{\mathbf{\lambda }_{%
\mathbf{j}}}\mathbf{B}_{3}^{\ast }\left( \mathbf{\lambda }_{\mathbf{j}};%
\mathbf{\theta }_{0}\right) =\mathbf{C}_{3}\left( 0\right) +O\left(
T^{-1}\right) ,
\end{equation*}%
for $\mathbf{C}_{3}\left( 0\right) :=\left( 2\pi \right) ^{-2}\int_{\Pi ^{2}}%
\mathbf{B}_{3}^{\ast }\left( \mathbf{\lambda };\mathbf{\theta }_{0}\right) d%
\mathbf{\lambda }=\mathbf{A}_{3}^{\prime }\left( \mathbf{\theta }_{0}\right)
\left( \mathbf{I}_{d^{3}}-\mathbf{S}_{3}\mathbf{S}_{3}^{\prime }\right) $,
while%
\begin{equation*}
\frac{1}{T^{2}}\sum_{\mathbf{\lambda }_{\mathbf{j}}}\mathbf{B}_{3}^{\ast
}\left( \mathbf{\lambda }_{\mathbf{j}};\mathbf{\theta }_{0}\right) A_{t,r}^{%
\left[ 3\right] }\left( \mathbf{\lambda }_{\mathbf{j}}\right) =\mathbf{C}%
_{3}\left( r-t\right) +O(T^{-1})
\end{equation*}%
as $T\rightarrow \infty ,$ where%
\begin{eqnarray*}
\mathbf{C}_{3}\left( j\right) &:= &\left[ \mathbf{C}_{\left[ 3,1\right]
}\left( j\right) \  \  \mathbf{C}_{\left[ 3,2\right] }\left( j\right) \  \ 
\mathbf{C}_{\left[ 3,3\right] }\left( j\right) \right] \\
&:= &\left[ \left( \mathbf{I}_{d}\otimes \mathbf{I}_{d}\otimes \mathbf{%
\mathbf{C}}\left( j\right) \right) \QATOP{{}}{{}}\left( \mathbf{I}%
_{d}\otimes \mathbf{\mathbf{C}}\left( j\right) \otimes \mathbf{I}_{d}\right) 
\QATOP{{}}{{}}\left( \mathbf{\mathbf{C}}\left( j\right) \otimes \mathbf{I}%
_{d}\otimes \mathbf{I}_{d}\right) \right] ,
\end{eqnarray*}%
and similarly, the $G_{t,r,s}^{\left[ 3\right] }$ terms only contribute when 
$r=s$ with%
\begin{equation*}
\frac{1}{T^{2}}\sum_{\mathbf{\lambda }_{\mathbf{j}}}\mathbf{B}_{3}^{\ast
}\left( \mathbf{\lambda }_{\mathbf{j}};\mathbf{\theta }_{0}\right)
G_{t,r,r}^{\left[ 3\right] }\left( \mathbf{\lambda }_{\mathbf{j}}\right) =%
\mathbf{C}_{3}\left( t-r\right) +O(T^{-1}).
\end{equation*}%

Then, using the same arguments as in VL, noting that the terms in $%
\mathbf{G}\ $for $r\neq s$ do not contribute, we can write 
\begin{equation*}
\mathbf{Z}_{3,T}=\sum_{t=1}^{T}Z_{3,t}^{0}+o_{p}\left( 1\right)
\end{equation*}%
where $Z_{3,t}^{0}$ is a MDS 
\begin{equation*}
Z_{3,t}^{0}=\frac{1}{T^{1/2}}\sum_{r=1}^{t-1}\left[ \mathbf{C}_{3}\left(
0\right) \  \  \mathbf{C}_{3}\left( r-t\right) \  \  \mathbf{C}_{3}\left(
t-r\right) \right] \left[ 
\begin{array}{c}
\frac{1}{t-1}\left( \mathbf{\varepsilon }_{t}^{\otimes 3}-E\left[ \mathbf{%
\varepsilon }_{t}^{\otimes 3}\right] \right) \\ 
\mathbf{\varepsilon }_{t,r}^{\left[ 3\right] } \\ 
\mathbf{\varepsilon }_{r,t}^{\left[ 3\right] }%
\end{array}%
\right]
\end{equation*}%
with $\mathbf{\varepsilon }_{t,r}^{\left[ 3\right] }$ defined in Section~5,
so that exploiting symmetries, $\mathbb{V}\left[ \mathbf{Z}_{3,T}\right] =
\sum_{t=1}^{T}\mathbb{V}\left[ Z_{3,t}^{0}\right] ${\small \ can be
approximated by%
\begin{eqnarray*}
&&\mathbf{C}_{3}\left( 0\right) \mathbb{V}\left[ \mathbf{\varepsilon }%
_{t}^{\otimes 3}\right] \mathbf{C}_{3}^{\prime }\left( 0\right) \\
&&+\frac{1}{T}\sum_{t=1}^{T}\sum_{r=1}^{t-1}\left[ \mathbf{C}^{\left[ 3%
\right] }\left( r-t\right) \  \  \mathbf{C}^{\left[ 3\right] }\left(
t-r\right) \right] \left[ 
\begin{array}{cc}
\mathbb{V}\left[ \mathbf{\varepsilon }_{t,r}^{\left[ 3\right] }\right] & 
\mathbb{C}\left[ \mathbf{\varepsilon }_{t,r}^{\left[ 3\right] },\mathbf{%
\varepsilon }_{r,t}^{\left[ 3\right] }\right] \\ 
\mathbb{C}\left[ \mathbf{\varepsilon }_{r,t}^{\left[ 3\right] },\mathbf{%
\varepsilon }_{t,r}^{\left[ 3\right] }\right] & \mathbb{V}\left[ \mathbf{%
\varepsilon }_{r,t}^{\left[ 3\right] }\right]%
\end{array}%
\right] \left[ 
\begin{array}{c}
\mathbf{C}^{\left[ 3\right] }\left( r-t\right) ^{\prime } \\ 
\mathbf{C}^{\left[ 3\right] }\left( t-r\right) ^{\prime }%
\end{array}%
\right] \\
&=&\mathbf{C}_{3}\left( 0\right) \mathbb{V}\left[ \mathbf{\varepsilon }%
_{t}^{\otimes 3}\right] \mathbf{C}_{3}^{\prime }\left( 0\right) +\frac{1}{T}%
\sum_{t=1}^{T}\sum_{r=1}^{t-1}\mathbf{C}^{\left[ 3\right] }\left( r-t\right) 
\mathbb{V}\left[ \mathbf{\varepsilon }_{t,r}^{\left[ 3\right] }\right] 
\mathbf{C}^{\left[ 3\right] }\left( r-t\right) ^{\prime }+\mathbf{C}^{\left[
3\right] }\left( t-r\right) \mathbb{V}\left[ \mathbf{\varepsilon }%
_{r,t}^{\left[ 3\right] }\right] \mathbf{C}^{\left[ 3\right] }\left(
t-r\right) ^{\prime } \\
&&+\mathbf{C}^{\left[ 3\right] }\left( r-t\right) \mathbb{C}\left[ \mathbf{%
\varepsilon }_{t,r}^{\left[ 3\right] },\mathbf{\varepsilon }_{r,t}^{\left(
3\right) }\right] \mathbf{C}^{\left[ 3\right] }\left( t-r\right) ^{\prime }+%
\mathbf{C}^{\left[ 3\right] }\left( t-r\right) \mathbb{C}\left[ \mathbf{%
\varepsilon }_{r,t}^{\left[ 3\right] },\mathbf{\varepsilon }_{t,r}^{\left(
3\right) }\right] \mathbf{C}^{\left[ 3\right] }\left( r-t\right) ^{\prime }
\\
&\rightarrow &\mathbf{C}_{3}\left( 0\right) \mathbb{V}\left[ \mathbf{%
\varepsilon }_{t}^{\otimes 3}\right] \mathbf{C}_{3}^{\prime }\left( 0\right)
+\sum_{j=-\infty }^{\infty }\mathbf{C}^{\left[ 3\right] }\left( j\right) 
\mathbb{V}\left[ \mathbf{\varepsilon }_{t,r}^{\left[ 3\right] }\right] 
\mathbf{C}^{\left[ 3\right] }\left( j\right) ^{\prime }+\sum_{j=-\infty
}^{\infty }\mathbf{C}^{\left[ 3\right] }\left( -j\right) \mathbb{C}\left[ 
\mathbf{\varepsilon }_{t,r}^{\left[ 3\right] },\mathbf{\varepsilon }%
_{r,t}^{\left[ 3\right] }\right] \mathbf{C}^{\left[ 3\right] }\left(
j\right) ^{\prime } \\
&=&\mathbf{\Phi }_{33}^{0}+\mathbf{\Phi }_{33}+\mathbf{\Phi }_{33}^{\dag }=%
\mathbf{\Phi }_{33}^{0}+\mathbf{\Phi }_{33}\left( \mathbf{\kappa }%
_{4}^{0}\right) +\mathbf{\Phi }_{33}^{\dag }\left( \mathbf{\kappa }%
_{3}^{0}\right) ,\  \text{say,}
\end{eqnarray*}%
where $\mathbf{\Phi }_{33}^{0}:=\mathbf{C}_{3}\left( 0\right) \mathbb{V}%
\left[ \mathbf{\varepsilon }_{t}^{\otimes 3}\right] \mathbf{C}_{3}^{\prime
}\left( 0\right) \ $depends on $\mathbf{\kappa }_{2},\mathbf{\kappa }_{3}$
and $\mathbf{\kappa }_{4},$ but not on $\mathbf{\kappa }_{6}$ because the
right factor $\mathbf{I}_{d^{3}}-\mathbf{S}_{3}\mathbf{S}_{3}^{\prime }=%
\mathbf{I}_{d^{3}}-\sum_{a}\mathbf{e}_{a}^{\otimes 3}\mathbf{e}_{a}^{\otimes
3\prime }$ of $\mathbf{C}_{3}\left( 0\right) $ is orthogonal to the
contribution of $\mathbf{\kappa }_{6}$ to $\mathbb{V}\left[ \mathbf{%
\varepsilon }_{t}^{\otimes 3}\right] ,$ namely$\  \sum_{a}\mathbf{\kappa }%
_{6,a}\mathbf{e}_{a}\mathbf{e}_{a}^{\prime }\otimes \mathbf{e}_{a}\mathbf{e}%
_{a}^{\prime }\otimes \mathbf{e}_{a}\mathbf{e}_{a}^{\prime }=\sum_{a}\mathbf{%
\kappa }_{6,a}\mathbf{e}_{a}^{\otimes 3}\mathbf{e}_{a}^{\otimes 3\prime }$,
and $\mathbb{V}\left[ \mathbf{\varepsilon }_{t,r}^{\left[ 3\right] }\right] =%
\mathbb{V}\left[ \mathbf{\varepsilon }_{r,t}^{\left[ 3\right] }\right] ,$ $%
t\neq r,$ where under Assumption 3$(h)$, $h\in \{3,4\}$, $\mathbb{V}\left[ \mathbf{\varepsilon }_{t,r}^{\left[ 3\right] }\right]$ is equal to
\begin{eqnarray*}
&&\sum_{a,b,c}\left[ 
\begin{array}{ccc}
\left \{ 
\begin{array}{c}
\mathbf{e}_{a}\mathbf{e}_{b}^{\prime }\otimes \mathbf{e}_{a}\mathbf{e}%
_{b}^{\prime }\otimes \mathbf{e}_{c}\mathbf{e}_{c}^{\prime } \\ 
+\mathbf{e}_{a}\mathbf{e}_{b}^{\prime }\otimes \mathbf{e}_{b}\mathbf{e}%
_{a}^{\prime }\otimes \mathbf{e}_{c}\mathbf{e}_{c}^{\prime }%
\end{array}%
\right \} & \left \{ 
\begin{array}{c}
\mathbf{e}_{a}\mathbf{e}_{a}^{\prime }\otimes \mathbf{e}_{b}\mathbf{e}%
_{c}^{\prime }\otimes \mathbf{e}_{c}\mathbf{e}_{b}^{\prime } \\ 
+\mathbf{e}_{a}\mathbf{e}_{b}^{\prime }\otimes \mathbf{e}_{b}\mathbf{e}%
_{c}^{\prime }\otimes \mathbf{e}_{c}\mathbf{e}_{a}^{\prime }%
\end{array}%
\right \} & \left \{ 
\begin{array}{c}
\mathbf{e}_{a}\mathbf{e}_{b}^{\prime }\otimes \mathbf{e}_{c}\mathbf{e}%
_{a}^{\prime }\otimes \mathbf{e}_{b}\mathbf{e}_{c}^{\prime } \\ 
+\mathbf{e}_{a}\mathbf{e}_{b}^{\prime }\otimes \mathbf{e}_{c}\mathbf{e}%
_{c}^{\prime }\otimes \mathbf{e}_{b}\mathbf{e}_{a}^{\prime }%
\end{array}%
\right \} \\ 
& \left \{ 
\begin{array}{c}
\mathbf{e}_{a}\mathbf{e}_{a}^{\prime }\otimes \mathbf{e}_{b}\mathbf{e}%
_{b}^{\prime }\otimes \mathbf{e}_{c}\mathbf{e}_{c}^{\prime } \\ 
+\mathbf{e}_{a}\mathbf{e}_{b}^{\prime }\otimes \mathbf{e}_{c}\mathbf{e}%
_{c}^{\prime }\otimes \mathbf{e}_{b}\mathbf{e}_{a}^{\prime }%
\end{array}%
\right \} & \left \{ 
\begin{array}{c}
\mathbf{e}_{a}\mathbf{e}_{b}^{\prime }\otimes \mathbf{e}_{b}\mathbf{e}%
_{a}^{\prime }\otimes \mathbf{e}_{c}\mathbf{e}_{c}^{\prime } \\ 
+\mathbf{e}_{a}\mathbf{e}_{b}^{\prime }\otimes \mathbf{e}_{b}\mathbf{e}%
_{c}^{\prime }\otimes \mathbf{e}_{c}\mathbf{e}_{a}^{\prime }%
\end{array}%
\right \} \\ 
&  & \left \{ 
\begin{array}{c}
\mathbf{e}_{c}\mathbf{e}_{c}^{\prime }\otimes \mathbf{e}_{a}\mathbf{e}%
_{b}^{\prime }\otimes \mathbf{e}_{a}\mathbf{e}_{b}^{\prime } \\ 
+\mathbf{e}_{c}\mathbf{e}_{c}^{\prime }\otimes \mathbf{e}_{a}\mathbf{e}%
_{b}^{\prime }\otimes \mathbf{e}_{b}\mathbf{e}_{a}^{\prime }%
\end{array}%
\right \}%
\end{array}%
\right] \\
&&+\sum_{a,b}\mathbf{\kappa }_{4,a}\left[ 
\begin{array}{ccc}
\mathbf{e}_{a}\mathbf{e}_{a}^{\prime }\otimes \mathbf{e}_{a}\mathbf{e}%
_{a}^{\prime }\otimes \mathbf{e}_{b}\mathbf{e}_{b}^{\prime } & \mathbf{e}_{a}%
\mathbf{e}_{a}^{\prime }\otimes \mathbf{e}_{a}\mathbf{e}_{b}^{\prime
}\otimes \mathbf{e}_{b}\mathbf{e}_{a}^{\prime } & \mathbf{e}_{a}\mathbf{e}%
_{b}^{\prime }\otimes \mathbf{e}_{a}\mathbf{e}_{a}^{\prime }\otimes \mathbf{e%
}_{b}\mathbf{e}_{a}^{\prime } \\ 
& \mathbf{e}_{a}\mathbf{e}_{a}^{\prime }\otimes \mathbf{e}_{b}\mathbf{e}%
_{b}^{\prime }\otimes \mathbf{e}_{a}\mathbf{e}_{a}^{\prime } & \mathbf{e}_{a}%
\mathbf{e}_{b}^{\prime }\otimes \mathbf{e}_{b}\mathbf{e}_{a}^{\prime
}\otimes \mathbf{e}_{a}\mathbf{e}_{a}^{\prime } \\ 
&  & \mathbf{e}_{b}\mathbf{e}_{b}^{\prime }\otimes \mathbf{e}_{a}\mathbf{e}%
_{a}^{\prime }\otimes \mathbf{e}_{a}\mathbf{e}_{a}^{\prime }%
\end{array}%
\right]
\end{eqnarray*}%
and $\mathbb{C}\left[ \mathbf{\varepsilon }_{t,r}^{\left[ 3\right] },\mathbf{%
\varepsilon }_{r,t}^{\left[ 3\right] }\right] =\mathbb{C}\left[ \mathbf{%
\varepsilon }_{r,t}^{\left[ 3\right] },\mathbf{\varepsilon }_{t,r}^{\left(
3\right) }\right] ,$ with
\begin{equation*}
\mathbb{C}\left[ \mathbf{\varepsilon }_{t,r}^{\left[ 3\right] },\mathbf{%
\varepsilon }_{r,t}^{\left[ 3\right] }\right] =\sum_{a,b}\mathbf{\kappa }%
_{3,a}\mathbf{\kappa }_{3,b}\left[ 
\begin{array}{ccc}
\mathbf{e}_{a}\mathbf{e}_{b}^{\prime }\otimes \mathbf{e}_{a}\mathbf{e}%
_{b}^{\prime }\otimes \mathbf{e}_{b}\mathbf{e}_{b}^{\prime } & \mathbf{e}_{a}%
\mathbf{e}_{b}^{\prime }\otimes \mathbf{e}_{a}\mathbf{e}_{a}^{\prime
}\otimes \mathbf{e}_{b}\mathbf{e}_{b}^{\prime } & \mathbf{e}_{a}\mathbf{e}%
_{a}^{\prime }\otimes \mathbf{e}_{a}\mathbf{e}_{b}^{\prime }\otimes \mathbf{e%
}_{b}\mathbf{e}_{b}^{\prime } \\ 
\mathbf{e}_{a}\mathbf{e}_{b}^{\prime }\otimes \mathbf{e}_{b}\mathbf{e}%
_{b}^{\prime }\otimes \mathbf{e}_{a}\mathbf{e}_{a}^{\prime } & \mathbf{e}_{a}%
\mathbf{e}_{b}^{\prime }\otimes \mathbf{e}_{b}\mathbf{e}_{a}^{\prime
}\otimes \mathbf{e}_{a}\mathbf{e}_{b}^{\prime } & \mathbf{e}_{a}\mathbf{e}%
_{a}^{\prime }\otimes \mathbf{e}_{b}\mathbf{e}_{b}^{\prime }\otimes \mathbf{e%
}_{a}\mathbf{e}_{b}^{\prime } \\ 
\mathbf{e}_{a}\mathbf{e}_{a}^{\prime }\otimes \mathbf{e}_{b}\mathbf{e}%
_{a}^{\prime }\otimes \mathbf{e}_{b}\mathbf{e}_{b}^{\prime } & \mathbf{e}_{a}%
\mathbf{e}_{a}^{\prime }\otimes \mathbf{e}_{b}\mathbf{e}_{b}^{\prime
}\otimes \mathbf{e}_{b}\mathbf{e}_{a}^{\prime } & \mathbf{e}_{a}\mathbf{e}%
_{b}^{\prime }\otimes \mathbf{e}_{b}\mathbf{e}_{a}^{\prime }\otimes \mathbf{e%
}_{b}\mathbf{e}_{a}^{\prime }%
\end{array}%
\right] .
\end{equation*}%

For $k=4,$%
\begin{eqnarray*}
\mathbf{B}_{4,1}^{\left( \ell \right) }\left( \mathbf{\lambda };\mathbf{%
\theta }\right) &=&\mathbf{B}_{4,1}^{\left( \ell \right) }\left( \lambda
_{1};\mathbf{\theta }\right) =\mathbf{I}_{d}\otimes \mathbf{I}_{d}\otimes 
\mathbf{I}_{d}\otimes \left \{ \mathbf{\dot{\Lambda}}^{(\ell )}\left(
e^{-i\lambda _{1}}\mathbf{;\theta }\right) -\mathbf{S}_{4}\mathbf{S}%
_{4}^{\prime }\mathbf{\bar{\Lambda}}^{(\ell )}\left( \mathbf{\theta }\right)
\right \} \\
\mathbf{B}_{4,2}^{\left( \ell \right) }\left( \mathbf{\lambda };\mathbf{%
\theta }\right) &=&\mathbf{B}_{4,2}^{\left( \ell \right) }\left( \lambda
_{2};\mathbf{\theta }\right) =\mathbf{I}_{d}\otimes \mathbf{I}_{d}\otimes
\left \{ \mathbf{\dot{\Lambda}}^{(\ell )}\left( e^{-i\lambda _{2}}\mathbf{%
;\theta }\right) -\mathbf{S}_{4}\mathbf{S}_{4}^{\prime }\mathbf{\bar{\Lambda}%
}^{(\ell )}\left( \mathbf{\theta }\right) \right \} \otimes \mathbf{I}_{d} \\
\mathbf{B}_{4,2}^{\left( \ell \right) }\left( \mathbf{\lambda };\mathbf{%
\theta }\right) &=&\mathbf{B}_{4,1}^{\left( \ell \right) }\left( \lambda
_{3};\mathbf{\theta }\right) =\mathbf{I}_{d}\otimes \left \{ \mathbf{\dot{%
\Lambda}}^{(\ell )}\left( e^{-i\lambda _{3}}\mathbf{;\theta }\right) -%
\mathbf{S}_{4}\mathbf{S}_{4}^{\prime }\mathbf{\bar{\Lambda}}^{(\ell )}\left( 
\mathbf{\theta }\right) \right \} \otimes \mathbf{I}_{d}\otimes \mathbf{I}%
_{d} \\
\mathbf{B}_{4,2}^{\left( \ell \right) }\left( \mathbf{\lambda };\mathbf{%
\theta }\right) &=&\mathbf{B}_{4,1}^{\left( \ell \right) }\left( \lambda
_{4};\mathbf{\theta }\right) =\left \{ \mathbf{\dot{\Lambda}}^{(\ell
)}\left( e^{-i\lambda _{4}}\mathbf{;\theta }\right) -\mathbf{S}_{4}\mathbf{S}%
_{4}^{\prime }\mathbf{\bar{\Lambda}}^{(\ell )}\left( \mathbf{\theta }\right)
\right \} \otimes \mathbf{I}_{d}\otimes \mathbf{I}_{d}\otimes \mathbf{I}_{d},
\end{eqnarray*}%
because 
\begin{eqnarray*}
\mathbf{\Psi }_{4}^{-1}(\mathbf{\lambda };\mathbf{\theta })\mathbf{\dot{\Psi}%
}_{4}^{\left( \ell \right) }(\mathbf{\lambda };\mathbf{\theta }) &=&\left( 
\mathbf{I}_{d}\otimes \mathbf{I}_{d}\otimes \mathbf{I}_{d}\otimes \mathbf{%
\dot{\Lambda}}^{(\ell )}\left( \lambda _{1}\mathbf{;\theta }\right) \right)
+\left( \mathbf{I}_{d}\otimes \mathbf{I}_{d}\otimes \mathbf{\dot{\Lambda}}%
^{(\ell )}\left( \lambda _{2}\mathbf{;\theta }\right) \otimes \mathbf{I}%
_{d}\right) \\
&+&\left( \mathbf{I}_{d}\otimes \mathbf{\dot{\Lambda}}^{(\ell )}\left(
\lambda _{3}\mathbf{;\theta }\right) \otimes \mathbf{I}_{d}\otimes \mathbf{I}%
_{d}\right) +\left( \mathbf{\dot{\Lambda}}^{(\ell )}\left( -\lambda
_{1}-\lambda _{2}-\lambda _{3}\mathbf{;\theta }\right) \otimes \mathbf{I}%
_{d}\otimes \mathbf{I}_{d}\otimes \mathbf{I}_{d}\right) ,
\end{eqnarray*}%
and $\mathbb{I}_{4}^{\mathbf{\varepsilon }}(\mathbf{\lambda }) = T^{-1}
w_{T}^{\mathbf{\varepsilon }}\left( -\lambda _{1}-\lambda _{2}-\lambda
_{3}\right) \otimes w_{T}^{\mathbf{\varepsilon }}\left( \lambda _{3}\right)
\otimes w_{T}^{\mathbf{\varepsilon }}\left( \lambda _{2}\right) \otimes
w_{T}^{\mathbf{\varepsilon }}\left( \lambda _{1}\right) $ is equal to
\begin{eqnarray*}
&&\frac{1}{T}\sum_{t,r,s,u=1}^{T}\exp \left( iu\left( \lambda _{1}+\lambda
_{2}+\lambda _{3}\right) -is\lambda _{3}-ir\lambda _{2}-it\lambda
_{1}\right) \left( \mathbf{\varepsilon }_{u}\otimes \mathbf{\varepsilon }%
_{s}\otimes \mathbf{\varepsilon }_{r}\otimes \mathbf{\varepsilon }_{t}\right)
\end{eqnarray*}
or 
\begin{eqnarray*}
&&\frac{1}{T}\sum_{t=1}^{T}\left( \mathbf{\varepsilon }_{t}\otimes \mathbf{%
\varepsilon }_{t}\otimes \mathbf{\varepsilon }_{t}\otimes \mathbf{%
\varepsilon }_{t}\right) +\frac{1}{T}\sum_{t=1}^{T}\sum_{r=1}^{t-1}A_{t,r}^{%
\left[ 4\right] }\left( \mathbf{\lambda }_{\mathbf{j}}\right) \left[ 
\begin{array}{c}
\mathbf{\varepsilon }_{t}\otimes \mathbf{\varepsilon }_{t}\otimes \mathbf{%
\varepsilon }_{t}\otimes \mathbf{\varepsilon }_{r} \\ 
\mathbf{\varepsilon }_{t}\otimes \mathbf{\varepsilon }_{t}\otimes \mathbf{%
\varepsilon }_{r}\otimes \mathbf{\varepsilon }_{t} \\ 
\mathbf{\varepsilon }_{t}\otimes \mathbf{\varepsilon }_{r}\otimes \mathbf{%
\varepsilon }_{t}\otimes \mathbf{\varepsilon }_{t} \\ 
\mathbf{\varepsilon }_{r}\otimes \mathbf{\varepsilon }_{t}\otimes \mathbf{%
\varepsilon }_{t}\otimes \mathbf{\varepsilon }_{t}%
\end{array}%
\right] \\
&+&\frac{1}{T}\sum_{t=1}^{T}\sum_{r,s=1}^{t-1}G_{t,r,s}^{\left[ 4\right]
}\left( \mathbf{\lambda }_{\mathbf{j}}\right) \left[ 
\begin{array}{c}
\mathbf{\varepsilon }_{s}\otimes \mathbf{\varepsilon }_{r}\otimes \mathbf{%
\varepsilon }_{t}\otimes \mathbf{\varepsilon }_{t} \\ 
\mathbf{\varepsilon }_{s}\otimes \mathbf{\varepsilon }_{t}\otimes \mathbf{%
\varepsilon }_{r}\otimes \mathbf{\varepsilon }_{t} \\ 
\mathbf{\varepsilon }_{t}\otimes \mathbf{\varepsilon }_{s}\otimes \mathbf{%
\varepsilon }_{r}\otimes \mathbf{\varepsilon }_{t} \\ 
\mathbf{\varepsilon }_{s}\otimes \mathbf{\varepsilon }_{t}\otimes \mathbf{%
\varepsilon }_{t}\otimes \mathbf{\varepsilon }_{r} \\ 
\mathbf{\varepsilon }_{t}\otimes \mathbf{\varepsilon }_{s}\otimes \mathbf{%
\varepsilon }_{t}\otimes \mathbf{\varepsilon }_{r} \\ 
\mathbf{\varepsilon }_{t}\otimes \mathbf{\varepsilon }_{t}\otimes \mathbf{%
\varepsilon }_{s}\otimes \mathbf{\varepsilon }_{r}%
\end{array}%
\right] +\frac{1}{T}\sum_{t=1}^{T}\sum_{r,s,u=1}^{t-1}F_{t,r,s,u}^{\left[ 4%
\right] }\left( \mathbf{\lambda }_{\mathbf{j}}\right) \left[ 
\begin{array}{c}
\mathbf{\varepsilon }_{u}\otimes \mathbf{\varepsilon }_{s}\otimes \mathbf{%
\varepsilon }_{r}\otimes \mathbf{\varepsilon }_{t} \\ 
\mathbf{\varepsilon }_{u}\otimes \mathbf{\varepsilon }_{s}\otimes \mathbf{%
\varepsilon }_{t}\otimes \mathbf{\varepsilon }_{r} \\ 
\mathbf{\varepsilon }_{u}\otimes \mathbf{\varepsilon }_{t}\otimes \mathbf{%
\varepsilon }_{r}\otimes \mathbf{\varepsilon }_{s} \\ 
\mathbf{\varepsilon }_{t}\otimes \mathbf{\varepsilon }_{s}\otimes \mathbf{%
\varepsilon }_{r}\otimes \mathbf{\varepsilon }_{u}%
\end{array}%
\right]
\end{eqnarray*}%
where%
\begin{equation*}
A_{t,r}^{\left[ 4\right] }\left( \mathbf{\lambda }_{\mathbf{j}}\right)
=\left \{ A_{t,r}^{\left[ 4,n\right] }\left( \mathbf{\lambda }_{\mathbf{j}%
}\right) \right \} :=\left[ 
\begin{array}{c}
\exp \left \{ -i\left( r-t\right) \lambda _{j_{1}}\right \} \\ 
\exp \left \{ -i\left( r-t\right) \lambda _{j_{2}}\right \} \\ 
\exp \left \{ -i\left( r-t\right) \lambda _{j_{3}}\right \} \\ 
\exp \left \{ -i\left( t-r\right) \left( \lambda _{j_{1}}+\lambda
_{j_{2}}+\lambda _{j_{3}}\right) \right \}%
\end{array}%
\right] ^{\prime }\ 
\end{equation*}%
\begin{equation*}
G_{t,r,s}^{\left[ 4\right] }\left( \mathbf{\lambda }_{\mathbf{j}}\right)
=\left \{ G_{t,r,s}^{\left[ 4,n\right] }\left( \mathbf{\lambda }_{\mathbf{j}%
}\right) \right \} :=\left[ 
\begin{array}{c}
\exp \left \{ -i\left( \left( t-s\right) \lambda _{j_{1}}+\left( t-s\right)
\lambda _{j_{2}}+\left( r-s\right) \lambda _{j_{3}}\right) \right \} \\ 
\exp \left \{ -i\left( \left( t-s\right) \lambda _{j_{1}}+\left( r-s\right)
\lambda _{j_{2}}+\left( t-s\right) \lambda _{j_{3}}\right) \right \} \\ 
\exp \left \{ -i\left( \left( r-t\right) \lambda _{j_{2}}+\left( s-t\right)
\lambda _{j_{3}}\right) \right \} \\ 
\exp \left \{ -i\left( \left( r-s\right) \lambda _{j_{1}}+\left( t-s\right)
\lambda _{j_{2}}+\left( t-s\right) \lambda _{j_{3}}\right) \right \} \\ 
\exp \left \{ -i\left( \left( r-t\right) \lambda _{j_{1}}+\left( s-t\right)
\lambda _{j_{3}}\right) \right \} \\ 
\exp \left \{ -i\left( \left( r-t\right) \lambda _{j_{1}}+\left( s-t\right)
\lambda _{j_{2}}\right) \right \}%
\end{array}%
\right] ^{\prime }
\end{equation*}%
\begin{equation*}
F_{t,r,s,u}^{\left[ 4\right] }\left( \mathbf{\lambda }_{\mathbf{j}}\right)
=\left \{ F_{t,r,s,u}^{\left[ 4,n\right] }\left( \mathbf{\lambda }_{\mathbf{j%
}}\right) \right \} :=\left[ 
\begin{array}{c}
\exp \left \{ -i\left( \left( t-u\right) \lambda _{j_{1}}+\left( r-u\right)
\lambda _{j_{2}}+\left( s-u\right) \lambda _{j_{3}}\right) \right \} \\ 
\exp \left \{ -i\left( \left( r-u\right) \lambda _{j_{1}}+\left( t-u\right)
\lambda _{j_{2}}+\left( s-u\right) \lambda _{j_{3}}\right) \right \} \\ 
\exp \left \{ -i\left( \left( s-u\right) \lambda _{j_{1}}+\left( r-u\right)
\lambda _{j_{2}}+\left( t-u\right) \lambda _{j_{3}}\right) \right \} \\ 
\exp \left \{ -i\left( \left( u-t\right) \lambda _{j_{1}}+\left( r-t\right)
\lambda _{j_{2}}+\left( s-t\right) \lambda _{j_{3}}\right) \right \}%
\end{array}%
\right] ^{\prime }.
\end{equation*}%

Then 
\begin{equation*}
\mathbf{Z}_{4,T}=\frac{T^{1/2}}{T^{2}}\sum_{\mathbf{\lambda }_{\mathbf{j}}}%
\func{Re}\left \{ \mathbf{B}_{4}^{\ast }\left( \mathbf{\lambda }_{\mathbf{j}%
}\right) \left( \mathbb{I}_{4}^{\varepsilon }(\mathbf{\lambda }_{\mathbf{j}%
})-E\left[ \mathbb{I}_{4}^{\varepsilon }(\mathbf{\lambda }_{\mathbf{j}})%
\right] \right) \right \} =\sum_{t=1}^{T}Z_{4,t} =\sum_{t=1}^{T} \sum_{n=0}^{3} Z^{(n)}_{4,t},
\end{equation*}%
where, noting that vec$\left( \mathbf{I}_{d^{4}}\right)
=\sum_{a=1}^{d}\left( \mathbf{e}_{a}\otimes \mathbf{e}_{a}\otimes \mathbf{e}%
_{a}\otimes \mathbf{e}_{a}\right) ,$ $\mathbf{e}_{a}$ equal to the $a$-th
column of $\mathbf{I}_{d},$}%
{\small 
\begin{equation*}
Z_{4,t}^{\left( 0\right) }:=\frac{1}{T^{1/2}}\mathbf{C}_{T}^{\left[ 4,0%
\right] }\left \{ \left( \mathbf{\varepsilon }_{t}\otimes \mathbf{\varepsilon 
}_{t}\otimes \mathbf{\varepsilon }_{t}\otimes \mathbf{\varepsilon }%
_{t}\right) -\sum_{a=1}^{d}\mathbf{\kappa }_{4,a}\left( \mathbf{e}%
_{a}\otimes \mathbf{e}_{a}\otimes \mathbf{e}_{a}\otimes \mathbf{e}%
_{a}\right) \right \} 
\end{equation*}%
\begin{eqnarray*}
Z_{4,t}^{\left( 1\right) } &:=&\frac{1}{T^{1/2}}\sum_{r=1}^{t-1}\left \{ 
\mathbf{C}_{T}^{\left[ 4,1\right] }\left( t-r\right) \left( \mathbf{%
\varepsilon }_{t}\otimes \mathbf{\varepsilon }_{t}\otimes \mathbf{%
\varepsilon }_{t}\otimes \mathbf{\varepsilon }_{r}\right) +\mathbf{C}_{T}^{%
\left[ 4,2\right] }\left( t-r\right) \left( \mathbf{\varepsilon }_{t}\otimes 
\mathbf{\varepsilon }_{t}\otimes \mathbf{\varepsilon }_{r}\otimes \mathbf{%
\varepsilon }_{t}\right) \right \}  \\
&&+\frac{1}{T^{1/2}}\sum_{r=1}^{t-1}\left \{ \mathbf{C}_{T}^{\left[ 4,3\right]
}\left( t-r\right) \left( \mathbf{\varepsilon }_{t}\otimes \mathbf{%
\varepsilon }_{r}\otimes \mathbf{\varepsilon }_{t}\otimes \mathbf{%
\varepsilon }_{t}\right) +\mathbf{C}_{T}^{\left[ 4,4\right] }\left(
t-r\right) \left( \mathbf{\varepsilon }_{r}\otimes \mathbf{\varepsilon }%
_{t}\otimes \mathbf{\varepsilon }_{t}\otimes \mathbf{\varepsilon }%
_{t}\right) \right \} 
\end{eqnarray*}%
\begin{eqnarray*}
Z_{4,t}^{\left( 2\right) } &:=&\frac{1}{T^{1/2}}\sum_{r,s=1}^{t-1}\left \{ 
\mathbf{G}_{T}^{\left[ 4,1\right] }\left( t-r,t-s\right) \left( \mathbf{%
\varepsilon }_{s}\otimes \mathbf{\varepsilon }_{r}\otimes \mathbf{%
\varepsilon }_{t}\otimes \mathbf{\varepsilon }_{t}\right) +\mathbf{G}_{T}^{%
\left[ 4,2\right] }\left( t-r,t-s\right) \left( \mathbf{\varepsilon }%
_{s}\otimes \mathbf{\varepsilon }_{t}\otimes \mathbf{\varepsilon }%
_{r}\otimes \mathbf{\varepsilon }_{t}\right) \right \}  \\
&&+\frac{1}{T^{1/2}}\sum_{r,s=1}^{t-1}\left \{ \mathbf{G}_{T}^{\left[ 4,3%
\right] }\left( t-r,t-s\right) \left( \mathbf{\varepsilon }_{t}\otimes 
\mathbf{\varepsilon }_{s}\otimes \mathbf{\varepsilon }_{r}\otimes \mathbf{%
\varepsilon }_{t}\right) +\mathbf{G}_{T}^{\left[ 4,4\right] }\left(
t-r,t-s\right) \left( \mathbf{\varepsilon }_{s}\otimes \mathbf{\varepsilon }%
_{t}\otimes \mathbf{\varepsilon }_{t}\otimes \mathbf{\varepsilon }%
_{r}\right) \right \}  \\
&&+\frac{1}{T^{1/2}}\sum_{r,s=1}^{t-1}\left \{ \mathbf{G}_{T}^{\left[ 4,5%
\right] }\left( t-r,t-s\right) \left( \mathbf{\varepsilon }_{t}\otimes 
\mathbf{\varepsilon }_{s}\otimes \mathbf{\varepsilon }_{t}\otimes \mathbf{%
\varepsilon }_{r}\right) +\mathbf{G}_{T}^{\left[ 4,6\right] }\left(
t-r,t-s\right) \left( \mathbf{\varepsilon }_{t}\otimes \mathbf{\varepsilon }%
_{t}\otimes \mathbf{\varepsilon }_{s}\otimes \mathbf{\varepsilon }%
_{r}\right) \right \} 
\end{eqnarray*}%
\begin{eqnarray*}
Z_{4,t}^{\left( 3\right) } &:=&\frac{1}{T^{1/2}}\sum_{r,s,u=1}^{t-1}\mathbf{F}%
_{T}^{\left[ 4,1\right] }\left( t-r,t-s,t-u\right) \left( \mathbf{%
\varepsilon }_{u}\otimes \mathbf{\varepsilon }_{s}\otimes \mathbf{%
\varepsilon }_{r}\otimes \mathbf{\varepsilon }_{t}\right)  \\
&&+\frac{1}{T^{1/2}}\sum_{r,s,u=1}^{t-1}\mathbf{F}_{T}^{\left[ 4,2\right]
}\left( t-r,t-s,t-u\right) \left( \mathbf{\varepsilon }_{u}\otimes \mathbf{%
\varepsilon }_{s}\otimes \mathbf{\varepsilon }_{t}\otimes \mathbf{%
\varepsilon }_{r}\right)  \\
&&+\frac{1}{T^{1/2}}\sum_{r,s,u=1}^{t-1}\mathbf{F}_{T}^{\left[ 4,3\right]
}\left( t-r,t-s,t-u\right) \left( \mathbf{\varepsilon }_{u}\otimes \mathbf{%
\varepsilon }_{t}\otimes \mathbf{\varepsilon }_{r}\otimes \mathbf{%
\varepsilon }_{s}\right)  \\
&&+\frac{1}{T^{1/2}}\sum_{r,s,u=1}^{t-1}\mathbf{F}_{T}^{\left[ 4,4\right]
}\left( t-r,t-s,t-u\right) \left( \mathbf{\varepsilon }_{t}\otimes \mathbf{%
\varepsilon }_{s}\otimes \mathbf{\varepsilon }_{r}\otimes \mathbf{%
\varepsilon }_{u}\right) .
\end{eqnarray*}%
Next,%
\begin{equation*}
\mathbf{C}_{T}^{\left[ 4,0\right] }:=\frac{1}{T^{3}}\sum_{\mathbf{\lambda }_{%
\mathbf{j}}}\mathbf{B}_{4}^{\ast }\left( \mathbf{\lambda }_{\mathbf{j}};%
\mathbf{\theta }_{0}\right) =\mathbf{C}_{4}\left( 0\right) +O\left(
T^{-1}\right) ,
\end{equation*}%
where $\mathbf{C}_{4}\left( 0\right) :=\left( 2\pi \right) ^{-3}\int_{\Pi
^{3}}\mathbf{B}_{4}^{\ast }\left( \mathbf{\lambda };\mathbf{\theta }%
_{0}\right) d\mathbf{\lambda =A}_{4}^{\prime }\left( \mathbf{\theta }%
_{0}\right) \left( \mathbf{I}_{d^{4}}-\mathbf{S}_{4}\mathbf{S}_{4}^{\prime
}\right) ,$ while the third block of terms does not contribute because for $%
t>r,s,$ as $T\rightarrow \infty ,$ 
\begin{eqnarray*}
\mathbf{C}_{T}^{\left[ 4,n\right] }\left( r-t\right)  &=&\frac{1}{T^{3}}%
\sum_{\mathbf{\lambda }_{\mathbf{j}}}\mathbf{B}_{[4,n]}^{\ast }\left( 
\mathbf{\lambda }_{\mathbf{j}\left( n\right) };\mathbf{\theta }_{0}\right)
A_{t,r}^{\left[ 4,n\right] }\left( \mathbf{\lambda }_{\mathbf{j}\left(
n\right) }\right) =\frac{1}{T}\sum_{j=1}^{T-1}\mathbf{B}_{[4,n]}^{\ast
}\left( \lambda _{j};\mathbf{\theta }_{0}\right) A_{t,r}^{\left[ 4,n\right]
}\left( \mathbf{\lambda }_{\mathbf{j}}\right)  \\
&\rightarrow &\mathbf{\mathbf{C}}^{\left[ 4,n\right] }\left( r-t\right) :=%
\frac{1}{2\pi }\int_{-\pi }^{\pi }\mathbf{B}_{[4,n]}^{\ast }\left( \lambda ;%
\mathbf{\theta }_{0}\right) \exp \left( i\left( t-r\right) \lambda \right)
d\lambda 
\end{eqnarray*}%
\begin{eqnarray*}
\mathbf{G}_{T}^{\left[ 4,n\right] }\left( r-t,s-t\right)  &=&\frac{1}{T^{3}}%
\sum_{\mathbf{\lambda }_{\mathbf{j}}}\mathbf{B}_{[4,n]}^{\ast }\left( 
\mathbf{\lambda }_{\mathbf{j}};\mathbf{\theta }_{0}\right) G_{t,r,s}^{\left[
4,n\right] }\left( \mathbf{\lambda }_{\mathbf{j}}\right) \rightarrow 0\  \  \\
\mathbf{F}_{T}^{\left[ 4,n\right] }\left( r-t,s-t,u-t\right)  &=&\frac{1}{%
T^{3}}\sum_{\mathbf{\lambda }_{\mathbf{j}}}\mathbf{B}_{[4,n]}^{\ast }\left( 
\mathbf{\lambda }_{\mathbf{j}};\mathbf{\theta }_{0}\right) F_{t,r,s,u}^{%
\left[ 4,n\right] }\left( \mathbf{\lambda }_{\mathbf{j}}\right) \rightarrow 
\mathbf{\mathbf{C}}^{\left[ 4,n\right] }\left( r-t\right) 1_{\left \{
r=s=u\right \} }.\  \  \  \  \  \  \  \  \  \  \  \  \  \  \  \  \ 
\end{eqnarray*}%

Then, proceeding as in VL, we can write%
\begin{equation*}
\mathbf{Z}_{4,T}=\sum_{t=1}^{T}Z_{4,t}^{0}+o_{p}\left( 1\right) 
\end{equation*}%
where $Z_{4,t}^{0}$ is a MDS 
\begin{equation*}
Z_{4,t}^{0}=\frac{1}{T^{1/2}}\sum_{r=1}^{t-1}\left[ \mathbf{C}_{4}\left(
0\right) \  \  \mathbf{C}_{4}\left( r-t\right) \  \  \mathbf{C}_{4}\left(
t-r\right) \right] \left[ 
\begin{array}{c}
\frac{1}{t-1}\left( \mathbf{\varepsilon }_{t}^{\otimes 4}-E\left[ \mathbf{%
\varepsilon }_{t}^{\otimes 4}\right] \right)  \\ 
\mathbf{\varepsilon }_{t,r}^{\left[ 4\right] } \\ 
\mathbf{\varepsilon }_{r,t}^{\left[ 4\right] }%
\end{array}%
\right] 
\end{equation*}%
with $\mathbf{C}_{4}\left( r-t\right) =\left[ \mathbf{C}_{\left[ 4,1\right]
}\left( r-t\right) \  \  \mathbf{C}_{\left[ 4,2\right] }\left( r-t\right) \  \ 
\mathbf{C}_{\left[ 4,3\right] }\left( r-t\right) \  \  \mathbf{C}_{\left[ 4,4%
\right] }\left( r-t\right) \right] $ and
$\mathbf{\varepsilon }_{t,r}^{\left[ 4\right] }$ defined in Section~5.
}

{\small Then, using that $\mathbb{V}\left[ \mathbf{\varepsilon }_{t,r}^{%
\left[ 4\right] }\right] =\mathbb{V}\left[ \mathbf{\varepsilon }_{r,t}^{%
\left[ 4\right] }\right] $ and $\mathbb{C}\left[ \mathbf{\varepsilon }%
_{t,r}^{\left[ 4\right] },\mathbf{\varepsilon }_{r,t}^{\left[ 4\right] }%
\right] =\mathbb{C}\left[ \mathbf{\varepsilon }_{r,t}^{\left[ 4\right] },%
\mathbf{\varepsilon }_{t,r}^{\left[ 4\right] }\right] $ and Assumption 3$(h)$%
, $h\in \{3,4,5,6\}$, 
\begin{eqnarray*}
\mathbb{V}\left[ \mathbf{Z}_{4,T}\right]  &\rightarrow &\mathbf{C}_{4}\left(
0\right) \mathbb{V}\left[ \mathbf{\varepsilon }_{t}^{\otimes 4}\right] 
\mathbf{C}_{4}^{\prime }\left( 0\right) +\sum_{j=-\infty }^{\infty }\mathbf{C%
}_{4}\left( j\right) \mathbb{V}\left[ \mathbf{\varepsilon }_{t,r}^{\left[ 4%
\right] }\right] \mathbf{C}_{4}^{\prime }\left( j\right) +\sum_{j=-\infty
}^{\infty }\mathbf{C}_{4}\left( -j\right) \mathbb{C}\left[ \mathbf{%
\varepsilon }_{t,r}^{\left[ 4\right] },\mathbf{\varepsilon }_{r,t}^{\left[ 4%
\right] }\right] \mathbf{C}_{4}^{\prime }\left( j\right)  \\
&=&\mathbf{\Phi }_{44}^{0}\left( \mathbf{\kappa }_{2}^{0},\mathbf{\kappa }%
_{3}^{0},\mathbf{\kappa }_{4}^{0},\mathbf{\kappa }_{5}^{0},\mathbf{\kappa }%
_{6}^{0}\right) +\mathbf{\Phi }_{44}\left( \mathbf{\mu }_{6}^{0},\mathbf{%
\kappa }_{3}^{0}\right) +\mathbf{\Phi }_{44}^{\dag }\left( \mathbf{\mu }%
_{4}^{0}\right) ,
\end{eqnarray*}%
where $\mathbf{\Phi }_{44}^{0}$ does not depend on $\mathbf{\kappa }_{8}^{0},
$ and $\mathbb{V}\left[ \mathbf{\varepsilon }_{t,r}^{\left[ 4\right] }\right]
$ is 
\begin{eqnarray*}
&&\sum_{ab} \! \mathbf{\mu }_{6a}\left[ \!
\begin{array}{cccc}
\mathbf{e}_{a}\mathbf{e}_{a}^{\prime }\otimes \mathbf{e}_{a}\mathbf{e}%
_{a}^{\prime }\otimes \mathbf{e}_{a}\mathbf{e}_{a}^{\prime }\otimes \mathbf{e%
}_{b}\mathbf{e}_{b}^{\prime } & \mathbf{e}_{a}\mathbf{e}_{a}^{\prime
}\otimes \mathbf{e}_{a}\mathbf{e}_{a}^{\prime }\otimes \mathbf{e}_{a}\mathbf{%
e}_{b}^{\prime }\otimes \mathbf{e}_{b}\mathbf{e}_{a}^{\prime } & \cdots  & 
\mathbf{e}_{a}\mathbf{e}_{b}^{\prime }\otimes \mathbf{e}_{a}\mathbf{e}%
_{a}^{\prime }\otimes \mathbf{e}_{a}\mathbf{e}_{a}^{\prime }\otimes \mathbf{e%
}_{b}\mathbf{e}_{a}^{\prime } \\ 
& \mathbf{e}_{a}\mathbf{e}_{a}^{\prime }\otimes \mathbf{e}_{a}\mathbf{e}%
_{a}^{\prime }\otimes \mathbf{e}_{b}\mathbf{e}_{b}^{\prime }\otimes \mathbf{e%
}_{a}\mathbf{e}_{a}^{\prime } &  & \vdots  \\ 
&  & \ddots  &  \\ 
&  &  & \mathbf{e}_{b}\mathbf{e}_{b}^{\prime }\otimes \mathbf{e}_{a}\mathbf{e%
}_{a}^{\prime }\otimes \mathbf{e}_{a}\mathbf{e}_{a}^{\prime }\otimes \mathbf{%
e}_{a}\mathbf{e}_{a}^{\prime }%
\end{array}%
\! \right]  \\
&-&\!\! \sum_{abc} \! \mathbf{\kappa }_{3a}\mathbf{\kappa }_{3a} \! \left[ \!\!
\begin{array}{cccc}
\mathbf{e}_{a}\mathbf{e}_{b}^{\prime }\otimes \mathbf{e}_{a}\mathbf{e}%
_{b}^{\prime }\otimes \mathbf{e}_{a}\mathbf{e}_{b}^{\prime }\otimes \mathbf{e%
}_{c}\mathbf{e}_{c}^{\prime } & \mathbf{e}_{a}\mathbf{e}_{b}^{\prime
}\otimes \mathbf{e}_{a}\mathbf{e}_{b}^{\prime }\otimes \mathbf{e}_{a}\mathbf{%
e}_{c}^{\prime }\otimes \mathbf{e}_{c}\mathbf{e}_{b}^{\prime } & \cdots  & 
\mathbf{e}_{a}\mathbf{e}_{c}^{\prime }\otimes \mathbf{e}_{a}\mathbf{e}%
_{b}^{\prime }\otimes \mathbf{e}_{a}\mathbf{e}_{b}^{\prime }\otimes \mathbf{e%
}_{c}\mathbf{e}_{b}^{\prime } \\ 
& \mathbf{e}_{a}\mathbf{e}_{b}^{\prime }\otimes \mathbf{e}_{a}\mathbf{e}%
_{b}^{\prime }\otimes \mathbf{e}_{c}\mathbf{e}_{c}^{\prime }\otimes \mathbf{e%
}_{a}\mathbf{e}_{b}^{\prime } &  & \vdots  \\ 
&  & \ddots  &  \\ 
&  &  & \mathbf{e}_{c}\mathbf{e}_{c}^{\prime }\otimes \mathbf{e}_{a}\mathbf{e%
}_{b}^{\prime }\otimes \mathbf{e}_{a}\mathbf{e}_{b}^{\prime }\otimes \mathbf{%
e}_{a}\mathbf{e}_{b}^{\prime }%
\end{array}%
\!\! \right] ,\  \  \  \  \  
\end{eqnarray*}%
and $\mathbb{C}\left[ \mathbf{\varepsilon }_{t,r}^{\left[ 4\right] },\mathbf{%
\varepsilon }_{r,t}^{\left[ 4\right] }\right] $ is%
\begin{equation*}
\sum_{ab}\mathbf{\mu }_{4a}\mathbf{\mu }_{4b}\left[ 
\begin{array}{cccc}
\mathbf{e}_{a}\mathbf{e}_{b}^{\prime }\otimes \mathbf{e}_{a}\mathbf{e}%
_{b}^{\prime }\otimes \mathbf{e}_{a}\mathbf{e}_{b}^{\prime }\otimes \mathbf{e%
}_{b}\mathbf{e}_{a}^{\prime } & \mathbf{e}_{a}\mathbf{e}_{b}^{\prime
}\otimes \mathbf{e}_{a}\mathbf{e}_{b}^{\prime }\otimes \mathbf{e}_{a}\mathbf{%
e}_{a}^{\prime }\otimes \mathbf{e}_{b}\mathbf{e}_{b}^{\prime } & \cdots  & 
\mathbf{e}_{a}\mathbf{e}_{a}^{\prime }\otimes \mathbf{e}_{a}\mathbf{e}%
_{b}^{\prime }\otimes \mathbf{e}_{a}\mathbf{e}_{b}^{\prime }\otimes \mathbf{e%
}_{b}\mathbf{e}_{b}^{\prime } \\ 
\vdots  & \mathbf{e}_{a}\mathbf{e}_{b}^{\prime }\otimes \mathbf{e}_{a}%
\mathbf{e}_{b}^{\prime }\otimes \mathbf{e}_{b}\mathbf{e}_{a}^{\prime
}\otimes \mathbf{e}_{a}\mathbf{e}_{b}^{\prime } &  & \vdots  \\ 
&  & \ddots  &  \\ 
\mathbf{e}_{b}\mathbf{e}_{b}^{\prime }\otimes \mathbf{e}_{a}\mathbf{e}%
_{b}^{\prime }\otimes \mathbf{e}_{a}\mathbf{e}_{b}^{\prime }\otimes \mathbf{e%
}_{a}\mathbf{e}_{a}^{\prime } & \cdots  &  & \mathbf{e}_{b}\mathbf{e}%
_{a}^{\prime }\otimes \mathbf{e}_{a}\mathbf{e}_{b}^{\prime }\otimes \mathbf{e%
}_{a}\mathbf{e}_{b}^{\prime }\otimes \mathbf{e}_{a}\mathbf{e}_{b}^{\prime }%
\end{array}%
\right] .
\end{equation*}%
$\bigskip $ }

{\small \noindent \textbf{Covariance terms. }For $k=2,3,$ using that $%
\mathbb{C}\left( \mathbf{\varepsilon }_{t}^{\left( 0\right) },\mathbf{%
\varepsilon }_{t,r}^{\left[ 3\right] }\right) =0,$ we find that%
\begin{eqnarray*}
\mathbb{C}\left( \mathbf{Z}_{2,T},\mathbf{Z}_{3,T}\right)  &\rightarrow &%
\mathbf{C}_{2}\left( 0\right) \mathbb{C}\left( \mathbf{\varepsilon }%
_{t}^{\otimes 3},\mathbf{\varepsilon }_{t}^{\otimes 3}\right) \mathbf{C}%
_{3}^{\prime }\left( 0\right) +\sum_{j=1}^{\infty }\left[ 
\begin{array}{c}
\mathbf{C}_{2}^{\prime }\left( -j\right)  \\ 
\mathbf{C}_{2}^{\prime }\left( j\right) 
\end{array}%
\right] ^{\prime }\mathbb{C}\left( \left[ 
\begin{array}{c}
\mathbf{\varepsilon }_{t,r}^{\left[ 2\right] } \\ 
\mathbf{\varepsilon }_{r,t}^{\left[ 2\right] }%
\end{array}%
\right] ,\left[ 
\begin{array}{c}
\mathbf{\varepsilon }_{t,r}^{\left[ 3\right] } \\ 
\mathbf{\varepsilon }_{r,t}^{\left[ 3\right] }%
\end{array}%
\right] \right) \left[ 
\begin{array}{c}
\mathbf{C}_{3}^{\prime }\left( -j\right)  \\ 
\mathbf{C}_{3}^{\prime }\left( j\right) 
\end{array}%
\right]  \\
&=&\mathbf{C}_{2}\left( 0\right) \  \mathbb{C}\left( \mathbf{\varepsilon }%
_{t}^{\otimes 2},\mathbf{\varepsilon }_{t}^{\otimes 3}\right) \mathbf{C}%
_{3}^{\prime }\left( 0\right) +\sum_{j=-\infty ,\neq 0}^{\infty }\mathbf{C}%
_{2}\left( j\right) \  \mathbb{C}\left( \mathbf{\varepsilon }_{r,t}^{\left[ 2%
\right] },\mathbf{\varepsilon }_{r,t}^{\left[ 3\right] }\right) \mathbf{C}%
_{3}^{\prime }\left( j\right)  \\
&&+\sum_{j=-\infty ,\neq 0}^{\infty }\mathbf{C}_{2}\left( j\right) \  \mathbb{%
C}\left( \mathbf{\varepsilon }_{r,t}^{\left[ 2\right] },\mathbf{\varepsilon }%
_{t,r}^{\left[ 3\right] }\right) \mathbf{C}_{3}^{\prime }\left( -j\right)  \\
&=&\Phi _{2,3}^{0}+\Phi _{2,3}+\Phi _{2,3}^{\dag }
\end{eqnarray*}%
because $\mathbb{C}\left( \mathbf{\varepsilon }_{t}^{\otimes 2},\  \  \mathbf{%
\varepsilon }_{t,r}^{\left[ 3\right] }\right) =\mathbb{C}\left( \mathbf{%
\varepsilon }_{t}^{\otimes 2},\  \  \mathbf{\varepsilon }_{r,t}^{\left[ 3%
\right] }\right) =0$, $\mathbb{C}\left( \mathbf{\varepsilon }_{r,t}^{\left[ 3%
\right] },\  \  \mathbf{\varepsilon }_{t}^{\otimes 2}\right) =\mathbb{C}\left( 
\mathbf{\varepsilon }_{t,r}^{\left[ 3\right] },\  \  \mathbf{\varepsilon }%
_{t}^{\otimes 2}\right) =0,$ and $\mathbb{C}\left( \mathbf{\varepsilon }%
_{t,r}^{\left[ 2\right] },\mathbf{\varepsilon }_{t,r}^{\left[ 3\right]
}\right) =\mathbb{C}\left( \mathbf{\varepsilon }_{r,t}^{\left[ 2\right] },%
\mathbf{\varepsilon }_{r,t}^{\left[ 3\right] }\right) $ with%
\begin{equation*}
\mathbb{C}\left( \mathbf{\varepsilon }_{t,r}^{\left[ 2\right] },\mathbf{%
\varepsilon }_{t,r}^{\left[ 3\right] }\right) =\sum_{a,b}\mathbf{\kappa }%
_{3,b}^{0}\left[ 
\begin{array}{ccc}
\mathbf{e}_{b}\mathbf{e}_{b}^{\prime }\otimes \mathbf{e}_{a}\mathbf{e}%
_{b}^{\prime }\otimes \mathbf{e}_{a}^{\prime } & \mathbf{e}_{b}\mathbf{e}%
_{b}^{\prime }\otimes \mathbf{e}_{a}\mathbf{e}_{a}^{\prime }\otimes \mathbf{e%
}_{b}^{\prime } & \mathbf{e}_{b}\mathbf{e}_{a}^{\prime }\otimes \mathbf{e}%
_{a}\mathbf{e}_{b}^{\prime }\otimes \mathbf{e}_{b}^{\prime } \\ 
\mathbf{e}_{a}\mathbf{e}_{b}^{\prime }\otimes \mathbf{e}_{b}\mathbf{e}%
_{b}^{\prime }\otimes \mathbf{e}_{a}^{\prime } & \mathbf{e}_{a}\mathbf{e}%
_{b}^{\prime }\otimes \mathbf{e}_{b}\mathbf{e}_{a}^{\prime }\otimes \mathbf{e%
}_{b}^{\prime } & \mathbf{e}_{a}\mathbf{e}_{a}^{\prime }\otimes \mathbf{e}%
_{b}\mathbf{e}_{b}^{\prime }\otimes \mathbf{e}_{b}^{\prime }%
\end{array}%
\right] 
\end{equation*}%
and $\mathbb{C}\left( \mathbf{\varepsilon }_{t,r}^{\left[ 2\right] },\mathbf{%
\varepsilon }_{r,t}^{\left[ 3\right] }\right) =\mathbb{C}\left( \mathbf{%
\varepsilon }_{r,t}^{\left[ 2\right] },\mathbf{\varepsilon }_{t,r}^{\left[ 3%
\right] }\right) $ with%
\begin{equation*}
\mathbb{C}\left( \mathbf{\varepsilon }_{t,r}^{\left[ 2\right] },\mathbf{%
\varepsilon }_{r,t}^{\left[ 3\right] }\right) =\sum_{a,b}\mathbf{\kappa }%
_{3,b}^{0}\left[ 
\begin{array}{ccc}
\mathbf{e}_{a}\mathbf{e}_{b}^{\prime }\otimes \mathbf{e}_{b}\mathbf{e}%
_{b}^{\prime }\otimes \mathbf{e}_{a}^{\prime } & \mathbf{e}_{a}\mathbf{e}%
_{b}^{\prime }\otimes \mathbf{e}_{b}\mathbf{e}_{a}^{\prime }\otimes \mathbf{e%
}_{b}^{\prime } & \mathbf{e}_{a}\mathbf{e}_{a}^{\prime }\otimes \mathbf{e}%
_{b}\mathbf{e}_{b}^{\prime }\otimes \mathbf{e}_{b}^{\prime } \\ 
\mathbf{e}_{b}\mathbf{e}_{b}^{\prime }\otimes \mathbf{e}_{a}\mathbf{e}%
_{b}^{\prime }\otimes \mathbf{e}_{a}^{\prime } & \mathbf{e}_{b}\mathbf{e}%
_{b}^{\prime }\otimes \mathbf{e}_{a}\mathbf{e}_{a}^{\prime }\otimes \mathbf{e%
}_{b}^{\prime } & \mathbf{e}_{b}\mathbf{e}_{a}^{\prime }\otimes \mathbf{e}%
_{a}\mathbf{e}_{b}^{\prime }\otimes \mathbf{e}_{b}^{\prime }%
\end{array}%
\right] .
\end{equation*}%

For $k=2,4,$ using now that $\mathbb{C}\left( \mathbf{\varepsilon }%
_{t}^{\left( 0\right) },\mathbf{\varepsilon }_{t,r}^{\left[ 4\right]
}\right) =0,$%
\begin{eqnarray*}
\mathbb{C}\left( \mathbf{Z}_{2,T},\mathbf{Z}_{4,T}\right)  &\rightarrow &%
\mathbf{C}_{2}\left( 0\right) \mathbb{C}\left( \mathbf{\varepsilon }%
_{t}^{\otimes 2},\mathbf{\varepsilon }_{t}^{\otimes 4}\right) \mathbf{C}%
_{4}^{\prime }\left( 0\right) +\sum_{j=1}^{\infty }\left[ 
\begin{array}{c}
\mathbf{C}_{2}^{\prime }\left( -j\right)  \\ 
\mathbf{C}_{2}^{\prime }\left( j\right) 
\end{array}%
\right] ^{\prime }\mathbb{C}\left( \left[ 
\begin{array}{c}
\mathbf{\varepsilon }_{t,r}^{\left[ 2\right] } \\ 
\mathbf{\varepsilon }_{r,t}^{\left[ 2\right] }%
\end{array}%
\right] ,\left[ 
\begin{array}{c}
\mathbf{\varepsilon }_{t,r}^{\left[ 4\right] } \\ 
\mathbf{\varepsilon }_{r,t}^{\left[ 4\right] }%
\end{array}%
\right] \right) \left[ 
\begin{array}{c}
\mathbf{C}_{4}^{\prime }\left( -j\right)  \\ 
\mathbf{C}_{4}^{\prime }\left( j\right) 
\end{array}%
\right]  \\
&=&\Phi _{2,4}^{0}+\Phi _{2,4}+\Phi _{2,4}^{\dag }
\end{eqnarray*}%
because $\mathbb{C}\left( \mathbf{\varepsilon }_{t}^{\otimes 2},\  \  \mathbf{%
\varepsilon }_{t,r}^{\left[ 4\right] }\right) =\mathbb{C}\left( \mathbf{%
\varepsilon }_{t}^{\otimes 2},\  \  \mathbf{\varepsilon }_{r-t}^{\left[ 4%
\right] }\right) =0$, $\mathbb{C}\left( \mathbf{\varepsilon }_{r-t}^{\left[ 2%
\right] },\  \  \mathbf{\varepsilon }_{t}^{\otimes 4}\right) =\mathbb{C}\left( 
\mathbf{\varepsilon }_{t,r}^{\left[ 2\right] },\  \  \mathbf{\varepsilon }%
_{t}^{\otimes 4}\right) =0,$ while 
\begin{eqnarray*}
\mathbb{C}\left( \left[ 
\begin{array}{c}
\mathbf{\varepsilon }_{t,r}^{\left[ 2\right] } \\ 
\mathbf{\varepsilon }_{r,t}^{\left[ 2\right] }%
\end{array}%
\right] ,\  \  \mathbf{\varepsilon }_{t,r}^{\left[ 4\right] }\right)  &=&%
\sum_{a,b}\mathbf{\mu }_{4,a}^{0}\left[ 
\begin{array}{cccc}
\mathbf{e}_{a}\mathbf{e}_{a}^{\prime }\otimes \mathbf{e}_{b}\mathbf{e}%
_{a}^{\prime }\otimes \mathbf{e}_{a}^{\prime }\otimes \mathbf{e}_{b}^{\prime
} & \mathbf{e}_{a}\mathbf{e}_{a}^{\prime }\otimes \mathbf{e}_{b}\mathbf{e}%
_{a}^{\prime }\otimes \mathbf{e}_{b}^{\prime }\otimes \mathbf{e}_{a}^{\prime
} & \ddots  &  \\ 
\mathbf{e}_{b}\mathbf{e}_{a}^{\prime }\otimes \mathbf{e}_{a}\mathbf{e}%
_{a}^{\prime }\otimes \mathbf{e}_{a}^{\prime }\otimes \mathbf{e}_{b}^{\prime
} & \mathbf{e}_{b}\mathbf{e}_{a}^{\prime }\otimes \mathbf{e}_{a}\mathbf{e}%
_{a}^{\prime }\otimes \mathbf{e}_{b}^{\prime }\otimes \mathbf{e}_{a}^{\prime
} &  & \ddots 
\end{array}%
\right] ,
\end{eqnarray*}%
and 
\begin{eqnarray*}
\mathbb{C}\left( \left[ 
\begin{array}{c}
\mathbf{\varepsilon }_{t,r}^{\left[ 2\right] } \\ 
\mathbf{\varepsilon }_{r,t}^{\left[ 2\right] }%
\end{array}%
\right] ,\  \  \mathbf{\varepsilon }_{r,t}^{\left[ 4\right] }\right)  &=&%
\sum_{a,b}\mathbf{\mu }_{4,a}^{0}\left[ 
\begin{array}{cccc}
\mathbf{e}_{b}\mathbf{e}_{a}^{\prime }\otimes \mathbf{e}_{a}\mathbf{e}%
_{a}^{\prime }\otimes \mathbf{e}_{a}^{\prime }\otimes \mathbf{e}_{b}^{\prime
} & \mathbf{e}_{b}\mathbf{e}_{a}^{\prime }\otimes \mathbf{e}_{a}\mathbf{e}%
_{a}^{\prime }\otimes \mathbf{e}_{b}^{\prime }\otimes \mathbf{e}_{a}^{\prime
} & \ddots  &  \\ 
\mathbf{e}_{a}\mathbf{e}_{a}^{\prime }\otimes \mathbf{e}_{b}\mathbf{e}%
_{a}^{\prime }\otimes \mathbf{e}_{a}^{\prime }\otimes \mathbf{e}_{b}^{\prime
} & \mathbf{e}_{a}\mathbf{e}_{a}^{\prime }\otimes \mathbf{e}_{b}\mathbf{e}%
_{a}^{\prime }\otimes \mathbf{e}_{b}^{\prime }\otimes \mathbf{e}_{a}^{\prime
} &  & \ddots 
\end{array}%
\right] .
\end{eqnarray*}%

Finally, for $k=3,4$. 
\begin{eqnarray*}
\mathbb{C}\left( \mathbf{Z}_{3,T},\mathbf{Z}_{4,T}\right)  &\rightarrow &%
\mathbf{C}_{3}\left( 0\right) \mathbb{C}\left( \mathbf{\varepsilon }%
_{t}^{\otimes 3},\mathbf{\varepsilon }_{t}^{\otimes 4}\right) \mathbf{C}%
_{4}^{\prime }\left( 0\right) +\sum_{j=0}^{\infty }\left[ 
\begin{array}{c}
\mathbf{C}_{3}^{\prime }\left( -j\right)  \\ 
\mathbf{C}_{3}^{\prime }\left( j\right) 
\end{array}%
\right] ^{\prime }\mathbb{C}\left( \left[ 
\begin{array}{c}
\mathbf{\varepsilon }_{t,r}^{\left[ 3\right] } \\ 
\mathbf{\varepsilon }_{r,t}^{\left[ 3\right] }%
\end{array}%
\right] ,\left[ 
\begin{array}{c}
\mathbf{\varepsilon }_{t,r}^{\left[ 4\right] } \\ 
\mathbf{\varepsilon }_{r,t}^{\left[ 4\right] }%
\end{array}%
\right] \right) \left[ 
\begin{array}{c}
\mathbf{C}_{4}\left( -j\right) ^{\prime } \\ 
\mathbf{C}_{4}\left( j\right) ^{\prime }%
\end{array}%
\right]  \\
&=&\Phi _{3,4}^{0}+\Phi _{3,4}+\Phi _{3,4}^{\dag }
\end{eqnarray*}%
because $\mathbb{C}\left( \mathbf{\varepsilon }_{t}^{\otimes 3},\  \  \mathbf{%
\varepsilon }_{t,r}^{\left[ 4\right] }\right) =\mathbb{C}\left( \mathbf{%
\varepsilon }_{t}^{\otimes 3},\  \  \mathbf{\varepsilon }_{r,t}^{\left[ 4%
\right] }\right) =0$, $\mathbb{C}\left( \mathbf{\varepsilon }_{t,r}^{\left[ 3%
\right] },\  \  \mathbf{\varepsilon }_{t}^{\otimes 4}\right) =\mathbb{C}\left( 
\mathbf{\varepsilon }_{r,t}^{\left[ 3\right] },\  \  \mathbf{\varepsilon }%
_{t}^{\otimes 4}\right) =0\ $and $\mathbb{C}\left( \mathbf{\varepsilon }%
_{t,r}^{\left[ 3\right] },\mathbf{\varepsilon }_{t,r}^{\left[ 4\right]
}\right) =\mathbb{C}\left( \mathbf{\varepsilon }_{r,t}^{\left[ 3\right] },%
\mathbf{\varepsilon }_{r,t}^{\left[ 4\right] }\right) $ is%
\begin{eqnarray*}
&&
\sum_{a,b}\mathbf{\mu }_{5a}^{0}\left( 
\begin{array}{cccc}
\mathbf{e}_{a}\mathbf{e}_{a}^{\prime }\otimes \mathbf{e}_{a}\mathbf{e}%
_{a}^{\prime }\otimes \mathbf{e}_{b}\mathbf{e}_{a}^{\prime }\otimes \mathbf{e%
}_{b}^{\prime } & \cdots  &  &  \\ 
\vdots  & \mathbf{e}_{a}\mathbf{e}_{a}^{\prime }\otimes \mathbf{e}_{b}%
\mathbf{e}_{a}^{\prime }\otimes \mathbf{e}_{a}\mathbf{e}_{b}^{\prime
}\otimes \mathbf{e}_{a}^{\prime } &  &  \\ 
&  & \mathbf{e}_{b}\mathbf{e}_{b}^{\prime }\otimes \mathbf{e}_{a}\mathbf{e}%
_{a}^{\prime }\otimes \mathbf{e}_{a}\mathbf{e}_{a}^{\prime }\otimes \mathbf{e%
}_{a}^{\prime } & \ddots 
\end{array}%
\right) ,
\end{eqnarray*}%
while $\mathbb{C}\left( \mathbf{\varepsilon }_{t,r}^{\left[ 3\right] },%
\mathbf{\varepsilon }_{r,t}^{\left[ 4\right] }\right) =\mathbb{C}\left( 
\mathbf{\varepsilon }_{r,t}^{\left[ 3\right] },\mathbf{\varepsilon }_{t,r}^{%
\left[ 4\right] }\right) $ is 
\begin{eqnarray*}
&&
\sum_{a,b}\mathbf{\mu }_{3a}^{0}\mathbf{\mu }_{4b}^{0}\left( 
\begin{array}{cccc}
\mathbf{e}_{a}\mathbf{e}_{b}^{\prime }\otimes \mathbf{e}_{a}\mathbf{e}%
_{b}^{\prime }\otimes \mathbf{e}_{b}\mathbf{e}_{b}^{\prime }\otimes \mathbf{e%
}_{a}^{\prime } & \cdots  &  &  \\ 
\vdots  & \mathbf{e}_{a}\mathbf{e}_{b}^{\prime }\otimes \mathbf{e}_{b}%
\mathbf{e}_{b}^{\prime }\otimes \mathbf{e}_{a}\mathbf{e}_{a}^{\prime
}\otimes \mathbf{e}_{b}^{\prime } &  &  \\ 
&  & \mathbf{e}_{b}\mathbf{e}_{b}^{\prime }\otimes \mathbf{e}_{a}\mathbf{e}%
_{a}^{\prime }\otimes \mathbf{e}_{a}\mathbf{e}_{b}^{\prime }\otimes \mathbf{e%
}_{b}^{\prime } & \ddots 
\end{array}%
\right) .
\end{eqnarray*}%

Then the joint CLT for averages of periodograms $\mathbf{Z}_{k,T}$
of orders $k=2,3,4$ follows as in VL using the Cramer-Wold device, noticing
that $\left( Z_{2,t}^{0\prime },Z_{3,t}^{0},Z_{4,t}^{0}\right) ^{\prime }$
is a martingale difference under independence of order $4$ with the given
asymptotic variance for independence of order $8,$ while the proof of the
convergence of conditional variances and Lindeberg Feller condition follows
as in the univariate case as they depend on the rate of decay of the norm of
the matrices $\mathbf{C}\left( j\right) $ of scores, which are similar to
the sequences $c_{j}$ in VL. \bigskip }

\section{\protect \small Appendix E: An alternative parameterization}

{\small We can formulate the $d$-dimensional VARMA$\left( p,q\right) $
process with no standardized errors by 
\begin{equation}
\Phi _{\mathbf{\theta }}\left( L\right) Y_{t}=\mu + \bar{\Theta}_{\mathbf{%
\theta }}\left( L\right) \mathbf{\varepsilon }_{t},\  \  \  \mathbf{\varepsilon 
}_{t}\sim iid_k\left( \mathbf{0},\text{v}\mathbf{\kappa }_{k}^{\text{IC}%
}\left( \mathbf{\alpha }_{k}\right) ,k=2,3,4\right)  \label{Repnew}
\end{equation}%
where 
\begin{equation*}
\bar{\Theta}_{\mathbf{\theta }}\left( L\right) =\bar{\Theta}_{0}\left( 
\mathbf{\theta }\right) +\bar{\Theta}_{1}\left( \mathbf{\theta }\right)
L+\cdots +\bar{\Theta}_{q}\left( \mathbf{\theta }\right) L^{q}
\end{equation*}%
with diag$\left \{ \bar{\Theta}_{0}\left( \mathbf{\theta }\right) \right \}
=\left( 1,\ldots ,1\right) ^{\prime }$ for all $\mathbf{\theta }\in \mathcal{%
S}\subset \mathbb{R}^{m}$ imposing a normalization on top of the
positiveness of Assumptions~6A and B on the diagonal values of $\Theta _{0}$%
, and v$\mathbf{\kappa }_{k}^{\text{IC}}\left( \mathbf{\alpha }_{k}\right) $
given in (\ref{vk3}) and (\ref{vk4}) for $k=3$ and 4, respectively, it is
also imposed for $k=2$ for a vector $\mathbf{\alpha }_{2}$ which now
includes the (positive) variance of all innovations, which are not
normalized to 1. Conversely, $\mathbf{\alpha }_{k}$ for $k=3,4,$ are now
interpreted as the level third and fourth cumulants (and not as skewness or
kurtosis coefficients), as $\bar{\Theta}_{\mathbf{\theta }}$ is now scaling
free. Likewise, $\mathbf{\theta }$ now, as $\bar{\Theta}_{\mathbf{\theta }},$
does not includes scale parameterization and the asymptotic properties of
spectral estimations are slightly simpler because of the symmetry for all $%
k=2,3,4$. Thus Assumption~3$\left( 2\right) $ implies uncorrelation, i.e.
independence of order $k=2$, so that the $k$-th order spectral density
parametric model for each time series combination $\mathbf{a}$ with
representation $\left( \ref{Repnew}\right) $ is given for all $k=2,3,4,$ by 
\begin{equation*}
f_{\mathbf{a},k}(\boldsymbol{\lambda };\mathbf{\theta },\mathbf{\alpha }%
_{k})=\left( \Phi _{\mathbf{\theta }}^{-1}\bar{\Theta}_{\mathbf{\theta }%
}\right) _{\mathbf{a}}^{\otimes k}\left( \mathbf{\lambda }\right) \text{vec}%
\left( \text{v}\mathbf{\kappa }_{k}^{\text{IC}}\left( \mathbf{\alpha }%
_{k}\right) \right) =\left( \Phi _{\mathbf{\theta }}^{-1}\bar{\Theta}_{%
\mathbf{\theta }}\right) _{\mathbf{a}}^{\otimes k}\left( \mathbf{\lambda }%
\right) \mathbf{S}_{k}\mathbf{\alpha }_{k}\mathbf{.}
\end{equation*}

Following with the same arguments as in Section~5, we can show the
consistency of $\mathbf{\hat{\alpha}}_{2,T}\left( \mathbf{\theta }\right) $
and set the simpler estimate $\mathbf{\hat{\alpha}}_{2,T}^{\dag }\left( 
\mathbf{\theta }\right) $ of the vector of variances, which is equivalent to
the usual prewhitened estimate in the frequency domain of the diagonal of
the variance covariance matrix of $\mathbf{\varepsilon }_{t},$ 
\begin{eqnarray*}
\mathbf{\hat{\alpha}}_{2,T}^{\dag }\left( \mathbf{\theta }\right) &:=&\frac{1%
}{T}\sum_{j=1}^{T-1}\func{Re}\left \{ \mathbf{S}_{2}^{\prime }\mathbf{\Psi }%
_{2}^{-1}(\lambda _{j};\mathbf{\theta })\mathbb{I}_{2}(\lambda _{j})\right \}
\\
&=&\mathbf{S}_{2}^{\prime }\frac{1}{T}\sum_{j=1}^{T-1}\func{Re}\left \{ 
\mathbf{\Psi }^{-1}(-\lambda _{j};\mathbf{\theta })w_{T}(-\lambda
_{j})\otimes \mathbf{\Psi }^{-1}(\lambda _{j};\mathbf{\theta })w_{T}(\lambda
_{j})\right \} \\
&=&\frac{1}{T}\sum_{j=1}^{T-1}\func{Re}\left \{ \mathbf{e}_{a}^{\prime }%
\mathbf{\Psi }^{-1}(-\lambda _{j};\mathbf{\theta })w_{T}(-\lambda
_{j})\otimes \mathbf{e}_{a}^{\prime }\mathbf{\Psi }^{-1}(\lambda _{j};%
\mathbf{\theta })w_{T}(\lambda _{j})\right \} _{a=1,\ldots ,d} \\
&=&\frac{1}{T}\sum_{j=1}^{T-1}\text{diag}\left \{ \mathbf{\Psi }%
^{-1}(\lambda _{j};\mathbf{\theta })I_{T}\left( \lambda _{j}\right) \mathbf{%
\Psi }^{-1}(-\lambda _{j};\mathbf{\theta })^{\prime }\right \}
\end{eqnarray*}%
where $I_{T}\left( \lambda _{j}\right) =w_{T}(\lambda _{j})w_{T}^{\prime
}(-\lambda _{j})$ is the usual (second order) periodogram matrix. Then,
following as in Theorem~$\ref{Th5},$ we can show the consistency of $\mathbf{%
\hat{\alpha}}_{2,T}^{\dag },$ and in fact, the consistency of $\mathbf{\hat{%
\theta}}_{k,T}$ without using $\mathcal{L}_{2,T}$ in the aggregated loss
function, as scaling is now excluded from $\mathbf{\theta }$ as $\mathbf{%
\bar{\Theta}}$ is normalized. Then, the three sets of cumulants (no
coefficients) can be identified independently and $\mathbf{\theta }$ can be
identified using a unique $\mathcal{L}_{k,T},$ $k=3,4.$ }

{\small To describe the asymptotic distribution of estimates, we also update
for $k=2 $ 
\begin{equation*}
\mathbf{B}_{2}\left( \mathbf{\lambda };\mathbf{\theta }\right) :=\mathbf{%
\Psi }_{2}^{-1}\left( \mathbf{\lambda };\mathbf{\theta }\right) \mathbf{\dot{%
\Psi}}_{2}(\mathbf{\lambda };\mathbf{\theta })-\mathbf{S}_{2}\mathbf{S}%
_{2}^{\prime }\bar{\mathbf{\Lambda}}_{2}\left( \mathbf{\theta }\right) ,\  \
\  \bar{\mathbf{\Lambda}}_{2}\left( \mathbf{\theta }\right) :=\left( 2\pi
\right) ^{-1}\int_{\Pi }\func{Re}\left \{ \mathbf{\Psi }_{2}^{-1}\left(
\lambda ;\mathbf{\theta }\right) \mathbf{\dot{\Psi}}_{2}(\lambda ;\mathbf{%
\theta })\right \} d\lambda
\end{equation*}%
and redefine for $\mathbf{\alpha }=\left( \mathbf{\alpha }_{2}^{\prime }%
\mathbf{,\alpha }_{3}^{\prime }\mathbf{,\alpha }_{4}^{\prime }\right)
^{\prime }$%
\begin{equation*}
\mathbf{\Sigma }\left( \mathbf{\theta ,\alpha }\right) :=\sum _{k\in 
\mathcal{K}}w_{k}\left( \mathbf{I}_{m}\otimes \mathbf{S}_{k}\mathbf{\alpha }%
_{k}\right) ^{\prime }\mathbf{H}_{k}\left( \mathbf{\theta }\right) \left( 
\mathbf{I}_{m}\otimes \mathbf{S}_{k}\mathbf{\alpha }_{k}\right)
\end{equation*}%
where now $w_{2}\geq 0$ and $\mathbf{H}_{2}$ is updated with the new $%
\mathbf{B}_{2}.$ }

{\small Then, Theorem~$\ref{Th6}$ holds for the new parameterization under
the same regularity conditions with%
\begin{equation*}
\mathbf{\delta }\left( \mathbf{\alpha }_{0}\right) :=\left[ \  \ w_{2} \left( 
\mathbf{I}_{m}\otimes \mathbf{S}_{2}\mathbf{\alpha }_{2}^{0} \right)
^{\prime }\  \  \left \vert \  \ w_{3}\left( \mathbf{I}_{m}\otimes \mathbf{S}%
_{3}\mathbf{\alpha }_{3}^{0}\right) ^{\prime }\  \  \right \vert \  \
w_{4}\left( \mathbf{I}_{m}\otimes \mathbf{S}_{4}\mathbf{\alpha }%
_{4}^{0}\right) ^{\prime }\  \  \right] .
\end{equation*}
and where now $\mathbf{C}_{2}\left( 0\right) =\left( \mathbf{I}_{d}\otimes 
\mathbf{\mathbf{C}}\left( 0\right) +\mathbf{\mathbf{C}}\left( 0\right)
\otimes \mathbf{I}_{d}\right) \left( \mathbf{I}_{d^{2}}-\mathbf{S}_{2}%
\mathbf{S}_{2}^{\prime }\right) $ neatly incorporates the effect of the
scaling estimation with $\mathbf{\Phi }_{22}^{0}\left( \mathbf{\theta }_{0};%
\mathbf{C}\right) =\mathbf{C}_{2}\left( 0\right) \mathbb{V}\left[ \mathbf{%
\varepsilon }_{t}^{\otimes 2}\right] \mathbf{C}_{2}^{\prime }\left( 0\right) 
$ not depending on fourth order cumulants for both fundamental or
non-fundamental models as in VL, unlike in the original parameterization, in
the same way as $\mathbf{\Phi }_{kk}^{0}\left( \mathbf{\theta }_{0};\mathbf{C%
}\right) $ does not depend on $\mathbf{\alpha }_{2k}^{0}$ for either
parametrization. }

{\small The properties of $\mathbf{\hat{\alpha}}_{2,T}^{\dag }\left( \mathbf{%
\hat{\theta}}_{w,T}^{\dag }\right) $ can be described in an extended version
of Theorem~$\ref{Th7},$ which covers $\mathbf{\hat{\alpha}}_{k,T}^{\dag }$
for all $k=2,3,4$ with $\mathbf{D}_{k,2}$ defined also by the general
formulation of $\mathbf{D}_{k,h},$ i.e.,%
\begin{equation*}
\mathbf{D}_{k,2}\left( 0\right) :=\mathbf{I}_{d^{k}}1_{\left \{ k=2\right \}
}-w_{2}\mathbf{\bar{\Lambda}}_{k}\left( \mathbf{\theta }_{0}\right) \left( 
\mathbf{I}_{m}\otimes \mathbf{S}_{k}\mathbf{\alpha }_{k}^{0}\right) \mathbf{%
\Sigma }^{-1}\left( \mathbf{\theta }_{0},\mathbf{\alpha }_{0}\right) \left( 
\mathbf{I}_{m}\otimes \mathbf{S}_{2}\mathbf{\alpha }_{2}^{0}\right) ^{\prime
}\mathbf{\mathbf{C}}_{2}\left( 0\right) .
\end{equation*}
}
\vspace{-.8cm}

\section{\protect \normalsize References}

\setlength{\parindent}{-.7cm}

{\normalsize Alessi L., M. Barigozzi and M. Capasso (2011).
Non-Fundamentalness in Structural Econometric Models: A Review, \textit{%
International Statistical Review}, 79, 16-47. }

{\normalsize Andrews, B., R.A. Davis and J. Breidt (2007). Rank-based
estimation for all-pass time series models, \textit{Annals of Statistics},
35, 844-869. }

{\normalsize Blanchard, O. and D. Quah (1989). The dynamic effects of
aggregate demand and supply disturbances. \textit{The American Economic
Review}, 79, 655-673. }

{\normalsize Boubacar Mainassara, Y.B. and C. Francq (2011). Estimating
structural VARMA models with uncorrelated but non-independent error terms. 
\textit{Journal of Multivariate Analysis}, 102, 496-505. }

{\normalsize Breidt, F.J., R.A. Davis and A.A. Trindade (2001). Least
absolute deviation estimation for all-pass time series models, \textit{%
Annals of Statistics}, 29, 919-946. }

{\normalsize Brillinger, D.R. (1975). \textit{Time Series: Data Analysis and
Theory}, Holden Day, San Francisco. }

{\normalsize Brillinger, D.R. (1985). Fourier inference: some methods for
the analysis of array and nongaussian series data, \textit{Water Resources
Bulletin}, 21, 744-756. }

{\normalsize Chan, K.-S., and L.-H. Ho (2004). On the unique representation
of non-Gaussian multivariate linear processes, Technical Report No 341,
Department of Statistics and Actuarial Science, The University of Iowa. }

{\normalsize Chan, K.-S., L.-H. Ho and H. Tong (2006). A Note on
Time-Reversibility of Multivariate Linear Processes, \textit{Biometrika},
93, 221-227. }

{\normalsize Comon, P. (1994). Independent Component Analysis, A New
Concept?, \textit{Signal Processing}, 36, 287-314. }

{\normalsize Gospodinov, N. and S. Ng (2015). Minimum Distance Estimation of
Possibly Noninvertible Moving Average Models.\textit{\ Journal of Business
and Economic Statistics}, 33, 403-417. }

{\normalsize Gouri\'{e}roux, C., A. Monfort and J.-P. Renne (2017).
Statistical inference for independent component analysis: Application to
structural VAR models. \textit{Journal of Econometrics}, 196, 111-126. }

{\normalsize Gouri\'{e}roux, C., A. Monfort and J.-P. Renne (2019).
Identification and Estimation in Non-Fundamental Structural VARMA Models. 
\textit{The Review of Economic Studies}, forthcoming. }

{\normalsize Granziera, E., H.R. Moon and F. Schorfheide (2018). Inference
for VARs identified with sign restrictions. \textit{Quantitative Economics}
9, 1087-1121. }

{\normalsize Hannan, E.J. (1970). \textit{Multiple Time Series}, John Wiley,
New York. }

{\normalsize Herwartz, H. and H. L\"{u}tkepohl (2014). Structural vector
autoregressions with Markov switching: Combining conventional with
statistical identification of shocks.\textit{\ Journsal of Econometrics},
183, 104-116. }

{\normalsize Hyv\"{a}rinen, A., K. Zhang, S. Shimizu and P.O. Hoyer (2010).
Estimation of a Structural Vector Autoregression Model Using
Non-Gaussianity, \textit{Journal of Machine Learning Research}, 11,
1709-1731. }

{\normalsize Jammalamadaka, S. R., T.S. Rao and G. Terdik (2006). Higher
order cumulants of random vectors and applications to statistical inference
and time series. \textit{Sankhya: The Indian Journal of Statistics}, 68,
326-356. }

{\normalsize Kumon, M. (1992). Identification of non-minimum phase transfer
function using higher-order spectrum", \textit{Annals of the Institute of
Statistical Mathematics}, 44, 239-260. }

{\normalsize 
}

{\normalsize Lanne, M. and J. Luoto (2019). GMM Estimation of Non-Gaussian
Structural Vector Autoregression. \textit{Journal of Business and Economic
Statistics}, forthcoming. }

{\normalsize Lanne, M. and H. L\"{u}tkepohl (2008). Identifying monetary
policy shocks via changes in volatility. \textit{Journal of Money, Credit
and Banking}, 40, 131-1149. }

{\normalsize Lanne, M. and H. L\"{u}tkepohl (2010). Structural Vector
Autoregressions With Nonnormal Residuals. \textit{Journal of Business and
Economic Statistics}, 28, 159-168. }

{\normalsize Lanne, M., H. L\"{u}tkepohl, H. and K. Maciejowska (2010).
Structural vector autoregressions with Markov switching, \textit{Journal of
Economic Dynamics and Control}, 34, 121-131. }

{\normalsize Lanne, M., M. Meitz and P. Saikkonen (2017). Identification and
estimation of non-Gaussian structural vector autoregressions. \textit{%
Journal of Econometrics}, 196, 288-304. }

{\normalsize Lanne, M. and P. Saikkonen (2013). Non causal Vector
Autoregression. \textit{Econometric Theory}, 29, 447-481. }

{\normalsize Lii, K.-S. and M. Rosenblatt (1982). Deconvolution and
estimation of transfer function phase and coefficients for nonGaussian
linear processes. \textit{The Annals of Statistics}, 4, 1196-1208. }

{\normalsize Lii, K.-S. and M. Rosenblatt (1992). An approximate maximum
likelihodd estimation for nonGaussian non-minimum phase moving average
processes. \textit{Journal of Multivariate Analysis}, 43, 272-299. }

{\normalsize Lii, K.-S. and M. Rosenblatt (1996). Maximum likelihodd
estimation for nonGaussian nonminimum phase ARMA sequences. \textit{%
Statistica Sinica}, 6, 1-22. }

{\normalsize Lippi, M. and Reichlin, L. (1994). VAR Analysis, Nonfundamental
Representations, Blaschke Matrices, \textit{Journal of Econometrics}, 63,
307-325. }

{\normalsize Lobato, I.N. and C. Velasco. (2018). Efficiency Improvements
for Minimum Distance Estimation of Causal and Invertible ARMA Models, 
\textit{Economics Letters}, 162, 150-152. }

{\normalsize L\"{u}tkepohl, H. and A. Net\^{s}unajev (2017). Structural
vector autoregressions with heteroskedasticity: A review of different
volatility models, \textit{Econometrics and Statistics}, 1, 2-18. }

{\normalsize Normandin, M. and L. Phnaneuf (2004). Monetary policy shocks:
Testing identification conditions under time-varying conditional volatility. 
\textit{Journal of Monetary Economics}, 51, 1217-1243. }

{\normalsize Pham, D.T. and P. Garat (1997). Blind separation of mixture of
independent source through a quasi-maximum likelihood approach, \textit{IEEE
Transactions on Signal Processing}, 45, 1712-1725. }

{\normalsize Rigobon, R. (2003). Identification through heteroskedasticity. 
\textit{Review of Economics and Statistics}, 85, 777-792. }

{\normalsize Rosenblatt, M. (1985), \textit{Stationary Sequences and Random
Fields}, Springer, New York. }

{\normalsize Rubio-Ramirez, J.F., D.F. Waggoner and T.A. Zha (2010).
Structural Vector Autoregressions: Theory of Identification and Algorithms
for Inference, \textit{Review of Economic Studies}, 77, 665-696. }

{\normalsize Velasco, C. and I.N. Lobato. (2018). Frequency Domain Minimum
Distance Inference for Possibly Noninvertible and Noncausal ARMA models, 
\textit{Annals of Statistics}, 46, 555-579. }

\end{document}